\documentclass[a4paper,DIVcalc,12pt]{amsart}
\usepackage[a4paper, left=3cm, right=3cm, top=3cm, bottom=2.5cm]{geometry} %% Seitenr?nder
\usepackage[onehalfspacing]{setspace}

\usepackage[T1]{fontenc}
\usepackage[latin1]{inputenc}
\usepackage{lmodern}
\usepackage[english]{babel}

\usepackage{amsmath}
\usepackage{amssymb}
\usepackage{amsthm}
\usepackage{mathrsfs}

\usepackage{braket}
\usepackage{dsfont}
\usepackage{faktor}
\usepackage{rotating}
\usepackage[arrow, matrix, curve]{xy}
\usepackage{mathtools}
\usepackage{hhline}
\usepackage{colortbl}
\usepackage{bbm}
\usepackage{stmaryrd}
\usepackage{paralist}
\usepackage{setspace}
\usepackage{hyperref}

\usepackage{tikz}

\usetikzlibrary{arrows,patterns,decorations.pathmorphing}

\tikzset{>=latex}
\usepackage{tkz-euclide}
\usetikzlibrary{shapes.misc}
\tikzset{cross/.style={cross out, draw=black, fill=none, minimum size=2*(#1-\pgflinewidth), inner sep=0pt, outer sep=0pt}, cross/.default={2pt}}

\DeclareFontFamily{OT1}{pzc}{}
\DeclareFontShape{OT1}{pzc}{m}{it}{<-> s * [1.10] pzcmi7t}{}
\DeclareMathAlphabet{\mathpzc}{OT1}{pzc}{m}{it}

\vfuzz \hfuzz
\newtheorem{prop}{Proposition}[section]
\newtheorem{thm}[prop]{Theorem}
\newtheorem{cor}[prop]{Corollary}
\newtheorem{lem}[prop]{Lemma}

\newtheorem{thmintro}{Theorem}

\theoremstyle{definition}
\newtheorem{defi}[prop]{Definition}
\newtheorem{con}[prop]{Construction}
\newtheorem{expl}[prop]{Example}

\newtheorem{notation}[prop]{Notation}

\theoremstyle{remark}
\newtheorem{rem}[prop]{Remark}

\newtheorem*{remintro}{Remark}

\newcommand{\iso}{\xrightarrow{\raisebox{-0.7ex}[0ex][0ex]{$\sim$}}}

\begin{document}
\numberwithin{equation}{section}
\numberwithin{figure}{section}
\numberwithin{table}{section}

\title[Tropical correspondence for del Pezzo log Calabi-Yau pairs]{Tropical correspondence for smooth del Pezzo log Calabi-Yau pairs}
\author{Tim Graefnitz}
\address{University of Hamburg \\ Department of Mathematics \\ Germany} 
\email{tim.graefnitz@gmx.de}

\begin{abstract}
Consider a log Calabi-Yau pair $(X,D)$ consisting of a smooth del Pezzo surface $X$ of degree $\geq 3$ and a smooth anticanonical divisor $D$. We prove a correspondence between genus zero logarithmic Gromov-Witten invariants of $X$ intersecting $D$ in a single point with maximal tangency and the consistent wall structure appearing in the dual intersection complex of $(X,D)$ from the Gross-Siebert reconstruction algorithm. More precisely, the logarithm of the product of functions attached to unbounded walls in the consistent wall structure gives a generating function for these invariants.
\end{abstract}

\maketitle

\setcounter{tocdepth}{1}
\tableofcontents

\section*{Introduction}								%%%

\thispagestyle{empty}

A smooth projective surface $X$ over the complex numbers $\mathbb{C}$ together with a reduced effective anticanonical divisor $D$ forms a \textit{log Calabi-Yau pair} $(X,D)$, meaning that $K_X+D$ is numerically trivial. The case where $D=D_1+\ldots+D_m$ is a cycle of smooth rational curves (\textit{maximal boundary}) has been studied in \cite{GHK1}. In \cite{GPS} it was shown that generating functions of logarithmic Gromov-Witten invariants of $X$ with maximal tangency at a single point on $D$ in this case can be read off from a certain \textit{scattering diagram}. The statement of \cite{GPS} was generalized in \cite{Bou2} to $q$-refined scattering diagrams and generating functions of higher genus Gromov-Witten invariants.

In this work we consider the somewhat complementary case with $D$ a smooth irreducible divisor. We restrict to the case where $X$ has very ample anticanonical bundle $-K_X$, i.e., is a smooth del Pezzo surface of degree $\geq 3$. In this case there is a Fano polytope $Q\subset\mathbb{R}^2$, that is, a convex lattice polytope containing the origin and with all vertices being primitive integral vectors, from which one can construct (see \S\ref{S:smoothing}) a family $(\mathfrak{X}_Q\rightarrow\mathbb{A}^2,\mathfrak{D}_Q)$ such that fixing one coordinate $s\neq 0$ on $\mathbb{A}^2$ one obtains a toric degeneration of $(X,D)$ and fixing $s=0$ gives a toric degeneration of $(X^0,D^0)$, where $X^0$ is a smooth nef toric surface admitting a $\mathbb{Q}$-Gorenstein smoothing to $X$ and $D^0=\partial X^0$ is the toric boundary.

Let $(B,\mathscr{P},\varphi)$ be the \textit{dual intersection complex} of the toric degeneration of $(X,D)$. The affine manifold with singularities $B$ is non-compact without boundary. In \cite{CPS} it is described how to construct a tropical superpotential from such a triple $(B,\mathscr{P},\varphi)$, leading to a Landau-Ginzburg model\footnote{In fact, there is additional information captured in so called \textit{gluing data} (\cite{DataI}, Definition 2.25). In this paper we always make the trivial choice, setting $s_e=1$ for any inclusion $e : \omega \rightarrow \tau$ of cells $\omega,\tau\in\mathscr{P}$, and will not mention gluing data.}. This perfectly fits into the picture, since the idea of the \textit{Gross-Siebert program} is that toric degenerations constructed from Legendre dual polarized polyhedral affine manifolds are mirror to each other, and in fact the mirror of a Fano variety together with a choice of anticanonical divisor is believed to be a Landau-Ginzburg model. The construction involves the scattering calculations described in \cite{GS11}, leading to a \textit{consistent wall structure} $\mathscr{S}_\infty$ on $(B,\mathscr{P},\varphi)$. This is a collection of codimension $1$ polyhedral subsets of $B$ (\textit{slabs} and \textit{walls}) with attached functions describing the gluing of canonical thickenings of affine pieces necessary to obtain a toric degeneration with intersection complex $(B,\mathscr{P},\varphi)$. See \cite{Inv} for an overview of this construction.

\subsection*{The statement}								%%

\label{defi:beta1}
For an effective curve class $\underline{\beta}\in H_2^+(X,\mathbb{Z})$ let $\beta$ be the class of $1$-marked stable log maps to $X$ of genus $0$, class $\underline{\beta}$ and maximal tangency with $D$ at a single unspecified point. Let $\mathscr{M}(X,\beta)$ be the moduli space of basic stable log maps of class $\beta$ (see \cite{LogGW}). By the results of \cite{LogGW} this is a proper Deligne-Mumford stack and admits a virtual fundamental class $\llbracket\mathscr{M}(X,\beta)\rrbracket$. It has virtual dimension zero. The corresponding logarithmic Gromov-Witten invariant is defined by integration, i.e., proper pushforward to a point:
\[ N_\beta = \int_{\llbracket\mathscr{M}(X,\beta)\rrbracket} 1. \]

\begin{thmintro}
\label{thm:trop}
Let $\mathfrak{H}_\beta$ be the set of tropical curves that arise as tropicalizations of stable log maps in $\mathscr{M}(X,\beta)$, see Definiton \ref{defi:H} for the precise definition. Denote by $\text{Mult}(h)$ the multiplicity of a tropical curve $h$ (Definition \ref{defi:mult}). Then
\[ N_\beta = \sum_{h\in\mathfrak{H}_\beta}\textup{Mult}(h). \]
\end{thmintro}

\begin{figure}[h!]
\renewcommand\thefigure{1}
\centering
\begin{tikzpicture}[xscale=3.6,yscale=0.9,define rgb/.code={\definecolor{mycolor}{RGB}{#1}}, rgb color/.style={define rgb={#1},mycolor},rotate=90]
\clip (-3,-1) rectangle (15,2);
\draw[->,rgb color={255,0,0}] (0.000,0.500) -- (0.000,6.50);
\draw[->,rgb color={255,0,0}] (0.000,0.500) -- (0.000,-5.50);
\draw[->,rgb color={255,0,0}] (-1.50,-0.500) -- (15.0,5.00);
\draw[->,rgb color={255,0,0}] (-1.50,1.50) -- (15.0,-4.00);
\draw[->,rgb color={255,0,0}] (-1.50,1.50) -- (-16.5,6.50);
\draw[->,rgb color={255,0,0}] (-1.50,-0.500) -- (-16.5,-5.50);
\draw[->,rgb color={255,0,0}] (-6.00,-1.50) -- (15.0,2.00);
\draw[->,rgb color={255,0,0}] (-6.00,2.50) -- (15.0,-1.00);
\draw[->,rgb color={255,0,0}] (-6.00,2.50) -- (-30.0,6.50);
\draw[->,rgb color={255,0,0}] (-6.00,-1.50) -- (-30.0,-5.50);
\draw[->,rgb color={255,0,0}] (-13.5,-2.50) -- (15.0,0.667);
\draw[->,rgb color={255,0,0}] (-13.5,3.50) -- (15.0,0.333);
\draw[->,rgb color={255,0,0}] (-13.5,3.50) -- (-40.5,6.50);
\draw[->,rgb color={255,0,0}] (-13.5,-2.50) -- (-40.5,-5.50);
\draw[->,rgb color={255,0,0}] (-24.0,-3.50) -- (15.0,-0.250);
\draw[->,rgb color={255,0,0}] (-24.0,4.50) -- (15.0,1.25);
\draw[->,rgb color={255,0,0}] (-24.0,4.50) -- (-48.0,6.50);
\draw[->,rgb color={255,0,0}] (-24.0,-3.50) -- (-48.0,-5.50);
\draw[->,rgb color={255,0,0}] (-37.5,-4.50) -- (15.0,-1.00);
\draw[->,rgb color={255,0,0}] (-37.5,5.50) -- (15.0,2.00);
\draw[->,rgb color={255,0,0}] (-37.5,5.50) -- (-52.5,6.50);
\draw[->,rgb color={255,0,0}] (-37.5,-4.50) -- (-52.5,-5.50);
\draw[->,rgb color={255,0,0}] (-54.0,-5.50) -- (15.0,-1.67);
\draw[->,rgb color={255,0,0}] (-54.0,6.50) -- (15.0,2.67);
\draw[->,red] (0.000,0.000) -- (15.0,0.000);
\draw[->,red] (-3.00,-1.00) -- (15.0,-1.00);
\draw[->,red] (0.000,1.00) -- (15.0,1.00);
\draw[->,red] (-3.00,2.00) -- (15.0,2.00);
\draw[->,red] (-45.0,-5.00) -- (15.0,-5.00);
\draw[->,red] (-30.0,-4.00) -- (15.0,-4.00);
\draw[->,red] (-18.0,-3.00) -- (15.0,-3.00);
\draw[->,red] (-9.00,-2.00) -- (15.0,-2.00);
\draw[->,red] (-9.00,3.00) -- (15.0,3.00);
\draw[->,red] (-18.0,4.00) -- (15.0,4.00);
\draw[->,red] (-30.0,5.00) -- (15.0,5.00);
\draw[->,red] (-45.0,6.00) -- (15.0,6.00);
\draw[->,black] (1.50,0.500) -- (15.0,0.500);
\draw[->,black] (-4.50,-1.50) -- (15.0,-1.50);
\draw[->,black] (0.000,-0.500) -- (15.0,-0.500);
\draw[->,black] (0.000,1.50) -- (15.0,1.50);
\draw[->,black] (-4.50,2.50) -- (15.0,2.50);
\draw[->,black] (-22.5,-3.50) -- (15.0,-3.50);
\draw[->,black] (-12.0,-2.50) -- (15.0,-2.50);
\draw[->,black] (-12.0,3.50) -- (15.0,3.50);
\draw[->,black] (-22.5,4.50) -- (15.0,4.50);
\draw[->,rgb color={255,169,0}] (0.000,0.000) -- (15.0,2.50);
\draw[->,rgb color={255,169,0}] (-18.0,-3.00) -- (15.0,-0.800);
\draw[->,rgb color={255,152,0}] (-9.00,-2.00) -- (15.0,0.000);
\draw[->,rgb color={255,132,0}] (-3.00,-1.00) -- (15.0,1.00);
\draw[->,rgb color={255,169,0}] (0.000,1.00) -- (15.0,6.00);
\draw[->,rgb color={255,132,0}] (-3.00,2.00) -- (-3.00,6.50);
\draw[->,rgb color={255,152,0}] (-9.00,3.00) -- (-19.5,6.50);
\draw[->,rgb color={255,169,0}] (-18.0,4.00) -- (-33.0,6.50);
\draw[->,rgb color={255,169,0}] (0.000,0.000) -- (15.0,-5.00);
\draw[->,rgb color={255,169,0}] (-18.0,-3.00) -- (-33.0,-5.50);
\draw[->,rgb color={255,152,0}] (-9.00,-2.00) -- (-19.5,-5.50);
\draw[->,rgb color={255,132,0}] (-3.00,-1.00) -- (-3.00,-5.50);
\draw[->,rgb color={255,169,0}] (0.000,1.00) -- (15.0,-1.50);
\draw[->,rgb color={255,132,0}] (-3.00,2.00) -- (15.0,0.000);
\draw[->,rgb color={255,152,0}] (-9.00,3.00) -- (15.0,1.00);
\draw[->,rgb color={255,169,0}] (-18.0,4.00) -- (15.0,1.80);
\draw[->,rgb color={255,152,0}] (3.00,0.000) -- (15.0,1.33);
\draw[->,rgb color={255,169,0}] (0.000,-1.00) -- (15.0,0.250);
\draw[->,rgb color={255,152,0}] (3.00,1.00) -- (15.0,3.00);
\draw[->,rgb color={255,169,0}] (0.000,2.00) -- (13.5,6.50);
\draw[->,rgb color={255,152,0}] (3.00,0.000) -- (15.0,-2.00);
\draw[->,rgb color={255,169,0}] (0.000,-1.00) -- (13.5,-5.50);
\draw[->,rgb color={255,152,0}] (3.00,1.00) -- (15.0,-0.333);
\draw[->,rgb color={255,169,0}] (0.000,2.00) -- (15.0,0.750);
\draw[->,rgb color={255,169,0}] (3.00,0.000) -- (15.0,-1.33);
\draw[->,rgb color={255,169,0}] (3.00,1.00) -- (15.0,0.000);
\draw[->,rgb color={255,152,0}] (1.50,0.500) -- (15.0,-1.00);
\draw[->,rgb color={255,169,0}] (0.000,-0.500) -- (15.0,-3.00);
\draw[->,rgb color={255,152,0}] (0.000,1.50) -- (15.0,0.250);
\draw[->,rgb color={255,152,0}] (1.50,0.500) -- (15.0,2.00);
\draw[->,rgb color={255,169,0}] (0.000,-0.500) -- (15.0,0.750);
\draw[->,rgb color={255,152,0}] (0.000,1.50) -- (15.0,4.00);
\draw[->,blue] (2.00,0.333) -- (15.0,0.333);
\draw[->,blue] (2.00,0.667) -- (15.0,0.667);
\draw[->,brown] (2.40,0.200) -- (15.0,0.200);
\draw[->,brown] (0.000,-0.800) -- (15.0,-0.800);
\draw[->,brown] (1.80,1.20) -- (15.0,1.20);
\draw[->,brown] (2.40,0.800) -- (15.0,0.800);
\draw[->,brown] (1.80,-0.200) -- (15.0,-0.200);
\draw[->,brown] (0.000,1.80) -- (15.0,1.80);
\draw[->,rgb color={255,169,0}] (3.00,0.000) -- (15.0,1.00);
\draw[->,rgb color={255,169,0}] (3.00,1.00) -- (15.0,2.33);
\draw[->,gray] (2.50,0.167) -- (15.0,0.167);
\draw[->,gray] (2.50,0.833) -- (15.0,0.833);
\draw[->,blue] (-13.0,-2.67) -- (15.0,-2.67);
\draw[->,blue] (-5.00,-1.67) -- (15.0,-1.67);
\draw[->,blue] (0.000,-0.667) -- (15.0,-0.667);
\draw[->,blue] (1.00,1.33) -- (15.0,1.33);
\draw[->,blue] (-3.00,2.33) -- (15.0,2.33);
\draw[->,blue] (-10.0,3.33) -- (15.0,3.33);
\draw[->,rgb color={255,169,0}] (2.00,0.333) -- (15.0,-0.750);
\draw[->,blue] (-10.0,-2.33) -- (15.0,-2.33);
\draw[->,blue] (-3.00,-1.33) -- (15.0,-1.33);
\draw[->,blue] (1.00,-0.333) -- (15.0,-0.333);
\draw[->,blue] (0.000,1.67) -- (15.0,1.67);
\draw[->,blue] (-5.00,2.67) -- (15.0,2.67);
\draw[->,blue] (-13.0,3.67) -- (15.0,3.67);
\draw[->,rgb color={255,169,0}] (2.00,0.667) -- (15.0,1.75);
\draw[->,rgb color={255,152,0}] (0.000,0.000) -- (15.0,3.33);
\draw[->,rgb color={255,169,0}] (-3.00,-1.00) -- (15.0,1.40);
\draw[->,rgb color={255,152,0}] (0.000,1.00) -- (8.25,6.50);
\draw[->,rgb color={255,169,0}] (-3.00,2.00) -- (-9.75,6.50);
\draw[->,rgb color={255,132,0}] (0.000,0.000) -- (15.0,1.67);
\draw[->,rgb color={255,169,0}] (-9.00,-2.00) -- (15.0,-0.400);
\draw[->,rgb color={255,152,0}] (-3.00,-1.00) -- (15.0,0.500);
\draw[->,rgb color={255,132,0}] (0.000,1.00) -- (15.0,3.50);
\draw[->,rgb color={255,152,0}] (-3.00,2.00) -- (10.5,6.50);
\draw[->,rgb color={255,169,0}] (-9.00,3.00) -- (-9.00,6.50);
\draw[->,rgb color={255,152,0}] (0.000,0.000) -- (8.25,-5.50);
\draw[->,rgb color={255,169,0}] (-3.00,-1.00) -- (-9.75,-5.50);
\draw[->,rgb color={255,152,0}] (0.000,1.00) -- (15.0,-2.33);
\draw[->,rgb color={255,169,0}] (-3.00,2.00) -- (15.0,-0.400);
\draw[->,rgb color={255,132,0}] (0.000,0.000) -- (15.0,-2.50);
\draw[->,rgb color={255,169,0}] (-9.00,-2.00) -- (-9.00,-5.50);
\draw[->,rgb color={255,152,0}] (-3.00,-1.00) -- (10.5,-5.50);
\draw[->,rgb color={255,132,0}] (0.000,1.00) -- (15.0,-0.667);
\draw[->,rgb color={255,152,0}] (-3.00,2.00) -- (15.0,0.500);
\draw[->,rgb color={255,169,0}] (-9.00,3.00) -- (15.0,1.40);
\draw[->,rgb color={255,169,0}] (1.80,0.400) -- (15.0,1.50);
\draw[->,rgb color={255,169,0}] (1.80,0.600) -- (15.0,-0.500);
\draw[->,green] (2.25,0.250) -- (15.0,0.250);
\draw[->,green] (-5.25,-1.75) -- (15.0,-1.75);
\draw[->,green] (0.000,-0.750) -- (15.0,-0.750);
\draw[->,green] (1.50,1.25) -- (15.0,1.25);
\draw[->,green] (-2.25,2.25) -- (15.0,2.25);
\draw[->,green] (2.25,0.750) -- (15.0,0.750);
\draw[->,green] (-2.25,-1.25) -- (15.0,-1.25);
\draw[->,green] (1.50,-0.250) -- (15.0,-0.250);
\draw[->,green] (0.000,1.75) -- (15.0,1.75);
\draw[->,green] (-5.25,2.75) -- (15.0,2.75);
\draw[->,brown] (1.80,0.400) -- (15.0,0.400);
\draw[->,brown] (0.000,-0.600) -- (15.0,-0.600);
\draw[->,brown] (0.600,1.40) -- (15.0,1.40);
\draw[->,brown] (1.80,0.600) -- (15.0,0.600);
\draw[->,brown] (0.600,-0.400) -- (15.0,-0.400);
\draw[->,brown] (0.000,1.60) -- (15.0,1.60);
\draw[->,rgb color={255,152,0}] (0.000,0.000) -- (15.0,1.25);
\draw[->,rgb color={255,169,0}] (-3.00,-1.00) -- (15.0,0.200);
\draw[->,rgb color={255,152,0}] (0.000,1.00) -- (15.0,2.67);
\draw[->,rgb color={255,169,0}] (-3.00,2.00) -- (15.0,5.00);
\draw[->,rgb color={255,152,0}] (0.000,0.000) -- (15.0,-1.67);
\draw[->,rgb color={255,169,0}] (-3.00,-1.00) -- (15.0,-4.00);
\draw[->,rgb color={255,152,0}] (0.000,1.00) -- (15.0,-0.250);
\draw[->,rgb color={255,169,0}] (-3.00,2.00) -- (15.0,0.800);
\draw[->,rgb color={255,169,0}] (0.000,0.000) -- (15.0,1.00);
\draw[->,rgb color={255,169,0}] (0.000,1.00) -- (15.0,2.25);
\draw[->,rgb color={255,169,0}] (0.000,0.000) -- (15.0,-1.25);
\draw[->,rgb color={255,169,0}] (0.000,1.00) -- (15.0,0.000);
\end{tikzpicture}
\caption{The wall structure of $(\mathbb{P}^2,E)$ consistent to order $6$. For the functions attached to unbounded walls see \S\ref{S:calcP2}.}
\label{fig:main}
\end{figure}
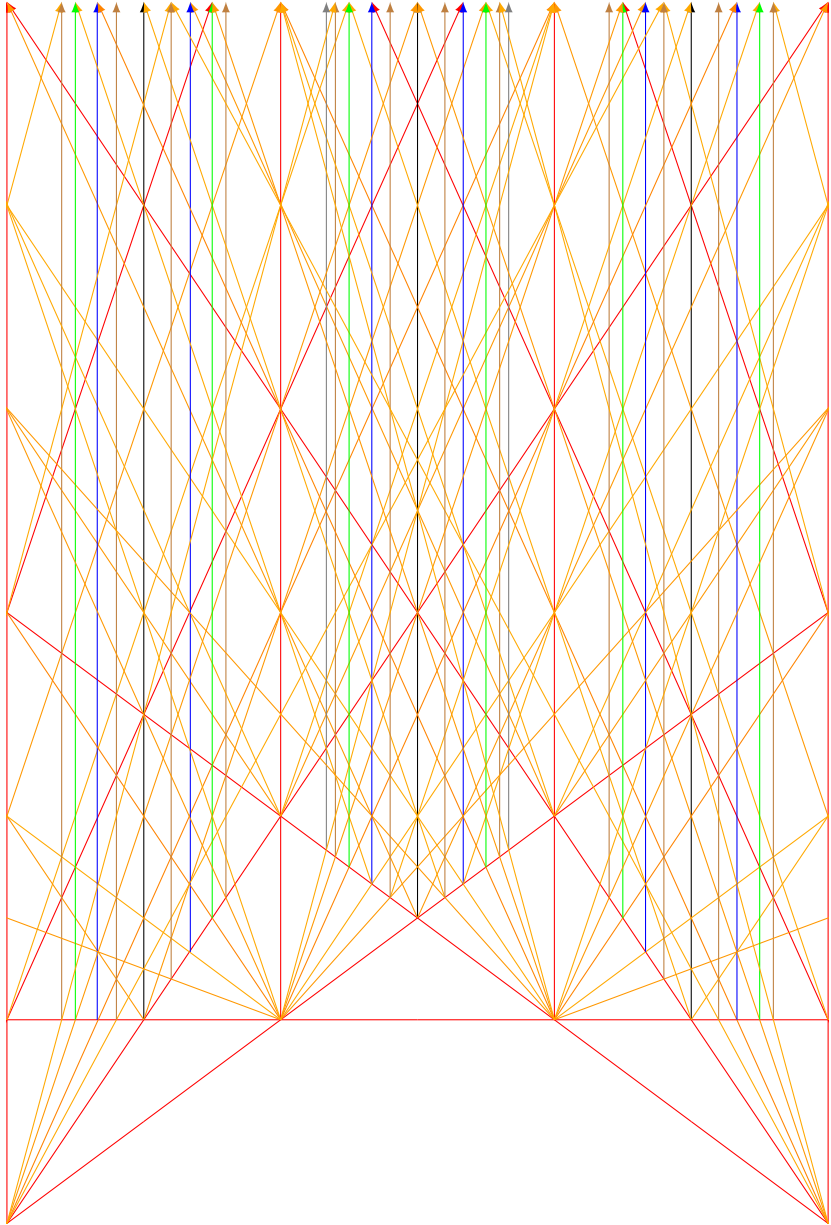

Let $\mathscr{S}_\infty$ be the consistent wall structure defined by the dual intersection complex $(B,\mathscr{P},\varphi)$ of $(X,D)$ via the Gross-Siebert algorithm \cite{GS11}. Figure \ref{fig:main} shows $\mathscr{S}_\infty$ for $(\mathbb{P}^2,E)$ up to order $6$.

The unbounded walls in $\mathscr{S}_\infty$ are all parallel in direction $m_{\textup{out}}\in \Lambda_B$. Here $\Lambda_B$ is the sheaf of integral tangent vectors on $B$ and $m_{\textup{out}}$ is the primitive vector in the unique unbounded direction of $B$ (the upward direction in Figure \ref{fig:main}). Let $f_{\textup{out}}$ be the product of all functions attached to unbounded walls in $\mathscr{S}_\infty$, regarded as elements of $\mathbb{C}\llbracket x\rrbracket$ for $x:=z^{(-m_{\textup{out}},0)}\in\mathbb{C}[\Lambda_B\oplus\mathbb{Z}]$. Then the main theorem is the following. It can be interpreted as a \textit{tropical correspondence theorem}, since the wall structure $\mathscr{S}_\infty$ is combinatorial in nature and supported on the dual intersection complex $(B,\mathscr{P},\varphi)$ of $(X,D)$. Notably, it is a tropical correspondence theorem in a \textit{non-toric} setting, as the smooth divisor $D$ has genus $1$ and thus is non-toric. So far, most such theorems have been obtained only in toric cases, a remarkable exception being \cite{Arg}.

\begin{thmintro}
\label{thm:main}
\[ \textup{log }f_{\textup{out}} = \sum_{\underline{\beta}\in H_2^+(X,\mathbb{Z})} (D \cdot\underline{\beta}) \cdot N_\beta \cdot x^{D \cdot \underline{\beta}}. \]
\end{thmintro}

For $(\mathbb{P}^2,E)$, this correspondence respects the torsion points on $E$: Consider the group law on $E$ with identity a flex point of $E$. The $3d$-torsion points form a subgroup $T_d$ of $S^1\times S^1$ isomorphic to $\mathbb{Z}_{3d}\times\mathbb{Z}_{3d}$. The stable log maps contributing to $N_d$ meet $E$ in such a $3d$-torsion point (Lemma \ref{lem:torsion}). For $P\in\cup_{d\geq 1}T_d$, let $k(P)$ be the smallest integer such that $P\in T_{k(P)}$. Let $N_{d,k}$ be the logarithmic Gromov-Witten invariant of stable log maps contributing to $N_d$ and intersecting $E$ in a point $P$ with $k(P)=k$. In \S\ref{S:torsion} we will show that this is well-defined. 

Let $s_{k,l}$ be the number of points in $T_d \simeq \mathbb{Z}_{3d}\times\mathbb{Z}_{3d}$ with $k(P)=k$ that are fixed by $M_l=\left(\begin{smallmatrix}1& 3l\\  0 & 1\end{smallmatrix}\right)$, but not fixed by $M_{l'}$ for any $l'<l$. Let $r_l$ be the number of points on $S^1$ of order $3l$, defined recursively in Lemma \ref{lem:rl}. For an unbounded wall $\mathfrak{p}\in\mathscr{S}_\infty$ let $l(\mathfrak{p})$ be the smallest integer such that $\textup{log }f_{\mathfrak{p}}$ has non-trivial $x^{3l(\mathfrak{p})}$-coefficient. The number of walls in $\mathscr{S}_\infty$ with $l(\mathfrak{p})=l$ is $r_l$.

\begin{thmintro}
\label{thm:torsion}
Let $\mathfrak{p}$ be an unbounded wall in $\mathscr{S}_\infty$ with $l(\mathfrak{p})=l$. Then
\[  \textup{log }f_{\mathfrak{p}} = \sum_{d=1}^\infty 3d \left(\sum_{k: l\mid k\mid d} \frac{s_{k,l}}{r_l} N_{d,k}\right) x^{3d}. \]
\end{thmintro}

Subtracting multiple cover contributions of curves of smaller degree, one obtains log BPS numbers $n_d$ and $n_{d,k}$ (see \S\ref{S:BPS}). Some of the $n_{d,k}$ have been calculated in \cite{Ta1}. The logarithmic Gromov-Witten invariants $N_d$ and $N_{d,k}$ and the log BPS numbers $n_d$ and $n_{d,k}$ are calculated, among other invariants, for $d\leq 6$, in \S\ref{S:calc}. Some of these numbers are new:

\begin{center}
\begin{tabular}{|l|l|l|l|l|l|} 															\hline
$n_{4,1}=14$	& $n_{4,2}=14$	& $n_{4,4}=16$		& $n_{6,1}=927$		& $n_{6,2}=938$		& $n_{6,3}=936$ \\ \hline
\end{tabular}
\end{center}

\begin{remintro}
There is a generalization of the above theorems to $q$-refined wall structures and higher genus logarithmic Gromov-Witten invariants, similar to the maximally degenerated case \cite{Bou2}, since the main argument to obtain higher genus statements in \cite{Bou1} and \cite{Bou2}, the gluing and vanishing properties of $\lambda$-classes, are purely local. We briefly sketch these ideas in \S\ref{S:genus}. See also \cite{Bou4}, Theorem 5.2.1.

An extension of the above theorems to $2$-marked invariants and broken lines will be established in \cite{Gra}. Building on this, the author jointly with Helge Ruddat and Eric Zaslow is working on an equality between the proper Landau-Ginzburg potential, defined via broken lines, and the open mirror map \cite{GRZ}.

In \cite{Bou3} Pierrick Bousseau proves an equality of the consistent wall structure $\mathscr{S}_\infty$ for $\mathbb{P}^2$ and a wall structure $\mathscr{S}_{\text{stab}}$ describing the wall crossing behavior of stability conditions on $D^b\text{Coh}(\mathbb{P}^2)$, the bounded derived of coherent sheaves on $\mathbb{P}^2$. $\mathscr{S}_{\text{stab}}$ can be interpreted as describing wall crossing of counts of coherent sheaves (generalized Donaldson-Thomas invariants) on $\mathbb{P}^2$. In \cite{Bou4}, building on \cite{Bou3} and Theorem \ref{thm:torsion} above, he proves a conjecture of Takahashi (\cite{Ta2}, Conjecture 1.6.) relating the $N_{d,k}$ with the primitive invariants $N_{d',d'}$. Yu-Shen Lin \cite{Lin} worked out a symplectic analogue of the correspondence described in this paper.
\end{remintro}

\subsection*{Motivation}									%%

The reason for an enumerative meaning of wall structures is the following. By the Strominger-Yau-Zaslow conjecture \cite{SYZ}, mirror dual Calabi-Yau varieties admit mirror dual Lagrangian torus fibrations. To construct the mirror to a given Calabi-Yau, one first constructs the \textit{semi-flat} mirror by dualizing the non-singular torus fibers. Then one corrects the complex structure of the semi-flat mirror such that it extends across the locus of singular fibers. It is expected that these corrections are determined by counts of holomorphic discs in the original variety with boundary on torus fibers \cite{SYZ}\cite{Fuk}.

Kontsevich and Soibelman \cite{KS} showed that in dimension two and with at most nodal singular fibers in the torus fibration, corrections of the complex structure are determined by algebraic self-consistency constraints which can be encoded by trees of gradient flow lines in the fan picture (dual intersection complex) of the degeneration, with certain automorphisms attached to the edges of the trees. From this they constructed a rigid analytic space from $B$, in dimension two. 

Under the discrete Legendre transform (\cite{DataI}, {\S}1.4) the gradient flow lines in the fan picture become straight lines in the cone picture (intersection complex). This was used by Gross and Siebert to construct a toric degeneration from the cone picture in any dimension. In the cone picture the self-consistency calculations are described by scattering diagrams (locally) and wall structures (globally). The fact that wall structures are used to construct a complex manifold in the cone picture and at the same time give generating functions for holomorphic curve counts of the fan picture can be seen as an explicit explanation for the connection between deformations and holomorphic curves in mirror symmetry.

In \cite{GHK1} Gross, Hacking and Keel construct the mirror to a log Calabi-Yau surface $(X,D)$ with maximal boundary. They use the above correspondence to define a canonical consistent scattering diagram from the enumerative geometry of $(X,D)$. There is an affine singularity at the vertex of this scattering diagram. Hence, the scattering diagram only gives an open subscheme of the mirror. To obtain the whole mirror they use broken lines to construct theta functions -- certain canonical global sections of line bundles on $\check{X}^\circ$. This gives enough functions to define an embedding of $\check{X}^\circ$ into an affine space. Taking the closure gives the mirror to $(X,D)$ as a partial compactification of $\check{X}^\circ$. It can be defined explicitly as the spectrum of an explicit algebra generated by theta functions, and with multiplication rule defined by the enumerative geometry of $(X,D)$. This has led to the modern viewpoint of \textit{intrinsic mirror symmetry} \cite{Intr1}\cite{Intr2}. It circumvents the constructions of scattering diagrams and broken lines and directly defines the mirror to $(X,D)$ as the spectrum of an algebra with multiplication rule defined by certain \textit{punctured Gromov-Witten invariants} of $(X,D)$ \cite{Intr1}\cite{ACGS2}.

\subsection*{Plan of the paper}								%%

In \S\ref{S:smoothing} we describe how smoothing the boundary of a Fano polytope leads to a family $(\mathfrak{X}_Q\rightarrow\mathbb{A}^2,\mathfrak{D}_Q)$ as above. Fixing one parameter $s\neq 0$ gives a toric degeneration $(\mathfrak{X}\rightarrow\mathbb{A}^1,\mathfrak{D})$ of $(X,D)$. It contains logarithmic singularities lying on the central fiber, corresponding to affine singularities in the dual intersection complex $(B,\mathscr{P},\varphi)$. In \S\ref{S:resolution} we describe a small log resolution of these singularities, leading to a log smooth degeneration $(\tilde{\mathfrak{X}}\rightarrow\mathbb{A}^1,\tilde{\mathfrak{D}})$ of $(X,D)$. In \S\ref{S:tropmap} we describe tropicalizations of stable log maps to the central fiber of $(\tilde{\mathfrak{X}}\rightarrow\mathbb{A}^1,\tilde{\mathfrak{D}})$ and show that there is a finite number of them. The tropicalizations induce a refinement of $\mathscr{P}$ and hence a logarithmic modification. This is a degeneration $(\tilde{\mathfrak{X}}_d\rightarrow\mathbb{A}^1,\tilde{\mathfrak{D}}_d)$ of $(X,D)$ such that stable log maps to the central fiber are torically transverse. This enables us to use the degeneration formula of logarithmic Gromov-Witten theory in \S\ref{S:degformula}. It gives a description of $N_\beta$ in terms of invariants $N_V$ labeled by vertices of the tropical curves found in \S\ref{S:tropmap}. In \S\ref{S:scattering} we show that the scattering calculations of \cite{GS11} give a similar formula for the logarithm of functions attached to unbounded walls in the consistent wall structure $\mathscr{S}_\infty$. This ultimately leads to a proof of Theorem \ref{thm:main}. In \S\ref{S:torsion} we explain that this correspondence respects the torsion points on $E$, leading to Theorem \ref{thm:torsion}. In \S\ref{S:genus} we discuss higher genus versions of the above statements. In \S\ref{S:calc} we explicitly calculate some invariants for $\mathbb{P}^2$, $\mathbb{P}^1\times\mathbb{P}^1$ and the cubic surface. In Appendix \ref{A:artin} we give some background on logarithmic modifications.

\subsection*{Acknowledgements}							%%

I would like to thank my supervisor Bernd Siebert for suggesting this interesting topic, as well as for his continuous support and many helpful discussions and explanations. 

I thank Helge Ruddat for his support as well as for the opportunity to present this work at several seminars and at the workshop ``Logarithmic Enumerative Geometry and Mirror Symmetry'' \cite{Oberwolfach} in Oberwolfach he organized together with Michel van Garrel and Dan Abramovich. Thanks to Michel van Garrel for carefully reading a late draft of this paper and for his continuous support. I thank Pierrick Bousseau for many comments and corrections to this paper, as well as for the suggestion and explanation of the more refined statement presented in \S\ref{S:torsion}. Thanks to Yu-Shen Lin for his invitation to present this work at seminars in Boston University and Boston College. I thank Lawrence Jack Barrott, Navid Nabijou and H\"ulya Arg\"uz for several discussions. Thanks to Simon Felten for pointing out some small errors, and to Mark Gross for sending me a Maple code for local computations of scattering diagrams.

This work was supported by the DFG funded Research Training Group 1670 ``Mathematics inspired by String theory and Quantum Field Theory''.

\section{Deforming toric degenerations}						%%%
\label{S:smoothing}

\begin{defi}
A \textit{smooth very ample log Calabi-Yau pair} is a log Calabi-Yau pair $(X,D)$ consisting of a smooth del Pezzo surface $X$ of degree $d\geq 3$ and a smooth very ample anticanonical divisor $D$.
\end{defi}

\subsection{The cone picture}								%%

\begin{con}
\label{con:family1}
Let $M\simeq\mathbb{Z}^2$ be a lattice and let $M_{\mathbb{R}}=M\otimes_{\mathbb{Z}}\mathbb{R}$ be the corresponding vector space. Let $Q\subset M_{\mathbb{R}}$ be a Fano polytope, i.e., a convex lattice polytope containing the origin and with all vertices being primitive integral vectors. The polytope $Q$ can be seen as an affine manifold via its embedding into $M_{\mathbb{R}}\simeq\mathbb{R}^2$. Let $\check{\mathscr{P}}$ be the polyhedral decomposition of $Q$ obtained by inserting edges connecting the vertices of $Q$ to the origin. Let $\check{\varphi} : Q \rightarrow \mathbb{R}$ be the strictly convex piecewise affine function on $(Q,\check{\mathscr{P}})$ defined by $\check{\varphi}(0)=0$ and $\check{\varphi}(v)=1$ for all vertices $v$ of $Q$. This means $\check{\varphi}$ is affine on the maximal cells of $\check{\mathscr{P}}$ and locally at each vertex $v$ of $\check{\mathscr{P}}$ gives a strictly convex piecewise affine function on the fan $\Sigma_v$ describing $\check{\mathscr{P}}$ locally. The triple $(Q,\check{\mathscr{P}},\check{\varphi})$ is a \textit{polarized polyhedral affine manifold} (\cite{GHS}, Construction 1.1). From this one obtains a toric degeneration of a toric del Pezzo surface with cyclic quotient singularities via the construction of Mumford \cite{Mum} (see also \cite{Inv}, \S1) as follows. Let 
\[ Q_{\check{\varphi}} = \left\{(m,h) \in M_{\mathbb{R}} \times \mathbb{R} \ | \ h \geq \varphi(m), m \in Q\right\} \]
be the convex upper hull of $\check{\varphi}$ and let 
\[ C(Q_{\check{\varphi}}) = \text{cl}\left(\mathbb{R}_{\geq 0} \cdot (Q_{\check{\varphi}} \times \{1\})\right) \subset M_{\mathbb{R}}\times\mathbb{R}\times\mathbb{R} \]
be the cone over $Q_{\check{\varphi}}$. The ring $\mathbb{C}[C(Q_{\check{\varphi}})\cap(M\times\mathbb{Z}\times\mathbb{Z})]$ is graded by the last component and we can define
\[ \mathfrak{X}^0 := \text{Proj}\left(\mathbb{C}[C(Q_{\check{\varphi}})\cap (M\times\mathbb{Z}\times\mathbb{Z})]\right). \]
By construction $\mathfrak{X}^0$ comes with an embedding into $\mathbb{P}^{N-1}\times\textup{Spec }\mathbb{C}[t]$, where $N$ is the number of lattice points of $Q$ and $t=z^{(0,1,0)}$. Projection to the last coordinate gives a toric degeneration $\mathfrak{X}^0 \rightarrow \mathbb{A}^1$ of a toric del Pezzo surface $X_0$ with quotient singularities. The polytope $Q$ is the momentum polytope of $X_0$ and the Fano condition on $Q$ corresponds to the condition on $X$ having very ample anticanonical bundle. The divisor
\[ \mathfrak{D}^0=\{z^{(0,0,0,1)}=0\} \subset \mathfrak{X}^0 \]
defined by setting the coordinate corresponding to the origin in $Q$ to zero is a very ample anticanonical divisor, since it corresponds to the pullback of the line bundle $\mathcal{O}_{\mathbb{P}^{N-1}}(1)$, which is the anticanonical bundle on the general fiber of $\mathfrak{X}^0$. By construction, the polarized polyhedral affine manifold $(Q,\check{\mathscr{P}},\check{\varphi})$ is the \textit{intersection complex} (\cite{DataI}, {\S}4.2) of the toric degeneration $\mathfrak{X}^0$.

\begin{rem}
Let $M'\subset M$ be the sublattice generated by the vertices of $\check{\mathscr{P}}$. We naturally have an embedding $\mathfrak{X}^0 \subset \mathbb{P}_{\check{B},\check{\mathscr{P}}}\times\mathbb{A}^1$, where $\mathbb{P}_{\check{B},\check{\mathscr{P}}}$ is the weighted projective space of dimension $|Q\cap M'-1|$ and weights $(1,\ldots,1,d)$ for $d$ the index of $M'$ in $M$.
\end{rem}

One can deform $\mathfrak{X}^0$ by perturbing its defining equations. This means we add a term $t^lsf$ to each equation, where $l$ is the lowest non-trivial $t$-order in the defining equations of $\mathfrak{X}^0$, $s\in\mathbb{A}^1$ is the deformation parameter and $f$ is a general polynomial defining a section of the anticanonical bundle of the general fiber of $\mathfrak{X}^0$. We give some examples below. By \cite{Pri1}, Theorem 1.1, this leads to a flat $2$-parameter family 
\[ (\mathfrak{X}_Q \rightarrow \mathbb{A}^2,\mathfrak{D}_Q) \]
such that
\begin{compactenum}[(1)]
\item for $s=0$ we have a toric degeneration $(\mathfrak{X}^0\rightarrow\mathbb{A}^1,\mathfrak{D}^0)$ of a log Calabi-Yau pair $(X^0,D^0)$ consisting of a toric del Pezzo surface with quotient singularities $X^0$ and its toric boundary $D^0=\partial X^0$;
\item for $s\neq 0$ we have a toric degeneration $(\mathfrak{X}\rightarrow\mathbb{A}^1,\mathfrak{D})$ of a smooth log Calabi-Yau pair $(X,D)$ consisting of a $\mathbb{Q}$-Gorenstein smoothing $X$ of $X^0$, i.e., a smooth del Pezzo surface of the same degree, and a smooth anticanonical divisor $D$. For different choices of $s\neq 0$ these toric degenerations are related via smooth deformation. We only care about $(X,D)$ up to smooth deformation, since log Gromov-Witten invariants are invariant under such deformations (\cite{MR}, Appendix A).
\end{compactenum}
\end{con}

\begin{notation}
We write the fibers of $(\mathfrak{X}_Q\rightarrow\mathbb{A}^2,\mathfrak{D}_Q)$ as $(X_t^s,D_t^s)$, where $s$ is the deformation parameter and $t$ is the parameter for the toric degeneration. We denote the $1$-parameter families defined by fixing one parameter by $(\mathfrak{X}_t\rightarrow\mathbb{A}^1,\mathfrak{D}_t)$ and $(\mathfrak{X}^s\rightarrow\mathbb{A}^1,\mathfrak{D}^s)$, respectively. When we fix a parameter different from zero, we sometimes omit the index, e.g. $X=X_t^s$ for $s,t\neq 0$. When writing $(\mathfrak{X}\rightarrow\mathbb{A}^1,\mathfrak{D})$ we will always mean the toric degeneration $(\mathfrak{X}^s\rightarrow\mathbb{A}^1,\mathfrak{D}^s)$ for some $s\neq 0$. By (2) above this notation makes sense. Moreover, we often supress the divisor in the notation.
\end{notation}

Let $\check{B}$ be the affine manifold with singularities obtained from $Q$ by introducing affine singularities on the interior edges of $\check{\mathscr{P}}$ such that the boundary of $\check{B}$ is a straight line, and let $(\check{B},\check{\mathscr{P}},\check{\varphi})$ be the corresponding polarized polyhedral affine manifold. Of course, there is a choice of the exact position of the affine singularities along the interior edges. In fact, one could form families of affine manifolds as in \cite{Pri1}. However, we don't care about the exact position, as we only care about the degeneration $(\mathfrak{X}\rightarrow\mathbb{A}^1,\mathfrak{D})$ up to deformation. So we may place the affine singularities in the middle of the interior edges. Note that in \cite{GS11} the affine singularities are required to have irrational coordinates, since otherwise some walls in the induced wall structure may cross them. However, this doesn't happen in our special case, so the middle of the edges will be a valid choice.

\begin{prop}
The intersection complex of $\mathfrak{X}\rightarrow\mathbb{A}^1$ is $(\check{B},\check{\mathscr{P}},\check{\varphi})$, while the intersection complex of $\mathfrak{X}^0\rightarrow\mathbb{A}^1$ is $(Q,\check{\mathscr{P}},\check{\varphi})$.
\end{prop}

\begin{proof}
First note that $\check{B}$ as above exists, since by reflexivity for any vertex $v$ the integral tangent vectors of any adjacent vertex together with $v-v_0$ generate the full lattice (see \cite{CPS}, Construction 6.2). Here $v_0$ is the unique interior vertex. 

The $t$-constant term in the defining equation for $X^0$ is independent of the variable $z^{(0,0,0,1)}$, since $(0,0,0,1)$ is the only lattice point at which $\varphi=0$. So the central fiber of $\mathfrak{X}^s\rightarrow\mathbb{A}^1$ is independent of $s$. As a consequence, the maximal cells of the intersection complex of $\mathfrak{X}^s\rightarrow\mathbb{A}^1$ are the same for each $s$. So the parameter $s$ only changes the affine structure, given by the fan structures at vertices of the intersection complex. These fan structures are defined by local models for the family at zero-dimensional toric strata of the central fiber. 

For $s\neq 0$, locally at the zero-dimensional toric stratum of $X_0^s$ corresponding to a vertex $v$ on the boundary of $B$, the family $\mathfrak{X}^s$ is given by $\{xy=t^l\}\subset\mathbb{A}^4$ for some $l>0$. So the fan structure at $v$ is given by the fan of $\mathbb{P}^1\times\mathbb{A}^1$. This shows that the boundary is a straight line.

Note that for $s=0$ locally at a $0$-dimensional stratum corresponding to $v\in\partial B$ the family $\mathfrak{X}^0$ is given by $\{xy=t^lw\}\subset\mathbb{A}^4$ for some $l>0$. The fan structure at $v$ is given by the fan with ray generators $(1,0)$, $(0,1)$ and $(1,1)$. So the affine charts are compatible and there are no affine singularities.
\end{proof}

\begin{expl}
\label{expl:P2}
Figure \ref{fig:P2} shows the intersection complex $(\check{B},\check{\mathscr{P}},\check{\varphi})$ of a toric degeneration of the log Calabi-Yau pair $(\mathbb{P}^2,E)$, where $E\subset\mathbb{P}^2$ is a smooth anticanonical divisor, i.e., an elliptic curve. This is obtained by smoothing a toric degeneration of $(\mathbb{P}^2,\partial\mathbb{P}^2)$, where $\partial\mathbb{P}^2$ is the toric boundary of $\mathbb{P}^2$. One can write down such a smoothing explicitly as follows.
\begin{eqnarray*}
\mathfrak{X}_Q &=& V\left(XYZ-t^3(W+sf_3)\right) \subset \mathbb{P}(1,1,1,3) \times \mathbb{A}^2 \\
\mathfrak{D}_Q &=& V(W)\subset\mathfrak{X}_Q
\end{eqnarray*}
Here $X,Y,Z,W$ are the coordinates of $\mathbb{P}(1,1,1,3)$, as shown in Figure \ref{fig:P2}, and $f_3$ is a general homogeneous degree $3$ polynomial in $X,Y,Z$. 

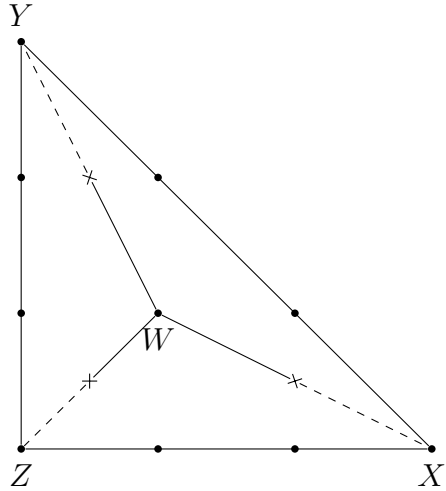
\begin{figure}[h!]
\centering
\begin{tikzpicture}[scale=1.8]
\coordinate[fill,circle,inner sep=1pt,label=below:${W}$] (0) at (0,0);
\coordinate[fill,circle,inner sep=1pt,label=below:${Z}$] (1) at (-1,-1);
\coordinate[fill,circle,inner sep=1pt,label=below:${X}$] (2) at (2,-1);
\coordinate[fill,circle,inner sep=1pt,label=above:${Y}$] (3) at (-1,2);
\coordinate[fill,cross,inner sep=2pt,rotate=45] (1a) at (-0.5,-0.5);
\coordinate[fill,cross,inner sep=2pt,rotate=63.43] (2a) at (1,-0.5);
\coordinate[fill,cross,inner sep=2pt,rotate=26.57] (3a) at (-0.5,1);
\draw (1) -- (2) -- (3) -- (1);
\draw (0) -- (1a);
\draw (0) -- (2a);
\draw (0) -- (3a);
\draw[dashed] (1a) -- (1);
\draw[dashed] (2a) -- (2);
\draw[dashed] (3a) -- (3);
\coordinate[fill,circle,inner sep=1pt] (a) at (0,-1);
\coordinate[fill,circle,inner sep=1pt] (b) at (1,-1);
\coordinate[fill,circle,inner sep=1pt] (c) at (-1,0);
\coordinate[fill,circle,inner sep=1pt] (d) at (-1,1);
\coordinate[fill,circle,inner sep=1pt] (e) at (1,0);
\coordinate[fill,circle,inner sep=1pt] (f) at (0,1);
\end{tikzpicture}
\caption{The intersection complex $(\check{B},\check{\mathscr{P}},\check{\varphi})$ of $(\mathbb{P}^2,E)$. The piecewise affine function $\check{\varphi}$ is $0$ at the center and $1$ on the boundary.}
\label{fig:P2}
\end{figure}

For $t\neq 0$ we have $X_t^s=\mathbb{P}^2$, since we can eliminate $W$ by $W=t^{-3}XYZ-sf_3$. For $s\neq 0$, $D_t^s\subset\mathbb{P}^2$ is defined by a general degree $3$ polynomial, so is an elliptic curve, and $D_t^0\subset\mathbb{P}^2$ is a cycle of three lines. For $t=0$ we have $X_0^s=V(XYZ)$ in $\mathbb{P}(1,1,1,3)$. This is a union of three $\mathbb{P}(1,1,3)$ glued along toric divisors as described by the combinatorics of Figure \ref{fig:P2}. $D_0^s$ is again a cycle of three lines.
\end{expl}

\begin{figure}[h!]
\centering
\begin{tikzpicture}[scale=1.8]
\coordinate[fill,circle,inner sep=1pt,label=below:${W}$] (0) at (0,0);
\coordinate[fill,circle,inner sep=1pt,label=right:${X}$] (1) at (1,0);
\coordinate[fill,circle,inner sep=1pt,label=above:${Y}$] (2) at (0,1);
\coordinate[fill,circle,inner sep=1pt,label=below:${Z}$] (3) at (-1,-1);
\coordinate[fill,cross,inner sep=2pt,rotate=0] (1a) at (0.5,0);
\coordinate[fill,cross,inner sep=2pt,rotate=0] (2a) at (0,0.5);
\coordinate[fill,cross,inner sep=2pt,rotate=45] (3a) at (-0.5,-0.5);
\draw (1) -- (2) -- (3) -- (1);
\draw (0) -- (1a);
\draw (0) -- (2a);
\draw (0) -- (3a);
\draw[dashed] (1a) -- (1);
\draw[dashed] (2a) -- (2);
\draw[dashed] (3a) -- (3);
\end{tikzpicture}
\caption{The intersection complex $(\check{B},\check{\mathscr{P}},\check{\varphi})$ of a smooth cubic surface, obtained by smoothing the Fano polytope of $\mathbb{P}^2/\mathbb{Z}_3$.}
\label{fig:cubic}
\end{figure}
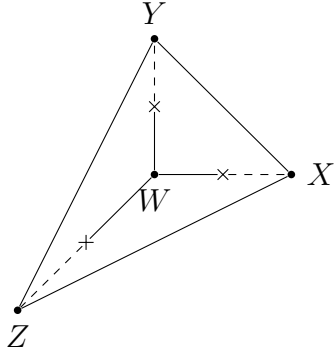

\begin{expl}
\label{expl:cubic}
Figure \ref{fig:cubic} shows the intersection complex of a toric degeneration of a smooth cubic surface $X$ (del Pezzo surface of degree $3$) obtained by smoothing the Fano polytope of the toric Gorenstein del Pezzo surface $X^0 = \mathbb{P}^2/\mathbb{Z}_3$, where $\mathbb{Z}_3$ acts by $(x,y,z) \mapsto (x,\zeta y,\zeta^{-1}z)$ for $\zeta$ a nontrivial third root of unity. This can be given explicitly as follows.
\begin{eqnarray*}
\mathfrak{X}_Q &=& V\left(XYZ-t^3(W^3+sf_3)\right) \subset \mathbb{P}^3\times\mathbb{A}^2 \\
\mathfrak{D}_Q &=& V(W) \subset \mathfrak{X}_Q
\end{eqnarray*}
Again, $X,Y,Z,W$ are the projective coordinates and $f_3$ is a general homogeneous degree $3$ polynomial in $X,Y,Z$. For $t,s\neq 0$, $X_t^s$ is a smooth cubic surface, and $D_t^s$ is a hyperplane section. For $t\neq 0,s=0$, $X_t^0$ is given by $V(XYZ-tW^3)\subset\mathbb{P}^3$, thus is a $\mathbb{Z}_3$-quotient of $\mathbb{P}^2$, and $D_t^0$ is a cycle of three lines. For $t=0$ we have $X_0^s=V(XYZ)\subset\mathbb{P}^3$. This is a union of three $\mathbb{P}^2$ glued as described by the combinatorics of Figure \ref{fig:cubic}, and again $D_0^s$ is a cycle of three lines.
\end{expl}

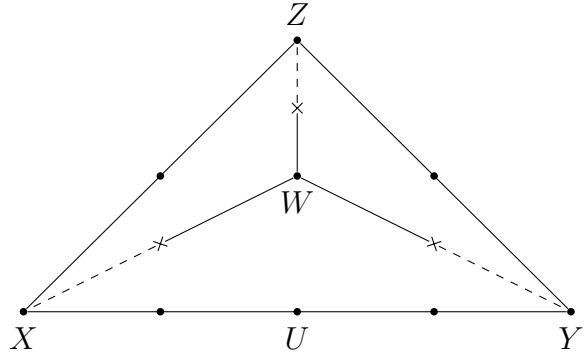
\begin{figure}[h!]
\centering
\begin{tikzpicture}[scale=1.8]
\coordinate[fill,circle,inner sep=1pt,label=below:${W}$] (0) at (0,0);
\coordinate[fill,circle,inner sep=1pt,label=above:${Z}$] (1) at (0,1);
\coordinate[fill,circle,inner sep=1pt,label=below:${Y}$] (2) at (2,-1);
\coordinate[fill,circle,inner sep=1pt,label=below:${X}$] (3) at (-2,-1);
\coordinate[fill,cross,inner sep=2pt,rotate=0] (1a) at (0,0.5);
\coordinate[fill,cross,inner sep=2pt,rotate=63.43] (2a) at (1,-0.5);
\coordinate[fill,cross,inner sep=2pt,rotate=26.57] (3a) at (-1,-0.5);
\draw (1) -- (2) -- (3) -- (1);
\draw (0) -- (1a);
\draw (0) -- (2a);
\draw (0) -- (3a);
\draw[dashed] (1a) -- (1);
\draw[dashed] (2a) -- (2);
\draw[dashed] (3a) -- (3);
\coordinate[fill,circle,inner sep=1pt] (a) at (-1,0);
\coordinate[fill,circle,inner sep=1pt] (b) at (1,0);
\coordinate[fill,circle,inner sep=1pt] (c) at (-1,-1);
\coordinate[fill,circle,inner sep=1pt,label=below:${U}$] (d) at (0,-1);
\coordinate[fill,circle,inner sep=1pt] (e) at (1,-1);
\end{tikzpicture}
\caption{The intersection complex $(\check{B},\check{\mathscr{P}},\check{\varphi})$ of $\mathbb{P}^1\times\mathbb{P}^1$, obtained by smoothing the Fano polytope of $\mathbb{P}(1,1,2)$.}
\label{fig:8'a}
\end{figure}

\begin{expl}
\label{expl:8'a}
Figure \ref{fig:8'a} shows the intersection complex of a toric degeneration of $\mathbb{P}^1\times\mathbb{P}^1$ obtained by smoothing the Fano polytope of $\mathbb{P}(1,1,2)$. This can be given explicitly as follows, with $f_2$ a general homogeneous degree $2$ polynomial in $X,Y,Z,U$ and $W$ the degree $2$ coordinate,
\begin{eqnarray*}
\mathfrak{X}_Q &=& V\left(XY-U^2+t^2sf_2,ZU-t^2(W+sf_2)\right)\subset\mathbb{P}(1,1,1,1,2)\times\mathbb{A}^2 \\
\mathfrak{D}_Q &=& V(W)\subset\mathfrak{X}_Q
\end{eqnarray*}
Indeed, $t=0$ implies $Z=0$ or $U=0$ which in turn implies $X=0$ or $Y=0$. We have $V(X)=V(Y)=\mathbb{P}(1,1,2)$ and $V(Z)=\{XY=U^2\}\subset\mathbb{P}(1,1,1,2)$ which is isomorphic to $\mathbb{P}(1,1,4)$. For $t\neq 0$ we have $X_t^0=\{XY=U^2+t^2sf_2\}\subset\mathbb{P}^3$ by elimination of $W$. For $s=0$ this is a singular quadric $X_t^0\simeq\mathbb{P}(1,1,2)$. For $s\neq 0$ it is a smooth quadric $X_t^s\simeq\mathbb{P}^1\times\mathbb{P}^1$. Again, $D_t^s$ is smooth if and only if $t,s\neq 0$.
\end{expl}

\begin{defi}
\label{defi:toricmodel}
Let $X$ be a smooth del Pezzo surface. A \textit{toric model} for $X$ is a toric del Pezzo surface with cyclic quotient singularities $X^0$ that admits a $\mathbb{Q}$-Gorenstein deformation to $X$.
\end{defi}

\begin{rem}
Note that there may be different toric models $X^0$ for the same smooth del Pezzo surface $X$. In fact, the Fano polytopes $Q$ of such $X^0$ are related via \textit{combinatorial mutations} (\cite{ACC+}, Theorem 3, see also \cite{CGG+}).
\end{rem}

\begin{prop}
For each smooth del Pezzo surface $X$ with very ample anticanonical class there exists a toric model $X^0$ with at most Gorenstein singularities.
\end{prop}

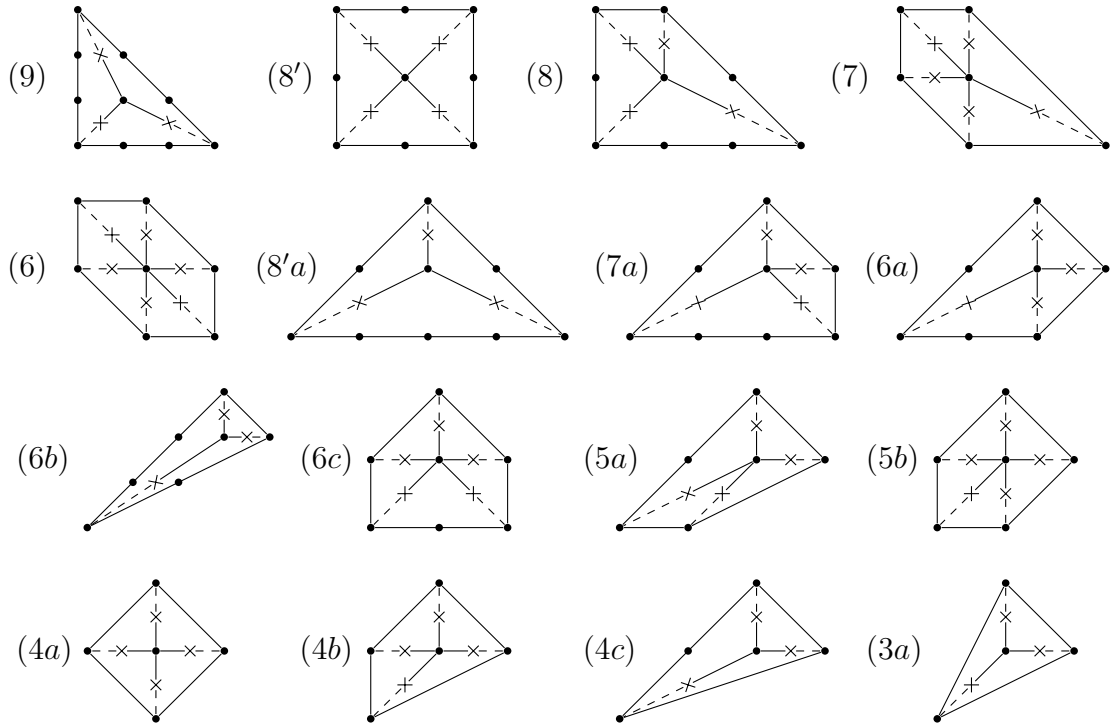
\begin{figure}[h!]
\centering
% (9)
\begin{minipage}[c]{0.05\textwidth}
$(9)$
\end{minipage}
\begin{minipage}[c]{0.16\textwidth}
\begin{tikzpicture}[scale=0.6]
\coordinate[fill,circle,inner sep=1pt] (0) at (0,0);
\coordinate[fill,circle,inner sep=1pt] (1) at (-1,-1);
\coordinate[fill,circle,inner sep=1pt] (2) at (2,-1);
\coordinate[fill,circle,inner sep=1pt] (3) at (-1,2);
\coordinate[fill,cross,inner sep=2pt,rotate=45] (1a) at (-0.5,-0.5);
\coordinate[fill,cross,inner sep=2pt,rotate=63.43] (2a) at (1,-0.5);
\coordinate[fill,cross,inner sep=2pt,rotate=26.57] (3a) at (-0.5,1);
\draw (1) -- (2) -- (3) -- (1);
\draw (0) -- (1a);
\draw (0) -- (2a);
\draw (0) -- (3a);
\draw[dashed] (1a) -- (1);
\draw[dashed] (2a) -- (2);
\draw[dashed] (3a) -- (3);
\coordinate[fill,circle,inner sep=1pt] (a) at (0,-1);
\coordinate[fill,circle,inner sep=1pt] (b) at (1,-1);
\coordinate[fill,circle,inner sep=1pt] (c) at (-1,0);
\coordinate[fill,circle,inner sep=1pt] (d) at (-1,1);
\coordinate[fill,circle,inner sep=1pt] (e) at (1,0);
\coordinate[fill,circle,inner sep=1pt] (f) at (0,1);
\end{tikzpicture}
\end{minipage}
% (8')
\begin{minipage}[c]{0.05\textwidth}
$(8')$
\end{minipage}
\begin{minipage}[c]{0.16\textwidth}
\begin{tikzpicture}[scale=0.9]
\coordinate[fill,circle,inner sep=1pt] (0) at (0,0);
\coordinate[fill,circle,inner sep=1pt] (1) at (-1,-1);
\coordinate[fill,circle,inner sep=1pt] (2) at (1,-1);
\coordinate[fill,circle,inner sep=1pt] (3) at (1,1);
\coordinate[fill,circle,inner sep=1pt] (4) at (-1,1);
\coordinate[fill,cross,inner sep=2pt,rotate=45] (1a) at (-0.5,-0.5);
\coordinate[fill,cross,inner sep=2pt,rotate=45] (2a) at (0.5,-0.5);
\coordinate[fill,cross,inner sep=2pt,rotate=45] (3a) at (0.5,0.5);
\coordinate[fill,cross,inner sep=2pt,rotate=45] (4a) at (-0.5,0.5);
\draw (1) -- (2) -- (3) -- (4) -- (1);
\draw (0) -- (1a);
\draw (0) -- (2a);
\draw (0) -- (3a);
\draw (0) -- (4a);
\draw[dashed] (1a) -- (1);
\draw[dashed] (2a) -- (2);
\draw[dashed] (3a) -- (3);
\draw[dashed] (4a) -- (4);
\coordinate[fill,circle,inner sep=1pt] (a) at (0,-1);
\coordinate[fill,circle,inner sep=1pt] (b) at (0,1);
\coordinate[fill,circle,inner sep=1pt] (c) at (-1,0);
\coordinate[fill,circle,inner sep=1pt] (d) at (1,0);
\end{tikzpicture}
\end{minipage}
% (8)
\begin{minipage}[c]{0.05\textwidth}
$(8)$
\end{minipage}
\begin{minipage}[c]{0.2\textwidth}
\begin{tikzpicture}[scale=0.9]
\coordinate[fill,circle,inner sep=1pt] (0) at (0,0);
\coordinate[fill,circle,inner sep=1pt] (1) at (-1,-1);
\coordinate[fill,circle,inner sep=1pt] (2) at (2,-1);
\coordinate[fill,circle,inner sep=1pt] (3) at (0,1);
\coordinate[fill,circle,inner sep=1pt] (4) at (-1,1);
\coordinate[fill,cross,inner sep=2pt,rotate=45] (1a) at (-0.5,-0.5);
\coordinate[fill,cross,inner sep=2pt,rotate=63.43] (2a) at (1,-0.5);
\coordinate[fill,cross,inner sep=2pt,rotate=0] (3a) at (0,0.5);
\coordinate[fill,cross,inner sep=2pt,rotate=45] (4a) at (-0.5,0.5);
\draw (1) -- (2) -- (3) -- (4) -- (1);
\draw (0) -- (1a);
\draw (0) -- (2a);
\draw (0) -- (3a);
\draw (0) -- (4a);
\draw[dashed] (1a) -- (1);
\draw[dashed] (2a) -- (2);
\draw[dashed] (3a) -- (3);
\draw[dashed] (4a) -- (4);
\coordinate[fill,circle,inner sep=1pt] (a) at (0,-1);
\coordinate[fill,circle,inner sep=1pt] (b) at (1,-1);
\coordinate[fill,circle,inner sep=1pt] (c) at (-1,0);
\coordinate[fill,circle,inner sep=1pt] (d) at (1,0);
\end{tikzpicture}
\end{minipage}
% (7)
\begin{minipage}[c]{0.05\textwidth}
$(7)$
\end{minipage}
\begin{minipage}[c]{0.22\textwidth}
\begin{tikzpicture}[scale=0.9]
\coordinate[fill,circle,inner sep=1pt] (0) at (0,0);
\coordinate[fill,circle,inner sep=1pt] (1) at (-1,0);
\coordinate[fill,circle,inner sep=1pt] (2) at (0,-1);
\coordinate[fill,circle,inner sep=1pt] (3) at (2,-1);
\coordinate[fill,circle,inner sep=1pt] (4) at (0,1);
\coordinate[fill,circle,inner sep=1pt] (5) at (-1,1);
\coordinate[fill,cross,inner sep=2pt,rotate=0] (1a) at (-0.5,0);
\coordinate[fill,cross,inner sep=2pt,rotate=0] (2a) at (0,-0.5);
\coordinate[fill,cross,inner sep=2pt,rotate=63.43] (3a) at (1,-0.5);
\coordinate[fill,cross,inner sep=2pt,rotate=0] (4a) at (0,0.5);
\coordinate[fill,cross,inner sep=2pt,rotate=45] (5a) at (-0.5,0.5);
\draw (1) -- (2) -- (3) -- (4) -- (5) -- (1);
\draw (0) -- (1a);
\draw (0) -- (2a);
\draw (0) -- (3a);
\draw (0) -- (4a);
\draw (0) -- (5a);
\draw[dashed] (1a) -- (1);
\draw[dashed] (2a) -- (2);
\draw[dashed] (3a) -- (3);
\draw[dashed] (4a) -- (4);
\draw[dashed] (5a) -- (5);
\end{tikzpicture}
\end{minipage} \\[6mm]
% (6)
\begin{minipage}[c]{0.05\textwidth}
$(6)$
\end{minipage}
\begin{minipage}[c]{0.15\textwidth}
\begin{tikzpicture}[scale=0.9]
\coordinate[fill,circle,inner sep=1pt] (0) at (0,0);
\coordinate[fill,circle,inner sep=1pt] (1) at (-1,0);
\coordinate[fill,circle,inner sep=1pt] (2) at (0,-1);
\coordinate[fill,circle,inner sep=1pt] (3) at (1,-1);
\coordinate[fill,circle,inner sep=1pt] (4) at (1,0);
\coordinate[fill,circle,inner sep=1pt] (5) at (0,1);
\coordinate[fill,circle,inner sep=1pt] (6) at (-1,1);
\coordinate[fill,cross,inner sep=2pt,rotate=0] (1a) at (-0.5,0);
\coordinate[fill,cross,inner sep=2pt,rotate=0] (2a) at (0,-0.5);
\coordinate[fill,cross,inner sep=2pt,rotate=45] (3a) at (0.5,-0.5);
\coordinate[fill,cross,inner sep=2pt,rotate=0] (4a) at (0.5,0);
\coordinate[fill,cross,inner sep=2pt,rotate=0] (5a) at (0,0.5);
\coordinate[fill,cross,inner sep=2pt,rotate=45] (6a) at (-0.5,0.5);
\draw (1) -- (2) -- (3) -- (4) -- (5) -- (6) -- (1);
\draw (0) -- (1a);
\draw (0) -- (2a);
\draw (0) -- (3a);
\draw (0) -- (4a);
\draw (0) -- (5a);
\draw (0) -- (6a);
\draw[dashed] (1a) -- (1);
\draw[dashed] (2a) -- (2);
\draw[dashed] (3a) -- (3);
\draw[dashed] (4a) -- (4);
\draw[dashed] (5a) -- (5);
\draw[dashed] (6a) -- (6);
\end{tikzpicture}
\end{minipage}
% (8a)
\begin{minipage}[c]{0.02\textwidth}
$(8'a)$
\end{minipage}
\begin{minipage}[c]{0.26\textwidth}
\begin{tikzpicture}[scale=0.9]
\coordinate[fill,circle,inner sep=1pt] (0) at (0,0);
\coordinate[fill,circle,inner sep=1pt] (1) at (0,1);
\coordinate[fill,circle,inner sep=1pt] (2) at (2,-1);
\coordinate[fill,circle,inner sep=1pt] (3) at (-2,-1);
\coordinate[fill,cross,inner sep=2pt,rotate=0] (1a) at (0,0.5);
\coordinate[fill,cross,inner sep=2pt,rotate=63.43] (2a) at (1,-0.5);
\coordinate[fill,cross,inner sep=2pt,rotate=26.57] (3a) at (-1,-0.5);
\draw (1) -- (2) -- (3) -- (1);
\draw (0) -- (1a);
\draw (0) -- (2a);
\draw (0) -- (3a);
\draw[dashed] (1a) -- (1);
\draw[dashed] (2a) -- (2);
\draw[dashed] (3a) -- (3);
\coordinate[fill,circle,inner sep=1pt] (a) at (-1,0);
\coordinate[fill,circle,inner sep=1pt] (b) at (1,0);
\coordinate[fill,circle,inner sep=1pt] (c) at (-1,-1);
\coordinate[fill,circle,inner sep=1pt] (d) at (0,-1);
\coordinate[fill,circle,inner sep=1pt] (e) at (1,-1);
\end{tikzpicture}
\end{minipage}
% (7a)
\begin{minipage}[c]{0.02\textwidth}
$(7a)$
\end{minipage}
\begin{minipage}[c]{0.2\textwidth}
\begin{tikzpicture}[scale=0.9]
\coordinate[fill,circle,inner sep=1pt] (0) at (0,0);
\coordinate[fill,circle,inner sep=1pt] (1) at (1,-1);
\coordinate[fill,circle,inner sep=1pt] (2) at (-2,-1);
\coordinate[fill,circle,inner sep=1pt] (3) at (0,1);
\coordinate[fill,circle,inner sep=1pt] (4) at (1,0);
\coordinate[fill,cross,inner sep=2pt,rotate=45] (1a) at (0.5,-0.5);
\coordinate[fill,cross,inner sep=2pt,rotate=63.43] (2a) at (-1,-0.5);
\coordinate[fill,cross,inner sep=2pt,rotate=0] (3a) at (0,0.5);
\coordinate[fill,cross,inner sep=2pt,rotate=0] (4a) at (0.5,0);
\draw (1) -- (2) -- (3) -- (4) -- (1);
\draw (0) -- (1a);
\draw (0) -- (2a);
\draw (0) -- (3a);
\draw (0) -- (4a);
\draw[dashed] (1a) -- (1);
\draw[dashed] (2a) -- (2);
\draw[dashed] (3a) -- (3);
\draw[dashed] (4a) -- (4);
\coordinate[fill,circle,inner sep=1pt] (a) at (0,-1);
\coordinate[fill,circle,inner sep=1pt] (b) at (-1,-1);
\coordinate[fill,circle,inner sep=1pt] (c) at (-1,0);
\end{tikzpicture}
\end{minipage}
% (6a)
\begin{minipage}[c]{0.02\textwidth}
$(6a)$
\end{minipage}
\begin{minipage}[c]{0.22\textwidth}
\begin{tikzpicture}[scale=0.9]
\coordinate[fill,circle,inner sep=1pt] (0) at (0,0);
\coordinate[fill,circle,inner sep=1pt] (1) at (0,1);
\coordinate[fill,circle,inner sep=1pt] (2) at (1,0);
\coordinate[fill,circle,inner sep=1pt] (3) at (0,-1);
\coordinate[fill,circle,inner sep=1pt] (4) at (-2,-1);
\coordinate[fill,cross,inner sep=2pt,rotate=0] (1a) at (0,0.5);
\coordinate[fill,cross,inner sep=2pt,rotate=0] (2a) at (0.5,0);
\coordinate[fill,cross,inner sep=2pt,rotate=0] (3a) at (0,-0.5);
\coordinate[fill,cross,inner sep=2pt,rotate=63.43] (4a) at (-1,-0.5);
\draw (1) -- (2) -- (3) -- (4) -- (1);
\draw (0) -- (1a);
\draw (0) -- (2a);
\draw (0) -- (3a);
\draw (0) -- (4a);
\draw[dashed] (1a) -- (1);
\draw[dashed] (2a) -- (2);
\draw[dashed] (3a) -- (3);
\draw[dashed] (4a) -- (4);
\coordinate[fill,circle,inner sep=1pt] (a) at (0,-1);
\coordinate[fill,circle,inner sep=1pt] (b) at (-1,-1);
\coordinate[fill,circle,inner sep=1pt] (c) at (-1,0);
\end{tikzpicture}
\end{minipage}\\[6mm]
% (6b)
\begin{minipage}[c]{0.05\textwidth}
$(6b)$
\end{minipage}
\begin{minipage}[c]{0.18\textwidth}
\begin{tikzpicture}[scale=0.6]
\coordinate[fill,circle,inner sep=1pt] (0) at (0,0);
\coordinate[fill,circle,inner sep=1pt] (1) at (0,1);
\coordinate[fill,circle,inner sep=1pt] (2) at (1,0);
\coordinate[fill,circle,inner sep=1pt] (3) at (-3,-2);
\coordinate[fill,cross,inner sep=2pt,rotate=0] (1a) at (0,0.5);
\coordinate[fill,cross,inner sep=2pt,rotate=0] (2a) at (0.5,0);
\coordinate[fill,cross,inner sep=2pt,rotate=33.69] (3a) at (-1.5,-1);
\draw (1) -- (2) -- (3) -- (1);
\draw (0) -- (1a);
\draw (0) -- (2a);
\draw (0) -- (3a);
\draw[dashed] (1a) -- (1);
\draw[dashed] (2a) -- (2);
\draw[dashed] (3a) -- (3);
\coordinate[fill,circle,inner sep=1pt] (a) at (-1,0);
\coordinate[fill,circle,inner sep=1pt] (b) at (-1,-1);
\coordinate[fill,circle,inner sep=1pt] (c) at (-2,-1);
\end{tikzpicture}
\end{minipage}
% (6c)
\begin{minipage}[c]{0.05\textwidth}
$(6c)$
\end{minipage}
\begin{minipage}[c]{0.18\textwidth}
\begin{tikzpicture}[scale=0.9]
\coordinate[fill,circle,inner sep=1pt] (0) at (0,0);
\coordinate[fill,circle,inner sep=1pt] (1) at (-1,0);
\coordinate[fill,circle,inner sep=1pt] (2) at (0,1);
\coordinate[fill,circle,inner sep=1pt] (3) at (1,0);
\coordinate[fill,circle,inner sep=1pt] (4) at (1,-1);
\coordinate[fill,circle,inner sep=1pt] (5) at (-1,-1);
\coordinate[fill,cross,inner sep=2pt,rotate=0] (1a) at (-0.5,0);
\coordinate[fill,cross,inner sep=2pt,rotate=0] (2a) at (0,0.5);
\coordinate[fill,cross,inner sep=2pt,rotate=0] (3a) at (0.5,0);
\coordinate[fill,cross,inner sep=2pt,rotate=45] (4a) at (0.5,-0.5);
\coordinate[fill,cross,inner sep=2pt,rotate=45] (5a) at (-0.5,-0.5);
\draw (1) -- (2) -- (3) -- (4) -- (5) -- (1);
\draw (0) -- (1a);
\draw (0) -- (2a);
\draw (0) -- (3a);
\draw (0) -- (4a);
\draw (0) -- (5a);
\draw[dashed] (1a) -- (1);
\draw[dashed] (2a) -- (2);
\draw[dashed] (3a) -- (3);
\draw[dashed] (4a) -- (4);
\draw[dashed] (5a) -- (5);
\coordinate[fill,circle,inner sep=1pt] (a) at (0,-1);
\end{tikzpicture}
\end{minipage}
% (5a)
\begin{minipage}[c]{0.02\textwidth}
$(5a)$
\end{minipage}
\begin{minipage}[c]{0.21\textwidth}
\begin{tikzpicture}[scale=0.9]
\coordinate[fill,circle,inner sep=1pt] (0) at (0,0);
\coordinate[fill,circle,inner sep=1pt] (1) at (0,1);
\coordinate[fill,circle,inner sep=1pt] (2) at (1,0);
\coordinate[fill,circle,inner sep=1pt] (3) at (-1,-1);
\coordinate[fill,circle,inner sep=1pt] (4) at (-2,-1);
\coordinate[fill,cross,inner sep=2pt,rotate=0] (1a) at (0,0.5);
\coordinate[fill,cross,inner sep=2pt,rotate=0] (2a) at (0.5,0);
\coordinate[fill,cross,inner sep=2pt,rotate=45] (3a) at (-0.5,-0.5);
\coordinate[fill,cross,inner sep=2pt,rotate=26.57] (4a) at (-1,-0.5);
\draw (1) -- (2) -- (3) -- (4) -- (1);
\draw (0) -- (1a);
\draw (0) -- (2a);
\draw (0) -- (3a);
\draw (0) -- (4a);
\draw[dashed] (1a) -- (1);
\draw[dashed] (2a) -- (2);
\draw[dashed] (3a) -- (3);
\draw[dashed] (4a) -- (4);
\coordinate[fill,circle,inner sep=1pt] (a) at (-1,0);
\end{tikzpicture}
\end{minipage}
% (5b)
\begin{minipage}[c]{0.05\textwidth}
$(5b)$
\end{minipage}
\begin{minipage}[c]{0.18\textwidth}
\begin{tikzpicture}[scale=0.9]
\coordinate[fill,circle,inner sep=1pt] (0) at (0,0);
\coordinate[fill,circle,inner sep=1pt] (1) at (0,1);
\coordinate[fill,circle,inner sep=1pt] (2) at (1,0);
\coordinate[fill,circle,inner sep=1pt] (3) at (0,-1);
\coordinate[fill,circle,inner sep=1pt] (4) at (-1,-1);
\coordinate[fill,circle,inner sep=1pt] (5) at (-1,0);
\coordinate[fill,cross,inner sep=2pt,rotate=0] (1a) at (0,0.5);
\coordinate[fill,cross,inner sep=2pt,rotate=0] (2a) at (0.5,0);
\coordinate[fill,cross,inner sep=2pt,rotate=0] (3a) at (0,-0.5);
\coordinate[fill,cross,inner sep=2pt,rotate=45] (4a) at (-0.5,-0.5);
\coordinate[fill,cross,inner sep=2pt,rotate=0] (5a) at (-0.5,0);
\draw (1) -- (2) -- (3) -- (4) -- (5) -- (1);
\draw (0) -- (1a);
\draw (0) -- (2a);
\draw (0) -- (3a);
\draw (0) -- (4a);
\draw (0) -- (5a);
\draw[dashed] (1a) -- (1);
\draw[dashed] (2a) -- (2);
\draw[dashed] (3a) -- (3);
\draw[dashed] (4a) -- (4);
\draw[dashed] (5a) -- (5);
\end{tikzpicture}
\end{minipage}\\[6mm]
% (4a)
\begin{minipage}[c]{0.05\textwidth}
$(4a)$
\end{minipage}
\begin{minipage}[c]{0.18\textwidth}
\begin{tikzpicture}[scale=0.9]
\coordinate[fill,circle,inner sep=1pt] (0) at (0,0);
\coordinate[fill,circle,inner sep=1pt] (1) at (0,1);
\coordinate[fill,circle,inner sep=1pt] (2) at (-1,0);
\coordinate[fill,circle,inner sep=1pt] (3) at (0,-1);
\coordinate[fill,circle,inner sep=1pt] (4) at (1,0);
\coordinate[fill,cross,inner sep=2pt,rotate=0] (1a) at (0,0.5);
\coordinate[fill,cross,inner sep=2pt,rotate=0] (2a) at (-0.5,0);
\coordinate[fill,cross,inner sep=2pt,rotate=0] (3a) at (0,-0.5);
\coordinate[fill,cross,inner sep=2pt,rotate=0] (4a) at (0.5,0);
\draw (1) -- (2) -- (3) -- (4) -- (1);
\draw (0) -- (1a);
\draw (0) -- (2a);
\draw (0) -- (3a);
\draw (0) -- (4a);
\draw[dashed] (1a) -- (1);
\draw[dashed] (2a) -- (2);
\draw[dashed] (3a) -- (3);
\draw[dashed] (4a) -- (4);
\end{tikzpicture}
\end{minipage}
% (4b)
\begin{minipage}[c]{0.05\textwidth}
$(4b)$
\end{minipage}
\begin{minipage}[c]{0.18\textwidth}
\begin{tikzpicture}[scale=0.9]
\coordinate[fill,circle,inner sep=1pt] (0) at (0,0);
\coordinate[fill,circle,inner sep=1pt] (1) at (0,1);
\coordinate[fill,circle,inner sep=1pt] (2) at (1,0);
\coordinate[fill,circle,inner sep=1pt] (3) at (-1,-1);
\coordinate[fill,circle,inner sep=1pt] (4) at (-1,0);
\coordinate[fill,cross,inner sep=2pt,rotate=0] (1a) at (0,0.5);
\coordinate[fill,cross,inner sep=2pt,rotate=0] (2a) at (0.5,0);
\coordinate[fill,cross,inner sep=2pt,rotate=45] (3a) at (-0.5,-0.5);
\coordinate[fill,cross,inner sep=2pt,rotate=0] (4a) at (-0.5,0);
\draw (1) -- (2) -- (3) -- (4) -- (1);
\draw (0) -- (1a);
\draw (0) -- (2a);
\draw (0) -- (3a);
\draw (0) -- (4a);
\draw[dashed] (1a) -- (1);
\draw[dashed] (2a) -- (2);
\draw[dashed] (3a) -- (3);
\draw[dashed] (4a) -- (4);
\end{tikzpicture}
\end{minipage}
% (4c)
\begin{minipage}[c]{0.02\textwidth}
$(4c)$
\end{minipage}
\begin{minipage}[c]{0.21\textwidth}
\begin{tikzpicture}[scale=0.9]
\coordinate[fill,circle,inner sep=1pt] (0) at (0,0);
\coordinate[fill,circle,inner sep=1pt] (1) at (0,1);
\coordinate[fill,circle,inner sep=1pt] (2) at (1,0);
\coordinate[fill,circle,inner sep=1pt] (3) at (-2,-1);
\coordinate[fill,cross,inner sep=2pt,rotate=0] (1a) at (0,0.5);
\coordinate[fill,cross,inner sep=2pt,rotate=0] (2a) at (0.5,0);
\coordinate[fill,cross,inner sep=2pt,rotate=26.57] (3a) at (-1,-0.5);
\draw (1) -- (2) -- (3) -- (1);
\draw (0) -- (1a);
\draw (0) -- (2a);
\draw (0) -- (3a);
\draw[dashed] (1a) -- (1);
\draw[dashed] (2a) -- (2);
\draw[dashed] (3a) -- (3);
\coordinate[fill,circle,inner sep=1pt] (a) at (-1,0);
\end{tikzpicture}
\end{minipage}
% (3a)
\begin{minipage}[c]{0.05\textwidth}
$(3a)$
\end{minipage}
\begin{minipage}[c]{0.18\textwidth}
\begin{tikzpicture}[scale=0.9]
\coordinate[fill,circle,inner sep=1pt] (0) at (0,0);
\coordinate[fill,circle,inner sep=1pt] (1) at (0,1);
\coordinate[fill,circle,inner sep=1pt] (2) at (1,0);
\coordinate[fill,circle,inner sep=1pt] (3) at (-1,-1);
\coordinate[fill,cross,inner sep=2pt,rotate=0] (1a) at (0,0.5);
\coordinate[fill,cross,inner sep=2pt,rotate=0] (2a) at (0.5,0);
\coordinate[fill,cross,inner sep=2pt,rotate=45] (3a) at (-0.5,-0.5);
\draw (1) -- (2) -- (3) -- (1);
\draw (0) -- (1a);
\draw (0) -- (2a);
\draw (0) -- (3a);
\draw[dashed] (1a) -- (1);
\draw[dashed] (2a) -- (2);
\draw[dashed] (3a) -- (3);
\end{tikzpicture}
\end{minipage}
\caption{Intersection complexes $(\check{B},\check{\mathscr{P}},\check{\varphi})$ of smoothings $X$ of toric Gorenstein del Pezzo surfaces $X^0$. The number in the labelling is the degree of $X$ and $X^0$. In the first five cases $X^0$ is smooth.}
\label{fig:list}
\end{figure}

\begin{proof}
For any $\mathbb{Q}$-Gorenstein deformation $\mathfrak{X}\rightarrow\mathbb{A}^1$, the relative canonical class $K_{\mathfrak{X}/\mathbb{A}^1}$ is $\mathbb{Q}$-Cartier. By definition, the degree of $X^0$ is the self-intersection of its (anti)canonical class. Hence, the degree of $X$ equals the degree of any of its toric models $X^0$. We need to show that the degrees of the given toric models are the ones shown in Figure \ref{fig:list}. The Fano polytope $Q$ of $X^0$ is exactly the Newton polytope of its anticanonical class. By the duality between subdivisions of Newton polytopes and tropical curves, the self intersection can easily be computed by intersecting two tropical curves dual to the Fano polytope $Q$. For (3a), i.e., $X^0=\mathbb{P}^2/\mathbb{Z}_3$ as in Example \ref{expl:cubic}, the intersection of tropical curves is the following:
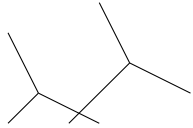
\begin{center}
\begin{tikzpicture}[scale=0.4]
\draw (0,0) -- (-1,-1);
\draw (0,0) -- (-1,2);
\draw (0,0) -- (2,-1);
\draw (3,1) -- (1,-1);
\draw (3,1) -- (2,3);
\draw (3,1) -- (5,0);
\end{tikzpicture}
\end{center}
The determinant of primitive tangent vectors at the intersection point is $|\textup{det}\left(\begin{smallmatrix}\textup{ }\ 1&2\\-1&1\end{smallmatrix}\right)|=3$. Indeed, this is the degree of $\mathbb{P}^2/\mathbb{Z}_3$. Similarly one computes the degrees of the other cases in Figure \ref{fig:list}. Alternatively, one can use the fact that the degree of a del Pezzo surface equals $|\check{B}\cap M|+1$ (see \cite{CPS}, {\S}6).

There are two smooth del Pezzo surfaces of degree $8$, the blow up of $\mathbb{P}^2$ at a point and $\mathbb{P}^1\times\mathbb{P}^1$. The del Pezzo surface $X^0$ in case (8'a) is a toric model for $\mathbb{P}^1\times\mathbb{P}^1$ . The other cases are determined, up to smooth deformations, by the degree, since del Pezzos of degree $\neq 8$ have a connected moduli space.
\end{proof}

\subsection{Fan picture and refinement}						%%

Let $Q$ be a Fano polytope and let $\mathfrak{X}_Q\rightarrow\mathbb{A}^2$ be the family from Construction \ref{con:family1}. Let $(\check{B},\check{\mathscr{P}},\check{\varphi})$ be the intersection complex of the toric degeneration $\mathfrak{X}:=\mathfrak{X}^{s\neq 0}\rightarrow\mathbb{A}^1$, i.e., one of the polarized polyhedral affine manifolds in Figure \ref{fig:list}. Performing the discrete Legendre transform (\cite{DataI}, \S1.4) we obtain another polarized polyhedral affine manifold that is the dual intersection complex (\cite{DataI}, \S4.1) of $\mathfrak{X}\rightarrow\mathbb{A}^1$.

\begin{defi}
\label{defi:sigma0}
Let $\sigma_0$ be the unique bounded maximal cell of the dual intersection complex of $\mathfrak{X}\rightarrow\mathbb{A}^1$.
\end{defi}

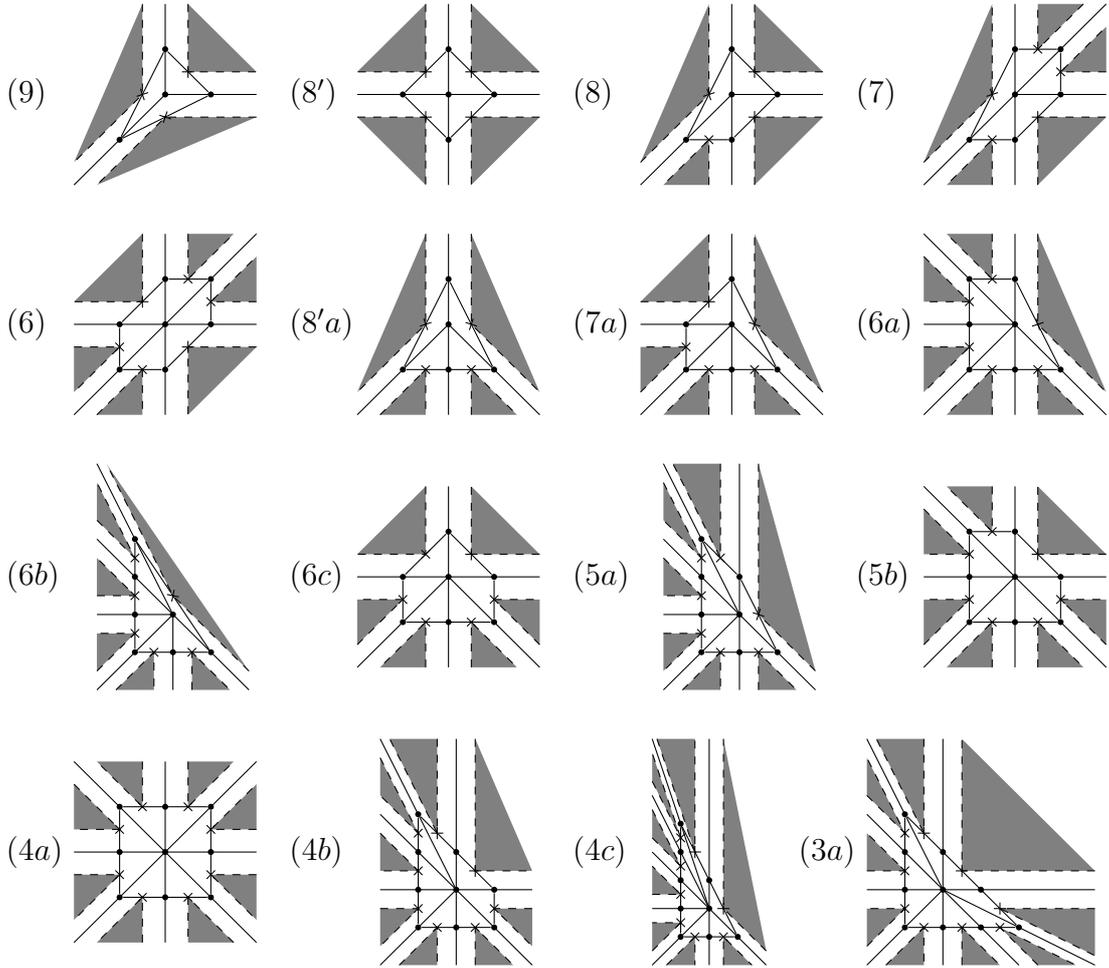
\begin{figure}[h!]
\centering
% (9)
\begin{minipage}[c]{0.05\textwidth}
$(9)$
\end{minipage}
\begin{minipage}[c]{0.18\textwidth}
\begin{tikzpicture}[scale=0.6]
\coordinate[fill,circle,inner sep=0.8pt] (0) at (0,0);
\coordinate[fill,circle,inner sep=0.8pt] (1) at (1,0);
\coordinate[fill,circle,inner sep=0.8pt] (2) at (0,1);
\coordinate[fill,circle,inner sep=0.8pt] (3) at (-1,-1);
\coordinate (1a) at (2,0);
\coordinate (2a) at (0,2);
\coordinate (3a) at (-2,-2);
\coordinate (12) at (0.5,0.5);
\coordinate (12a) at (0.5,2);
\coordinate (12b) at (2,0.5);
\coordinate (23) at (-0.5,0);
\coordinate (23a) at (-0.5,2);
\coordinate (23b) at (-2,-1.5);
\coordinate (31) at (0,-0.5);
\coordinate (31a) at (2,-0.5);
\coordinate (31b) at (-1.5,-2);
\draw (1) -- (2) -- (3) -- (1);
\draw (0,0) -- (1a);
\draw (0,0) -- (2a);
\draw (0,0) -- (3a);
\fill[color=gray] (12) -- (12a) -- (12b) -- (12);
\fill[color=gray] (23) -- (23a) -- (23b) -- (23);
\fill[color=gray] (31) -- (31a) -- (31b) -- (31);
\draw[dashed] (12a) -- (12) -- (12b);
\draw[dashed] (23a) -- (23) -- (23b);
\draw[dashed] (31a) -- (31) -- (31b);
\draw (12) node[fill,cross,inner sep=1pt,rotate=45]{};
\draw (23) node[fill,cross,inner sep=1pt,rotate=63.43]{};
\draw (31) node[fill,cross,inner sep=1pt,rotate=26.57]{};
\end{tikzpicture}
\end{minipage}
% (8')
\begin{minipage}[c]{0.05\textwidth}
$(8')$
\end{minipage}
\begin{minipage}[c]{0.18\textwidth}
\begin{tikzpicture}[scale=0.6]
\coordinate[fill,circle,inner sep=0.8pt] (0) at (0,0);
\coordinate[fill,circle,inner sep=0.8pt] (1) at (1,0);
\coordinate[fill,circle,inner sep=0.8pt] (2) at (0,1);
\coordinate[fill,circle,inner sep=0.8pt] (3) at (-1,0);
\coordinate[fill,circle,inner sep=0.8pt] (4) at (0,-1);
\coordinate (1a) at (2,0);
\coordinate (2a) at (0,2);
\coordinate (3a) at (-2,0);
\coordinate (4a) at (0,-2);
\coordinate (12) at (0.5,0.5);
\coordinate (12a) at (0.5,2);
\coordinate (12b) at (2,0.5);
\coordinate (23) at (-0.5,0.5);
\coordinate (23a) at (-0.5,2);
\coordinate (23b) at (-2,0.5);
\coordinate (34) at (-0.5,-0.5);
\coordinate (34a) at (-2,-0.5);
\coordinate (34b) at (-0.5,-2);
\coordinate (41) at (0.5,-0.5);
\coordinate (41a) at (2,-0.5);
\coordinate (41b) at (0.5,-2);
\draw (1) -- (2) -- (3) -- (4) -- (1);
\draw (0,0) -- (1a);
\draw (0,0) -- (2a);
\draw (0,0) -- (3a);
\draw (0,0) -- (4a);
\fill[color=gray] (12) -- (12a) -- (12b) -- (12);
\fill[color=gray] (23) -- (23a) -- (23b) -- (23);
\fill[color=gray] (34) -- (34a) -- (34b) -- (34);
\fill[color=gray] (41) -- (41a) -- (41b) -- (41);
\draw[dashed] (12a) -- (12) -- (12b);
\draw[dashed] (23a) -- (23) -- (23b);
\draw[dashed] (34a) -- (34) -- (34b);
\draw[dashed] (41a) -- (41) -- (41b);
\draw (12) node[fill,cross,inner sep=1pt,rotate=45]{};
\draw (23) node[fill,cross,inner sep=1pt,rotate=45]{};
\draw (34) node[fill,cross,inner sep=1pt,rotate=45]{};
\draw (41) node[fill,cross,inner sep=1pt,rotate=45]{};
\end{tikzpicture}
\end{minipage}
% (8)
\begin{minipage}[c]{0.05\textwidth}
$(8)$
\end{minipage}
\begin{minipage}[c]{0.18\textwidth}
\begin{tikzpicture}[scale=0.6]
\coordinate[fill,circle,inner sep=0.8pt] (0) at (0,0);
\coordinate[fill,circle,inner sep=0.8pt] (1) at (-1,-1);
\coordinate[fill,circle,inner sep=0.8pt] (2) at (0,-1);
\coordinate[fill,circle,inner sep=0.8pt] (3) at (1,0);
\coordinate[fill,circle,inner sep=0.8pt] (4) at (0,1);
\coordinate (1a) at (-2,-2);
\coordinate (2a) at (0,-2);
\coordinate (3a) at (2,0);
\coordinate (4a) at (0,2);
\coordinate (12) at (-0.5,-1);
\coordinate (12a) at (-0.5,-2);
\coordinate (12b) at (-1.5,-2);
\coordinate (23) at (0.5,-0.5);
\coordinate (23a) at (0.5,-2);
\coordinate (23b) at (2,-0.5);
\coordinate (34) at (0.5,0.5);
\coordinate (34a) at (2,0.5);
\coordinate (34b) at (0.5,2);
\coordinate (41) at (-0.5,-0);
\coordinate (41a) at (-2,-1.5);
\coordinate (41b) at (-0.5,2);
\draw (1) -- (2) -- (3) -- (4) -- (1);
\draw (0,0) -- (1a);
\draw (0,0) -- (2a);
\draw (0,0) -- (3a);
\draw (0,0) -- (4a);
\fill[color=gray] (12) -- (12a) -- (12b) -- (12);
\fill[color=gray] (23) -- (23a) -- (23b) -- (23);
\fill[color=gray] (34) -- (34a) -- (34b) -- (34);
\fill[color=gray] (41) -- (41a) -- (41b) -- (41);
\draw[dashed] (12a) -- (12) -- (12b);
\draw[dashed] (23a) -- (23) -- (23b);
\draw[dashed] (34a) -- (34) -- (34b);
\draw[dashed] (41a) -- (41) -- (41b);
\draw (12) node[fill,cross,inner sep=1pt,rotate=0]{};
\draw (23) node[fill,cross,inner sep=1pt,rotate=45]{};
\draw (34) node[fill,cross,inner sep=1pt,rotate=45]{};
\draw (41) node[fill,cross,inner sep=1pt,rotate=63.43]{};
\end{tikzpicture}
\end{minipage}
% (7)
\begin{minipage}[c]{0.05\textwidth}
$(7)$
\end{minipage}
\begin{minipage}[c]{0.18\textwidth}
\begin{tikzpicture}[scale=0.6]
\coordinate[fill,circle,inner sep=0.8pt] (0) at (0,0);
\coordinate[fill,circle,inner sep=0.8pt] (1) at (-1,-1);
\coordinate[fill,circle,inner sep=0.8pt] (2) at (0,-1);
\coordinate[fill,circle,inner sep=0.8pt] (3) at (1,0);
\coordinate[fill,circle,inner sep=0.8pt] (4) at (1,1);
\coordinate[fill,circle,inner sep=0.8pt] (5) at (0,1);
\coordinate (1a) at (-2,-2);
\coordinate (2a) at (0,-2);
\coordinate (3a) at (2,0);
\coordinate (4a) at (2,2);
\coordinate (5a) at (0,2);
\coordinate (12) at (-0.5,-1);
\coordinate (12a) at (-0.5,-2);
\coordinate (12b) at (-1.5,-2);
\coordinate (23) at (0.5,-0.5);
\coordinate (23a) at (0.5,-2);
\coordinate (23b) at (2,-0.5);
\coordinate (34) at (1,0.5);
\coordinate (34a) at (2,0.5);
\coordinate (34b) at (2,1.5);
\coordinate (45) at (0.5,1);
\coordinate (45a) at (1.5,2);
\coordinate (45b) at (0.5,2);
\coordinate (51) at (-0.5,0);
\coordinate (51a) at (-2,-1.5);
\coordinate (51b) at (-0.5,2);
\draw (1) -- (2) -- (3) -- (4) -- (5) -- (1);
\draw (0,0) -- (1a);
\draw (0,0) -- (2a);
\draw (0,0) -- (3a);
\draw (0,0) -- (4a);
\draw (0,0) -- (5a);
\fill[color=gray] (12) -- (12a) -- (12b) -- (12);
\fill[color=gray] (23) -- (23a) -- (23b) -- (23);
\fill[color=gray] (34) -- (34a) -- (34b) -- (34);
\fill[color=gray] (45) -- (45a) -- (45b) -- (45);
\fill[color=gray] (51) -- (51a) -- (51b) -- (51);
\draw[dashed] (12a) -- (12) -- (12b);
\draw[dashed] (23a) -- (23) -- (23b);
\draw[dashed] (34a) -- (34) -- (34b);
\draw[dashed] (45a) -- (45) -- (45b);
\draw[dashed] (51a) -- (51) -- (51b);
\draw (12) node[fill,cross,inner sep=1pt,rotate=0]{};
\draw (23) node[fill,cross,inner sep=1pt,rotate=45]{};
\draw (34) node[fill,cross,inner sep=1pt,rotate=0]{};
\draw (45) node[fill,cross,inner sep=1pt,rotate=0]{};
\draw (51) node[fill,cross,inner sep=1pt,rotate=63.43]{};
\end{tikzpicture}
\end{minipage}\\[6mm]
% (6)
\begin{minipage}[c]{0.05\textwidth}
$(6)$
\end{minipage}
\begin{minipage}[c]{0.18\textwidth}
\begin{tikzpicture}[scale=0.6]
\coordinate[fill,circle,inner sep=0.8pt] (0) at (0,0);
\coordinate[fill,circle,inner sep=0.8pt] (1) at (1,1);
\coordinate[fill,circle,inner sep=0.8pt] (2) at (0,1);
\coordinate[fill,circle,inner sep=0.8pt] (3) at (-1,0);
\coordinate[fill,circle,inner sep=0.8pt] (4) at (-1,-1);
\coordinate[fill,circle,inner sep=0.8pt] (5) at (0,-1);
\coordinate[fill,circle,inner sep=0.8pt] (6) at (1,0);
\coordinate (1a) at (2,2);
\coordinate (2a) at (0,2);
\coordinate (3a) at (-2,0);
\coordinate (4a) at (-2,-2);
\coordinate (5a) at (0,-2);
\coordinate (6a) at (2,0);
\coordinate (12) at (0.5,1);
\coordinate (12a) at (0.5,2);
\coordinate (12b) at (1.5,2);
\coordinate (23) at (-0.5,0.5);
\coordinate (23a) at (-0.5,2);
\coordinate (23b) at (-2,0.5);
\coordinate (34) at (-1,-0.5);
\coordinate (34a) at (-2,-0.5);
\coordinate (34b) at (-2,-1.5);
\coordinate (45) at (-0.5,-1);
\coordinate (45a) at (-1.5,-2);
\coordinate (45b) at (-0.5,-2);
\coordinate (56) at (0.5,-0.5);
\coordinate (56a) at (2,-0.5);
\coordinate (56b) at (0.5,-2);
\coordinate (61) at (1,0.5);
\coordinate (61a) at (2,0.5);
\coordinate (61b) at (2,1.5);
\draw (1) -- (2) -- (3) -- (4) -- (5) -- (6) -- (1);
\draw (0,0) -- (1a);
\draw (0,0) -- (2a);
\draw (0,0) -- (3a);
\draw (0,0) -- (4a);
\draw (0,0) -- (5a);
\draw (0,0) -- (6a);
\fill[color=gray] (12) -- (12a) -- (12b) -- (12);
\fill[color=gray] (23) -- (23a) -- (23b) -- (23);
\fill[color=gray] (34) -- (34a) -- (34b) -- (34);
\fill[color=gray] (45) -- (45a) -- (45b) -- (45);
\fill[color=gray] (56) -- (56a) -- (56b) -- (56);
\fill[color=gray] (61) -- (61a) -- (61b) -- (61);
\draw[dashed] (12a) -- (12) -- (12b);
\draw[dashed] (23a) -- (23) -- (23b);
\draw[dashed] (34a) -- (34) -- (34b);
\draw[dashed] (45a) -- (45) -- (45b);
\draw[dashed] (56a) -- (56) -- (56b);
\draw[dashed] (61a) -- (61) -- (61b);
\draw (12) node[fill,cross,inner sep=1pt,rotate=0]{};
\draw (23) node[fill,cross,inner sep=1pt,rotate=45]{};
\draw (34) node[fill,cross,inner sep=1pt,rotate=0]{};
\draw (45) node[fill,cross,inner sep=1pt,rotate=0]{};
\draw (56) node[fill,cross,inner sep=1pt,rotate=45]{};
\draw (61) node[fill,cross,inner sep=1pt,rotate=0]{};
\end{tikzpicture}
\end{minipage}
% (8a)
\begin{minipage}[c]{0.05\textwidth}
$(8'a)$
\end{minipage}
\begin{minipage}[c]{0.18\textwidth}
\begin{tikzpicture}[scale=0.6]
\coordinate[fill,circle,inner sep=0.8pt] (0) at (0,0);
\coordinate[fill,circle,inner sep=0.8pt] (1) at (-1,-1);
\coordinate[fill,circle,inner sep=0.8pt] (2) at (1,-1);
\coordinate[fill,circle,inner sep=0.8pt] (3) at (0,1);
\coordinate (1a) at (-2,-2);
\coordinate (2a) at (2,-2);
\coordinate (3a) at (0,2);
\coordinate (12) at (-0.5,-1);
\coordinate (12a) at (-1.5,-2);
\coordinate (12b) at (-0.5,-2);
\coordinate (12') at (0.5,-1);
\coordinate (12'a) at (0.5,-2);
\coordinate (12'b) at (1.5,-2);
\coordinate (23) at (0.5,0);
\coordinate (23a) at (0.5,2);
\coordinate (23b) at (2,-1.5);
\coordinate (31) at (-0.5,0);
\coordinate (31a) at (-0.5,2);
\coordinate (31b) at (-2,-1.5);
\draw (1) -- (2) -- (3) -- (1);
\draw (0,0) -- (1a);
\draw (0,0) -- (2a);
\draw (0,0) -- (3a);
\draw (0,0) -- (0,-2);
\fill[color=gray] (12) -- (12a) -- (12b) -- (12);
\fill[color=gray] (12') -- (12'a) -- (12'b) -- (12');
\fill[color=gray] (23) -- (23a) -- (23b) -- (23);
\fill[color=gray] (31) -- (31a) -- (31b) -- (31);
\draw[dashed] (12a) -- (12) -- (12b);
\draw[dashed] (12'a) -- (12') -- (12'b);
\draw[dashed] (23a) -- (23) -- (23b);
\draw[dashed] (31a) -- (31) -- (31b);
\draw (12) node[fill,cross,inner sep=1pt,rotate=0]{};
\draw (12') node[fill,cross,inner sep=1pt,rotate=0]{};
\draw (23) node[fill,cross,inner sep=1pt,rotate=26.57]{};
\draw (31) node[fill,cross,inner sep=1pt,rotate=63.43]{};
\coordinate[fill,circle,inner sep=0.8pt] (b) at (0,-1);
\end{tikzpicture}
\end{minipage}
% (7a)
\begin{minipage}[c]{0.05\textwidth}
$(7a)$
\end{minipage}
\begin{minipage}[c]{0.18\textwidth}
\begin{tikzpicture}[scale=0.6]
\coordinate[fill,circle,inner sep=0.8pt] (0) at (0,0);
\coordinate[fill,circle,inner sep=0.8pt] (1) at (1,-1);
\coordinate[fill,circle,inner sep=0.8pt] (2) at (-1,-1);
\coordinate[fill,circle,inner sep=0.8pt] (3) at (-1,0);
\coordinate[fill,circle,inner sep=0.8pt] (4) at (0,1);
\coordinate (1a) at (2,-2);
\coordinate (2a) at (-2,-2);
\coordinate (3a) at (-2,0);
\coordinate (4a) at (0,2);
\coordinate (12) at (0.5,-1);
\coordinate (12a) at (1.5,-2);
\coordinate (12b) at (0.5,-2);
\coordinate (12') at (-0.5,-1);
\coordinate (12'a) at (-0.5,-2);
\coordinate (12'b) at (-1.5,-2);
\coordinate (23) at (-1,-0.5);
\coordinate (23a) at (-2,-1.5);
\coordinate (23b) at (-2,-0.5);
\coordinate (34) at (-0.5,0.5);
\coordinate (34a) at (-2,0.5);
\coordinate (34b) at (-0.5,2);
\coordinate (41) at (0.5,0);
\coordinate (41a) at (2,-1.5);
\coordinate (41b) at (0.5,2);
\draw (1) -- (2) -- (3) -- (4) -- (1);
\draw (0,0) -- (1a);
\draw (0,0) -- (2a);
\draw (0,0) -- (3a);
\draw (0,0) -- (4a);
\draw (0,0) -- (0,-2);
\fill[color=gray] (12) -- (12a) -- (12b) -- (12);
\fill[color=gray] (12') -- (12'a) -- (12'b) -- (12');
\fill[color=gray] (23) -- (23a) -- (23b) -- (23);
\fill[color=gray] (34) -- (34a) -- (34b) -- (34);
\fill[color=gray] (41) -- (41a) -- (41b) -- (41);
\draw[dashed] (12a) -- (12) -- (12b);
\draw[dashed] (12'a) -- (12') -- (12'b);
\draw[dashed] (23a) -- (23) -- (23b);
\draw[dashed] (34a) -- (34) -- (34b);
\draw[dashed] (41a) -- (41) -- (41b);
\draw (12) node[fill,cross,inner sep=1pt,rotate=0]{};
\draw (12') node[fill,cross,inner sep=1pt,rotate=0]{};
\draw (23) node[fill,cross,inner sep=1pt,rotate=0]{};
\draw (34) node[fill,cross,inner sep=1pt,rotate=45]{};
\draw (41) node[fill,cross,inner sep=1pt,rotate=26.57]{};
\coordinate[fill,circle,inner sep=0.8pt] (b) at (0,-1);
\end{tikzpicture}
\end{minipage}
% (6a)
\begin{minipage}[c]{0.05\textwidth}
$(6a)$
\end{minipage}
\begin{minipage}[c]{0.18\textwidth}
\begin{tikzpicture}[scale=0.6]
\coordinate[fill,circle,inner sep=0.8pt] (0) at (0,0);
\coordinate[fill,circle,inner sep=0.8pt] (1) at (1,-1);
\coordinate[fill,circle,inner sep=0.8pt] (2) at (-1,-1);
\coordinate[fill,circle,inner sep=0.8pt] (3) at (-1,1);
\coordinate[fill,circle,inner sep=0.8pt] (4) at (0,1);
\coordinate (1a) at (2,-2);
\coordinate (2a) at (-2,-2);
\coordinate (3a) at (-2,2);
\coordinate (4a) at (0,2);
\coordinate (12a) at (1.5,-2);
\coordinate (12b) at (0.5,-2);
\coordinate (12') at (-0.5,-1);
\coordinate (12'a) at (-0.5,-2);
\coordinate (12'b) at (-1.5,-2);
\coordinate (23) at (-1,0.5);
\coordinate (23a) at (-2,1.5);
\coordinate (23b) at (-2,0.5);
\coordinate (23') at (-1,-0.5);
\coordinate (23'a) at (-2,-0.5);
\coordinate (23'b) at (-2,-1.5);
\coordinate (34) at (-0.5,1);
\coordinate (34a) at (-0.5,2);
\coordinate (34b) at (-1.5,2);
\coordinate (41) at (0.5,0);
\coordinate (41a) at (0.5,2);
\coordinate (41b) at (2,-1.5);
\draw (1) -- (2) -- (3) -- (4) -- (1);
\draw (0,0) -- (1a);
\draw (0,0) -- (2a);
\draw (0,0) -- (3a);
\draw (0,0) -- (4a);
\draw (0,0) -- (0,-2);
\draw (0,0) -- (-2,0);
\fill[color=gray] (12) -- (12a) -- (12b) -- (12);
\fill[color=gray] (12') -- (12'a) -- (12'b) -- (12');
\fill[color=gray] (23) -- (23a) -- (23b) -- (23);
\fill[color=gray] (23') -- (23'a) -- (23'b) -- (23');
\fill[color=gray] (34) -- (34a) -- (34b) -- (34);
\fill[color=gray] (41) -- (41a) -- (41b) -- (41);
\draw[dashed] (12a) -- (12) -- (12b);
\draw[dashed] (12'a) -- (12') -- (12'b);
\draw[dashed] (23a) -- (23) -- (23b);
\draw[dashed] (23'a) -- (23') -- (23'b);
\draw[dashed] (34a) -- (34) -- (34b);
\draw[dashed] (41a) -- (41) -- (41b);
\draw (12) node[fill,cross,inner sep=1pt,rotate=0]{};
\draw (12') node[fill,cross,inner sep=1pt,rotate=0]{};
\draw (23) node[fill,cross,inner sep=1pt,rotate=0]{};
\draw (23') node[fill,cross,inner sep=1pt,rotate=0]{};
\draw (34) node[fill,cross,inner sep=1pt,rotate=0]{};
\draw (41) node[fill,cross,inner sep=1pt,rotate=63.43]{};
\coordinate[fill,circle,inner sep=0.8pt] (b) at (-1,0);
\coordinate[fill,circle,inner sep=0.8pt] (c) at (0,-1);
\end{tikzpicture}
\end{minipage}\\[6mm]
% (6b)
\begin{minipage}[c]{0.07\textwidth}
$(6b)$
\end{minipage}
\begin{minipage}[c]{0.16\textwidth}
\begin{tikzpicture}[scale=0.5]
\coordinate[fill,circle,inner sep=0.8pt] (0) at (0,0);
\coordinate[fill,circle,inner sep=0.8pt] (1) at (-1,-1);
\coordinate[fill,circle,inner sep=0.8pt] (2) at (1,-1);
\coordinate[fill,circle,inner sep=0.8pt] (3) at (-1,2);
\coordinate (1a) at (-2,-2);
\coordinate (2a) at (2,-2);
\coordinate (3a) at (-2,4);
\coordinate (12) at (-0.5,-1);
\coordinate (12a) at (-1.5,-2);
\coordinate (12b) at (-0.5,-2);
\coordinate (12') at (0.5,-1);
\coordinate (12'a) at (0.5,-2);
\coordinate (12'b) at (1.5,-2);
\coordinate (23) at (0,0.5);
\coordinate (23a) at (2,-1.5);
\coordinate (23b) at (-1.75,4);
\coordinate (31) at (-1,-0.5);
\coordinate (31a) at (-2,-1.5);
\coordinate (31b) at (-2,-0.5);
\coordinate (31') at (-1,0.5);
\coordinate (31'a) at (-2,0.5);
\coordinate (31'b) at (-2,1.5);
\coordinate (31'') at (-1,1.5);
\coordinate (31''a) at (-2,2.5);
\coordinate (31''b) at (-2,3.5);
\draw (1) -- (2) -- (3) -- (1);
\draw (0,0) -- (1a);
\draw (0,0) -- (2a);
\draw (0,0) -- (3a);
\draw (0,0) -- (0,-2);
\draw (0,0) -- (-2,0);
\draw (0,0) -- (-2,2);
\fill[color=gray] (12) -- (12a) -- (12b) -- (12);
\fill[color=gray] (12') -- (12'a) -- (12'b) -- (12');
\fill[color=gray] (23) -- (23a) -- (23b) -- (23);
\fill[color=gray] (31) -- (31a) -- (31b) -- (31);
\fill[color=gray] (31') -- (31'a) -- (31'b) -- (31');
\fill[color=gray] (31'') -- (31''a) -- (31''b) -- (31'');
\draw[dashed] (12a) -- (12) -- (12b);
\draw[dashed] (12'a) -- (12') -- (12'b);
\draw[dashed] (23a) -- (23) -- (23b);
\draw[dashed] (31a) -- (31) -- (31b);
\draw[dashed] (31'a) -- (31') -- (31'b);
\draw[dashed] (31''a) -- (31'') -- (31''b);
\draw (12) node[fill,cross,inner sep=1pt,rotate=0]{};
\draw (12') node[fill,cross,inner sep=1pt,rotate=0]{};
\draw (23) node[fill,cross,inner sep=1pt,rotate=33.69]{};
\draw (31) node[fill,cross,inner sep=1pt,rotate=0]{};
\draw (31') node[fill,cross,inner sep=1pt,rotate=0]{};
\draw (31'') node[fill,cross,inner sep=1pt,rotate=0]{};
\coordinate[fill,circle,inner sep=0.8pt] (b) at (-1,0);
\coordinate[fill,circle,inner sep=0.8pt] (c) at (-1,1);
\coordinate[fill,circle,inner sep=0.8pt] (d) at (0,-1);
\end{tikzpicture}
\end{minipage}
% (6c)
\begin{minipage}[c]{0.05\textwidth}
$(6c)$
\end{minipage}
\begin{minipage}[c]{0.18\textwidth}
\begin{tikzpicture}[scale=0.6]
\coordinate[fill,circle,inner sep=0.8pt] (0) at (0,0);
\coordinate[fill,circle,inner sep=0.8pt] (1) at (-1,-1);
\coordinate[fill,circle,inner sep=0.8pt] (2) at (1,-1);
\coordinate[fill,circle,inner sep=0.8pt] (3) at (1,0);
\coordinate[fill,circle,inner sep=0.8pt] (4) at (0,1);
\coordinate[fill,circle,inner sep=0.8pt] (5) at (-1,0);
\coordinate (1a) at (-2,-2);
\coordinate (2a) at (2,-2);
\coordinate (3a) at (2,0);
\coordinate (4a) at (0,2);
\coordinate (5a) at (-2,0);
\coordinate (12) at (-0.5,-1);
\coordinate (12a) at (-0.5,-2);
\coordinate (12b) at (-1.5,-2);
\coordinate (12') at (0.5,-1);
\coordinate (12'a) at (0.5,-2);
\coordinate (12'b) at (1.5,-2);
\coordinate (23) at (1,-0.5);
\coordinate (23a) at (2,-1.5);
\coordinate (23b) at (2,-0.5);
\coordinate (34) at (0.5,0.5);
\coordinate (34a) at (2,0.5);
\coordinate (34b) at (0.5,2);
\coordinate (45) at (-0.5,0.5);
\coordinate (45a) at (-0.5,2);
\coordinate (45b) at (-2,0.5);
\coordinate (51) at (-1,-0.5);
\coordinate (51a) at (-2,-1.5);
\coordinate (51b) at (-2,-0.5);
\draw (1) -- (2) -- (3) -- (4) -- (5) -- (1);
\draw (0,0) -- (1a);
\draw (0,0) -- (2a);
\draw (0,0) -- (3a);
\draw (0,0) -- (4a);
\draw (0,0) -- (5a);
\draw (0,0) -- (0,-2);
\fill[color=gray] (12) -- (12a) -- (12b) -- (12);
\fill[color=gray] (12') -- (12'a) -- (12'b) -- (12');
\fill[color=gray] (23) -- (23a) -- (23b) -- (23);
\fill[color=gray] (34) -- (34a) -- (34b) -- (34);
\fill[color=gray] (45) -- (45a) -- (45b) -- (45);
\fill[color=gray] (51) -- (51a) -- (51b) -- (51);
\draw[dashed] (12a) -- (12) -- (12b);
\draw[dashed] (12'a) -- (12') -- (12'b);
\draw[dashed] (23a) -- (23) -- (23b);
\draw[dashed] (34a) -- (34) -- (34b);
\draw[dashed] (45a) -- (45) -- (45b);
\draw[dashed] (51a) -- (51) -- (51b);
\draw (12) node[fill,cross,inner sep=1pt,rotate=0]{};
\draw (12') node[fill,cross,inner sep=1pt,rotate=0]{};
\draw (23) node[fill,cross,inner sep=1pt,rotate=0]{};
\draw (34) node[fill,cross,inner sep=1pt,rotate=45]{};
\draw (45) node[fill,cross,inner sep=1pt,rotate=0]{};
\draw (51) node[fill,cross,inner sep=1pt,rotate=0]{};
\coordinate[fill,circle,inner sep=0.8pt] (b) at (0,-1);
\end{tikzpicture}
\end{minipage}
% (5a)
\begin{minipage}[c]{0.07\textwidth}
$(5a)$
\end{minipage}
\begin{minipage}[c]{0.16\textwidth}
\begin{tikzpicture}[scale=0.5]
\coordinate[fill,circle,inner sep=0.8pt] (0) at (0,0);
\coordinate[fill,circle,inner sep=0.8pt] (1) at (1,-1);
\coordinate[fill,circle,inner sep=0.8pt] (2) at (-1,-1);
\coordinate[fill,circle,inner sep=0.8pt] (3) at (-1,2);
\coordinate[fill,circle,inner sep=0.8pt] (4) at (0,1);
\coordinate (1a) at (2,-2);
\coordinate (2a) at (-2,-2);
\coordinate (3a) at (-2,4);
\coordinate (4a) at (0,4);
\coordinate (12) at (-0.5,-1);
\coordinate (12a) at (-1.5,-2);
\coordinate (12b) at (-0.5,-2);
\coordinate (12') at (0.5,-1);
\coordinate (12'a) at (0.5,-2);
\coordinate (12'b) at (1.5,-2);
\coordinate (23) at (-1,1.5);
\coordinate (23a) at (-2,3.5);
\coordinate (23b) at (-2,2.5);
\coordinate (23') at (-1,0.5);
\coordinate (23'a) at (-2,1.5);
\coordinate (23'b) at (-2,0.5);
\coordinate (23'') at (-1,-0.5);
\coordinate (23''a) at (-2,-0.5);
\coordinate (23''b) at (-2,-1.5);
\coordinate (34) at (-0.5,1.5);
\coordinate (34a) at (-0.5,4);
\coordinate (34b) at (-1.75,4);
\coordinate (41) at (0.5,0);
\coordinate (41a) at (0.5,4);
\coordinate (41b) at (2,-1.5);
\draw (1) -- (2) -- (3) -- (4) -- (1);
\draw (0,0) -- (1a);
\draw (0,0) -- (2a);
\draw (0,0) -- (3a);
\draw (0,0) -- (4a);
\draw (0,0) -- (0,-2);
\draw (0,0) -- (-2,0);
\draw (0,0) -- (-2,2);
\fill[color=gray] (12) -- (12a) -- (12b) -- (12);
\fill[color=gray] (12') -- (12'a) -- (12'b) -- (12');
\fill[color=gray] (23) -- (23a) -- (23b) -- (23);
\fill[color=gray] (23') -- (23'a) -- (23'b) -- (23');
\fill[color=gray] (23'') -- (23''a) -- (23''b) -- (23'');
\fill[color=gray] (34) -- (34a) -- (34b) -- (34);
\fill[color=gray] (41) -- (41a) -- (41b) -- (41);
\draw[dashed] (12a) -- (12) -- (12b);
\draw[dashed] (12'a) -- (12') -- (12'b);
\draw[dashed] (23a) -- (23) -- (23b);
\draw[dashed] (23'a) -- (23') -- (23'b);
\draw[dashed] (23''a) -- (23'') -- (23''b);
\draw[dashed] (34a) -- (34) -- (34b);
\draw[dashed] (41a) -- (41) -- (41b);
\draw (12) node[fill,cross,inner sep=1pt,rotate=0]{};
\draw (12') node[fill,cross,inner sep=1pt,rotate=0]{};
\draw (23) node[fill,cross,inner sep=1pt,rotate=0]{};
\draw (23') node[fill,cross,inner sep=1pt,rotate=0]{};
\draw (23'') node[fill,cross,inner sep=1pt,rotate=0]{};
\draw (34) node[fill,cross,inner sep=1pt,rotate=0]{};
\draw (41) node[fill,cross,inner sep=1pt,rotate=26.57]{};
\coordinate[fill,circle,inner sep=0.8pt] (b) at (-1,1);
\coordinate[fill,circle,inner sep=0.8pt] (c) at (-1,0);
\coordinate[fill,circle,inner sep=0.8pt] (d) at (0,-1);
\end{tikzpicture}
\end{minipage}
% (5b)
\begin{minipage}[c]{0.05\textwidth}
$(5b)$
\end{minipage}
\begin{minipage}[c]{0.18\textwidth}
\begin{tikzpicture}[scale=0.6]
\coordinate[fill,circle,inner sep=0.8pt] (0) at (0,0);
\coordinate[fill,circle,inner sep=0.8pt] (1) at (-1,-1);
\coordinate[fill,circle,inner sep=0.8pt] (2) at (1,-1);
\coordinate[fill,circle,inner sep=0.8pt] (3) at (1,0);
\coordinate[fill,circle,inner sep=0.8pt] (4) at (0,1);
\coordinate[fill,circle,inner sep=0.8pt] (5) at (-1,1);
\coordinate (1a) at (-2,-2);
\coordinate (2a) at (2,-2);
\coordinate (3a) at (2,0);
\coordinate (4a) at (0,2);
\coordinate (5a) at (-2,2);
\coordinate (12) at (-0.5,-1);
\coordinate (12a) at (-0.5,-2);
\coordinate (12b) at (-1.5,-2);
\coordinate (12') at (0.5,-1);
\coordinate (12'a) at (0.5,-2);
\coordinate (12'b) at (1.5,-2);
\coordinate (23) at (1,-0.5);
\coordinate (23a) at (2,-1.5);
\coordinate (23b) at (2,-0.5);
\coordinate (34) at (0.5,0.5);
\coordinate (34a) at (2,0.5);
\coordinate (34b) at (0.5,2);
\coordinate (45) at (-0.5,1);
\coordinate (45a) at (-0.5,2);
\coordinate (45b) at (-1.5,2);
\coordinate (51) at (-1,0.5);
\coordinate (51a) at (-2,1.5);
\coordinate (51b) at (-2,0.5);
\coordinate (51') at (-1,-0.5);
\coordinate (51'a) at (-2,-1.5);
\coordinate (51'b) at (-2,-0.5);
\draw (1) -- (2) -- (3) -- (4) -- (5) -- (1);
\draw (0,0) -- (1a);
\draw (0,0) -- (2a);
\draw (0,0) -- (3a);
\draw (0,0) -- (4a);
\draw (0,0) -- (5a);
\draw (0,0) -- (0,-2);
\draw (0,0) -- (-2,0);
\fill[color=gray] (12) -- (12a) -- (12b) -- (12);
\fill[color=gray] (12') -- (12'a) -- (12'b) -- (12');
\fill[color=gray] (23) -- (23a) -- (23b) -- (23);
\fill[color=gray] (34) -- (34a) -- (34b) -- (34);
\fill[color=gray] (45) -- (45a) -- (45b) -- (45);
\fill[color=gray] (51) -- (51a) -- (51b) -- (51);
\fill[color=gray] (51') -- (51'a) -- (51'b) -- (51');
\draw[dashed] (12a) -- (12) -- (12b);
\draw[dashed] (12'a) -- (12') -- (12'b);
\draw[dashed] (23a) -- (23) -- (23b);
\draw[dashed] (34a) -- (34) -- (34b);
\draw[dashed] (45a) -- (45) -- (45b);
\draw[dashed] (51a) -- (51) -- (51b);
\draw[dashed] (51'a) -- (51') -- (51'b);
\draw (12) node[fill,cross,inner sep=1pt,rotate=0]{};
\draw (12') node[fill,cross,inner sep=1pt,rotate=0]{};
\draw (23) node[fill,cross,inner sep=1pt,rotate=0]{};
\draw (34) node[fill,cross,inner sep=1pt,rotate=45]{};
\draw (45) node[fill,cross,inner sep=1pt,rotate=0]{};
\draw (51) node[fill,cross,inner sep=1pt,rotate=0]{};
\draw (51') node[fill,cross,inner sep=1pt,rotate=0]{};
\coordinate[fill,circle,inner sep=0.8pt] (b) at (0,-1);
\coordinate[fill,circle,inner sep=0.8pt] (c) at (-1,0);
\end{tikzpicture}
\end{minipage}\\[6mm]
% (4a)
\begin{minipage}[c]{0.05\textwidth}
$(4a)$
\end{minipage}
\begin{minipage}[c]{0.18\textwidth}
\begin{tikzpicture}[scale=0.6]
\coordinate[fill,circle,inner sep=0.8pt] (0) at (0,0);
\coordinate[fill,circle,inner sep=0.8pt] (1) at (1,-1);
\coordinate[fill,circle,inner sep=0.8pt] (2) at (-1,-1);
\coordinate[fill,circle,inner sep=0.8pt] (3) at (-1,1);
\coordinate[fill,circle,inner sep=0.8pt] (4) at (1,1);
\coordinate (1a) at (2,-2);
\coordinate (2a) at (-2,-2);
\coordinate (3a) at (-2,2);
\coordinate (4a) at (2,2);
\coordinate (12) at (0.5,-1);
\coordinate (12a) at (1.5,-2);
\coordinate (12b) at (0.5,-2);
\coordinate (12') at (-0.5,-1);
\coordinate (12'a) at (-0.5,-2);
\coordinate (12'b) at (-1.5,-2);
\coordinate (23) at (-1,0.5);
\coordinate (23a) at (-2,1.5);
\coordinate (23b) at (-2,0.5);
\coordinate (23') at (-1,-0.5);
\coordinate (23'a) at (-2,-0.5);
\coordinate (23'b) at (-2,-1.5);
\coordinate (34) at (-0.5,1);
\coordinate (34a) at (-0.5,2);
\coordinate (34b) at (-1.5,2);
\coordinate (34') at (0.5,1);
\coordinate (34'a) at (0.5,2);
\coordinate (34'b) at (1.5,2);
\coordinate (41) at (1,0.5);
\coordinate (41a) at (2,0.5);
\coordinate (41b) at (2,1.5);
\coordinate (41') at (1,-0.5);
\coordinate (41'a) at (2,-0.5);
\coordinate (41'b) at (2,-1.5);
\draw (1) -- (2) -- (3) -- (4) -- (1);
\draw (0,0) -- (1a);
\draw (0,0) -- (2a);
\draw (0,0) -- (3a);
\draw (0,0) -- (4a);
\draw (0,0) -- (0,-2);
\draw (0,0) -- (-2,0);
\draw (0,0) -- (0,2);
\draw (0,0) -- (2,0);
\fill[color=gray] (12) -- (12a) -- (12b) -- (12);
\fill[color=gray] (12') -- (12'a) -- (12'b) -- (12');
\fill[color=gray] (23) -- (23a) -- (23b) -- (23);
\fill[color=gray] (23') -- (23'a) -- (23'b) -- (23');
\fill[color=gray] (34) -- (34a) -- (34b) -- (34);
\fill[color=gray] (34') -- (34'a) -- (34'b) -- (34');
\fill[color=gray] (41) -- (41a) -- (41b) -- (41);
\fill[color=gray] (41') -- (41'a) -- (41'b) -- (41');
\draw[dashed] (12a) -- (12) -- (12b);
\draw[dashed] (12'a) -- (12') -- (12'b);
\draw[dashed] (23a) -- (23) -- (23b);
\draw[dashed] (23'a) -- (23') -- (23'b);
\draw[dashed] (34a) -- (34) -- (34b);
\draw[dashed] (34'a) -- (34') -- (34'b);
\draw[dashed] (41a) -- (41) -- (41b);
\draw[dashed] (41'a) -- (41') -- (41'b);
\draw (12) node[fill,cross,inner sep=1pt,rotate=0]{};
\draw (12') node[fill,cross,inner sep=1pt,rotate=0]{};
\draw (23) node[fill,cross,inner sep=1pt,rotate=0]{};
\draw (23') node[fill,cross,inner sep=1pt,rotate=0]{};
\draw (34) node[fill,cross,inner sep=1pt,rotate=0]{};
\draw (34') node[fill,cross,inner sep=1pt,rotate=0]{};
\draw (41) node[fill,cross,inner sep=1pt,rotate=0]{};
\draw (41') node[fill,cross,inner sep=1pt,rotate=0]{};
\coordinate[fill,circle,inner sep=0.8pt] (b) at (-1,0);
\coordinate[fill,circle,inner sep=0.8pt] (c) at (0,-1);
\coordinate[fill,circle,inner sep=0.8pt] (d) at (1,0);
\coordinate[fill,circle,inner sep=0.8pt] (e) at (0,1);
\end{tikzpicture}
\end{minipage}
% (4b)
\begin{minipage}[c]{0.07\textwidth}
$(4b)$
\end{minipage}
\begin{minipage}[c]{0.16\textwidth}
\begin{tikzpicture}[scale=0.5]
\coordinate[fill,circle,inner sep=0.8pt] (0) at (0,0);
\coordinate[fill,circle,inner sep=0.8pt] (1) at (1,-1);
\coordinate[fill,circle,inner sep=0.8pt] (2) at (-1,-1);
\coordinate[fill,circle,inner sep=0.8pt] (3) at (-1,2);
\coordinate[fill,circle,inner sep=0.8pt] (4) at (1,0);
\coordinate (1a) at (2,-2);
\coordinate (2a) at (-2,-2);
\coordinate (3a) at (-2,4);
\coordinate (4a) at (2,0);
\coordinate (12) at (-0.5,-1);
\coordinate (12a) at (-1.5,-2);
\coordinate (12b) at (-0.5,-2);
\coordinate (12') at (0.5,-1);
\coordinate (12'a) at (0.5,-2);
\coordinate (12'b) at (1.5,-2);
\coordinate (23) at (-1,1.5);
\coordinate (23a) at (-2,3.5);
\coordinate (23b) at (-2,2.5);
\coordinate (23') at (-1,0.5);
\coordinate (23'a) at (-2,1.5);
\coordinate (23'b) at (-2,0.5);
\coordinate (23'') at (-1,-0.5);
\coordinate (23''a) at (-2,-0.5);
\coordinate (23''b) at (-2,-1.5);
\coordinate (34) at (-0.5,1.5);
\coordinate (34a) at (-0.5,4);
\coordinate (34b) at (-1.75,4);
\coordinate (34') at (0.5,0.5);
\coordinate (34'a) at (0.5,4);
\coordinate (34'b) at (2,0.5);
\coordinate (41) at (1,-0.5);
\coordinate (41a) at (2,-0.5);
\coordinate (41b) at (2,-1.5);
\draw (1) -- (2) -- (3) -- (4) -- (1);
\draw (0,0) -- (1a);
\draw (0,0) -- (2a);
\draw (0,0) -- (3a);
\draw (0,0) -- (4a);
\draw (0,0) -- (0,-2);
\draw (0,0) -- (-2,0);
\draw (0,0) -- (0,4);
\draw (0,0) -- (-2,2);
\fill[color=gray] (12) -- (12a) -- (12b) -- (12);
\fill[color=gray] (12') -- (12'a) -- (12'b) -- (12');
\fill[color=gray] (23) -- (23a) -- (23b) -- (23);
\fill[color=gray] (23') -- (23'a) -- (23'b) -- (23');
\fill[color=gray] (23'') -- (23''a) -- (23''b) -- (23'');
\fill[color=gray] (34) -- (34a) -- (34b) -- (34);
\fill[color=gray] (34') -- (34'a) -- (34'b) -- (34');
\fill[color=gray] (41) -- (41a) -- (41b) -- (41);
\draw[dashed] (12a) -- (12) -- (12b);
\draw[dashed] (12'a) -- (12') -- (12'b);
\draw[dashed] (23a) -- (23) -- (23b);
\draw[dashed] (23'a) -- (23') -- (23'b);
\draw[dashed] (23''a) -- (23'') -- (23''b);
\draw[dashed] (34a) -- (34) -- (34b);
\draw[dashed] (34'a) -- (34') -- (34'b);
\draw[dashed] (41a) -- (41) -- (41b);
\draw (12) node[fill,cross,inner sep=1pt,rotate=0]{};
\draw (12') node[fill,cross,inner sep=1pt,rotate=0]{};
\draw (23) node[fill,cross,inner sep=1pt,rotate=0]{};
\draw (23') node[fill,cross,inner sep=1pt,rotate=0]{};
\draw (23'') node[fill,cross,inner sep=1pt,rotate=0]{};
\draw (34) node[fill,cross,inner sep=1pt,rotate=45]{};
\draw (34') node[fill,cross,inner sep=1pt,rotate=45]{};
\draw (41) node[fill,cross,inner sep=1pt,rotate=0]{};
\coordinate[fill,circle,inner sep=0.8pt] (b) at (-1,1);
\coordinate[fill,circle,inner sep=0.8pt] (c) at (-1,0);
\coordinate[fill,circle,inner sep=0.8pt] (d) at (0,1);
\coordinate[fill,circle,inner sep=0.8pt] (e) at (0,-1);
\end{tikzpicture}
\end{minipage}
% (4c)
\begin{minipage}[c]{0.06\textwidth}
$(4c)$
\end{minipage}
\begin{minipage}[c]{0.12\textwidth}
\begin{tikzpicture}[scale=0.375]
\coordinate[fill,circle,inner sep=0.8pt] (0) at (0,0);
\coordinate[fill,circle,inner sep=0.8pt] (1) at (-1,-1);
\coordinate[fill,circle,inner sep=0.8pt] (2) at (1,-1);
\coordinate[fill,circle,inner sep=0.8pt] (3) at (-1,3);
\coordinate (1a) at (-2,-2);
\coordinate (2a) at (2,-2);
\coordinate (3a) at (-2,6);
\coordinate (12) at (-0.5,-1);
\coordinate (12a) at (-1.5,-2);
\coordinate (12b) at (-0.5,-2);
\coordinate (12') at (0.5,-1);
\coordinate (12'a) at (0.5,-2);
\coordinate (12'b) at (1.5,-2);
\coordinate (23) at (0.5,0);
\coordinate (23a) at (2,-1.5);
\coordinate (23b) at (0.5,6);
\coordinate (23') at (-0.5,2);
\coordinate (23'a) at (-1.833,6);
\coordinate (23'b) at (-0.5,6);
\coordinate (31) at (-1,-0.5);
\coordinate (31a) at (-2,-1.5);
\coordinate (31b) at (-2,-0.5);
\coordinate (31') at (-1,0.5);
\coordinate (31'a) at (-2,0.5);
\coordinate (31'b) at (-2,1.5);
\coordinate (31'') at (-1,1.5);
\coordinate (31''a) at (-2,2.5);
\coordinate (31''b) at (-2,3.5);
\coordinate (31''') at (-1,2.5);
\coordinate (31'''a) at (-2,4.5);
\coordinate (31'''b) at (-2,5.5);
\draw (1) -- (2) -- (3) -- (1);
\draw (0,0) -- (1a);
\draw (0,0) -- (2a);
\draw (0,0) -- (3a);
\draw (0,0) -- (0,-2);
\draw (0,0) -- (-2,0);
\draw (0,0) -- (0,6);
\draw (0,0) -- (-2,2);
\draw (0,0) -- (-2,4);
\fill[color=gray] (12) -- (12a) -- (12b) -- (12);
\fill[color=gray] (12') -- (12'a) -- (12'b) -- (12');
\fill[color=gray] (23) -- (23a) -- (23b) -- (23);
\fill[color=gray] (23') -- (23'a) -- (23'b) -- (23');
\fill[color=gray] (31) -- (31a) -- (31b) -- (31);
\fill[color=gray] (31') -- (31'a) -- (31'b) -- (31');
\fill[color=gray] (31'') -- (31''a) -- (31''b) -- (31'');
\fill[color=gray] (31''') -- (31'''a) -- (31'''b) -- (31''');
\draw[dashed] (12a) -- (12) -- (12b);
\draw[dashed] (12'a) -- (12') -- (12'b);
\draw[dashed] (23a) -- (23) -- (23b);
\draw[dashed] (23'a) -- (23') -- (23'b);
\draw[dashed] (31a) -- (31) -- (31b);
\draw[dashed] (31'a) -- (31') -- (31'b);
\draw[dashed] (31''a) -- (31'') -- (31''b);
\draw[dashed] (31'''a) -- (31''') -- (31'''b);
\draw (12) node[fill,cross,inner sep=1pt,rotate=0]{};
\draw (12') node[fill,cross,inner sep=1pt,rotate=0]{};
\draw (23) node[fill,cross,inner sep=1pt,rotate=45]{};
\draw (23') node[fill,cross,inner sep=1pt,rotate=45]{};
\draw (31) node[fill,cross,inner sep=1pt,rotate=0]{};
\draw (31') node[fill,cross,inner sep=1pt,rotate=0]{};
\draw (31'') node[fill,cross,inner sep=1pt,rotate=0]{};
\draw (31''') node[fill,cross,inner sep=1pt,rotate=0]{};
\coordinate[fill,circle,inner sep=0.8pt] (b) at (-1,0);
\coordinate[fill,circle,inner sep=0.8pt] (c) at (-1,1);
\coordinate[fill,circle,inner sep=0.8pt] (d) at (-1,2);
\coordinate[fill,circle,inner sep=0.8pt] (e) at (0,-1);
\coordinate[fill,circle,inner sep=0.8pt] (f) at (0,1);
\end{tikzpicture}
\end{minipage}
% (3)
\begin{minipage}[c]{0.05\textwidth}
$(3a)$
\end{minipage}
\begin{minipage}[c]{0.23\textwidth}
\begin{tikzpicture}[scale=0.5]
\coordinate[fill,circle,inner sep=0.8pt] (0) at (0,0);
\coordinate[fill,circle,inner sep=0.8pt] (1) at (-1,-1);
\coordinate[fill,circle,inner sep=0.8pt] (2) at (2,-1);
\coordinate[fill,circle,inner sep=0.8pt] (3) at (-1,2);
\coordinate (1a) at (-2,-2);
\coordinate (2a) at (4,-2);
\coordinate (3a) at (-2,4);
\coordinate (12) at (-0.5,-1);
\coordinate (12a) at (-1.5,-2);
\coordinate (12b) at (-0.5,-2);
\coordinate (12') at (0.5,-1);
\coordinate (12'a) at (0.5,-2);
\coordinate (12'b) at (1.5,-2);
\coordinate (12'') at (1.5,-1);
\coordinate (12''a) at (2.5,-2);
\coordinate (12''b) at (3.5,-2);
\coordinate (23) at (-0.5,1.5);
\coordinate (23a) at (-0.5,4);
\coordinate (23b) at (-1.75,4);
\coordinate (23') at (0.5,0.5);
\coordinate (23'a) at (4,0.5);
\coordinate (23'b) at (0.5,4);
\coordinate (23'') at (1.5,-0.5);
\coordinate (23''a) at (4,-0.5);
\coordinate (23''b) at (4,-1.75);
\coordinate (31) at (-1,-0.5);
\coordinate (31a) at (-2,-1.5);
\coordinate (31b) at (-2,-0.5);
\coordinate (31') at (-1,0.5);
\coordinate (31'a) at (-2,0.5);
\coordinate (31'b) at (-2,1.5);
\coordinate (31'') at (-1,1.5);
\coordinate (31''a) at (-2,2.5);
\coordinate (31''b) at (-2,3.5);
\draw (1) -- (2) -- (3) -- (1);
\draw (0,0) -- (1a);
\draw (0,0) -- (2a);
\draw (0,0) -- (3a);
\draw (0,0) -- (0,-2);
\draw (0,0) -- (-2,0);
\draw (0,0) -- (0,4);
\draw (0,0) -- (4,0);
\draw (0,0) -- (2,-2);
\draw (0,0) -- (-2,2);
\fill[color=gray] (12) -- (12a) -- (12b) -- (12);
\fill[color=gray] (12') -- (12'a) -- (12'b) -- (12');
\fill[color=gray] (12'') -- (12''a) -- (12''b) -- (12'');
\fill[color=gray] (23) -- (23a) -- (23b) -- (23);
\fill[color=gray] (23') -- (23'a) -- (23'b) -- (23');
\fill[color=gray] (23'') -- (23''a) -- (23''b) -- (23'');
\fill[color=gray] (31) -- (31a) -- (31b) -- (31);
\fill[color=gray] (31') -- (31'a) -- (31'b) -- (31');
\fill[color=gray] (31'') -- (31''a) -- (31''b) -- (31'');
\draw[dashed] (12a) -- (12) -- (12b);
\draw[dashed] (12'a) -- (12') -- (12'b);
\draw[dashed] (12''a) -- (12'') -- (12''b);
\draw[dashed] (23a) -- (23) -- (23b);
\draw[dashed] (23'a) -- (23') -- (23'b);
\draw[dashed] (23''a) -- (23'') -- (23''b);
\draw[dashed] (31a) -- (31) -- (31b);
\draw[dashed] (31'a) -- (31') -- (31'b);
\draw[dashed] (31''a) -- (31'') -- (31''b);
\draw (12) node[fill,cross,inner sep=1pt,rotate=0]{};
\draw (12') node[fill,cross,inner sep=1pt,rotate=0]{};
\draw (12'') node[fill,cross,inner sep=1pt,rotate=0]{};
\draw (23) node[fill,cross,inner sep=1pt,rotate=45]{};
\draw (23') node[fill,cross,inner sep=1pt,rotate=45]{};
\draw (23'') node[fill,cross,inner sep=1pt,rotate=45]{};
\draw (31) node[fill,cross,inner sep=1pt,rotate=0]{};
\draw (31') node[fill,cross,inner sep=1pt,rotate=0]{};
\draw (31'') node[fill,cross,inner sep=1pt,rotate=0]{};
\coordinate[fill,circle,inner sep=0.8pt] (b) at (-1,0);
\coordinate[fill,circle,inner sep=0.8pt] (c) at (-1,1);
\coordinate[fill,circle,inner sep=0.8pt] (d) at (0,-1);
\coordinate[fill,circle,inner sep=0.8pt] (e) at (1,-1);
\coordinate[fill,circle,inner sep=0.8pt] (f) at (1,0);
\coordinate[fill,circle,inner sep=0.8pt] (g) at (0,1);
\end{tikzpicture}
\end{minipage}
\caption{Dual intersection complexes $(B,\mathscr{P},\varphi)$ of smooth very ample log Calabi-Yau pairs. The shaded regions are cut out and the dashed lines are mutually identified. Compare this with \cite{KM}, Figure 2, and \cite{Pum}, Figure 5.15.}
\label{fig:listb}
\end{figure}

\begin{con}
\label{con:family}
Refine the dual intersection complex of $\mathfrak{X}\rightarrow\mathbb{A}^1$ by introducing rays starting at the origin and pointing to the integral points of the bounded maximal cell. This yields another polarized polyhedral affine manifold $(B,\mathscr{P},\varphi)$, as shown in Figure \ref{fig:listb}. A refinement of the dual intersection complex gives a logarithmic modification of $\mathfrak{X}_Q\rightarrow\mathbb{A}^2$ (see \S\ref{A:artin}). Since the deformation parameter $s$ is not part of the logarithmic data, the logarithmic modification does not change the general fiber $X=X_{t\neq0}^{s\neq0}$. It can be constructed as follows.
\begin{compactenum}[(1)]
\item Blow up $\mathfrak{X}_Q\subset\mathbb{P}_{\check{\mathscr{P}}}\times\mathbb{A}^2$ at $X_{\sigma_0}\times\{(0,0)\}$, where $X_{\sigma_0}$ is the point corresponding to $\sigma_0$. This corresponds to inserting edges from the origin to corners of $\sigma_0$.
\item Introducing the ray starting at the origin and pointing in the direction of an integral vector on the interior of a bounded edge $\omega$ of $\mathscr{P}$ corresponds to a blow up at $X_\omega\times\mathbb{A}^1\times\{s=0\}$, where $X_\omega=\mathbb{P}^1$ is the component corresponding to $\omega$, i.e., the line through the points corresponding to the bounded maximal cell and the unbounded maximal cell containing $\omega$, respectively.
\end{compactenum}
In cases where $X^0$ is not smooth we refine the asymptotic fan of $\mathscr{P}$. This corresponds to a toric blow up of the toric model $X^0$. This blow up is nef but not necessarily ample. By \cite{KM}, Proposition A.2, the deformation of such a nef toric model still is $(X,D)$ and has Picard group isomorphic to $\text{Pic}(X)$. Note that $\text{Pic}(X)$ is isomorphic to $H_2(X,\mathbb{Z})$ for the del Pezzo surface $X$ by the Kodaira vanishing theorem and Poincar\'e duality.
\end{con}

\begin{defi}
By abuse of notation, from now on $\mathfrak{X}_Q\rightarrow\mathbb{A}^2$ will denote the logarithmic modification from Construction \ref{con:family}. Note that $X^0$ is smooth, toric and nef, but not necessarily ample. We call it a \textit{smooth toric model} of $X$. If $X^0$ is ample it coincides with the toric model of $X$ (Definition \ref{defi:toricmodel}).
\end{defi}

\begin{rem}
Note that $(B,\mathscr{P})$ is simple (\cite{DataI}, Definition 1.60), since all affine singularities have monodromy $\left(\begin{smallmatrix}1& 1\\ 0& 1\end{smallmatrix}\right)$ in suitable coordinates. Thus we can apply the reconstruction theorem (\cite{GS11}, Proposition 2.42) together with the construction of a tropical superpotential from \cite{CPS} to construct the mirror Landau-Ginzburg model to $(X,D)$.
\end{rem}

\begin{expl}
\label{expl:8'a2}
Consider the smoothing of $\mathbb{P}(1,1,2)$ to $\mathbb{P}^1\times\mathbb{P}^1$ (case (8'a)) from Example \ref{expl:8'a}. The logarithmic modification from Construction \ref{con:family} is a $2$-parameter family $\mathfrak{X}_Q\rightarrow\mathbb{A}^2$ such that $\mathfrak{X}\rightarrow\mathbb{A}^1$ is a toric degeneration of $\mathbb{P}^1\times\mathbb{P}^1$ and $\mathfrak{X}^0\rightarrow\mathbb{A}^1$ is a toric degeneration of the Hirzebruch surface $\mathbb{F}_2$, the $\mathbb{P}^1$-bundle over $\mathbb{P}^1$ given by $\mathbb{F}_2=\mathbb{P}(\mathcal{O}_{\mathbb{P}^1}\oplus\mathcal{O}_{\mathbb{P}^1}(2))$. This is the smooth surface obtained by blowing up the singular point on $\mathbb{P}(1,1,2)$, corresponding the the subdivision of the asymptotic fan given in Figure \ref{fig:8'a2}.

The Picard group $\text{Pic}(\mathbb{F}_2)\simeq H_2(\mathbb{F}_2,\mathbb{Z})$ is generated by the class of a fiber $F$ and the class of a section, e.g., the exceptional divisor $E$ of the blow up. The intersection numbers are $F^2=0$, $E^2=-2$ and $E\cdot F=1$. The anticanonical bundle is $-K_{\mathbb{F}_2}=2F+S+E=4F+2E$, where $S=2F+E$ is the class of a section different from the exceptional divisor. The classes of curves corresponding to the rays in the fan of $\mathbb{F}_2$ are given in Figure \ref{fig:8'a2}.

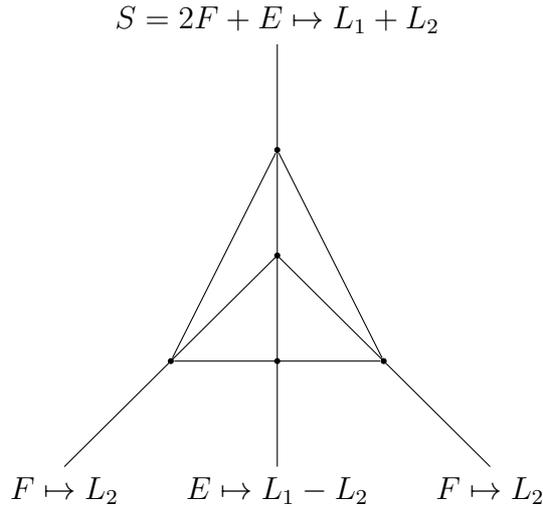
\begin{figure}[h!]
\centering
\begin{tikzpicture}[scale=1.4]
\coordinate[fill,circle,inner sep=0.8pt] (0) at (0,0);
\coordinate[fill,circle,inner sep=0.8pt] (1) at (-1,-1);
\coordinate[fill,circle,inner sep=0.8pt] (2) at (1,-1);
\coordinate[fill,circle,inner sep=0.8pt] (3) at (0,1);
\coordinate (1a) at (-2,-2);
\coordinate (2a) at (2,-2);
\coordinate (3a) at (0,2);
\coordinate (12) at (-0.5,-1);
\coordinate (12a) at (-1.5,-2);
\coordinate (12b) at (-0.5,-2);
\coordinate (12') at (0.5,-1);
\coordinate (12'a) at (0.5,-2);
\coordinate (12'b) at (1.5,-2);
\coordinate (23) at (0.5,0);
\coordinate (23a) at (0.5,2);
\coordinate (23b) at (2,-1.5);
\coordinate (31) at (-0.5,0);
\coordinate (31a) at (-0.5,2);
\coordinate (31b) at (-2,-1.5);
\draw (1) -- (2) -- (3) -- (1);
\draw (0,0) -- (1a) node[below]{$F\mapsto L_2$};
\draw (0,0) -- (2a) node[below]{$F\mapsto L_2$};
\draw (0,0) -- (3a) node[above]{$S=2F+E\mapsto L_1+L_2$};
\draw (0,0) -- (0,-2) node[below]{$E\mapsto L_1-L_2$};
\coordinate[fill,circle,inner sep=0.8pt] (b) at (0,-1);
\end{tikzpicture}
\caption{The dual intersection complex of a toric degeneration of $\mathbb{F}_2$. The classes of curves corresponding to the rays and their images under $\text{Pic}(\mathbb{F}_2) \iso \text{Pic}(\mathbb{P}^1\times\mathbb{P}^1)$ are given.}
\label{fig:8'a2}
\end{figure}

The Picard group of $\mathbb{P}^1\times\mathbb{P}^1$ is generated by the class of a bidegree $(1,0)$ curve $L_1$ and a bidegree $(0,1)$ curve $L_2$, with intersection numbers $L_1^2=0$, $L_2^2=0$ and $L_1\cdot L_2=1$. There is an isomorphism
\[ \text{Pic}(\mathbb{F}_2) \iso \text{Pic}(\mathbb{P}^1\times\mathbb{P}^1), \ F \mapsto L_2, \ E \mapsto L_1-L_2. \]
Note that there is another isomorphism by the symmetry $L_1\leftrightarrow L_2$ and we made a choice here, fixed by the deformation $\mathfrak{X}^0\hookrightarrow\mathfrak{X}_Q$. We will use this isomorphism in \S\ref{S:calc8'a} to calculate the logarithmic Gromov-Witten invariants of $\mathbb{P}^1\times\mathbb{P}^1$ in an alternative way.
\end{expl}

\subsection{Affine charts}									%%
\label{S:affinecharts}

Figure \ref{fig:listb} shows the dual intersection complexes $(B,\mathscr{P},\varphi)$ in the chart of $\sigma_0$ (Definition \ref{defi:sigma0}). The shaded regions are cut out and the dashed lines are mutually identified, so in fact all unbounded edges are parallel.

\begin{defi}
\label{defi:mout}
Let $m_{\text{out}}\in\Lambda_B$ denote the primitive integral tangent vector pointing in the unique unbounded direction on $B$.
\end{defi}

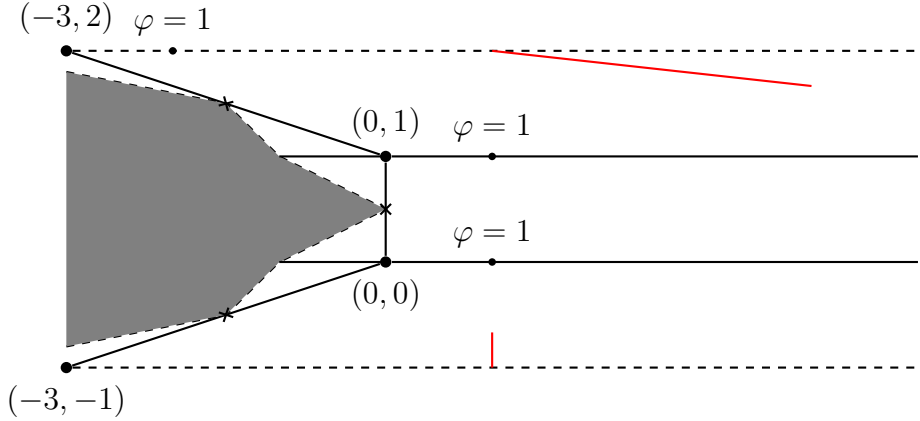
\begin{figure}[h!]
\centering
\begin{tikzpicture}[scale=1.4]
\coordinate[label=below:${(0,0)}$,fill,circle,inner sep = 1.5pt] (1) at (0,0);
\coordinate[label=above:${(0,1)}$,fill,circle,inner sep = 1.5pt] (2) at (0,1);
\coordinate[label=below:${(-3,-1)}$,fill,circle,inner sep = 1.5pt] (3) at (-3,-1);
\coordinate[label=above:${(-3,2)}$,fill,circle,inner sep = 1.5pt] (4) at (-3,2);
% phi
\coordinate[label=above:${\varphi=1}$,fill,circle,inner sep = 1pt] (0a) at (-2,2);
\coordinate[label=above:${\varphi=1}$,fill,circle,inner sep = 1pt] (0b) at (1,1);
\coordinate[label=above:${\varphi=1}$,fill,circle,inner sep = 1pt] (0c) at (1,0);
\draw[thick] (3) -- (1) -- (2) -- (4);
\draw[thick] (1) -- (5,0);
\draw[thick] (2) -- (5,1);
\draw[thick,dashed] (3) -- (5,-1);
\draw[thick,dashed] (4) -- (5,2);
% cut
\draw[thick,dashed] (-1.5,-0.25) -- (0,0.5) -- (-1.5,1.25);
\draw[thick,dashed] (-3,-0.8) -- (-1.5,-0.5) -- (-0.75,0.25);
\draw[thick,dashed] (-3,1.8) -- (-1.5,1.5) -- (-0.75,0.75);
\fill[color=gray] (-1.5,-0.25) -- (0,0.5) -- (-1.5,1.25) -- (-1.5,-0.25) -- cycle;
\fill[color=gray] (-3,1.8) -- (-1.5,1.5) -- (-0.75,0.75) -- (-0.75,0.25) -- (-1.5,-0.5) -- (-3,-0.8);
\coordinate[fill,cross,thick,inner sep=2pt] (5) at (0,0.5);
\coordinate[fill,cross,thick,inner sep=2pt,rotate=26.57] (6) at (-1.5,-0.5);
\coordinate[fill,cross,thick,inner sep=2pt,rotate=-30] (7) at (-1.5,1.5);
% refinement
\draw[thick] (0,0) -- (-1,0);
\draw[thick] (0,1) -- (-1,1);
% red
\draw[red,thick] (1,2) -- (4,1.667);
\draw[red,thick] (1,-1) -- (1,-0.667);
% shift
\draw[white] (6,0) -- (6,1);
\end{tikzpicture}
\caption{$(B,\mathscr{P},\varphi)$ for $(\mathbb{P}^2,E)$ in the chart of an unbounded maximal cell. The dark region is cut out and the dashed lines are mutually identified. A straight line is shown in red.}
\label{fig:dualb}
\end{figure}

\begin{expl}
Figure \ref{fig:dualb} shows the dual intersection complex $(B,\mathscr{P},\varphi)$ of $(\mathbb{P}^2,E)$ in the chart of an unbounded maximal cell. Intuitively, this picture can be obtained by mutually gluing the dashed lines in Figure \ref{fig:listb}, (9). The two horizontal dashed lines are identified. The monodromy transformation by passing across the upper horizontal dashed line is given by $\Lambda_B \rightarrow \Lambda_B, m \mapsto \left(\begin{smallmatrix}1&9\\ 0&1\end{smallmatrix}\right) \cdot m$.
\end{expl}

\begin{figure}[h!]
\centering
\begin{tikzpicture}[scale=0.8]
\coordinate[label=below:${(0,0)}$,fill,circle,inner sep = 1.5pt] (1) at (0,0);
\coordinate[label=above:${(0,1)}$,fill,circle,inner sep = 1.5pt] (2) at (0,1);
\coordinate[label=below:${(-3,-1)}$,fill,circle,inner sep = 1.5pt] (3) at (-3,-1);
\coordinate[label=above:${(-3,2)}$,fill,circle,inner sep = 1.5pt] (4) at (-3,2);
\draw[thick] (3) -- (1) -- (2) -- (4);
\draw[thick] (1) -- (4.5,0);
\draw[thick] (2) -- (4.5,1);
\draw[thick] (3) -- (4.5,-1);
\draw[thick] (4) -- (4.5,2);
\draw[thick] (3) -- (-9,-2) -- (-9.9,-2.1);
\draw[thick] (4) -- (-9,3) -- (-9.9,3.1);
\draw[thick] (-9,-2) -- (4.5,-2);
\draw[thick] (-9,3) -- (4.5,3);
\fill[opacity=0.1] (-9.9,-2.1) -- (4.5,-2.1) -- (4.5,-1) -- (-9.9,-1) -- (-9.9,-2.1);
\fill[opacity=0.1] (-9.9,3.1) -- (4.5,3.1) -- (4.5,2) -- (-9.9,2) -- (-9.9,3.1);
% cut
\draw[thick,dashed] (-1.5,-0.25) -- (0,0.5) -- (-1.5,1.25);
\draw[thick,dashed] (-6,-1.4) -- (-1.5,-0.5) -- (-0.75,0.25);
\draw[thick,dashed] (-6,2.4) -- (-1.5,1.5) -- (-0.75,0.75);
\draw[thick,dashed] (-3,-0.75) -- (-6,-1.5) -- (-9.9,-1.9875);
\draw[thick,dashed] (-3,1.75) -- (-6,2.5) -- (-9.9,2.9875);
\fill[color=gray] (-1.5,-0.25) -- (0,0.5) -- (-1.5,1.25) -- (-1.5,-0.25);
\fill[color=gray] (-6,2.4) -- (-1.5,1.5) -- (-0.75,0.75) -- (-0.75,0.25) -- (-1.5,-0.5) -- (-6,-1.4);
\fill[color=gray] (-9.9,-1.9875) -- (-6,-1.5) -- (-3,-0.75) -- (-3,1.75) -- (-6,2.5) -- (-9.9,2.9875);
\coordinate[fill,cross,thick,inner sep=2pt] (5) at (0,0.5);
\coordinate[fill,cross,thick,inner sep=2pt,rotate=26.57] (6) at (-1.5,-0.5);
\coordinate[fill,cross,thick,inner sep=2pt,rotate=-30] (7) at (-1.5,1.5);
\coordinate[fill,cross,thick,inner sep=2pt,rotate=15] (8) at (-6,-1.5);
\coordinate[fill,cross,thick,inner sep=2pt,rotate=-15] (9) at (-6,2.5);
% refinement
\draw[thick] (0,0) -- (-1,0);
\draw[thick] (0,1) -- (-1,1);
\draw[thick] (0,-1) -- (-4,-1);
\draw[thick] (0,2) -- (-4,2);
% red
\draw[red,thick] (4,1.667) -- (-2,2.333);
\draw[red,thick] (1,-0.667) -- (1,-1.333);
\end{tikzpicture}
\caption{A chart of a covering space $\bar{B}$ of $B$ with fundamental domain the white region (including one of the rays on its border). The preimage of the straight line from Figure \ref{fig:dualb} is shown in red.}
\label{fig:dualc}
\end{figure}
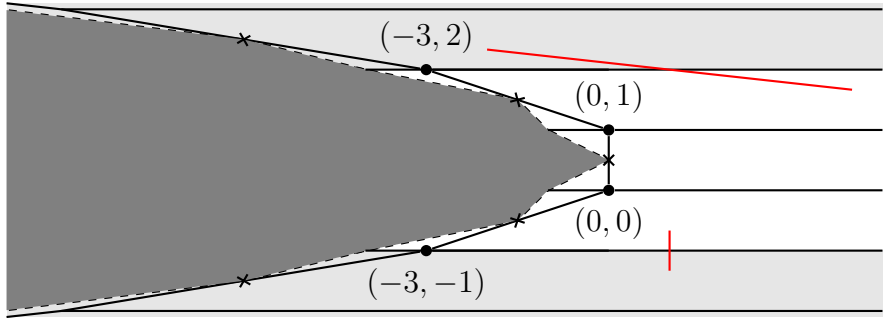

We can extend the description of the affine structure across the horizontal dashed line by giving a chart of a discrete covering space $\bar{B}$ of $B$ (Figure \ref{fig:dualc}). Passing from one fundamental domain to an adjacent one amounts to applying the monodromy transformation by crossing the horizontal dashed line in $B$.

This gives a trivialization $\Lambda_{\bar{B}} \simeq M = \mathbb{Z}^2$ on $\bar{B}\setminus(\textup{Int}(\bar{\sigma_0})\cup\bar{\Delta})$, where $\bar{\sigma_0}$ and $\bar{\Delta}$ are the preimages of the bounded maximal cell $\sigma_0$ and the discriminant locus $\Delta$, respectively. We will see in Lemma \ref{lem:global} that the consistent wall structure $\mathscr{S}_\infty$ defined by $(B,\mathscr{P},\varphi)$ has support disjoint from $\text{Int}(\sigma_0)$. Hence, the whole scattering procedure can be described in this affine chart of the covering space $\bar{B}$. This allows for a simple implementation of the scattering algorithm (see \S\ref{S:calc}).

\section{Resolution of log singularities}					%%
\label{S:resolution}

Let $Q$ be a Fano polytope and consider the family $\mathfrak{X}_Q\rightarrow\mathbb{A}^2$ from Construction \ref{con:family}. Equip $\mathbb{A}^2$ with the divisorial log structure defined by $V(t)\subset\mathbb{A}^2$ and $\mathfrak{X}_Q$ with the divisorial log structure defined by $\mathfrak{X}_0\cup\mathfrak{D}_Q\subset\mathfrak{X}_Q$, that is, the sheaf of monoids
\[ \mathcal{M}_{(\mathfrak{X}_Q,\mathfrak{X}_0\cup\mathfrak{D}_Q)} := (j_\star\mathcal{O}_{\mathfrak{X}_Q\setminus (\mathfrak{X}_0\cup\mathfrak{D}_Q)}^\times)\cap\mathcal{O}_{\mathfrak{X}_Q}, \quad j : \mathfrak{X}_Q\setminus (\mathfrak{X}_0\cup\mathfrak{D}_Q) \hookrightarrow \mathfrak{X}_Q. \]
For an introduction to log structures see e.g. \cite{Kat1} or \cite{Gr10}, {\S}3.

If we consider the fibers $X_t^s$ or the families $\mathfrak{X}_t\rightarrow\mathbb{A}^1$ or $\mathfrak{X}^s\rightarrow\mathbb{A}^1$ as log schemes, we always mean equipped with the log structure by restriction of the above log structure. Now for each $s\in\mathbb{A}^1$ the family $\mathfrak{X}^s\rightarrow\mathbb{A}^1$ is log smooth away from finitely many points on the central fiber, corresponding to the affine singularities of the dual intersection complex $(B,\mathscr{P},\varphi)$. At these points $\mathfrak{X}^s$ is locally given by $\textup{Spec }\mathbb{C}[x,y,w,t]/(xy-t^l(w+s))$ with log structure given by $V(t)\cup V(w)$. This is isomorphic to $\text{Spec }\mathbb{C}[x,y,\tilde{w},t]/(xy-t^l\tilde{w})$ with $\tilde{w}=w+s$. The log structure is given by $V(t)\cup V(\tilde{w})$ for $s=0$ and by $V(t)$ for $s\neq 0$. Arguments as in \cite{Gr10}, Example 3.20, show that for $s\neq 0$ this is not fine at the point given by $x=y=w=t=0$. For $s=0$ the log structure is fine saturated but not log smooth. Following \cite{DataII}, Lemma 2.12, we describe a small log resolution $\tilde{\mathfrak{X}}^s \rightarrow \mathfrak{X}^s$ such that $\tilde{\mathfrak{X}}^s$ is fine and log smooth over $\mathbb{A}^1$.

\subsection{The local picture}						%%
\label{S:resloc}

$\text{Spec }\mathbb{C}[x,y,\tilde{w},t]/(xy-t^l\tilde{w})$ is the affine toric variety defined by the cone $\sigma$ generated by $(0,0,1)$, $(0,1,0)$, $(1,0,1)$ and $(l,1,0)$. In fact,
\begin{eqnarray*}
\textup{Spec }\mathbb{C}[\sigma^\vee \cap \mathbb{Z}^3] &=& \textup{Spec }\mathbb{C}[z^{(1,0,0)},z^{(-1,l,1)},z^{(0,0,1)},z^{(0,1,0)}] \\
&=& \textup{Spec }\faktor{\mathbb{C}[x,y,\tilde{w},t]}{(xy-\tilde{w}t^l)}.
\end{eqnarray*}
We obtain a toric blow up by subdividing the fan consisting of the single cone $\sigma$. There are two ways of doing this and they are related by a \textit{flop}. We choose the subdivision $\Sigma$ as in \cite{DataII}, Lemma 2.12, with maximal cones $\sigma_1$ generated by $(0,0,1)$, $(1,0,1)$ and $(0,1,0)$, and $\sigma_2$ generated by $(l,1,0)$, $(1,0,1)$ and $(0,1,0)$.

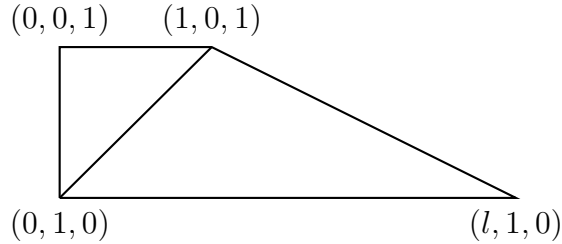
\begin{figure}[h!]
\centering
\begin{tikzpicture}[scale=2]
\draw[thick] (0,0) node[below]{$(0,1,0)$} -- (3,0) node[below]{$(l,1,0)$} -- (1,1) node[above]{$(1,0,1)$} -- (0,1) node[above]{$(0,0,1)$} -- (0,0);
\draw[thick] (0,0) -- (1,1);
\end{tikzpicture}
\caption{Generators of the cone defining a toric model of a log singularity and a choice of subdivision.}
\label{fig:cone}
\end{figure}

These cones define affine toric varieties
\begin{eqnarray*}
X_{\sigma_1} &=& \textup{Spec }\mathbb{C}[z^{(1,0,0)},z^{(0,1,0)},z^{(-1,0,1)}] = \textup{Spec }\mathbb{C}[x,t,u] = \mathbb{A}^3, \\
X_{\sigma_2} &=& \textup{Spec }\mathbb{C}[z^{(-1,l,1)},z^{(0,0,1)},z^{(0,1,0)},z^{(1,0,-1)}] = \textup{Spec }\mathbb{C}[y,\tilde{w},t,v]/(yv-t^l), \\
X_{\sigma_{12}} &=& \textup{Spec }\mathbb{C}[z^{(1,0,0)},z^{(0,1,0)},z^{\pm(-1,0,1)}] = \textup{Spec }\mathbb{C}[x,t,u^{\pm1}] = \mathbb{A}^2 \times \mathbb{G}_m.
\end{eqnarray*}
The toric variety $X_\Sigma$ defined by $\Sigma$ is obtained by gluing $X_{\sigma_1}$ and $X_{\sigma_2}$ along $X_{\sigma_{12}}$. This is the fibered coproduct (with $u=U/V$ and $v=V/U$)
\[ X_\Sigma = X_{\sigma_1} \amalg_{X_{\sigma_{12}}} X_{\sigma_2} = \textup{Proj }\faktor{\mathbb{C}[x,y,\tilde{w},t][U,V]}{(\tilde{w}V-xU,yV-t^lU)}. \]
Note that we take $\textup{Proj}$ of the polynomial ring with variables $U,V$ over the ring $\mathbb{C}[x,y,\tilde{w},t]$, so only $U,V$ are homogeneous coordinates, of degree $1$. The grading is given by degree in $U$ and $V$. The exceptional set of the resolution $X_\Sigma\rightarrow X_\sigma$ is a line contained in the irreducible component of the central fiber $X_{\Sigma,0}$ given by $y=0$.

\begin{figure}[h!]
\centering
\begin{tikzpicture}[scale=2]
\draw (-2,-0.6) -- (0,-0.6) -- (2.4,0.6);
\draw (0,-0.6) -- (0,0.6);
% exceptional line
\draw[thick,red] (0,0) -- (0.4,0) -- (0.8,0);
% arrows
\draw[thick,->] (0.8,0) -- (1,0.1) node[right]{$x$};
\draw[thick,->] (0,0) -- (-0.2,0) node[below]{$y$};
\draw[thick,->] (0,0) -- (0,0.2) node[right]{$w$};
\draw[thick,->] (0.8,0) -- (0.6,0) node[above]{$u$};
\draw[thick,->] (0,0) -- (0.2,0) node[below]{$v$};
\end{tikzpicture}
\caption{Local picture of the central fiber of the resolution, with exceptional line shown in red.}
\label{fig:local}
\end{figure}
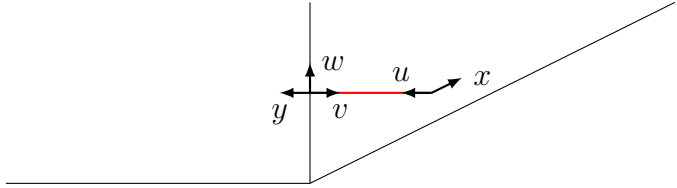

Equip $X_\Sigma$ with the divisorial log structure by its central fiber $X_{\Sigma,0}$ and pull this log structure back to $X_{\Sigma,0}$. Then $X_{\sigma_1}$ and $X_{\sigma_2}$ are log smooth with respect to the restriction of this log structure, since they are simple normal crossings. They form an affine cover of $X_{\Sigma,0}$, so $X_{\Sigma,0}$ is log smooth.

Similarly, if we make the opposite choice of subdivision, the exceptional line is contained in the irreducible component of the central fiber given by $x=0$.

\subsection{The global picture}							%%

There are two geometric descriptions of the toric blow up $X_\Sigma \rightarrow X_\sigma$ considered in \S\ref{S:resloc} above:
\begin{compactenum}[(1)]
\item $X_\Sigma\rightarrow X_\sigma$ is given by blowing up $X_\sigma$ along $\{y=t=0\}$. Indeed, this corresponds to subdividing $\sigma$ by cones connecting the face of $\sigma$ corresponding to $\{y=t=0\}$, in our case the ray generated by $(1,0,1)$, with all other faces of $\sigma$, leading to the fan $\Sigma$. 
\item Let $X_{\Sigma'}$ be the blow up of $X_\sigma$ along the origin. This corresponds to inserting a ray in the center of $\sigma$ and connecting all faces of $\sigma$ with this ray, leading to a fan $\Sigma'$. The exceptional set of this blow up is isomorphic to $\mathbb{P}^1\times\mathbb{P}^1$. Choose one of the $\mathbb{P}^1$-factors and partially contract the exceptional set in $X_{\Sigma'}$ by projecting to this factor. This corresponds to one of the two ways to pair off the four maximal cones in $\Sigma'$ into two cones. One choice leads to $\Sigma$, so we obtain $X_{\Sigma'} \rightarrow X_\Sigma$ by a partial contraction of the exceptional set in $X_{\Sigma'}$.
\end{compactenum}
These constructions can also be performed globally on $\mathfrak{X}_Q$. 

In (1) we blow up $\mathfrak{X}_Q$ along one of the irreducible components of $\mathfrak{X}_0$. We can do this for all irreducible components of $\mathfrak{X}_0$ and obtain a log smooth family over $\mathbb{A}^2$. However, this family will depend on the order of the blow ups and the irreducible components of its central fiber will contain different numbers of exceptional lines. 

In (2) we blow up $\mathfrak{X}_Q$ along curves on $\mathfrak{X}_0$ and then partially contract the exceptional sets. In each step we have two ways to choose the contraction. Making the right choices we obtain a more symmetric resolution.

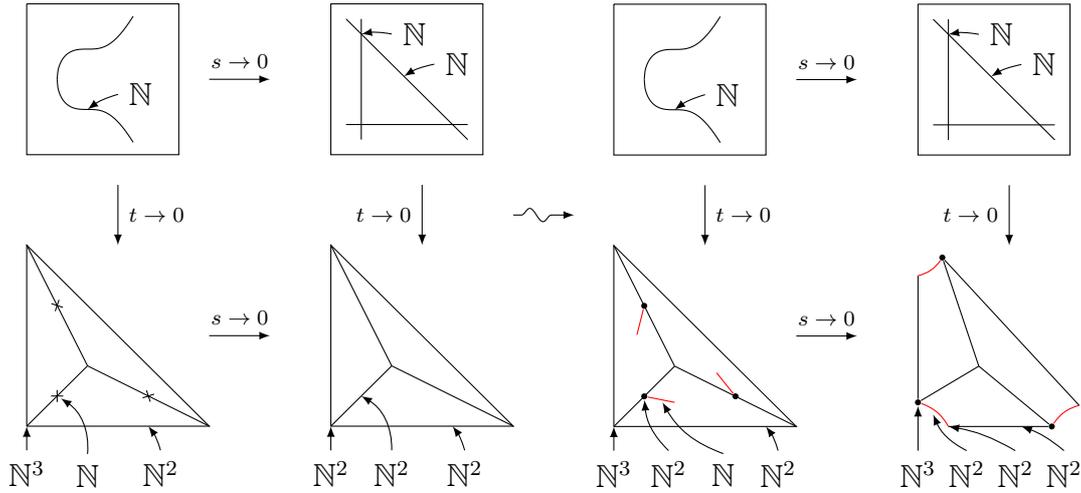
\begin{figure}[h!]
\centering
\begin{tikzpicture}[scale=0.8,decoration={snake,pre length=3pt,post length=7pt}]
% arrows
\draw[->] (2,0.5) -- (2.5,0.5) node[above]{\tiny$s\rightarrow 0$} -- (3,0.5);
\draw[->] (0.5,3) -- (0.5,2.5) node[right]{\tiny$t\rightarrow 0$} -- (0.5,2);
\draw[->] (2,4.75) -- (2.5,4.75) node[above]{\tiny$s\rightarrow 0$} -- (3,4.75);
\draw[->] (5.5,3) -- (5.5,2.5) node[left]{\tiny$t\rightarrow 0$} -- (5.5,2);
\draw[->,decorate] (7,2.5) -- (8,2.5);
% t=0 s!=0
\coordinate (0) at (0,0);
\coordinate (1) at (-1,-1);
\coordinate (2) at (2,-1);
\coordinate (3) at (-1,2);
\coordinate[fill,cross,inner sep=1pt,rotate=45] (1a) at (-0.5,-0.5);
\coordinate[fill,cross,inner sep=1pt,rotate=63.43] (2a) at (1,-0.5);
\coordinate[fill,cross,inner sep=1pt,rotate=26.57] (3a) at (-0.5,1);
\draw (1) -- (2) -- (3) -- (1);
\draw (0) -- (1);
\draw (0) -- (2);
\draw (0) -- (3);
\coordinate (4) at (1,-1);
\coordinate (5) at (1.2,-1.4);
\draw[->] (5) node[below]{$\mathbb{N}^2$} to [bend right=10] (4);
\coordinate (6) at (-1,-1.4);
\draw[->] (6) node[below]{$\mathbb{N}^3$} to [bend right=0] (1);
\coordinate (7) at (0,-1.46);
\draw[->] (7) node[below]{$\mathbb{N}$} to [bend right=30] (1a);
% t=0 s=0
\coordinate (0) at (5,0);
\coordinate (1) at (4,-1);
\coordinate (2) at (7,-1);
\coordinate (3) at (4,2);
\coordinate (1a) at (4.5,-0.5);
\coordinate (2a) at (6,-0.5);
\coordinate (3a) at (4.5,1);
\draw (1) -- (2) -- (3) -- (1);
\draw (0) -- (1);
\draw (0) -- (2);
\draw (0) -- (3);
\coordinate (4) at (6,-1);
\coordinate (5) at (6.2,-1.4);
\draw[->] (5) node[below]{$\mathbb{N}^2$} to [bend right=10] (4);
\coordinate (6) at (4,-1.4);
\draw[->] (6) node[below]{$\mathbb{N}^2$} to [bend right=0] (1);
\coordinate (7) at (5,-1.4);
\draw[->] (7) node[below]{$\mathbb{N}^2$} to [bend right=30] (1a);
% t!=0 s!=0
\draw (-1,6) -- (-1,3.5) -- (1.5,3.5) -- (1.5,6) -- (-1,6);
\draw[smooth,samples=100,domain=-0.5:0.75] plot (\x,{sqrt((2*\x)^3+1)/2+4.75});
\draw[smooth,samples=100,domain=-0.5:0.75] plot (\x,{-sqrt((2*\x)^3+1)/2+4.75});
\coordinate (4) at (0,4.25);
\coordinate (5) at (0.5,4.5);
\draw[->] (5) node[right]{$\mathbb{N}$} to [bend right=10] (4);
% t!=0 s=0
\draw (4,6) -- (4,3.5) -- (6.5,3.5) -- (6.5,6) -- (4,6);
\draw (4.5,3.75) -- (4.5,5.75);
\draw (4.25,4) -- (6.25,4);
\draw (4.25,5.75) -- (6.25,3.75);
\coordinate (4) at (4.5,5.5);
\coordinate (5) at (5,5.5);
\draw[->] (5) node[right]{$\mathbb{N}$} to [bend right=10] (4);
\coordinate (6) at (5.2,4.8);
\coordinate (7) at (5.7,5);
\draw[->] (7) node[right]{$\mathbb{N}$} to [bend right=10] (6);
\end{tikzpicture}
\begin{tikzpicture}[scale=0.8]
% arrows
\draw[->] (2,0.5) -- (2.5,0.5) node[above]{\tiny$s\rightarrow 0$} -- (3,0.5);
\draw[->] (0.5,3) -- (0.5,2.5) node[right]{\tiny$t\rightarrow 0$} -- (0.5,2);
\draw[->] (2,4.75) -- (2.5,4.75) node[above]{\tiny$s\rightarrow 0$} -- (3,4.75);
\draw[->] (5.5,3) -- (5.5,2.5) node[left]{\tiny$t\rightarrow 0$} -- (5.5,2);
% t=0 s!=0
\coordinate (0) at (0,0);
\coordinate (1) at (-1,-1);
\coordinate (2) at (2,-1);
\coordinate (3) at (-1,2);
\coordinate[fill,circle,inner sep=0.8pt] (1a) at (-0.5,-0.5);
\coordinate[fill,circle,inner sep=0.8pt] (2a) at (1,-0.5);
\coordinate[fill,circle,inner sep=0.8pt] (3a) at (-0.5,1);
\draw (1) -- (2) -- (3) -- (1);
\draw (0) -- (1);
\draw (0) -- (2);
\draw (0) -- (3);
\draw[red] (1a) -- (0,-0.6);
\draw[red] (2a) -- (0.69,-0.11);
\draw[red] (3a) -- (-0.62,0.52);
\coordinate (4) at (1.5,-1);
\coordinate (5) at (1.7,-1.4);
\draw[->] (5) node[below]{$\mathbb{N}^2$} to [bend right=10] (4);
\coordinate (6) at (-1,-1.4);
\draw[->] (6) node[below]{$\mathbb{N}^3$} to [bend right=0] (1);
\coordinate (7) at (-0.1,-1.4);
\draw[->] (7) node[below]{$\mathbb{N}^2$} to [bend left=20] (1a);
\coordinate (8) at (0.8,-1.45);
\coordinate (9) at (-0.2,-0.6);
\draw[->] (8) node[below]{$\mathbb{N}$} to [bend left=20] (9);
% t=0 s=0
\coordinate (0) at (5,0);
\coordinate[fill,circle,inner sep=0.8pt] (1) at (4,-0.6);
\coordinate[fill,circle,inner sep=0.8pt] (2) at (6.2,-1);
\coordinate[fill,circle,inner sep=0.8pt] (3) at (4.4,1.8);
\coordinate (1a) at (4.5,-1);
\coordinate (2a) at (6.65,-0.65);
\coordinate (3a) at (4,1.5);
\draw (1a) -- (2);
\draw (2a) -- (3);
\draw (3a) -- (1);
\draw (0) -- (1);
\draw (0) -- (2);
\draw (0) -- (3);
\draw[red] (1) to [bend left=20] (1a);
\draw[red] (2) to [bend left=20] (2a);
\draw[red] (3) to [bend left=20] (3a);
\coordinate (4) at (5.7,-1);
\coordinate (5) at (6.4,-1.4);
\draw[->] (5) node[below]{$\mathbb{N}^2$} to [bend right=10] (4);
\coordinate (6) at (4,-1.4);
\draw[->] (6) node[below]{$\mathbb{N}^3$} to [bend right=0] (1);
\coordinate (7) at (5.6,-1.4);
\draw[->] (7) node[below]{$\mathbb{N}^2$} to [bend right=10] (1a);
\coordinate (8) at (4.25,-0.8);
\coordinate (9) at (4.8,-1.4);
\draw[->] (9) node[below]{$\mathbb{N}^2$} to [bend left=20] (8);
% t!=0 s!=0
\draw (-1,6) -- (-1,3.5) -- (1.5,3.5) -- (1.5,6) -- (-1,6);
\draw[smooth,samples=100,domain=-0.5:0.75] plot (\x,{sqrt((2*\x)^3+1)/2+4.75});
\draw[smooth,samples=100,domain=-0.5:0.75] plot (\x,{-sqrt((2*\x)^3+1)/2+4.75});
\coordinate (4) at (0,4.25);
\coordinate (5) at (0.5,4.5);
\draw[->] (5) node[right]{$\mathbb{N}$} to [bend right=10] (4);
% t!=0 s=0
\draw (4,6) -- (4,3.5) -- (6.5,3.5) -- (6.5,6) -- (4,6);
\draw (4.5,3.75) -- (4.5,5.75);
\draw (4.25,4) -- (6.25,4);
\draw (4.25,5.75) -- (6.25,3.75);
\coordinate (4) at (4.5,5.5);
\coordinate (5) at (5,5.5);
\draw[->] (5) node[right]{$\mathbb{N}$} to [bend right=10] (4);
\coordinate (6) at (5.2,4.8);
\coordinate (7) at (5.7,5);
\draw[->] (7) node[right]{$\mathbb{N}$} to [bend right=10] (6);
\end{tikzpicture}
\caption{Fibers of the families $\mathfrak{X}_Q\rightarrow\mathbb{A}^2$ (left) and $\tilde{\mathfrak{X}}_Q\rightarrow\mathbb{A}^2$ (right) for $(\mathbb{P}^2,E)$. The exceptional lines are red. Some stalks of the ghost sheaves are given.}
\label{fig:bigpicture1}
\end{figure}

\begin{con}[The log smooth degeneration]
\label{con:deg2}
For each log singularity on $\mathfrak{X}_Q$ we have two choices of a small resolution as in (2) above, fixed by choosing which irreducible component of $\tilde{X}_0^{s\neq 0}$ contains the exceptional line. We make a symmetric choice such that we have one exceptional line on each irreducible component of $\tilde{X}_0$ (see Figure \ref{fig:bigpicture1}). The only reason for doing so is to avoid distinction of cases. We obtain a log smooth family $(\tilde{\mathfrak{X}}_Q,\tilde{\mathfrak{D}}_Q)\rightarrow \mathbb{A}^2$. Since we only change the fibers $X_0^s$, we still have that $(\tilde{\mathfrak{X}}^0,\mathfrak{D}^0)\rightarrow\mathbb{A}^1$ is a degeneration of $X^0$ and $\mathfrak{X}\rightarrow\mathbb{A}^1$ is a degeneration of $X$, but these are not toric degenerations.
\end{con}

\begin{figure}[h!]
\centering
\begin{tikzpicture}[scale=1.4]
\coordinate (1) at (0,0);
\coordinate[label=above:${v}$] (2) at (0,1);
\coordinate (3) at (-3,-1);
\coordinate (4) at (-3,2);
% cut
\fill[color=gray] (-3,-1) -- (0,0) -- (0,1) -- (-3,2);
\draw[thick] (3) -- (1) -- (2) -- (4);
\draw[thick] (1) -- (5,0);
\draw[thick] (2) -- (5,1);
\draw[thick,dashed] (3) -- (5,-1);
\draw[thick,dashed] (4) -- (5,2);
% directions
\draw[very thick,->] (0,1) -- (-0.6,1.2) node[above]{$m_{v,+}$};
\draw[very thick,->] (0,1) -- (0,0.5) node[right]{$m_{v,-}$};
\draw[very thick,->] (0,1) -- (0.5,1) node[above]{$m_{\textup{out}}$};
% red
\draw[very thick,red] (0,1) -- (-0.3,1.1);
\draw[very thick,red] (0,0) -- (0,0.316);
\draw[very thick,red] (-3,-1) -- (-2.7,-0.9);
\coordinate[fill,cross,thick,inner sep=2pt] (5) at (0,0);
\coordinate[fill,cross,thick,inner sep=2pt] (6) at (0,1);
\coordinate[fill,cross,thick,inner sep=2pt,rotate=26.57] (7) at (-3,-1);
\coordinate[fill,cross,thick,inner sep=2pt,rotate=-30] (8) at (-3,2);
\end{tikzpicture}
\caption{The dual intersection complex $(\tilde{B},\mathscr{P},\varphi)$ of $\tilde{X}_0$ for $(\mathbb{P}^2,E)$ away from $\sigma_0$, with choices of resolutions indicated.}
\label{fig:duald}
\end{figure}
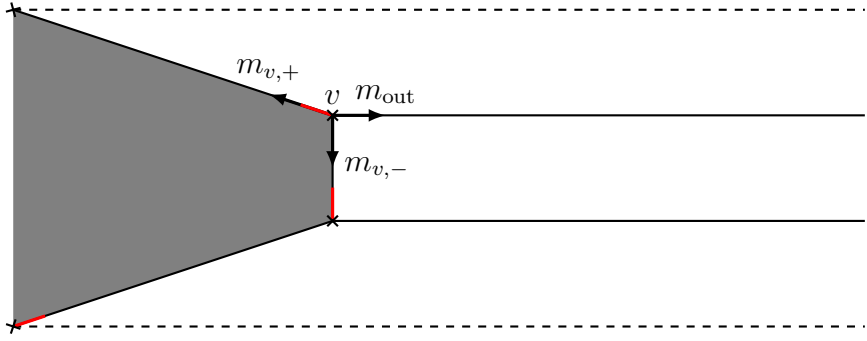

The small resolution does not change the local toric models at generic points of toric strata. As a consequence, the dual intersection complex $\tilde{B}$ of $\tilde{\mathfrak{X}}\rightarrow\mathbb{A}^1$ is homeomorphic to the dual intersection complex $B$ of $\mathfrak{X}\rightarrow\mathbb{A}^1$. But there is one difference here. The irreducible components of $\tilde{X}_0$ are non-toric, so there is no natural fan structure at the vertices. Further, $\tilde{X}_0$ has no log singularities. Hence, there is no focus-focus singularity on bounded edges in the dual intersection complex. However, the gluing is still in such a way that the unbounded edges are parallel, leading to affine singularities at the vertices of $\tilde{B}$, coming from the gluing. This gives a triple $(\tilde{B},\mathscr{P},\varphi)$ as in Figure \ref{fig:duald}. Affine manifolds with singularities at vertices have been considered by Gross-Hacking-Keel \cite{GHK1} to construct mirrors to log Calabi-Yau surfaces.

\begin{defi}
\label{defi:Lexc}
For a vertex $v$ of $\mathscr{P}$, let $L^{\textup{exc}}_v$ be the unique exceptional line contained in the irreducible component $\tilde{X}_v$ of $\tilde{X}_0$ corresponding to $v$.
\end{defi}

For later convenience we indicate the choices of small resolutions by red stubs attached to the vertices of $\mathscr{P}$. The stub at a vertex $v$ points in the direction corresponding to the toric divisor of $X_v$ intersecting $L^{\textup{exc}}_v$. Denote the primitive vector in the direction of the red stub adjacent to $v$ by $m_{v,+}\in\Lambda_{\tilde{B},v}$. Denote the primitive vector in the direction of the other edge of $\sigma_0$ adjacent to $v$ by $m_{v,-}\in\Lambda_{\tilde{B},v}$. Further, $m_{\textup{out}}$ is the unique unbounded direction (Definition \ref{defi:mout}).

\subsection{Logarithmic Gromov-Witten invariants}				%%

Logarithmic Gromov-Witten invariants have been defined in \cite{Che1}\cite{AC} and \cite{LogGW} as counts of stable log maps. A stable log map is a stable map defined in the category of log schemes with additional logarithmic data at the marked points, allowing for specification of contact orders. This leads to a generalization of Gromov-Witten theory in log smooth situations. For example, Gromov-Witten invariants relative to a (log) smooth divisor can be defined in this context, avoiding the target expansion of relative Gromov-Witten theory \cite{Li}\cite{Li2}. This is the case of interest to us. 

Let $\tilde{\mathfrak{X}}:=\mathfrak{X}^{s\neq 0}\rightarrow\mathbb{A}^1$ be the log smooth family from Construction \ref{con:deg2}. Note that $\tilde{X}_{t\neq 0}=X$. For the definition of stable log maps and their classes, see \cite{LogGW}.

The group of curve (= divisor) classes on $X$ is isomorphic to the singular homology group $H_2(X,\mathbb{Z})$ by Poincar\'e duality and since $H^1(X,\mathcal{O}_X)=0$ for del Pezzo surfaces by the Kodaira vanishing theorem. We write $H_2^+(X,\mathbb{Z})$ for the monoid of effective curve classes.

\begin{defi}
\label{defi:beta}
For an effective curve class $\underline{\beta}\in H_2^+(X,\mathbb{Z})\simeq H_2^+(X^0,\mathbb{Z})$ define a class $\beta$ of stable log maps to $\tilde{\mathfrak{X}}_Q \rightarrow \mathbb{A}^2$ as follows:
\begin{compactenum}[(1)]
\item genus $g=0$;
\item $k=1$ marked point $p$;
\item fibers have curve class $\underline{\beta}$;
\item contact data $u_p = (D \cdot \underline{\beta}) m_{\textup{out}}$, that is, full tangency with $D$ at the marked point. Here $m_{\textup{out}}\in\Lambda_{\tilde{B}}$ is the primitive integral tangent vector pointing in the unbounded direction on $\tilde{B}$ (Definition \ref{defi:mout}).
\end{compactenum}
A choice of $s\in\mathbb{A}^1$ gives an embedding $\gamma : \mathbb{A}^1 \rightarrow \mathbb{A}^2$. Let $\gamma^!$ be the corresponding refined Gysin homomorphism (\cite{Ful}, {\S}6.6). Then $\gamma^!\beta$ defines a class of stable log maps to $\mathfrak{X}^s\rightarrow\mathbb{A}^1$ that, by abuse of notation, we also write as $\beta$.
\end{defi}

\begin{rem}
One comment is in order about the space in which $u_p$ lives. By definition (\cite{LogGW}, Discussion 1.8(ii)), $u_p$ is an element of $P_p^\vee:=\textup{Hom}(f^\star\overline{\mathcal{M}}_{\tilde{\mathfrak{X}}_Q}|_p,\mathbb{N})$ and $m_{\textup{out}}$ is an element of $\Lambda_{\tilde{B}}$, the sheaf of integral tangent vectors on the dual intersection complex $\tilde{B}$. Let $\omega$ be the cell of $\mathscr{P}$ corresponding to the minimal stratum of $\tilde{X}_0^s$ to which the marked point is mapped. Then $\omega$ is an unbounded $1$- or $2$-dimensional cell and $m_{\textup{out}}$ defines an element of $\Lambda_{\tilde{B},\omega}$. Both, $P_p^\vee$ and $\Lambda_{\tilde{B},\omega}$ are submonoids of $\Lambda_{\Sigma(\tilde{X}_0^s),\omega}$ and their intersection is $\mathbb{N}\cdot m_{\textup{out}}\subseteq\Lambda_{\tilde{B},\omega}$. Thus $\mathbb{N}\cdot m_{\textup{out}}$ can be viewed as a submonoid of $P_p^\vee$, so the definition above makes sense.
\end{rem}

\begin{defi}
Let $\mathscr{M}(\tilde{\mathfrak{X}},\beta)$ be the moduli space of basic stable log maps to $\tilde{\mathfrak{X}}:=\mathfrak{X}^{s\neq 0}\rightarrow\mathbb{A}^1$ of class $\beta$. 
\end{defi}

By \cite{LogGW}, Theorems 0.2 and 0.3, $\mathscr{M}(\tilde{\mathfrak{X}},\beta)$ is a proper Deligne-Mumford stack and admits a virtual fundamental class $\llbracket\mathscr{M}(\tilde{\mathfrak{X}},\beta)\rrbracket$. Since $(X,D)$ is a log Calabi-Yau pair, the class $\beta$ is combinatorially finite (\cite{LogGW}, Definition 3.3). Hence, the virtual dimension of $\mathscr{M}(\tilde{\mathfrak{X}},\beta)$ is zero and the following definition makes sense.

\begin{defi}
\label{defi:Nb}
For $\beta$ as in Definition \ref{defi:beta} define the logarithmic Gromov-Witten invariant
\[ N_\beta = \int_{\llbracket\mathscr{M}(\tilde{\mathfrak{X}},\beta)\rrbracket} 1. \]
\end{defi}

\begin{defi}
\label{defi:Nd}
Define $w_{\textup{out}} = \textup{min}\{D\cdot\underline{\beta} \ | \ \underline{\beta}\in H_2^+(X,\mathbb{Z}) \}$ so e.g. for $(X,D)=(\mathbb{P}^2,E)$ we have $w_{\text{out}}=3$, since $E$ has degree $3$. For $d>0$ define
\[ N_d = \sum_{\substack{\underline{\beta} \in H_2^+(X,\mathbb{Z}) \\ D\cdot\underline{\beta}=dw_{\textup{out}}}} N_\beta. \]
\end{defi}

\begin{rem}
\label{rem:constant}
Logarithmic Gromov-Witten invariants are constant in log smooth families (\cite{MR}, Appendix A). This means the following. Let $\gamma : \{\text{pt}\} \rightarrow \mathbb{A}^1$ be a point and let $\gamma^!$ be the corresponding refined Gysin homomorphism . Then $\gamma^!\beta$ defines a class of stable log maps to the fiber $\tilde{X}_t$ that, by abuse of notation, we also write as $\beta$. We get a moduli space and a virtual fundamental class $\llbracket\mathscr{M}(\tilde{X}_t,\beta)\rrbracket$. Then, for all $t\in\mathbb{A}^1$,
\[ N_\beta = \int_{\llbracket\mathscr{M}(\tilde{X}_t,\beta)\rrbracket}1. \]
This shows that $N_\beta$ equals the logarithmic Gromov-Witten invariant $N_\beta$ defined in the introduction. Moreover, as shown in \cite{AMW}, $N_\beta$ equals the relative Gromov-Witten invariant of the smooth pair $(X,D)$ as defined in \cite{Li}.
\end{rem}

\subsection{Log BPS numbers}							%%
\label{S:BPS}

The logarithmic Gromov-Witten invariants $N_\beta$ are not integers but rather rational numbers. The fractional part comes from multiple cover contributions of curves with class $\beta'$ such that $\underline{\beta}=k\cdot\underline{\beta}'$. 

\begin{prop}[\cite{GPS}, Proposition 6.1]
The $k$-fold cover of an irreducible curve of class $\beta'$ contributes the following factor to $N_{k\cdot\beta'}$:
\[ M_{\beta'}[k] = \frac{1}{k^2}\binom{k(D\cdot\underline{\beta}'-1)-1}{k-1} \]
\end{prop}

We use the same formula for reducible curves, though it may be unclear how to interpret this as a multiple cover contribution. 

\begin{defi}
Define numbers $n_\beta$ by subtracting multiple cover contributions:
\[ N_\beta = \sum_{\beta' : \underline{\beta}=k\cdot\underline{\beta}'} M_{\beta'}[k] \cdot n_{\beta'}. \]
They are called \textit{Gopakumar-Vafa invariants} or \textit{log BPS numbers} as they are related to BPS state counts in string theory \cite{GV}.
\end{defi}

\begin{rem}
\label{rem:local}
The logarithmic Gromov-Witten invariants $N_\beta$ are related to local Gromov-Witten invariants $N_\beta^{\textup{loc}}$ of the total space of the canonical bundle $K_X$ of $X$ by the formula $N_\beta = (-1)^{D\cdot\underline{\beta}-1}(D\cdot\underline{\beta})N_\beta^{\textup{loc}}$. This was conjectured by Takahashi (\cite{Ta2}, Remark 1.11) and proved by Gathmann (\cite{Ga}, Example 2.2) and more generally by van Garrel, Graber and Ruddat \cite{vGGR}. The log BPS numbers $n_d$ were shown to be integers in \cite{vGWZ}, using integrality of local BPS numbers. 
\end{rem}

\section{Tropical curves and refinement}					%%%
\label{S:tropmap}

In this section we analyze what tropicalizations of stable log maps contributing to $N_d$ look like. We prove that for each $d$ there are only finitely many such tropical curves (Corollary \ref{cor:finite}). Choosing a subdivision of the dual intersection complex $(\tilde{B},\mathscr{P},\varphi)$ such that tropicalizations contributing to $N_d$ are contained in the $1$-skeleton of the polyhedral decomposition leads to a logarithmic modification $\tilde{\mathfrak{X}}_d$ of $\tilde{\mathfrak{X}}$ (Construction \ref{con:deg3}) with the property that stable log maps to the central fiber $Y$ of $\tilde{\mathfrak{X}}_d$ contributing to $N_d$ are torically transverse.

\subsection{Tropicalization of stable log maps}					%%

\begin{defi}[\cite{ACGS1}, 2.1.1, 2.1.2]
\label{defi:Cones}
Define $\textbf{Cones}$ to be the category whose objects are pairs $(\sigma_{\mathbb{R}},M)$ where $M\cong\mathbb{Z}^n$ is a lattice and $\sigma_{\mathbb{R}}\subseteq M_{\mathbb{R}}=M\otimes_{\mathbb{Z}}\mathbb{R}$ is a top-dimensional strictly convex rational polyhedral cone. A morphism of cones $\varphi : \sigma_1\rightarrow \sigma_2$ is a homomorphism $\varphi : M_1 \rightarrow M_2$ which takes $\sigma_{1\mathbb{R}}$ into $\sigma_{2\mathbb{R}}$. It is a \textit{face morphism} if it identifies $\sigma_{1\mathbb{R}}$ with a face of $\sigma_{2\mathbb{R}}$ and $M_1$ with a saturated sublattice of $M_2$. A \textit{generalized cone complex} is a topological space with a presentation as the colimit of an arbitrary finite diagram in the category $\textbf{Cones}$ with all morphisms being face morphisms.
\end{defi}

\begin{defi}[\cite{ACGS1}, 2.1.4]
\label{defi:trop}
Let $X$ be a fine saturated log scheme with log structure defined in the Zariski topology. For the generic point $\eta$ of a stratum of $X$, its characteristic monoid $\overline{\mathcal{M}}_{X,\eta}$ defines a dual monoid $(\overline{\mathcal{M}}_{X,\eta})^\vee := \textup{Hom}(\overline{\mathcal{M}}_{X,\eta},\mathbb{N})$ lying in the group $(\overline{\mathcal{M}}_{X,\eta})^\star := \textup{Hom}(\overline{\mathcal{M}}_{X,\eta},\mathbb{Z})$, hence a dual cone
\[ \sigma_\eta := \left((\overline{\mathcal{M}}_{X,\eta})_{\mathbb{R}}^\vee,(\overline{\mathcal{M}}_{X,\eta})^\star\right). \]
If $\eta$ is specialization of $\eta'$, there is a well-defined generization map $\overline{\mathcal{M}}_{X,\eta}\rightarrow \overline{\mathcal{M}}_{X,\eta'}$. Dualizing, we obtain a face morphism $\sigma_{\eta'}\rightarrow\sigma_\eta$. This gives a diagram of cones indexed by strata of $X$ with face morphisms, hence gives a generalized cone complex $\Sigma(X)$, the \textit{tropicalization} of $X$. This construction is functorial.
\end{defi}

Let $Q$ be a Fano polytope and let $\tilde{\mathfrak{X}}\rightarrow\mathbb{A}^1$ be the log smooth degeneration of the corresponding smooth very ample log Calabi-Yau pair $(X,D)$  from Construction \ref{con:deg2}. Let $\underline{\beta}\in H_2^+(X,\mathbb{Z})$ be an effective curve class and consider a basic stable log map $f:C/\textup{pt}_{Q_{\text{basic}}}\rightarrow\tilde{X}_0/\textup{pt}_{\mathbb{N}}$ of class $\beta$ (Definition \ref{defi:beta}). Here $Q_{\text{basic}}$ is the basic monoid\footnote{The basic monoid $Q$ has the property that $\Sigma(\textup{pt}_{Q_{\text{basic}}})=\textup{Hom}(Q_{\text{basic}},\mathbb{R}_{\geq 0})$ is the moduli space of deformations of $\Sigma(C)$ as a tropical curve preserving its combinatorial type (\cite{LogGW}, Remark 1.21).} of $f$ (\cite{LogGW}, {\S}1.5). We will see in Corollary \ref{cor:finite} that in our situation $Q_{\text{basic}}=\mathbb{N}$.
\begin{equation}
\label{eq:stablelog}
\begin{xy}
\xymatrix{
C \ar[r]^f \ar[d]^\gamma					& \tilde{X}_0 \ar[d]^{\tilde{\pi}_0} \\
\textup{pt}_{Q_{\text{basic}}} \ar[r]^g 		& \textup{pt}_{\mathbb{N}}
}
\end{xy}
\end{equation}
Tropicalization of \eqref{eq:stablelog} gives a diagram of generalized cone complexes. Note that $\Sigma(\textup{pt}_{\mathbb{N}})=\mathbb{R}_{\geq 0}$. The fiber $\Sigma(\tilde{\pi}_0)^{-1}(1)$ is homeomorphic to the dual intersection complex $\tilde{B}$ of $\tilde{X}_0$. Similarly, for a general element $b$ of the cone $\Sigma(\text{pt}_{Q_{\text{basic}}})$ the fiber $\Sigma(\gamma)^{-1}(b)$ is homeomorphic to the dual intersection graph $\Gamma_C$ of $C$. Hence, tropicalization of \eqref{eq:stablelog} and restriction to the fiber over $1\in\Sigma(\textup{pt}_{\mathbb{N}})=\mathbb{R}_{\geq 0}$ gives a map
\begin{equation}
\label{eq:Delta}
\tilde{h} : \Gamma_C \rightarrow \tilde{B}.
\end{equation}

There is additional data on $\Gamma_C$ making $\tilde{h} : \Gamma_C\rightarrow \tilde{B}$ into a tropical curve in the sense of \cite{ACGS1}. Note that such a tropical curve only fulfills a modified version of the balancing condition (\cite{LogGW}, Proposition 1.15). In \S\ref{S:balancing} we will see what this means in our case. 

To make the connection with scattering diagrams in \S\ref{S:scattering} it is useful to consider tropical curves on $B$ (not $\tilde{B}$) that are balanced in the usual sense but may have some bounded legs.

There are many slightly different definitions of parametrized tropical curves, depending on the context in which they are used. The following definition is a synthesis of the definition in \cite{ACGS1} and \cite{Gr10}, Definition 1.32. In \cite{ACGS1} only tropical curves with no bounded legs are considered, while in \cite{Gr10} tropical curves are required to be balanced.

\begin{defi}
\label{defi:tropical}
Let $B$ be a $2$-dimensional integral affine manifold with singularities. Let $\Delta\subset B$ be the discriminant locus and write $B_0:=B\setminus\Delta$. A \textit{(parametrized) tropical curve on $B$}, written $h : \Gamma \rightarrow B$, is a homogeneous map $h : \Gamma \rightarrow B_0$, where $\Gamma$ is the topological realization of a graph\footnote{The topological realization of a graph $\Gamma$ is a topological space which is the union of line segments corresponding to the edges. By abuse of notation, we also denote this by $\Gamma$. Whenever we talk about a map from a graph we mean a homogeneous map from its topological realization.}, possibly with some non-compact edges (\textit{legs}), together with
\begin{compactenum}[(1)]
\item a non-negative integer $g_V$ (\textit{genus}) for each vertex $V$;
\item a non-negative integer $\ell_E$ (\textit{length}) for each compact edge $E$;
\item an element $u_{(V,E)}\in i_\star\Lambda_{B_0,h(V)}$ (\textit{weight vector}) for every vertex $V$ and edge or leg $E$ adjacent to $V$. Here $\Lambda_{B_0}$ is the sheaf of integral affine tangent vectors on $B$ and $i : B_0 \hookrightarrow B$ is the inclusion. The index of $u_{(V,E)}$ in the lattice $\Lambda_{B,h(V)}$ is called the \textit{weight} $w_E$ of $E$;
\end{compactenum}
such that
\begin{compactenum}[(i)]
\item if $E$ is a compact edge with vertices $V_1$, $V_2$, then $h$ maps $E$ affine linearly\footnote{The affine structure on $\Gamma$ is given by the lengths $\ell_E$ of its edges.} to the line segment connecting $h(V_1)$ and $h(V_2)$, and $h(V_2)-h(V_1)=\ell_Eu_{(V_1,E)}$. In particular, $u_{(V_1,E)} = -u_{(V_2,E)}$;
\item if $E$ is a leg with vertex $V$, then $h$ maps $E$ affine linearly either to the ray $h(V)+\mathbb{R}_{\geq 0}u_{(V,E)}$ or to the line segment $[h(V),\delta)$ for $\delta$ an affine singularity of $B$ such that $\delta-h(V)\in\mathbb{R}_{>0} u_{(V,E)}$, i.e., $u_{(V,E)}$ points from $h(V)$ to $\delta$.
\end{compactenum}
We write the set of compact edges of $\Gamma$ as $E(\Gamma)$, the set of legs as $L(\Gamma)$, the set of legs mapped to a ray (\textit{unbounded legs}) as $L_\infty(\Gamma)$ and the set of legs mapped to an open line segment (\textit{bounded legs}) as $L_\Delta(\Gamma)$ (since such edges end at the singular locus $\Delta$ of $B$).

The \textit{genus} of a parametrized tropical curve $h : \Gamma \rightarrow B$ is defined by
\[ g_h := g_\Gamma + \sum_{V\in V(\Gamma)}g(V), \]
where $g_\Gamma$ is the genus (first Betti number) of the graph $\Gamma$.
\end{defi}

\begin{rem}
Note that if $E$ is a leg of $\Gamma$, then $h(E)$ must be parallel to the edge of $\mathscr{P}$ containing $\delta$, since there is only one tangent direction at $\delta$, i.e., $\Lambda_{B,\delta}\simeq\mathbb{Z}$.
\end{rem}

\begin{defi}
\label{defi:aut}
An \textit{isomorphism} of tropical curves $h:\Gamma\rightarrow B$ and $h':\Gamma'\rightarrow B$ is a homeomorphism $\phi : \Gamma \rightarrow \Gamma'$ such that $h=h'\circ\phi$, $g_{\phi(V)}=g_V$ and $u_{(\phi(V),\phi(E))}=u_{(V,E)}$. An \textit{automorphism} of a tropical curve $h$ is an isomorphism of $h$ with itself. Here we use the convention that an edge $E$ is a pair of orientations of $E$, so that the automorphism group of a graph with a single loop is $\mathbb{Z}/2\mathbb{Z}$.
\end{defi}

\begin{rem}
We will only consider tropical curves of genus $0$. In particular, our tropical curves will have no loops.
\end{rem}

\begin{defi}
Let $(B,\mathscr{P})$ be a $2$-dimensional polyhedral affine manifold. A tropical curve $h : \Gamma\rightarrow B$ is \textit{compatible with $\mathscr{P}$} if
\begin{compactenum}[(1)]
\item the edges of $\Gamma$ do not extend across several maximal cells of $\mathscr{P}$. In other words, we have a well-defined map $E(\Gamma) \cup L(\Gamma) \rightarrow \mathscr{P}$ associating to an edge or leg $E$ the minimal cell of $\mathscr{P}$ containing it.
\item there are no bivalent vertices in $\Gamma$ mapped to a maximal cell of $\mathscr{P}$.
\end{compactenum}
\end{defi}

\begin{con}
\label{con:tropical}
Let $\tilde{h} : \Gamma_C \rightarrow \tilde{B}$ be the continuous map from \eqref{eq:Delta}. We describe additional data making $\tilde{h}$ a tropical curve compatible with $\mathscr{P}$.
\begin{compactenum}[(1)]
\item For each vertex $V$, the genus is $g_V=0$.
\item Let $E\in E(\Gamma_C)$ be a compact edge with vertices $V_1,V_2$, corresponding to a node $q\in C$. Then $\overline{\mathcal{M}}_{C,q}$ is isomorphic to the submonoid $S_{e_q}$ of $\mathbb{N}^2$ generated by $(e_q,0)$, $(0,e_q)$ and $(1,1)$ for some $e_q\in\mathbb{N}_{>0}$ (\cite{Kat2}, 1.8). Moreover, there is an equation $\tilde{h}(V_2)-\tilde{h}(V_1)=\pm e_qu_q$ for some $u_q\in\Lambda_{\tilde{B}}$ (see \cite{LogGW}, Discussions 1.8, 1.13). Then the length of $E$ is $\ell_E=e_q$ and the weight vectors are $u_{(V_i,E)}=\pm u_q$, with sign chosen such that $u_{(V_i,E)}$ points away from $\tilde{h}(V_i)$.
\item Let $E\in L_\infty(\Gamma_C)$ be an unbounded leg with vertex $V$, corresponding to a marked point $p\in C$. Then $\overline{\mathcal{M}}_{C,p}$ is isomorphic to $\mathbb{N}\oplus\mathbb{N}$ and 
\[ f^\star\overline{\mathcal{M}}_{X}|_p\rightarrow\overline{\mathcal{M}}_{C,p}\overset{\text{pr}_2}{\rightarrow}\mathbb{N} \]
is determined by an element of $P_p^\vee=\textup{Hom}(f^\star\overline{\mathcal{M}}_{X}|_p,\mathbb{N})$, inducing an element $u_p\in\Lambda_{\tilde{B},\tilde{h}(V)}$. The weight vector is $u_{(V,E)}=u_p$.
\end{compactenum}
The properties for $\tilde{h} : \Gamma_C \rightarrow \tilde{B}$ to be compatible with $\mathscr{P}$ can be achieved by (1) inserting vertices at points mapping to vertices or edges of $\mathscr{P}$ and (2) removing bivalent vertices mapping to a maximal cell of $\mathscr{P}$, by replacing a chain of edges connected via bivalent vertices with a single edge. The latter is possible, since vertices of $\Gamma_C$ not mapping to vertices of $\mathscr{P}$ are balanced by Proposition \ref{prop:balancing}, (I), below.
\end{con}

\begin{defi}
Let $B$ be a $2$-dimensional integral affine manifold with singularities and $m\in\Lambda_B$ an integral tangent vector. A tropical curve $h:\Gamma\rightarrow B$ is called \textit{of degree $d$ relative to $m$} if there is exactly one unbounded leg $E_{\textup{out}}\in L_\infty(\Gamma)$, and its weight vector is $u_{(V_{\textup{out}},E_{\textup{out}})}=d \cdot m$. Here $V_{\textup{out}}$ is the unique vertex of $E_{\textup{out}}$.
\end{defi}

\begin{prop}
\label{prop:tropical}
Given a stable log map $f:C/\textup{pt}_{Q_{\text{basic}}}\rightarrow \tilde{X}_0/\textup{pt}_{\mathbb{N}}$ of class $\beta$ as in Definition \ref{defi:beta}, the continuous map $\tilde{h}:\Gamma_C\rightarrow \tilde{B}$ from \eqref{eq:Delta} together with the additional data defined in Construction \ref{con:tropical} is a tropical curve without bounded legs, of genus $0$, degree $D\cdot\underline{\beta}$ relative to $m_{\textup{out}}\in\Lambda_{\tilde{B}}$ (Definition \ref{defi:mout}) and compatible with the dual intersection complex $\mathscr{P}$. By abuse of notation, we call this tropical curve $\tilde{h}:\Gamma_C\rightarrow \tilde{B}$ the \textit{tropicalization} of $f$.
\end{prop}

\begin{proof}
The properties (i) and (ii) of Definition \ref{defi:tropical} follow by the structure of $f : C \rightarrow \tilde{X}_0$ on the level of ghost sheaves (see \cite{LogGW}, Discussions 1.8, 1.13). Hence, $\tilde{h}$ is a tropical curve. Moreover, these discussions show that $\tilde{h}$ has no bounded legs. $\tilde{h}$ is of degree $dw_{\textup{out}}$ relative to $m_{\textup{out}}$ by Definition \ref{defi:beta}, (4).
\end{proof}

\begin{rem}
There is one issue here, since we lost some information by smoothing $(X^0,D^0)$. An effective curve class $\beta\in H_2^+(X^0,\mathbb{Z})$ is determined by its intersection numbers $d_1,\ldots,d_k$ with the components of $D^0=D_1,\ldots,D_k$. After smoothing $D^0$ we only see the sum $d=d_1+\ldots+d_k$. In particular, if $X^0$ is a smooth toric del Pezzo surface with Picard number $>1$, i.e., different from $(\mathbb{P}^2,E)$, then we only see the total degree, not the multi-degree. One could solve this problem using non-trivial gluing data capturing information of the divisor $D^0$. We will give a more geometric solution in \S\ref{S:degdiv} by looking at the limit of curves under $s\rightarrow 0$, where $s\in\mathbb{A}^2$ is the deformation parameter of the family $\mathfrak{X}_Q\rightarrow\mathbb{A}^2$ from Construction \ref{con:family}.
\end{rem}

\subsection{Types of vertices}								%%
\label{S:balancing}

Let $f:C/\textup{pt}_{Q_{\text{basic}}}\rightarrow\tilde{X}_0/\textup{pt}_{\mathbb{N}}$ be a stable log map of class $\beta$ as in Definition \ref{defi:beta} and let $\tilde{h}:\Gamma_C\rightarrow \tilde{B}$ be the corresponding tropical curve.

\begin{prop}
\label{prop:balancing}
Let $C_V$ be an irreducible component of $C$, corresponding to a vertex $V$ of $\Gamma_C$. Then the following three cases can occur:
\begin{compactenum}[(I)]
\item If $C_V$  is mapped to a $0$- or $1$-dimensional toric stratum of $\tilde{X}_0$, i.e., if $V$ is not mapped to a vertex of $\mathscr{P}$, then the ordinary balancing condition holds:
\[ \sum_{E\ni V}u_{(V,E)} = 0. \]
The sum is over all edges or legs $E\in E(\Gamma_C)\cup L(\Gamma_C)$ containing $V$.
\item If $C_V$ is mapped onto an exceptional divisor $L^{\textup{exc}}_v$ on some component $\tilde{X}_v$ of $\tilde{X}_0$ (Definition \ref{defi:Lexc}), then $C_V$ is a $k$-fold multiple cover of $L^{\textup{exc}}_v\simeq\mathbb{P}^1$ for some $k>0$. It is fully ramified at the point $p=L^{\textup{exc}}_v\cap\partial\tilde{X}_v$, where $\partial\tilde{X}_v$ is the proper transform of the toric boundary $\partial X_v$ under the resolution from \S\ref{S:resolution}. The vertex $V$ is mapped to the vertex $v$ of $\mathscr{P}$. It is $1$-valent with adjacent edge $E$ mapped onto the edge of $\mathscr{P}$ containing the red stub adjacent to $v$. The balancing condition reads (with $m_{v,+}$ as in Figure \ref{fig:duald})
\[ u_{(V,E)} = km_{v,+}. \]
\item Otherwise, $V$ is mapped to a vertex $v$ of $\mathscr{P}$ and has exactly one adjacent edge or leg $E_{V,\textup{out}}$ that is not mapped onto a compact edge of $\mathscr{P}$. All other edges (possibly none) are compact with other vertex of type (II) above. In this case, for some $k\geq 0$, the following balancing condition holds:
\[ \sum_{E\ni V} u_{(V,E)} + km_{v,+} = 0. \]
\end{compactenum}
\end{prop}

\begin{proof}
If $C_V$ does not intersect an exceptional line, the log structure on $\tilde{X}_0$ along the image of $C_V$ is the toric one. Then by \cite{ACGS2}, Remark 2.26, the ordinary balancing condition holds. This proves (I).

If $C_V$ is mapped onto an exceptional line $L^{\textup{exc}}_v\simeq\mathbb{P}^1$ on some component $\tilde{X}_v$, it is a $k$-fold multiple cover for some $k>0$. Suppose it is not fully ramified at the point where $L^{\textup{exc}}_v$ meets $\partial \tilde{X}_v$, i.e., $V$ has valency $>1$. Let $E_1,E_2$ be two distinct edges adjacent to $V$. We have $V \neq V_{\textup{out}}$, since $C_V$ does not meet the toric divisor of $\tilde{X}_v$ belonging to $\tilde{D}_0$. By Proposition \ref{prop:tropical}, $\Gamma_C$ has only one leg, and this leg is attached to $V_{\textup{out}}$. Thus $E_1$ and $E_2$ are bounded. Let $V_1,V_2$ be the vertices of $E_1,E_2$ different from $V$, respectively. There is a chain of vertices and edges (possibly the trivial one) connecting $V_1$ to $V_{\textup{out}}$ and similarly for $V_2$. These two chains form a cycle of the graph $\Gamma_C$, so $g(\Gamma_C)>0$ in contradiction with rationality of $C$. Hence, there is a unique bounded edge $E$ adjacent to $V$. Let $V'$ be its other vertex. Then $h(V')-h(V)$ points in the direction of $m_{v,+}$, since the only special point (node) of $C_V$ is mapped to $L^{\textup{exc}}_v\cap\partial\tilde{X}_v$. It follows by Definition \ref{defi:tropical} that $u_{(V,E)}$ points in the direction of $m_{v,+}$. Its affine length, the weight $w_E$, is the multiplicity of the node which is the ramification order $k$. This proves (II).

Let $C_V$ be a component of $C$ that intersects an exceptional divisor $L^{\textup{exc}}_v$ on some component $\tilde{X}_v$ but is not mapped onto it. We first prove the balancing condition. Let $m\in\Gamma(C_V,f^\star\overline{\mathcal{M}}_X|_{C_V})$ be the generator of the submonoid $\mathbb{N}$ of $\Gamma(C_V,f^\star\overline{\mathcal{M}}_X|_{C_V})$ corresponding to the exponent of the degeneration parameter $t$. Then by \cite{LogGW}, Lemma 1.14, the line bundle associated to $m$ is $f^\star\nu^\star\mathcal{O}_{X_v}(-kD_{v,+})$, where $\nu : \tilde{X}_v \rightarrow X_v$ is the resolution and $D_{v,+}$ is the toric divisor of $X_v$ whose proper transform $\tilde{X}_v$ intersects $L_v^{\textup{exc}}$. This is because the condition that $C_V$ intersects $L_v^{\textup{exc}}$ is equivalent to the condition that $\nu(C_V)$ intersects $D_{v,+}$ in the point $\nu(L_v^{\text{exc}})$. The corresponding contact data is $km_{v,+}$. This shows $\sum_{E\ni V} u_{(V,E)} + km_{v,+} = 0$, where $k$ is the sum of the affine lengths of the additional $u_p$. Now we show uniqueness of an edge or leg $E_{V,\text{out}}$ as claimed. If all legs and edges adjacent to $V$ are mapped onto compact edges of $\mathscr{P}$, this balancing condition can not be achieved. So there is at least one such edge or leg. We show by contradiction that there is at most one. Assume that there are two edges or legs $E, E'$ adjacent to $V$ that are not mapped onto compact edges of $\mathscr{P}$. Then $E,E'$ are either unbounded legs or bounded edges with other vertex of type (I). If $E$ or $E'$ is unbounded, then, since vertices of type (I) fulfill the ordinary balancing condition, we would have at least two unbounded legs, contradicting the assumptions. If $E$ and $E'$ are bounded edges with other vertices of type (I), their paths to $V_{\text{out}}$ form a cycle, contradicting $g=0$. So there is a unique edge not mapped onto a compat edge of $\mathscr{P}$. This proves (III).
\end{proof}

\begin{defi}
\label{defi:types}
Denote the set of vertices of the given types in Proposition \ref{prop:balancing} by $V_{I}(\Gamma_C)$, $V_{II}(\Gamma_C)$ and $V_{III}(\Gamma_C)$, respectively.
\end{defi}

\begin{defi}
\label{defi:tildeH}
Let $\tilde{\mathfrak{H}}_d$ be the set of isomorphism classes of tropical curves $\tilde{h}:\tilde{\Gamma}\rightarrow \tilde{B}$ compatible with $\mathscr{P}$ of genus $0$ and degree $dw_{\textup{out}}$ relative to $m_{\textup{out}}\in\Lambda_{\tilde{B}}$, without bounded legs and with vertices of one of the types (I)-(III) above. 
\end{defi}

\begin{figure}[h!]
\centering
\begin{tikzpicture}[scale=1]
% SOURCE
\draw[very thick,blue] (-7,-1) node[fill,circle,inner sep=1.5pt]{}-- (-6,0) node[fill,circle,inner sep=1.5pt]{} -- (-5.8,0.67) node[fill,circle,inner sep=1.5pt]{} -- (-6,1) node[fill,circle,inner sep=1.5pt]{} -- (-6.5,0) node[fill,circle,inner sep=1.5pt]{};
\draw[very thick,blue] (-5.8,0.67) -- (-5.3,0.833) node[below]{$9$} -- (-4.8,1);
% ARROW
\draw[->,thick] (-4.5,0.5) -- (-3.5,0.5);
% TARGET
\coordinate (1) at (0,0);
\coordinate (2) at (0,1);
\coordinate (3) at (-3,-1);
\coordinate (4) at (-3,2);
% cut
\fill[color=gray] (-3,-1) -- (0,0) -- (0,1) -- (-3,2);
\draw[thick] (3) -- (1) -- (2) -- (4);
\draw[thick] (1) -- (5,0);
\draw[thick] (2) -- (5,1);
\draw[thick,dashed] (3) -- (5,-1);
\draw[thick,dashed] (4) -- (5,2);
% blue
\draw[very thick,blue] (0,0.05) node[fill,circle,inner sep=1.5pt,label=left:(II)]{} -- (0,1) node[fill,circle,inner sep=1.5pt,label=above:(III)]{} -- (2,0.67) node[fill,circle,inner sep=1.5pt,label=below:(I)]{} -- (0,-0) node[fill,circle,inner sep=1.5pt,label=below:(III)]{} -- (-3,-1) node[fill,circle,inner sep=1.5pt,label=below:(II)]{};
\draw[very thick,blue] (2,0.67) -- (3.5,0.67) node[below]{$9$} -- (5,0.67);
% red
\draw[very thick,red] (0,1) -- (-0.3,1.1);
\draw[very thick,red] (0,0) -- (0,0.316);
\draw[very thick,red] (-3,-1) -- (-2.7,-0.9);
% crosses
\coordinate[fill,cross,thick,inner sep=2pt] (5) at (0,0);
\coordinate[fill,cross,thick,inner sep=2pt] (6) at (0,1);
\coordinate[fill,cross,thick,inner sep=2pt,rotate=26.57] (7) at (-3,-1);
\coordinate[fill,cross,thick,inner sep=2pt,rotate=-30] (8) at (-3,2);
\end{tikzpicture}
\caption{A tropical curve $\tilde{h}:\tilde{\Gamma}\rightarrow \tilde{B}$ in $\tilde{\mathfrak{H}}_3$ for $(\mathbb{P}^2,E)$, with types of vertices indicated. The integer $9$ is the weight of the outgoing edge. All other edges have weight $1$. Two vertices are mapped to the same vertex of $\mathscr{P}$, but not connected by an edge.}
\label{fig:balancing}
\end{figure}
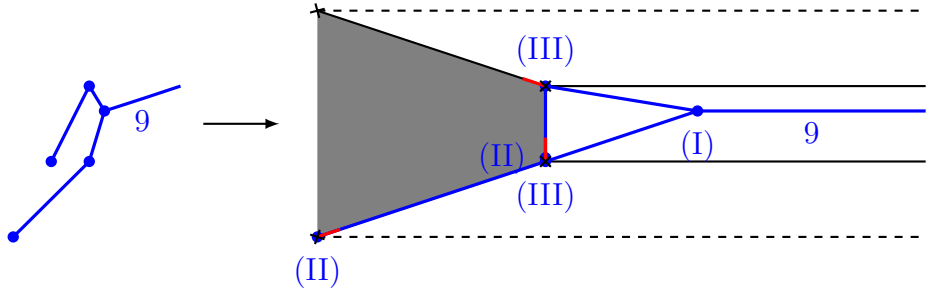

\begin{lem}
\label{lem:finite}
Let $\tilde{h}:\tilde{\Gamma}\rightarrow \tilde{B}$ be a tropical curve in $\tilde{\mathfrak{H}}_d$ for some $d>0$. Then $\tilde{h}(\tilde{\Gamma})$ is disjoint from the interior of $\sigma_0$ (Definition \ref{defi:sigma0}).
\end{lem}

\begin{proof}
Let $\tilde{h}:\tilde{\Gamma}\rightarrow \tilde{B}$ be a tropical curve in $\tilde{\mathfrak{H}}_d$. Give $\tilde{\Gamma}$ the structure of a rooted tree by defining the root vertex to be the vertex $V_{\textup{out}}$ of the unique unbounded leg $E_{\textup{out}}$. Let $V$ be a vertex of $\tilde{\Gamma}$ and let $E_{V,\textup{out}}$ be the edge connecting $V$ with its parent, or $E_{V,\textup{out}}=E_{\textup{out}}$ if $V$ is the root vertex $V_{\textup{out}}$. By Proposition \ref{prop:balancing}, if $V$ is mapped to a vertex $v$ of $\mathscr{P}$, hence of type (II) or (III), then $E_{V,\textup{out}}$ is mapped to the conical subset $\mathbb{R}_{\leq 0}m_{v,+}+\mathbb{R}_{\leq 0}m_{v,-}$ of $\tilde{B}$, and if $V$ is of type (I), then by induction $E_{V,\textup{out}}$ is mapped to the subset $\bigcup_{V'} \mathbb{R}_{\leq 0}m_{\tilde{h}(V'),+}+\mathbb{R}_{\leq 0}m_{\tilde{h}(V'),-}$ of $\tilde{B}$, where the union is over all vertices of type (II) and (III) in the subgraph of $\tilde{\Gamma}$ with root $V$. In particular, $\tilde{h}(\tilde{\Gamma})$ is contained in
\[\bigcup_{V'\in V_{II}(\tilde{\Gamma})\cup V_{III}(\tilde{\Gamma})} \mathbb{R}_{\leq 0}m_{\tilde{h}(V'),+}+\mathbb{R}_{\leq 0}m_{\tilde{h}(V'),-} \]
This is disjoint from the interior of $\sigma_0$.
\end{proof}

\subsection{Balanced tropical curves}							%%

We describe a procedure to obtain tropical curves in $\tilde{\mathfrak{H}}_d$ from tropical curves to $B$ (not $\tilde{B}!$) that are balanced in the usual sense. This makes the connection to scattering diagrams in \S\ref{S:scattering} more transparent. Moreover, the degeneration formula gets more symmetric when expressed in invariants labeled by balanced tropical curves (see Theorem \ref{thm:degmax}).

\begin{defi}
\label{defi:balancedtrop}
Let $\mathfrak{H}_d$ be the set of isomorphism classes of tropical curves $h:\Gamma\rightarrow B$ compatible with $\mathscr{P}$, possibly with bounded legs, of genus $0$ and degree $dw_{\textup{out}}$ relative to $m_{\textup{out}}$, satisfying the ordinary balancing condition at each vertex $V$ of $\Gamma$:
\[ \sum_{E\ni V}u_{(V,E)} = 0, \]
\end{defi}

\begin{con}
\label{con:tilde}
We construct a surjective map $\mathfrak{H}_d \rightarrow \tilde{\mathfrak{H}}_d$ as follows.

Let $h:\Gamma\rightarrow B$ be a tropical curve in $\mathfrak{H}_d$. Let $E\in L_\Delta(\Gamma)$ be a bounded leg with vertex $V$. Then $E$ is mapped to the line segment $[h(V),\delta)$ for $\delta$ an affine singularity on an edge $\omega$ of $\mathscr{P}$. Since $\Lambda_{B,\delta}$ is one-dimensional, $h(E)$ is parallel to $\omega$. Since $h$ is compatible with $\mathscr{P}$ and by the balancing condition, $h(V)$ must be a vertex $v$ of $\mathscr{P}$. Let $m_{v,\delta}$ be the primitive integral tangent vector pointing from $v$ to $\delta$ and let $m_{v,+},m_{v,-}$ be as in Figure \ref{fig:duald}. Two cases can occur.
\begin{compactenum}[(1)]
\item If $m_{v,\delta}=m_{v,+}$, i.e., if $E$ is mapped in the direction of the red stub attached to $v$, then remove $E$ from $\Gamma$.
\item Otherwise, $m_{v,\delta}=m_{v,-}$. Then add a vertex $V'$ to $E$ to obtain a compact edge $\tilde{E}$. Define $u_{(V',E)}=-u_{(V,E)}$ and $\tilde{h}(\tilde{E})=\omega$, such that $\tilde{h}(V')=v'$ is a vertex of $\mathscr{P}$. This determines the length $\ell_{\tilde{E}}$ by Definition \ref{defi:tropical}, (i).
\end{compactenum}
We show that the map $\mathfrak{H}_d\rightarrow\tilde{\mathfrak{H}}_d$ constructed this way is surjective. Let $\tilde{h} : \tilde{\Gamma} \rightarrow \tilde{B}$ be a tropical curve in $\tilde{\mathfrak{H}}_d$. We can construct a preimage of $\tilde{h}$ as follows. (1) For each vertex $V\in V_{III}(\tilde{\Gamma})$, add a bounded leg $E$ with vertex $V$ and weight vector $u_{(V,E)}=-\sum_{E'\ni V} u_{(V,E')}$. The image of $E$ is specified by Definition \ref{defi:tropical}, (ii). (2) For each vertex $V\in V_{II}(\tilde{\Gamma})$, let $E$ be the unique adjacent edge. It is a bounded edge and we remove the vertex $V$ from $E$ to obtain a bounded leg. This shows that the map $\mathfrak{H}_d\rightarrow\tilde{\mathfrak{H}}_d$ is surjective. Note that in step (1) we could also add several bounded legs with weights a partition of $\sum_{E'\ni V} u_{(V,E')}$, so the number of preimages of $\tilde{h}$ is the number of such partitions.
\end{con}

\begin{defi}
Let $(\bar{B},\bar{\mathscr{P}})$ be the covering space of $(B,\mathscr{P})$ described in \S\ref{S:affinecharts}. Let $\bar{\mathfrak{H}}_d$ be the set of isomorphism classes of balanced tropical curves $\bar{h}:\bar{\Gamma}\rightarrow\bar{B}$ compatible with $\bar{\mathscr{P}}$ of genus $0$ and degree $dw_{\textup{out}}$ relative to $m_{\textup{out}}$ satisfying the ordinary balancing condition and such that the image of $E_{\textup{out}}$ lies in a fixed fundamental domain.
\end{defi}

\begin{con}
\label{con:bij}
Define a map $\bar{\mathfrak{H}}_d \rightarrow \mathfrak{H}_d$ by sending $\bar{h} : \Gamma \rightarrow \bar{B}$ to $h : \Gamma \rightarrow B$, where $h$ is the composition of $\bar{h}$ with the covering map $\bar{B} \rightarrow B$. This map is bijective. The inverse map is given as follows. Let $h : B \rightarrow \Gamma$ be a tropical curve in $\mathfrak{H}_d$ and choose an unbounded maximal cell of $\mathscr{P}$. Choose a fundamental domain of $\bar{B} \rightarrow B$ and let $\bar{h}(E_{\textup{out}})$ be the preimage of $h(E_{\textup{out}})$ in that fundamental domain. Whenever the image $h(V)$ of an edge $V$ lies on the horizontal dashed line in Figure \ref{fig:dualb} with respect to the chart on the unbounded maximal cell chosen, we change the fundamental domain and apply the monodromy transformation.
\end{con}

\begin{figure}[h!]
\centering
\begin{tikzpicture}[scale=1]
% SOURCE
\draw[very thick,blue] (-6.25,-0.5) -- (-6,0) node[fill,circle,inner sep=1.5pt]{} -- (-5.8,0.67) node[fill,circle,inner sep=1.5pt]{} -- (-6,1) node[fill,circle,inner sep=1.5pt]{} -- (-6.25,0.5);
\draw[very thick,blue] (-6,1) -- (-6.25,1.5) node[right]{$2$};
\draw[very thick,blue] (-5.8,0.67) -- (-5.3,0.833) node[below]{$9$} -- (-4.8,1);
% ARROW
\draw[->,thick] (-4.5,0.5) -- (-3.5,0.5);
% TARGET
\coordinate (1) at (0,0);
\coordinate (2) at (0,1);
\coordinate (3) at (-3,-1);
\coordinate (4) at (-3,2);
\draw[thick] (3) -- (1) -- (2) -- (4);
\draw[thick] (1) -- (5,0);
\draw[thick] (2) -- (5,1);
\draw[thick,dashed] (3) -- (5,-1);
\draw[thick,dashed] (4) -- (5,2);
% cut
\draw[thick,dashed] (-1.5,-0.25) -- (0,0.5) -- (-1.5,1.25);
\draw[thick,dashed] (-3,-0.8) -- (-1.5,-0.5) -- (-0.75,0.25);
\draw[thick,dashed] (-3,1.8) -- (-1.5,1.5) -- (-0.75,0.75);
\fill[color=gray] (-1.5,-0.25) -- (0,0.5) -- (-1.5,1.25) -- (-1.5,-0.25) -- cycle;
\fill[color=gray] (-3,1.8) -- (-1.5,1.5) -- (-0.75,0.75) -- (-0.75,0.25) -- (-1.5,-0.5) -- (-3,-0.8);
\coordinate[fill,cross,thick,inner sep=2pt] (5) at (0,0.5);
\coordinate[fill,cross,thick,inner sep=2pt,rotate=26.57] (6) at (-1.5,-0.5);
\coordinate[fill,cross,thick,inner sep=2pt,rotate=-30] (7) at (-1.5,1.5);
% refinement
\draw[thick] (0,0) -- (-1,0);
\draw[thick] (0,1) -- (-1,1);
% blue
\fill[color=blue] (0,0) circle (1.5pt);
\fill[color=blue] (0,1) circle (1.5pt);
\fill[color=blue] (2,0.66) circle (1.5pt);
\draw[very thick,blue] (-1.5,1.5) -- (-0.75,1.25) node[above]{$2$} -- (0,1) -- (2,0.66) -- (0,0) -- (-1.5,-0.5);
\draw[very thick,blue] (2,0.66) -- (3.5,0.66) node[above]{$9$} -- (5,0.66);
\draw[very thick,blue] (0,1) -- (0,0.5);
\end{tikzpicture}
\caption{A balanced tropical curve $h:\Gamma\rightarrow B$ in $\mathfrak{H}_3$ for $(\mathbb{P}^2,E)$ giving the tropical curve in Figure \ref{fig:balancing} under the map from Construction \ref{con:tilde}. The integers are weights of edges $\neq 1$.}
\label{fig:balancing2}
\end{figure}
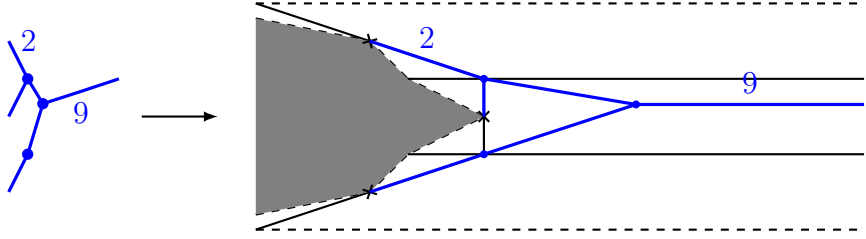

\begin{lem}
\label{lem:finite2}
The set $\bar{\mathfrak{H}}_d$ is finite.
\end{lem}

\begin{proof}
Let $\bar{h}:\bar{\Gamma}\rightarrow \bar{B}$ be a tropical curve in $\bar{\mathfrak{H}}_d$. The graph $\bar{\Gamma}$ together with the set of weight vectors of bounded legs $\{u_{(V,E)} \ | \ V\in E\in L_\Delta(\bar{\Gamma})\}$ determines the image of $\bar{h}$. Indeed, for a bounded leg $E$, the weight vector $u_{(V,E)}$ determines its image, since all edges containing affine singularities have different direction. The images of all other edges are determined by the balancing condition. If we only know the set $\{u_{(V,E)} \ | \ V\in E\in L_\Delta(\bar{\Gamma})\}$, there are finitely many possibilities for $\bar{\Gamma}$, since the number of leaves is specified. So we need to show that there are only finitely many possible sets $\{u_{(V,E)} \ | \ V\in E\in L_\Delta(\bar{\Gamma})\}$ for a tropical curve $\bar{h}:\bar{\Gamma}\rightarrow \bar{B}$ in $\bar{\mathfrak{H}}_d$.

Let $V_{\text{out}}$ be the vertex of the unique unbounded edge $E_{\text{out}}$ and let $\sigma_{\text{out}}$ be an unbounded maximal cell containing $V_{\text{out}}$. Let $\varphi_{\text{out}}$ be a representative of the piecewise affine function $\varphi$ on $\sigma_{\text{out}}$. Note that via an affine transformation of $\bar{B}$ we can achieve that $\varphi_{\text{out}}(m)=\braket{m,m_{\text{out}}}$. Let $E$ be a bounded leg of $\bar{\Gamma}$ with vertex $V$ and write $u_{(V,E)}=w_Em_E$ with $w_E\in\mathbb{Z}_{>0}$ and $m_E\in\Lambda_{\bar{h}(V)}\simeq \mathbb{Z}^2$ primitive. Let $E'$ be an edge on the path from $V$ to $V_{\text{out}}$ that is \textit{pointing towards} $V_{\text{out}}$, i.e., such that the ray $\bar{h}(V)+\mathbb{R}_{>0}m_E$ intersects the interior of the bounded maximal cell $\sigma_0$. Then $\varphi_{\text{out}}(m_{E'})>0$ and by convexity of the bounded maximal cell and of $\varphi_{\text{out}}$ we have $w_{E'}\varphi_{\text{out}}(m_{E'}) \geq w_E|\varphi_{\text{out}}(m_E)|$. In particular $w_{E_{\text{out}}}\varphi_{\text{out}}(m_{\text{out}})=dw_{\text{out}} \geq w_E|\varphi(m_E)|$. This gives a bound $w_E \leq \lfloor\frac{dw_{\text{out}}}{|\varphi_{\text{out}}(m_E)|}\rfloor$ and there are only finitely many bounded legs $E$ with $|\varphi_{\text{out}}(m_E)| \leq dw_{\text{out}}$, i.e., such that this bound is nonzero.
\end{proof}

\begin{cor}
\label{cor:finite}
The sets $\mathfrak{H}_d$ and $\tilde{\mathfrak{H}}_d$ are finite. In particular, tropical curves in $\tilde{\mathfrak{H}}_d$ are rigid and the basic monoid of stable log maps in $\tilde{\mathfrak{H}_d}$ is $Q_{\text{basic}}=\mathbb{N}$.
\end{cor}

\subsection{The limit $s\rightarrow 0$}						%%%
\label{S:degdiv}

Here we consider the $2$-parameter family $\tilde{\mathfrak{X}}_Q\rightarrow\mathbb{A}^2$ and describe limits of stable log maps to $\mathfrak{X}:=\mathfrak{X}^{s\neq 0}$ under $s\rightarrow 0$. This will enable us to read off the curve class $\underline{\beta}\in H_2^+(X,\mathbb{Z})$ of a stable log map from its tropicalization.

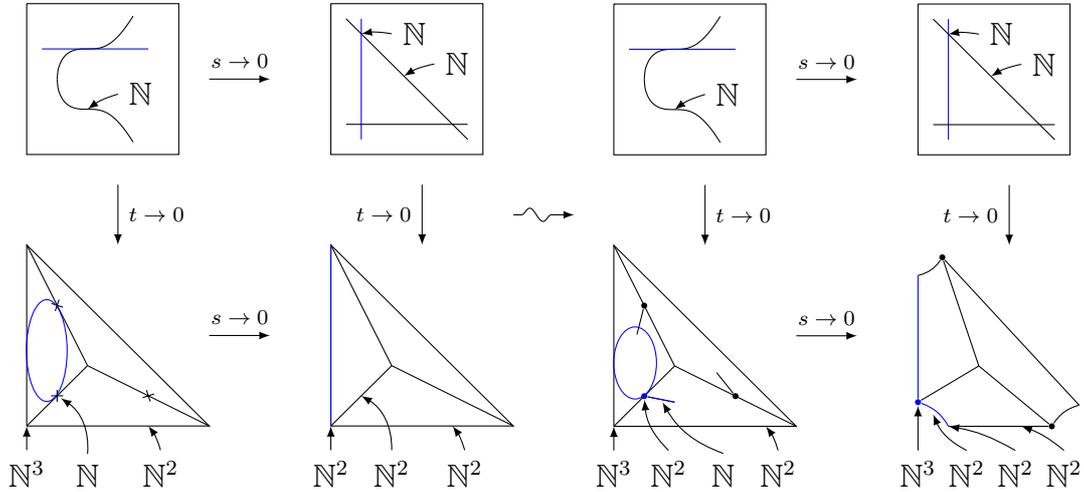
\begin{figure}[h!]
\centering
\begin{tikzpicture}[scale=0.8,decoration={snake,pre length=3pt,post length=7pt}]
% arrows
\draw[->] (2,0.5) -- (2.5,0.5) node[above]{\tiny$s\rightarrow 0$} -- (3,0.5);
\draw[->] (0.5,3) -- (0.5,2.5) node[right]{\tiny$t\rightarrow 0$} -- (0.5,2);
\draw[->] (2,4.75) -- (2.5,4.75) node[above]{\tiny$s\rightarrow 0$} -- (3,4.75);
\draw[->] (5.5,3) -- (5.5,2.5) node[left]{\tiny$t\rightarrow 0$} -- (5.5,2);
\draw[->,decorate] (7,2.5) -- (8,2.5);
% t=0 s!=0
\coordinate (0) at (0,0);
\coordinate (1) at (-1,-1);
\coordinate (2) at (2,-1);
\coordinate (3) at (-1,2);
\coordinate[fill,cross,inner sep=1pt,rotate=45] (1a) at (-0.5,-0.5);
\coordinate[fill,cross,inner sep=1pt,rotate=63.43] (2a) at (1,-0.5);
\coordinate[fill,cross,inner sep=1pt,rotate=26.57] (3a) at (-0.5,1);
\draw (1) -- (2) -- (3) -- (1);
\draw (0) -- (1);
\draw (0) -- (2);
\draw (0) -- (3);
\draw[blue] (-0.67,0.25) ellipse (9.5pt and 24pt);
\coordinate (4) at (1,-1);
\coordinate (5) at (1.2,-1.4);
\draw[->] (5) node[below]{$\mathbb{N}^2$} to [bend right=10] (4);
\coordinate (6) at (-1,-1.4);
\draw[->] (6) node[below]{$\mathbb{N}^3$} to [bend right=0] (1);
\coordinate (7) at (0,-1.46);
\draw[->] (7) node[below]{$\mathbb{N}$} to [bend right=30] (1a);
% t=0 s=0
\coordinate (0) at (5,0);
\coordinate (1) at (4,-1);
\coordinate (2) at (7,-1);
\coordinate (3) at (4,2);
\coordinate (1a) at (4.5,-0.5);
\coordinate (2a) at (6,-0.5);
\coordinate (3a) at (4.5,1);
\draw (1) -- (2) -- (3) -- (1);
\draw (0) -- (1);
\draw (0) -- (2);
\draw (0) -- (3);
\draw[blue] (1) -- (3);
\coordinate (4) at (6,-1);
\coordinate (5) at (6.2,-1.4);
\draw[->] (5) node[below]{$\mathbb{N}^2$} to [bend right=10] (4);
\coordinate (6) at (4,-1.4);
\draw[->] (6) node[below]{$\mathbb{N}^2$} to [bend right=0] (1);
\coordinate (7) at (5,-1.4);
\draw[->] (7) node[below]{$\mathbb{N}^2$} to [bend right=30] (1a);
% t!=0 s!=0
\draw (-1,6) -- (-1,3.5) -- (1.5,3.5) -- (1.5,6) -- (-1,6);
\draw[smooth,samples=100,domain=-0.5:0.75] plot (\x,{sqrt((2*\x)^3+1)/2+4.75});
\draw[smooth,samples=100,domain=-0.5:0.75] plot (\x,{-sqrt((2*\x)^3+1)/2+4.75});
\draw[blue] (-0.75,5.25) -- (1,5.25);
\coordinate (4) at (0,4.25);
\coordinate (5) at (0.5,4.5);
\draw[->] (5) node[right]{$\mathbb{N}$} to [bend right=10] (4);
% t!=0 s=0
\draw (4,6) -- (4,3.5) -- (6.5,3.5) -- (6.5,6) -- (4,6);
\draw[blue] (4.5,3.75) -- (4.5,5.75);
\draw (4.25,4) -- (6.25,4);
\draw (4.25,5.75) -- (6.25,3.75);
\coordinate (4) at (4.5,5.5);
\coordinate (5) at (5,5.5);
\draw[->] (5) node[right]{$\mathbb{N}$} to [bend right=10] (4);
\coordinate (6) at (5.2,4.8);
\coordinate (7) at (5.7,5);
\draw[->] (7) node[right]{$\mathbb{N}$} to [bend right=10] (6);
\end{tikzpicture}
\begin{tikzpicture}[scale=0.8]
% arrows
\draw[->] (2,0.5) -- (2.5,0.5) node[above]{\tiny$s\rightarrow 0$} -- (3,0.5);
\draw[->] (0.5,3) -- (0.5,2.5) node[right]{\tiny$t\rightarrow 0$} -- (0.5,2);
\draw[->] (2,4.75) -- (2.5,4.75) node[above]{\tiny$s\rightarrow 0$} -- (3,4.75);
\draw[->] (5.5,3) -- (5.5,2.5) node[left]{\tiny$t\rightarrow 0$} -- (5.5,2);
% t=0 s!=0
\coordinate (0) at (0,0);
\coordinate (1) at (-1,-1);
\coordinate (2) at (2,-1);
\coordinate (3) at (-1,2);
\coordinate[fill,circle,inner sep=0.8pt,blue] (1a) at (-0.5,-0.5);
\coordinate[fill,circle,inner sep=0.8pt] (2a) at (1,-0.5);
\coordinate[fill,circle,inner sep=0.8pt] (3a) at (-0.5,1);
\draw (1) -- (2) -- (3) -- (1);
\draw (0) -- (1);
\draw (0) -- (2);
\draw (0) -- (3);
\draw (1a) -- (0,-0.6);
\draw (2a) -- (0.69,-0.11);
\draw (3a) -- (-0.62,0.52);
\draw[blue] (-0.65,0.05) ellipse (10pt and 17pt);
\draw[blue] (1a) -- (0,-0.6);
\coordinate (4) at (1.5,-1);
\coordinate (5) at (1.7,-1.4);
\draw[->] (5) node[below]{$\mathbb{N}^2$} to [bend right=10] (4);
\coordinate (6) at (-1,-1.4);
\draw[->] (6) node[below]{$\mathbb{N}^3$} to [bend right=0] (1);
\coordinate (7) at (-0.1,-1.4);
\draw[->] (7) node[below]{$\mathbb{N}^2$} to [bend left=20] (1a);
\coordinate (8) at (0.8,-1.45);
\coordinate (9) at (-0.2,-0.6);
\draw[->] (8) node[below]{$\mathbb{N}$} to [bend left=20] (9);
% t=0 s=0
\coordinate (0) at (5,0);
\coordinate[fill,circle,inner sep=0.8pt,blue] (1) at (4,-0.6);
\coordinate[fill,circle,inner sep=0.8pt] (2) at (6.2,-1);
\coordinate[fill,circle,inner sep=0.8pt] (3) at (4.4,1.8);
\coordinate (1a) at (4.5,-1);
\coordinate (2a) at (6.65,-0.65);
\coordinate (3a) at (4,1.5);
\draw (1a) -- (2);
\draw (2a) -- (3);
\draw[blue] (3a) -- (1);
\draw (0) -- (1);
\draw (0) -- (2);
\draw (0) -- (3);
\draw[blue] (1) to [bend left=20] (1a);
\draw (2) to [bend left=20] (2a);
\draw (3) to [bend left=20] (3a);
\coordinate (4) at (5.7,-1);
\coordinate (5) at (6.4,-1.4);
\draw[->] (5) node[below]{$\mathbb{N}^2$} to [bend right=10] (4);
\coordinate (6) at (4,-1.4);
\draw[->] (6) node[below]{$\mathbb{N}^3$} to [bend right=0] (1);
\coordinate (7) at (5.6,-1.4);
\draw[->] (7) node[below]{$\mathbb{N}^2$} to [bend right=10] (1a);
\coordinate (8) at (4.25,-0.8);
\coordinate (9) at (4.8,-1.4);
\draw[->] (9) node[below]{$\mathbb{N}^2$} to [bend left=20] (8);
% t!=0 s!=0
\draw (-1,6) -- (-1,3.5) -- (1.5,3.5) -- (1.5,6) -- (-1,6);
\draw[smooth,samples=100,domain=-0.5:0.75] plot (\x,{sqrt((2*\x)^3+1)/2+4.75});
\draw[smooth,samples=100,domain=-0.5:0.75] plot (\x,{-sqrt((2*\x)^3+1)/2+4.75});
\draw[blue] (-0.75,5.25) -- (1,5.25);
\coordinate (4) at (0,4.25);
\coordinate (5) at (0.5,4.5);
\draw[->] (5) node[right]{$\mathbb{N}$} to [bend right=10] (4);
% t!=0 s=0
\draw (4,6) -- (4,3.5) -- (6.5,3.5) -- (6.5,6) -- (4,6);
\draw[blue] (4.5,3.75) -- (4.5,5.75);
\draw (4.25,4) -- (6.25,4);
\draw (4.25,5.75) -- (6.25,3.75);
\coordinate (4) at (4.5,5.5);
\coordinate (5) at (5,5.5);
\draw[->] (5) node[right]{$\mathbb{N}$} to [bend right=10] (4);
\coordinate (6) at (5.2,4.8);
\coordinate (7) at (5.7,5);
\draw[->] (7) node[right]{$\mathbb{N}$} to [bend right=10] (6);
\end{tikzpicture}
\caption{Fibers of a stable log map of degree $1$ to $\tilde{\mathfrak{X}}_Q\rightarrow\mathbb{A}^2$ for $(\mathbb{P}^2,E)$ (right) and their image on $\mathfrak{X}_Q\rightarrow\mathbb{A}^2$ (left) under the resolution from \S\ref{S:resolution}. Some stalks of the ghost sheaves are given.}
\label{fig:bigpicture}
\end{figure}

\begin{defi}
For an effective curve class $\underline{\beta}\in H_2^+(X,\mathbb{Z})$ let $\mathscr{M}(\tilde{\mathfrak{X}}_Q,\beta)$ be the moduli space of basic stable log maps to $\tilde{\mathfrak{X}}_Q\rightarrow\mathbb{A}^2$ of class $\beta$ (Definition \ref{defi:beta}). Since $\tilde{\mathfrak{X}}_Q$ is projective over $\mathbb{A}^1$ by projection to $s$, the moduli space $\mathscr{M}(\tilde{\mathfrak{X}}_Q,\beta)$ is proper over $\mathbb{A}^1$ by \cite{LogGW}, Theorem 0.2. Figure \ref{fig:bigpicture} shows the fibers of a stable log map of degree $1$ in $\mathscr{M}(\tilde{\mathfrak{X}}_Q,\beta)$ for $(\mathbb{P}^2,E)$.
\end{defi}

\begin{lem}
Let $f : \mathfrak{C} \rightarrow \mathfrak{X}$ be a stable log map in $\mathscr{M}(\tilde{\mathfrak{X}}_Q,\beta)$. Then the fibers $f_t^0 : C_t^0 \rightarrow X_t^0$ map entirely to the divisor $D_t^0$.
\end{lem}

\begin{proof}
Suppose there is an irreducible component of $C_t^0$ not mapped to $D_t^0$. Then, since the marked point is mapped to $D_t^0$, there exists an irreducible component of $C_t^0$ that is not mapped onto $D_t^0$ and is not contracted to a point. But then the tropicalization of $f_t^0$ must have at least two legs, since the balancing condition implies balancing of the legs. This means $C_t^0$ must have at least two marked points, in contradiction with the definition of $\beta$ (Definition \ref{defi:beta}).
\end{proof}

\begin{defi}
\label{defi:G}
Let $\Gamma_{D_t^0}$ be the dual intersection graph of $D_t^0$. This is a cycle with $r$ vertices. Let $\mathfrak{G}_d$ be the set of graph morphisms $g : \Gamma \rightarrow \Gamma_{D_t^0}$ where $\Gamma$ is a tree (genus $0$ graph) with vertices $V$ decorated by $d_V\in\mathbb{N}_{>0}$ such that $\sum_V d_V = d$.
\end{defi}

\begin{con}
\label{con:HG}
Note that projection to the unique unbounded direction defines a map $B \rightarrow \Gamma_{D_t^0}$, where vertices of $\Gamma_{D_t^0}$ correspond to unbounded edges of $B$. Define a surjective map
\[ \mathfrak{H}_d \rightarrow \mathfrak{G}_d \]
by composing $h : \Gamma \rightarrow B$ with this projection and defining the label $d_V$ at a vertex $V$ as follows. Let $h:\Gamma\rightarrow B$ be a tropical curve in $\mathfrak{H}_d$. Give $\Gamma$ the structure of a rooted tree by defining the root vertex to be the vertex $V_{\textup{out}}$ of the unique unbounded leg $E_{\textup{out}}$. Let $V$ be a vertex of $\Gamma$ and let $E_{V,\textup{out}}$ be the edge connecting $V$ with its parent, or $E_{V,\textup{out}}=E_{\textup{out}}$ if $V$ is the root vertex $V_{\textup{out}}$. Let $E_1,\ldots,E_r$ be the other edges of $\Gamma$ adjacent to $V$. Then define
\[ d_V = \varphi(u_{(V,E_{V,\text{out}}}) - \sum_{i=1}^r \varphi(-u_{(V,E_i)}) \]
\end{con}

Let $f : C \rightarrow \mathfrak{X} := \mathfrak{X}^{s\neq 0}$ be a stable log map in $\mathscr{M}(\mathfrak{X},\beta)$. Since all the fibers $\mathfrak{X}^s$ are isomorphic for $s\neq 0$ this gives a family of stable log maps over $\mathbb{A}^1 \times (\mathbb{A}^1 \setminus \{0\})$. Since $\mathscr{M}(\tilde{\mathfrak{X}}_Q,\beta)$ is proper this family can be uniquely completed to a family over $\mathbb{A}^2$. In other words, the limit of a stable log map in $\mathscr{M}(\mathfrak{X},\beta)$ under $s\rightarrow 0$ is well defined.

\begin{prop}
\label{prop:limit}
Let $f : \mathfrak{C} \rightarrow \mathfrak{X}_Q$ be a stable log map with tropicalization $\tilde{h}$ mapping to $h\in\mathfrak{H}_d$ under the map from Construction \ref{con:tilde}. The limit of $f$ with respect to the family $\mathfrak{X}_t^0\rightarrow\mathbb{A}^1$ has dual graph that is given by the image of $h$ under the map from Construction \ref{con:HG}.
\end{prop}

\begin{proof}
Let $f : \mathfrak{C} \rightarrow \mathfrak{X}_Q$ be a stable log map with tropicalization $\tilde{h}$ mapping to $h : \Gamma \rightarrow B$. Consider the fiber $f_0^0 : C_0^0 \rightarrow X_0^0$. If a vertex $V$ of $\Gamma$ is mapped to a vertex $v$ of $\mathscr{P}$ or the unbounded edge adjacent to $v$, then the corresponding irreducible component $C_V$ of $C_0^0$ is mapped to the irreducible component $X_v$ of $X_0^0$ corresponding to $v$. But then, for $t\neq 0$, the corresponding irreducible component $C_V$ of $C_t^0$ is mapped to the irreducible component $D_v$ of $D_t^0$ corresponding to $v$. This is the image of $V$ under the map from Construction \ref{con:HG}. The map $f_t^0 : C_t^0 \rightarrow D_v = \mathbb{P}^1$ is a multiple cover of a line. Its degree is precisely the label $d_V$ of $V$ as above.
\end{proof}

\begin{figure}[h!]
\centering
\begin{tikzpicture}[scale=1.2]
\coordinate (1) at (0,0);
\coordinate (2) at (0,1);
\coordinate (3) at (-3,-1);
\coordinate (4) at (-3,2);
\draw[thick] (3) -- (1) -- (2) -- (4);
\draw[thick] (1) -- (5,0);
\draw[thick] (2) -- (5,1);
\draw[thick,dashed] (3) -- (5,-1);
\draw[thick] (4) -- (5,2);
% cut
\draw[thick,dashed] (-1.5,-0.25) -- (0,0.5) -- (-1.5,1.25);
\draw[thick,dashed] (-3,-0.8) -- (-1.5,-0.5) -- (-0.75,0.25);
\draw[thick,dashed] (-3,1.8) -- (-1.5,1.5) -- (-0.75,0.75);
\fill[color=gray] (-1.5,-0.25) -- (0,0.5) -- (-1.5,1.25) -- (-1.5,-0.25) -- cycle;
\fill[color=gray] (-3,1.8) -- (-1.5,1.5) -- (-0.75,0.75) -- (-0.75,0.25) -- (-1.5,-0.5) -- (-3,-0.8);
\coordinate[fill,cross,thick,inner sep=2pt] (5) at (0,0.5);
\coordinate[fill,cross,thick,inner sep=2pt,rotate=26.57] (6) at (-1.5,-0.5);
\coordinate[fill,cross,thick,inner sep=2pt,rotate=-30] (7) at (-1.5,1.5);
% refinement
\draw[thick] (0,0) -- (-1,0);
\draw[thick] (0,1) -- (-1,1);
% blue
\fill[color=blue] (0,0) circle (1.5pt);
\fill[color=blue] (0,1) circle (1.5pt);
\fill[color=blue] (2,0.66) circle (1.5pt);
\draw[very thick,blue] (-1.5,1.5) -- (-0.75,1.25) node[above]{$2$} -- (0,1) -- (2,0.66) -- (0,0) -- (-1.5,-0.5);
\draw[very thick,blue] (2,0.66) -- (3.5,0.66) node[below]{$9$} -- (5,0.66);
\draw[very thick,blue] (0,1) -- (0,0.5);
% projection
\draw[->] (6,0.5) -- (7,0.5);
\draw[thick] (8,-0.9) -- (8,2);
\draw[very thick, blue] (8,0) node[right]{$1$} -- (8,1) node[right]{$2$};
\coordinate[draw,dashed,circle,inner sep=2pt] (a) at (8,-1);
\coordinate[fill,circle,inner sep=2pt,blue] (b) at (8,0);
\coordinate[fill,circle,inner sep=2pt,blue] (c) at (8,1);
\coordinate[fill,circle,inner sep=2pt] (d) at (8,2);
\end{tikzpicture}
\caption{The tropical curve from Figure \ref{fig:balancing2} under $\mathfrak{H}_d \rightarrow \mathfrak{G}_d$ gives $[2,1]$ for some choice of $D_i$ and some choice of ordering.}
\label{fig:balancing2}
\end{figure}
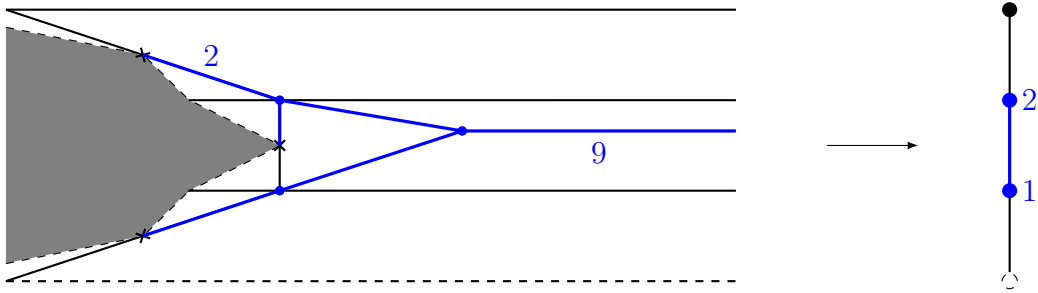

For $s=0$ fix a cyclic labelling of the cycle of lines $D_t^0 = D_1 + \ldots + D_k$ and write $v_1,\ldots,v_r$ for the corresponding vertices of $\Gamma_{D_t^0}$.

\begin{defi}
\label{defi:splitting}
Given $g : \Gamma \rightarrow \Gamma_{D_t^0}$ in $\mathfrak{G}_d$ write
\[ d_i := \sum_{\substack{V \in \Gamma^{[0]} \\ g(V) = v_i}} d_V. \]
Then the collection $[d_1,\ldots,d_k]$ gives the degrees of the curve over the lines $D_i=\mathbb{P}^1$. If all components of $D_t^0$ are isomorphic as divisors of $X_t^0$, e.g. for $X=\mathbb{P}^2$ or $X=\mathbb{P}^1\times\mathbb{P}^1$, we omit all zeros in this collection. We call this collection the \textit{degree splitting} corresponding to $g$ or any tropical curve $h$ mapping to $g$ under the map from Construction \ref{con:HG}. For example, the tropical curve for $(\mathbb{P}^2,E)$ in Figure \ref{fig:balancing2} has degree splitting $[2,1,0]$ for any choice of labelling of $D_t^0=D_1+D_2+D_3$, and we simply write this as $[2,1]$.
\end{defi}

\begin{con}
\label{con:GH}
Define a surjective map
\[ \mathfrak{G}_d \rightarrow H_2^+(X,\mathbb{Z}) \]
by sending an element $g : \Gamma \rightarrow \Gamma_{D_t^0}$ of $\mathfrak{G}_d$ with degree splitting $[d_1,\ldots,d_k]$ to the curve class $\underline{\beta}$ defined by 
\[ D_i \cdot \underline{\beta} = d_i. \]
This is well-defined by the balancing condition and since $H_2^+(X,\mathbb{Z})\simeq H_2^+(X^0,\mathbb{Z})$, where $X^0$ is a toric variety.
\end{con}

\begin{defi}
\label{defi:H}
Let $\mathfrak{H}_\beta$ be the set of tropical curves mapping to $\underline{\beta}$ under the decomposition of the maps from Constructions \ref{con:HG} and \ref{con:GH}. Let $\tilde{\mathfrak{H}}_\beta$ be the preimage of $\mathfrak{H}_\beta$ under the map from Construction \ref{con:tilde}.
\end{defi}

A direct consequence of Proposition \ref{prop:limit} is the following.

\begin{cor}
The tropicalization of a stable log map of class $\beta$ is in $\tilde{\mathfrak{H}}_\beta$.
\end{cor}

\subsection{Refinement and logarithmic modification}				%%
\label{S:refinement}

To apply the degeneration formula in \S\ref{S:degformula} we need a degeneration of $(X,D)$ such that all stable log maps to the central fiber are torically transverse. We achieve this as follows. 

\begin{con}[The refined degeneration]
\label{con:deg3}
Let $\mathscr{P}_d$ be a refinement of $\mathscr{P}$ such that each tropical curve in $\mathfrak{H}_{\leq d}=\cup_{d'\leq d}\mathfrak{H}_{d'}$ (or equivalently in $\tilde{\mathfrak{H}}_{\leq d}$) is contained in the $1$-skeleton of $\mathscr{P}_d$. This is well-defined by finiteness of $\mathfrak{H}_d$ (Corollary \ref{cor:finite}) and defines a refinement of the generalized cone complex $\Sigma(\tilde{X}_0)$ by taking cones over cells of $\mathscr{P}_d$. In turn, $\mathscr{P}_d$ induces a logarithmic modification $\tilde{\mathfrak{X}}_d\rightarrow\mathbb{A}^1$ of $\tilde{\mathfrak{X}}\rightarrow\mathbb{A}^1$ (see \S\ref{A:artin}) without changing the generic fiber. By making a base change $\mathbb{A}^1 \rightarrow \mathbb{A}^1, t \mapsto t^e$ we can scale $\mathscr{P}_d$ and thus assume it has integral vertices (c.f. \cite{NS}, Proposition 6.3).
\end{con}

\begin{rem}
\label{rem:transverse}
The dual intersection complex of the central fiber $Y$ of $\tilde{\mathfrak{X}}_d$ is given by $(\tilde{B},\mathscr{P}_d,\varphi)$. Hence, all stable log maps to $Y\rightarrow\textup{pt}_{\mathbb{N}}$ of class $\beta$ as in Definition \ref{defi:beta} are torically transverse, since their tropicalizations map onto the $1$-skeleton of $\mathscr{P}_d$, with vertices mapping to vertices of $\mathscr{P}_d$ (see \cite{MR}, Proposition 4.6).
\end{rem}

Gromov-Witten invariants are invariant under logarithmic modifications \cite{AW}. Hence,
\[ N_\beta = \int_{\llbracket\mathscr{M}(Y,\beta)\rrbracket}1. \]
In the next section we will apply the degeneration formula of logarithmic Gromov-Witten theory to get a formula for $N_d$ in terms of logarithmic Gromov-Witten invariants of irreducible components of $Y$.

\section{The degeneration formula}							%%%
\label{S:degformula}

Consider a projective semi-stable degeneration $\pi : \mathfrak{X} \rightarrow T = \textup{Spec }R$, for $R$ a discrete valuation ring. This is a projective surjection such that the generic fiber is smooth and the fiber $X=\pi^{-1}(0)$ over the closed point $0\in T$ is simple normal crossings with two smooth connected (hence irreducible) components $X_1,X_2$ meeting in a smooth connected divisor $D$. In the logarithmic language this means that $\pi$ is log smooth when $T$ and $\mathfrak{X}$ carry the divisorial log structures given by $0\in T$ and $X\subseteq\mathfrak{X}$, respectively. The degeneration formula relates invariants (relative or logarithmic Gromov-Witten invariants) on the generic fiber of $\pi : \mathfrak{X} \rightarrow T$ to invariants on the components $X_1,X_2$ of the special fiber. 

The degeneration formula was proved for stable relative maps, in symplectic geometry \cite{LiRu}\cite{IoPa} and in algebraic geometry \cite{Li2}\cite{AF}, as well as for stable log maps using expanded degenerations \cite{Che2}. A pure log-geometric version avoiding the target expansions of relative Gromov-Witten theory was worked out by Kim, Lho and Ruddat \cite{KLR} using logarithmic Gromov-Witten theory \cite{LogGW}. 

The formula in \cite{KLR} is stated in the setup above, with two smooth irreducible components $X_1$ and $X_2$. In our setup we have several irreducible components, indexed by the vertices of a tropical curve. One could generalize the formula in \cite{KLR} to the case where $X_1,X_2$ are only log smooth, in particular they might be reducible, and then apply it repeatedly to get a formula for several components. For the sake of self-containedness we will take a different approach and prove the formula more explicitly in our setting, only referring to \cite{KLR} for some general statements.

Fix an integer $d>0$ and let $\tilde{\mathfrak{X}}_d\rightarrow\mathbb{A}^1$ be the refined log smooth degeneration of $(X,D)$ from Construction \ref{con:deg3}. Write the central fiber as $Y$ and let $\beta$ be a class of stable log maps as in Definition \ref{defi:Nd} with $D\cdot\underline{\beta}\leq dw_{\text{out}}$. Let $Y^\circ$ be the complement of zero-dimensional toric strata in $Y$ and write 
\[ \mathscr{M}_\beta := \mathscr{M}(Y,\beta). \]

Since $\mathscr{M}(Y^\circ,\beta)$ is canonically isomorphic to the moduli space of torically transverse stable log maps to $Y$ of class $\beta$ and all such maps are torically transverse (Remark \ref{rem:transverse}), the canonical inclusion gives an isomorphism $\mathscr{M}_\beta \cong \mathscr{M}(Y^\circ,\beta)$.

For a vertex $v\in\mathscr{P}_d^{[0]}$, let $Y_v^\circ$ be the complement of the $0$-dimensional toric strata of the irreducible component $Y_v$ of $Y$ corresponding to $v$. Then $Y^\circ$ is a union of finitely many log smooth schemes $Y_v^\circ$ over $\textup{pt}_{\mathbb{N}}$, with $Y_v^\circ \cap Y_{v'}^\circ = \emptyset$ if there is no edge connecting $v$ and $v'$, and $D_E^\circ:=Y_v^\circ \cap Y_{v'}^\circ$ log smooth and a divisor of both $Y_v^\circ$ and $Y_{v'}^\circ$ if there is an edge $E$ connecting $v$ and $v'$. The intersection of any triple of components is empty. Hence, we can apply the degeneration formula.

\subsection{Toric invariants}							%%
\label{S:toricinv}

We introduce logarithmic Gromov-Witten invariants of toric varieties with point conditions on the toric boundary, following \cite{GPS}.

Let $M\simeq\mathbb{Z}^2$ be a lattice and let $M_{\mathbb{R}}=M\otimes_{\mathbb{Z}}\mathbb{R}$ be the associated vector space. Let $(m_1,\ldots,m_n)$ be an $n$-tuple of distinct nonzero primitive vectors in $M$ and let $\textbf{w}=(\textbf{w}_1,\ldots,\textbf{w}_n)$ be an $n$-tuple of weight vectors $\textbf{w}_i=(w_{i1},\ldots,w_{il_i})$ with $l_i>0$, $w_{ij}\in \mathbb{N}$ such that
\[ \sum_{i=1}^n|\textbf{w}_i|m_i = w_{\textup{out}} m_{\textup{out}} \]
for $0\neq m_{\text{out}}\in M$ primitive and $w_{\textup{out}}> 0$. Here $|\textbf{w}_i| := \sum_{j=1}^{l_i} w_{ij}$. Let $\Sigma$ be the complete rational fan in $M_{\mathbb{R}}$ whose rays are generated by $-m_1,\ldots,-m_n,m_{\textup{out}}$ and let $X_\Sigma$ be the corresponding toric surface over $\mathbb{C}$. By refining $\Sigma$ if necessary, we can assume that $X_\Sigma$ is nonsingular. Let $D_1,\ldots,D_n,D_{\textup{out}}\subseteq X_\Sigma$ be the toric divisors corresponding to the given rays. Let $X_\Sigma^\circ$ be the complement of the $0$-dimensional torus orbits in $X_\Sigma$, and let $D_i^\circ=D_i\cap X_\Sigma^\circ$, $D_{\textup{out}}^\circ = D_{\textup{out}}\cap X_\Sigma^\circ$. Then define a class $\beta_{\textbf{w}}$ of stable log maps to $X_\Sigma$ as follows.
\begin{compactenum}[(1)]
\item genus $g=0$;
\item $k=l_1+\ldots+l_n+1$ marked points $p_{ij}, i=1,\ldots,n, j=1,\ldots,l_i$ and $p$;
\item $\underline{\beta}_{\textbf{w}}\in H_2(X_\Sigma,\mathbb{Z})$ defined by intersection numbers with toric divisors,
\[ D_i \cdot \underline{\beta}_{\textbf{w}} = |\textbf{w}_i|, \quad D_{\textup{out}}\cdot\underline{\beta}_{\textbf{w}} = w_{\textup{out}}; \]
\item contact data $u_{p_{ij}}=w_{ij}m_i$ and $u_p=w_{\textup{out}}m_{\textup{out}}$.
\end{compactenum}
By restriction we get a class of stable log maps to $X_\Sigma^\circ$ that we also denote by $\beta_{\textbf{w}}$. The moduli space $\mathscr{M}(X_\Sigma^\circ,\beta_{\textbf{w}})$ in general is not proper, since $X_\Sigma^\circ$ is not proper. However, the evaluation map 
\[ \textup{ev}^\circ : \mathscr{M}(X_\Sigma^\circ,\beta_{\textbf{w}}) \rightarrow \prod_{i=1}^n(D_i^\circ)^{l_i} \]
is proper (\cite{GPS}, Proposition 4.2) and we obtain a proper moduli space via base change to a point. To be precise, let $\gamma : \textup{Spec }\mathbb{C} \rightarrow \prod_{i=1}^n(D_i^\circ)^{l_i}$ be a point. Then 
\[ \mathscr{M}_\gamma := \textup{Spec }\mathbb{C} \times_{\prod_{i=1}^n(D_i^\circ)^{l_i}} \mathscr{M}(X_\Sigma^\circ,\beta_{\textbf{w}}) \]
is a proper Deligne-Mumford stack admitting a virtual fundamental class, and we can define the logarithmic Gromov-Witten invariant
\begin{equation}
\label{eq:Ntor}
N_{\textbf{m}}(\textbf{w}) := \int_{\mathscr{M}_\gamma}\gamma^!\llbracket\mathscr{M}(X_\Sigma^\circ,\beta_{\textbf{w}})\rrbracket.
\end{equation}
Since the codimension of $\gamma$ equals the virtual dimension of $\mathscr{M}(X_\Sigma^\circ,\beta_{\textbf{w}})$, this definition makes sense. Note that we may add further primitive vectors $m_i$ to $\textbf{m}$, with weight vectors $\textbf{w}_i=0$. This leads to a subdivision of $\Sigma$, hence to a toric blow up of $X_\Sigma$, but the logarithmic Gromov-Witten invariants do not change.

\subsection{The decomposition formula}					%%
\label{S:decomposition}

By the decomposition formula for stable log maps (\cite{ACGS1}, Theorem 1.2), the moduli space $\mathscr{M}_\beta$ decomposes into moduli spaces indexed by certain decorated tropical curves. Here decorated means that there are classes of stable log maps $\beta_V$ attached to the vertices. In this section we show that a tropical curve in $\tilde{\mathfrak{H}}_d$ automatically carries such decorations.

\begin{prop}
\label{prop:decorations}
Let $f : C/\textup{pt}_{\mathbb{N}} \rightarrow Y/\textup{pt}_{\mathbb{N}}$ be a stable log map in $\mathscr{M}_\beta$ with tropicalization $\tilde{h} : \tilde{\Gamma} \rightarrow \tilde{B}$. For each vertex $V\in\tilde{\mathscr{P}}_d^{[0]}$, the class $[C_V]\in H_2^+(Y_{\tilde{h}(V)},\mathbb{Z})$ is uniquely determined by the intersection numbers of $C_V$ with components of $\partial Y_{\tilde{h}(V)}$, i.e., by $\tilde{h}$.
\end{prop}

\begin{proof}
If $V$ is of type (I) as in Definition \ref{defi:types}, then $Y_{\tilde{h}(V)}$ is a toric variety, so the statement is true. If $V$ is of type (II), then $C_V$ is a multiple cover of some exceptional line $L^{\textup{exc}}_v$ (Definition \ref{defi:Lexc}). Its intersection with $\partial Y_v$ determines the degree $d$ of the multiple cover, hence the curve class $[C_V] = d[L^{\textup{exc}}_v]\in H_2^+(Y_v,\mathbb{Z})$. Let V be a vertex of type (III). It is mapped to a vertex $v$ of $\mathscr{P}$. Let $X_v$ be the corresponding component of $X_0$. This is a toric variety. By Proposition \ref{prop:balancing}, (III), we know the intersection of the image of $C_V$ under the resolution $\nu : \tilde{\mathfrak{X}} \rightarrow \mathfrak{X}$ from \S\ref{S:resolution} with the toric divisors of $X_v$, hence the curve class $[\nu(C_V)]\in H_2^+(X_v,\mathbb{Z})$. But this determines $[C_V] = [\nu(C_V)] - k[L^{\textup{exc}}_v] \in H_2^+(Y_v,\mathbb{Z})$, where $k$ is as in Proposition \ref{prop:balancing}, (III).
\end{proof}

\begin{defi}
\label{defi:Mh}
For $\tilde{h} : \tilde{\Gamma} \rightarrow \tilde{B}$ in $\tilde{\mathfrak{H}}_d$, let $\mathscr{M}_{\tilde{h}}$ be the moduli space of stable log maps with tropicalization $\tilde{h}$. This is proper by \cite{ACGS1}, Proposition 2.34. 
\end{defi}

\begin{rem}
\label{rem:tau}
In fact \cite{ACGS1} deals with moduli spaces $\mathscr{M}_{\tilde{\tau}}$ of stable log maps \textit{marked by} $\tilde{\tau}=(\tau,\textbf{A})$, where $\tau$ is a type of tropical maps and $\textbf{A}$ is a vertex decoration by curve classes.  Since the virtual dimension of $\mathscr{M}_\beta$ is zero and tropical curves in $\tilde{\mathfrak{H}}_d$ are rigid, such $\tilde{\tau}$ are in bijection with vertex decorated tropical curves. We showed in Proposition \ref{prop:decorations} that tropical curves in $\tilde{\mathfrak{H}}_d$ carry unique vertex decorations. So $\tilde{\tau}$ uniquely defines a tropical curve $\tilde{h}$ and $\mathscr{M}_{\tilde{\tau}}$ equals $\mathscr{M}_{\tilde{h}}$.
\end{rem}

\begin{prop}[Decomposition formula]
\label{prop:dec}
\[ \llbracket\mathscr{M}_\beta\rrbracket = \sum_{\tilde{h}\in\tilde{\mathfrak{H}}_\beta}\frac{l_{\tilde{\Gamma}}}{|\textup{Aut}(\tilde{h})|} F_\star\llbracket\mathscr{M}_{\tilde{h}}\rrbracket, \]
where $l_{\tilde{\Gamma}} := \textup{lcm}\{w_E \ | \ E\in E(\tilde{\Gamma})\}$ and $F:\mathscr{M}_{\tilde{h}}\rightarrow\mathscr{M}_\beta$ is the forgetful map. Here $\text{Aut}(\tilde{h})$ is the group of automorphisms of $\tilde{h}$ (Definition \ref{defi:aut}).
\end{prop}

\begin{proof}
The decomposition formula (\cite{ACGS1}, Theorem 1.2) gives $\llbracket\mathscr{M}_\beta\rrbracket$ as a sum over decorated types of tropical maps $\tilde{\tau}$. By Remark \ref{rem:tau} this is a summation over $\tilde{\mathfrak{H}}_\beta$. The multiplicity $m_\tau$ in \cite{ACGS1}, Theorem 1.2, is defined as the index of the image of the lattice $\Sigma(\textup{pt}_{\mathbb{N}})=\mathbb{N}$ inside the lattice $\Sigma(\textup{pt}_\mathbb{N})=\mathbb{N}$. Here the first $\text{pt}_{\mathbb{N}}$ is the base of the curve while the second $\text{pt}_{\mathbb{N}}$ is the base of $Y$. In other words, $m_\tau$ is the smallest integer such that scaling $\tilde{B}$ by $m_\tau$ leads to a tropical curve with integral vertices and edge lengths. By  Construction \ref{con:deg3} $\mathscr{P}_d$ has integral vertices (by the  base change $t\mapsto t^e$) and tropical curves in $\tilde{\mathfrak{H}}_\beta$ are contained in the $1$-skeleton of $\mathscr{P}_d$ with vertices mapping to vertices of $\mathscr{P}_d$. The affine length of the image of an edge $E\in E(\tilde{\Gamma})$ is $\ell_Ew_E$ by Definition \ref{defi:tropical}, (i). So the scaling necessary to obtain integral edge lengths is $m_\tau = l_{\tilde{\Gamma}} := \textup{lcm}\{w_E \ | \ E\in E(\tilde{\Gamma})\}$. Moreover $\text{Aut}(\tau)=\text{Aut}(\tilde{h})$ by our definition of automorphisms (see Definition \ref{defi:aut}). Then \cite{ACGS1}, Theorem 1.2, gives the formula above.
\end{proof}

\subsection{Contributions of the vertices}						%%
\label{S:contributions}

Let $\tilde{h} : \tilde{\Gamma} \rightarrow \tilde{B}$ be a tropical curve in $\tilde{\mathfrak{H}}_d$. Define
\[ \mathscr{M}_V^\circ := \mathscr{M}(Y_{\tilde{h}(V)}^\circ,\beta_V), \]
where $Y_{\tilde{h}(V)}^\circ$ is the complement of the $0$-dimensional toric strata in $Y_{\tilde{h}(V)}$.

For $V\in V_{II}(\tilde{\Gamma})$ (Definition \ref{defi:types}) with adjacent edge $E$, the moduli space $\mathscr{M}_V^\circ$ is proper, since it is isomorphic to the moduli space of $w_E$-fold multiple covers of $\mathbb{P}^1$ totally ramified at a point.

For $V\in V(\tilde{\Gamma})\setminus V_{II}(\tilde{\Gamma})$ we obtain a proper moduli space as follows. Again, $\tilde{\Gamma}$ is a rooted tree with root vertex $V_{\textup{out}}$. There is a natural orientation of the edges of $\tilde{\Gamma}$ by choosing edges to point from a vertex to its parent. For each vertex $V\in V(\tilde{\Gamma})\setminus V_{II}(\tilde{\Gamma})$ there is an evaluation map
\[ \textup{ev}_{V,-}^\circ : \mathscr{M}_V^\circ \rightarrow \prod_{E\rightarrow V}D_E^\circ, \]
where the product is over all edges of $\tilde{\Gamma}$ adjacent to $V$ and pointing towards $V$.

\begin{lem}
\label{lem:proper}
The evaluation map $\textup{ev}_{V,-}^\circ$ is proper.
\end{lem}

\begin{proof}
For $V\in V_I(\tilde{\Gamma})$ this is \cite{GPS}, Proposition 4.2. For $V\in V_{III}(\tilde{\Gamma})$ it is similar to \cite{GPS}, Proposition 5.1. Let us carry this out. Let $V\in V(\tilde{\Gamma})$ be a vertex. We use the valuative criterion for properness, so let $R$ be a valuation ring with residue field $K$, and suppose we are given a diagram
\begin{equation*}
\begin{xy}
\xymatrix{
T=\text{Spec }K \ar[r]\ar[d] 	& \mathscr{M}_V^\circ \ar[r]\ar[d]^{\textup{ev}_{V,-}^\circ} 	& \mathscr{M}_V \ar[d]^{\textup{ev}_{V,-}} \\
S=\text{Spec }R \ar[r]		& \prod_{E\rightarrow V}D_E^\circ \ar[r]					& \prod_{E\rightarrow V}D_E
}
\end{xy}
\end{equation*}
Since $\mathscr{M}_V$ is proper, $\textup{ev}_{V,-}$ is proper and we obtain a unique family of stable log maps
\begin{equation*}
\begin{xy}
\xymatrix{
\mathcal{C} \ar[r]\ar[d]			& Y_{\tilde{h}(V)}\times S \ar[d] \\
S \ar[r]^{=}					& S
}
\end{xy}
\end{equation*}
We will show that $f$ is a family of stable log maps to $Y_{\tilde{h}(V)}^\circ$. Let $0\in S$ be the closed point and consider $f_0 : C_0 \rightarrow Y_{\tilde{h}(V)}$. The marked points of $C_0$ map to $Y_{\tilde{h}(V)}^\circ$. 

Suppose $f_0(C_0)$ intersects a toric divisor $D\subset Y_{\tilde{h}(V)}$ at a point of $D\setminus Y_{\tilde{h}(V)}^\circ$. The intersection number of $C_0$ with $D$ is accounted for in $Y_{\tilde{h}(V)}^\circ$. For $V\in V_I(\tilde{\Gamma})$ this is clear and for $V\in V_{III}(\tilde{\Gamma})$ this follows since the intersection number is accounted for after composing with the resolution $\tilde{\mathfrak{X}}\rightarrow\mathfrak{X}$ from \S\ref{S:resolution}. Hence, there must be an irreducible component $C$ of $C_0$ dominating $D$. Let $D_1$ and $D_2$ be the two distinct toric divisors of $Y_{\tilde{h}(V)}$ intersecting $D$ only at two distinct torus fixed points of $D$. It was shown in \cite{GPS}, Proposition 4.2, that there are irreducible components $C_1,C_2\subset C_0$ intersecting $C$ and dominating $E_1$ and $E_2$, respectively. By applying this statement repeatedly, replacing $C$ with $C_1$ or $C_2$, we find that $C_0$ contains a cycle of components dominating the union of toric divisors of $Y_{\tilde{h}(V)}$. But then $C_0$ would have genus $g>0$, contradicting the assumptions. We have shown by contradiction that $f : \mathcal{C} \rightarrow Y_{\tilde{h}(V)} \times S$ is a family of stable log maps to $Y_{\tilde{h}(V)}^\circ$. Hence, $\textup{ev}_{V,-}^\circ$ is proper by the valuative criterion for properness.
\end{proof}

Since properness of morphisms is stable under base change, we obtain a proper moduli space by base change to a point
\[ \gamma_V : \textup{Spec }\mathbb{C} \rightarrow \prod_{E\rightarrow V}D_E^\circ, \]
that is,
\[ \mathscr{M}_{\gamma_V} := \textup{Spec }\mathbb{C} \times_{\prod_{E\rightarrow V}D_E^\circ} \mathscr{M}_V^\circ \]
is a proper Deligne-Mumford stack.

\begin{lem}
\label{lem:vdim}
For $V\in V_{II}(\tilde{\Gamma})$ the virtual dimension of $\mathscr{M}_V$ is zero. Otherwise the virtual dimension of $\mathscr{M}_V$ equals the codimension of $\gamma_V$.
\end{lem}

\begin{proof}
For $V\in V_{II}(\tilde{\Gamma})$ with adjacent edge $E$ the moduli space $\mathscr{M}_V$ is isomorphic to the moduli space of $w_E$-fold multiple covers of $\mathbb{P}^1$ totally ramified at a point. However, the two moduli spaces carry obstruction theories which differ by $H^1(C,f^\star\mathcal{O}_{\mathbb{P}^1}(-1))$ at a moduli point $[f:C\rightarrow L_V^{\text{exc}}]$ (c.f. \cite{GPS}, {\S}5.3). The rank of $H^1(C,f^\star\mathcal{O}_{\mathbb{P}^1}(-1))$ is $w_E-1$ and so is the virtual dimension of $\mathscr{M}(\mathbb{P}^1/\infty,w_E)$. Hence, the virtual dimension of $\mathscr{M}_V$ is zero.

Otherwise, the virtual dimension is easily seen to be the number of edges of $\tilde{\Gamma}$ pointing towards $V$, with orientation of $\tilde{\Gamma}$ as given above. By definition this is the codimension of $\gamma_V$.
\end{proof}

\begin{defi}
\label{defi:NV}
For a vertex $V$ of $\tilde{\Gamma}$ define
\[ N_V := \begin{cases} \int_{\llbracket\mathscr{M}_V^\circ\rrbracket}1, & V \in V_{II}(\tilde{\Gamma}); \\ \int_{\mathscr{M}_{\gamma_V}}\gamma_V^!\llbracket\mathscr{M}_V^\circ\rrbracket, & V \in V_I(\tilde{\Gamma})\cup V_{III}(\tilde{\Gamma}). \end{cases} \]
This is a finite number by Lemma \ref{lem:vdim} and independent of $\gamma_V$ by Lemma \ref{lem:proper}.
\end{defi}

\begin{prop}
\label{prop:N}
We give $N_V$ for the different types of vertices (Definition \ref{defi:types}).
\begin{compactenum}[(I)]
\item For $V\in V_I(\tilde{\Gamma})$ let $e_1,\ldots,e_n$ be the edges of $\mathscr{P}_d$ adjacent to $\tilde{h}(V)$ and let $m_1,\ldots,m_n$ be the corresponding primitive vectors. Let $\textbf{w}_{i}=(w_{i1},\ldots,w_{il_i})$ be the weights of edges of $\tilde{\Gamma}$ mapping to $e_i$ and write $\textbf{w}=(\textbf{w}_1,\ldots,\textbf{w}_n)$. Then $N_V$ is the toric invariant from \eqref{eq:Ntor}
\[ N_V = N_\textbf{m}(\textbf{w}), \]
\item If $V\in V_{II}(\tilde{\Gamma})$, then
\[ N_V = \frac{(-1)^{w_E-1}}{w_E^2}, \]
where $E$ is the unique edge adjacent to $V$.
\item If $V\in V_{III}(\tilde{\Gamma})$, then
\[ N_V = \sum_{\textbf{w}_{V,+}} \frac{N_{\textbf{m}}(\textbf{w})}{|\textup{Aut}(\textbf{w}_{V,+})|} \prod_{i=1}^l\frac{(-1)^{w_{V,i}-1}}{w_{V,i}}. \]
The sum is over all weight vectors $\textbf{w}_{V,+}=(w_{V,1},\ldots,w_{V,l_V})$ such that $|\textbf{w}_{V,+}| := \sum_{i=1}^{l_V} w_{V,i} = k$, with $k$ as in Proposition \ref{prop:balancing}, (III). Further, $N_{\textbf{m}}(\textbf{w})$ is as in \eqref{eq:Ntor} with $\textbf{m}=(m_{v,-},m_{v,+})$ and $\textbf{w}=((w_E)_{E\in E_{V,-}(\tilde{\Gamma})},\textbf{w}_{V,+})$, where $E_{V,-}$ is the set of edges adjacent to $V$ and mapped to direction $m_{v,-}$.
\end{compactenum}
\end{prop}

\begin{proof}
(I) is by the definition of $N_{\textbf{m}}(\textbf{w})$ in \eqref{eq:Ntor}. For (II) recall from the proof of Lemma \ref{lem:vdim} that $\mathscr{M}_V$ is isomorphic to $\mathscr{M}(\mathbb{P}^1/\infty,w_E)$ with obstruction theory differing by $H^1(C,f^\star\mathcal{O}_{\mathbb{P}^1}(-1))$. Hence, 
\[ N_V=\int_{\llbracket\mathscr{M}(\mathbb{P}^1/\infty,w_E)\rrbracket}e(H^1(C,f^\star\mathcal{O}_{\mathbb{P}^1}(-1))) \]
which is equal to $(-1)^{w_E-1}/w_E^2$ by the genus zero part of \cite{BP}, Theorem 5.1, see also \cite{GPS}, Propositions 5.2 and 6.1. For (III) we apply \cite{GPS}, Proposition 5.3. We only blow up one point on the divisor $D_{v,+}$, so in the notation of \cite{GPS}, Proposition 5.3, we have $\textbf{P}=(P_+)$ with $P_+ = k$ as in Proposition \ref{prop:balancing}, (III). Note that our $i$ is called $j$ in \cite{GPS} and the $i$ of \cite{GPS} is equal to $1$ here. Further, $R_{\textbf{P}_+|\textbf{w}_+} = \prod_{i=1}^{l_i}\frac{(-1)^{w_i-1}}{w_i^2}$ by \cite{GPS}, Proposition 5.2, and the discussion thereafter. Then \cite{GPS}, Proposition 5.3, gives (III).
\end{proof}

\subsection{Gluing}										%%

Define $\bigtimes_{V\in V(\tilde{\Gamma})}\mathscr{M}_V$ to be the moduli space of stable log maps in $\prod_V\mathscr{M}_V$ matching over the divisors $D_E$, $E\in E(\tilde{\Gamma})$, i.e., the fiber product
\begin{equation*}
\begin{xy}
\xymatrix{
\displaystyle\bigtimes_V\mathscr{M}_V \ar[rr]\ar[d]						&& \displaystyle\prod_V\mathscr{M}_V \ar[d]^{\textup{ev}} \\
\displaystyle\prod_{E\in E(\tilde{\Gamma})} D_E \ar[rr]^\delta				&& \displaystyle\prod_V\prod_{\substack{E\in E(\tilde{\Gamma}) \\ V \in E}} D_E
}
\end{xy}
\end{equation*}
Here $\textup{ev}$ is the product of evaluation maps to common divisors (labeled by compact edges) and $\delta$ is the diagonal map. Similarly define $\bigtimes_V\mathscr{M}_V^\circ$.

\begin{defi}
Let $\textup{cut} : \mathscr{M}_{\tilde{h}} \rightarrow \bigtimes_V\mathscr{M}_V$ be the morphism defined by cutting a curve along its gluing nodes. For a precise definition see \cite{Bou1}, {\S}7.1. Here $\mathscr{M}_{\tilde{h}}$ denotes the moduli space of stable log maps with tropicalization $\tilde{h}$ (Definition \ref{defi:Mh}). Since every stable log map in $\mathscr{M}_{\tilde{h}}$ is torically transverse (Remark \ref{rem:transverse}) this is in fact a morphism
\[ \text{cut} : \mathscr{M}_{\tilde{h}} \rightarrow \bigtimes_V\mathscr{M}_V^\circ. \]
\end{defi}

\begin{lem}
\label{lem:degcut}
The morphism $\text{cut}$ is \'etale of degree
\[ \textup{deg}(\textup{cut}) = \frac{\prod_{E\in E(\tilde{\Gamma})}w_E}{l_{\tilde{\Gamma}}}, \]
where $l_{\tilde{\Gamma}} = \textup{lcm}\{w_E\}$.
\end{lem}

\begin{proof}
Since stable log maps in $\mathscr{M}_{\tilde{h}}$ are torically transverse by Construction \ref{con:deg3}, locally we are gluing along a smooth divisor as in \cite{KLR}. Then the statement is \cite{KLR}, Lemma 9.2 (4), and the degree is computed by \cite{KLR}, (6.13). We will briefly explain how to arrive at this expression.

For each edge $E$ we have a choice of $w_E$-th root of unity in the log structure of $C$ at the corresponding node, contributing a factor of $w_E$ to $\text{deg}(\text{cut})$. This was computed e.g. in \cite{NS}, Proposition 7.1, \cite{Gr10}, Proposition 4.23, and \cite{Bou1}, Proposition 18. In fact its is a bit more involved: In the definition of tropicalization we removed some of the bivalent vertices from the dual intersection graph $\Gamma_C$. For each compact edge $E$ in $\Gamma_C$ we have a choice of $w_E$-th root of unity as follows. Locally at a node $q$ we have $C=\text{Spec }\mathbb{C}[u,v]/(uv)$ and $\tilde{\mathfrak{X}}_d=\text{Spec }\mathbb{C}[x,y,w^{\pm 1},t]/(xy-t^\ell)$ for some $\ell\in\mathbb{Z}_{>0}$, so $Y:=\tilde{X}_{d,0}=\text{Spec }\mathbb{C}[x,y,w^{\pm 1}]/(xy)$. Locally at $q$ a log structure on $C$ is given by a commutative diagram
\begin{equation*}
\begin{xy}
\xymatrix{
f^{-1}\mathcal{M}_Y \ar[r]^{f^\#}\ar[d]^{\alpha_Y} 	& \mathcal{M}_C \ar[d]^{\alpha_C} \\
f^{-1}\mathcal{O}_Y \ar[r]^{f^\star} 				& \mathcal{O}_C
}
\end{xy}
\end{equation*}
For any $w_E$-th root of unity $\zeta$ there is a chart for $\mathcal{M}_C$ locally at $q$ given by
\[ S_\ell \rightarrow \mathcal{O}_C, \left((a,b,c)\right) \mapsto \begin{cases} (\zeta^{-1}u)^av^b & c = 0, \\ 0 & c \neq 0. \end{cases} \]
Here $S_\ell = \mathbb{N}^2 \oplus_{\mathbb{N}} \mathbb{N}$ with $\mathbb{N} \rightarrow \mathbb{N}^2$ the diagonal embedding and $\mathbb{N} \rightarrow \mathbb{N}, 1 \mapsto \ell$ (see Construction \ref{con:tropical}, (2)). None of these choices are identified via a scheme theoretically trivial isomorphism and all possible extensions are of the above form (see \cite{Gr10}, Proposition 4.23, Step 2). 

Now consider a chain of edges of $\Gamma_C$ connected by bivalent vertices not mapping to vertices of $\mathscr{P}_d$. Then these bivalent vertices get removed by producing the tropicalization and the chain of edges is replaced by a single edge $E$. In this case there are some isomorphisms between the above stable log maps that are not scheme theoretically trivial. Up to such isomorphisms there are exactly $w_E$ stable log maps (see \cite{Gr10}, Proposition 4.23, Step 3). So we really only get one factor of $w_E$ for each edge in the tropicalization. The log structure at general points and marked points is uniquely determined (see \cite{Gr10}, Proposition 4.23, Step 1). This gives the nominator of \textup{deg}(\textup{cut}) as above.

There are further isomorphisms of stable log maps given by the action of $l_\Gamma$-th roots of unity on the base of the curves (see \cite{KLR}, discussion before (6.13). This gives the denominator of \textup{deg}(\textup{cut}) as above.
\end{proof}

\begin{prop}[Gluing formula]
\label{prop:gluing}
We have
\[ \textup{cut}_\star\llbracket\mathscr{M}_{\tilde{h}}\rrbracket = \frac{1}{\ell_{\tilde{\Gamma}}}\prod_{E\in E(\tilde{\Gamma})}w_E \cdot \delta^!\prod_{V\in V(\tilde{\Gamma})}\llbracket\mathscr{M}_V^\circ\rrbracket. \]
\end{prop}

\begin{proof}
By compatibility of obstruction theories (see \cite{KLR}, {\S}9, and \cite{Bou1}, {\S}7.3) we have
\[ \llbracket\mathscr{M}_{\tilde{h}}\rrbracket = \text{cut}^\star \delta^!\prod_{V\in V(\tilde{\Gamma})}\llbracket\mathscr{M}_V^\circ\rrbracket. \]
By the projection formula, $\text{cut}_\star\text{cut}^\star$ is multiplication with $\text{deg}(\text{cut})$ which is $\frac{1}{\ell_{\tilde{\Gamma}}}\prod_{E\in E(\tilde{\Gamma})}w_E$ by Lemma \ref{lem:degcut}.
\end{proof}

\begin{prop}
\label{prop:loop}
We have
\[ \int_{\llbracket\mathscr{M}_{\tilde{h}}\rrbracket}1 = \frac{1}{\ell_{\tilde{\Gamma}}}\prod_{E\in E(\tilde{\Gamma})}w_E \cdot \int_{\delta^!\prod_V\llbracket\mathscr{M}_V\rrbracket}1. \]
\end{prop}

\begin{proof}
By Proposition \ref{prop:gluing} the cycles $\text{cut}_\star\llbracket\mathscr{M}_{\tilde{h}}\rrbracket$ and $\frac{1}{\ell_{\tilde{\Gamma}}}\prod_{E\in E(\tilde{\Gamma})}w_E \cdot \delta^!\prod_V\llbracket\mathscr{M}_V\rrbracket$ have the same restriction to the open substack $\bigtimes_V\mathscr{M}_V^\circ$ of $\bigtimes_V\mathscr{M}_V$. Hence by \cite{Ful}, Proposition 1.8, their difference is rationally equivalent to a cycle supported on the closed substack $Z:=(\bigtimes_V\mathscr{M}_V)\setminus(\bigtimes_V\mathscr{M}_V^\circ)$. Suppose there exists an element $(f_V : C_V \rightarrow Y_{\tilde{h}(V)})_{V\in V(\tilde{\Gamma})}\in Z$. Then by the loop construction in the proof of Lemma \ref{lem:proper} at least one of the source curves $C_V$ would contain a nontrivial cycle of components, contradicting $g=0$. So $Z$ is empty, completing the proof.
\end{proof}

\begin{prop}[Identifying the pieces]
\label{prop:pieces}
We have
\[ \int_{\delta^!\prod_{V}\llbracket\mathscr{M}_V\rrbracket} 1 = \prod_{V\in V(\tilde{\Gamma})}N_V. \]
\end{prop}

\begin{proof}
This is similar to the proof of \cite{Bou1}, Proposition 22. By definition of $\delta$ we have
\[ \int_{\delta^!\prod_{V}\llbracket\mathscr{M}_V\rrbracket} 1 = \int_{\prod_{V}\llbracket\mathscr{M}_V\rrbracket} \textup{ev}^\star[\delta], \]
where $[\delta]$ is the class of the diagonal $\prod_E D_E$. Since each $D_E$ is a projective line, we have
\[ [\delta] = \prod_{E\in E(\tilde{\Gamma})}(\textup{pt}_E \times 1 + 1 \times \textup{pt}_E). \]
As before we give $\tilde{\Gamma}$ the structure of a rooted tree by choosing the root vertex to be the vertex $V_{\textup{out}}$ of the unique unbounded leg $E_{\textup{out}}$. For a bounded edge $E$ let $V_{E,+}$ and $V_{E,-}$ be the vertices of $E$ such that $V_{E,+}$ is the parent of $V_{E,-}$. 

We will show by dimensional arguments that the only term of
\[ \textup{ev}^\star[\delta] = \prod_{E\in E(\tilde{\Gamma})} \left((\textup{ev}_{V_{E,-}})^\star[\textup{pt}_E] + (\textup{ev}_{V_{E,+}})^\star[\textup{pt}_E]\right) \]
giving a nonzero contribution after integration over $\prod_{V}\llbracket\mathscr{M}_V\rrbracket$ is $\prod(\textup{ev}_{V_{E,+}})^\star[\textup{pt}_E]$. In other words:

\underline{Claim:} For each compact edge $E$, a term of $\textup{ev}^\star[\delta]$ giving a nonzero contribution after integration over $\prod_{V}\llbracket\mathscr{M}_V\rrbracket$ does not contain a factor $(\textup{ev}_{V_{E,-}})^\star[\textup{pt}_E]$.

Let $E$ be a compact edge with $V_{E,-}$ a vertex of type (II) as in Definition \ref{defi:types}. By Proposition \ref{prop:N}, (II), the virtual dimension of $\mathscr{M}_{V_{E,-}}$ is zero. Hence, $(\textup{ev}_{V_{E,-}})^\star[\textup{pt}_E]=0$, since its insertion over $\mathscr{M}_{V_{E,-}}$ defines an enumerative problem of virtual dimension $-1$.

Now consider a compact edge $E$ with $V_{E,-}$ of type (III). Let $E_i, i\in I$ be the edges adjacent to $V_{E,-}$ and different from $E$ (possibly $I=\emptyset$). By Proposition \ref{prop:balancing}, (III), the edges $E_i$ connect $V_{E,-}$ with a vertex $V_i$ of type (II). The terms in $\textup{ev}^\star[\delta]$ containing a factor $(\textup{ev}_{V_{E_i,-}})^\star[\textup{pt}_{E_i}]$ give zero after integration over $\llbracket\mathscr{M}_{V_i}\rrbracket$ by the dimensional argument above. Hence, to give a nonzero contribution, a term of $\textup{ev}^\star[\delta]$ must contain the factor $\prod_{i\in I}(\textup{ev}_{V_{E_i,+}})^\star[\textup{pt}_{E_i}]$. By Proposition \ref{prop:N}, (III), the virtual dimension of $\mathscr{M}_V$ is $|I|$, so the insertion of $\prod_{i\in I}(\textup{ev}_{V_{E_i,+}})^\star[\textup{pt}_{E_i}]$ in $\mathscr{M}_V$ defines an enumerative problem of virtual dimension $0$. Any further insertion would reduce the virtual dimension to $-1$, so a term of $\textup{ev}^\star[\delta]$ giving a nonzero contribution does not contain the factor $(\textup{ev}_{V_{E,-}})^\star[\textup{pt}_E]$.

We will show the claim for compact edges $E$ with $V_{E,-}$ a vertex of type (I) by induction on the height of $V_{E,-}$, that is, the maximal length of chains connecting $V_{E,-}$ with a leaf of $\tilde{\Gamma}$. By Proposition \ref{prop:balancing}, (I), a vertex of type (I) fulfills the ordinary balancing condition. In particular, it must have more than one adjacent edge, hence cannot be a leaf. This shows the set of leaves of $\tilde{\Gamma}$ is contained in $V_{II}(\tilde{\Gamma})\cup V_{III}(\tilde{\Gamma})$. Thus we have already shown that the claim is true for compact edges $E$ with $V_{E,-}$ of height $0$. This is the base case. For the induction step assume that the claim is true for compact edges $E$ with $V_{E,-}$ of height $\leq k$ for some $k\in\mathbb{N}$ and consider a compact edge $E$ with $V_{E,-}$ of height $k+1$. Assume that $V_{E,+}$ is of type (I), since otherwise the claim is true by the above arguments. Let $E_i,i\in I$ be the edges connecting $V_{E,-}$ with its childs. By Proposition \ref{prop:N}, (I), the virtual dimension of $\mathscr{M}_{V_{E,-}}$ is $|I|$. By the induction hypothesis, a term of $\textup{ev}^\star[\delta]$ giving a nonzero contribution must contain the factor $\prod_{i\in I}(\textup{ev}_{V_{E_i,+}})^\star[\textup{pt}_{E_i}]$. Inserting this factor over $\mathscr{M}_{V_{E,-}}$ gives an enumerative problem of virtual dimension $0$. Again, for dimensional reasons, a term of $\textup{ev}^\star[\delta]$ giving a nonzero contribution cannot contain the factor $(\textup{ev}_{V_{E,-}})^\star[\textup{pt}_E]$, hence it must contain the factor $(\textup{ev}_{V_{E,+}})^\star[\textup{pt}_E]$. This proves the claim. Now
\[ \int_{\prod_{V\in V(\tilde{\Gamma})}\llbracket\mathscr{M}_V\rrbracket} \prod_{E\in E(\tilde{\Gamma})}(\textup{ev}_{V_{E,+}})^\star[\textup{pt}_E] = \prod_{V\in V(\tilde{\Gamma})} \int_{\llbracket\mathscr{M}_V\rrbracket}(\textup{ev}_{E\rightarrow V})^\star [\textup{pt}] = \prod_{V\in V(\tilde{\Gamma})}N_V, \]
completing the proof.
\end{proof}

\subsection{The degeneration formula}							%%%

Combining the decomposition formula and the gluing formula, we obtain the \textit{degeneration formula}, expressing $N_\beta$ in terms of logarithmic Gromov-Witten invariants $N_V$ labeled by vertices of tropical curves. 

\begin{prop}[Degeneration formula]
\label{prop:deg}
\[ N_\beta = \sum_{\tilde{h}\in\tilde{\mathfrak{H}}_\beta} \frac{1}{|\textup{Aut}(\tilde{h})|} \cdot \prod_{E\in E(\tilde{\Gamma})}w_E \cdot \prod_{V\in V(\tilde{\Gamma})} N_V. \]
\end{prop}

\begin{proof}
Since the virtual dimension of $\mathscr{M}_\beta$ is zero, integration (i.e., proper pushforward to a point) of the decomposition formula (Proposition \ref{prop:dec}) gives
\[ N_\beta = \sum_{\tilde{h}\in\tilde{\mathfrak{H}}_\beta} \frac{1}{|\textup{Aut}(\tilde{h})|}\int_{\llbracket\mathscr{M}_{\tilde{h}}\rrbracket} 1. \]
Using Propositions \ref{prop:loop} and \ref{prop:pieces} we get the above formula.
\end{proof}

As mentioned earlier, summation over balanced tropical curves in $\mathfrak{H}_\beta$ will give a more symmetric version of the above formula:

\begin{defi}
\label{defi:Ntor}
Let $h : \Gamma \rightarrow B$ be a tropical curve in $\mathfrak{H}_\beta$ and let $V$ be a vertex of $\Gamma$. Then the image of $V$ under the map from Construction \ref{con:tilde} is a vertex of $\tilde{\Gamma}$ of type (I) or (III). Let $\textbf{m}$ and $\textbf{w}$ be as in the respective case of Proposition \ref{prop:N} and define 
$ N_V^{\textup{tor}} := N_{\textbf{m}}(\textbf{w}). $
Note that $N_V^{\textup{tor}}=N_V$ for vertices of type (I).
\end{defi}

\begin{defi}
\label{defi:Nh}
For a tropical curve $h : \Gamma \rightarrow B$ in $\mathfrak{H}_d$ for some $d$ define
\[ N_h := \left(\frac{1}{|\textup{Aut}(h)|} \cdot \prod_{E\in E(\Gamma)}w_E\cdot \prod_{E\in L_\Delta(\Gamma)}\frac{(-1)^{w_E-1}}{w_E}\cdot\prod_{V\in V(\Gamma)} N_V^{\textup{tor}}\right), \]
where $L_\Delta(\Gamma)$ is the set of bounded legs (see Definition \ref{defi:tropical}).
\end{defi}

\begin{thm}[Symmetric version of the degeneration formula]
\label{thm:degmax}
\[ N_\beta = \sum_{h\in\mathfrak{H}_\beta} N_h. \]
\end{thm}

\begin{proof}
Using Propositions \ref{prop:deg} and \ref{prop:N}, we have
\begin{eqnarray*}
N_d &=& \sum_{\tilde{h}\in\tilde{\mathfrak{H}}_\beta}\left( \frac{1}{|\textup{Aut}(\tilde{h})|} \cdot \prod_{E\in E(\tilde{\Gamma})} w_E \cdot \prod_{V\in V_I(\tilde{\Gamma})}N_V^{\textup{tor}} \cdot \prod_{V\in V_{II}(\tilde{\Gamma})}\frac{(-1)^{w_{E_V}-1}}{w_{E_V}^2}\right. \\
&& \textup{ } \hspace{3cm} \left.\cdot \prod_{V\in V_{III}(\tilde{\Gamma})}\left(\sum_{\textbf{w}_{V,+}} N_V^{\textup{tor}} \frac{1}{|\textup{Aut}(\textbf{w}_{V,+})|}\prod_{i=1}^{l_V}\frac{(-1)^{w_{V,i}-1}}{w_{V,i}}\right)\right).
\end{eqnarray*}
Canceling the $w_{E_V}$ for vertices of type (II) in the first product against the ones in the denominator of the third product and factoring out the second sum we get
\begin{eqnarray*}
N_d &=& \sum_{\tilde{h}\in\tilde{\mathfrak{H}}_\beta}\sum_{(\textbf{w}_{V,+})_{V\in V_{III}(\tilde{\Gamma})}} \left(\frac{1}{|\textup{Aut}(\tilde{h})||\textup{Aut}(\textbf{w}_{V,+})|} \cdot \prod_{E\in E(\tilde{\Gamma})\setminus\cup_{V\in V_{II}(\tilde{\Gamma})}\{E_V\}}w_E  \right. \\
&& \textup{ } \hspace{0cm} \left. \cdot \prod_{V\in V_I(\tilde{\Gamma})\cup V_{III}(\tilde{\Gamma})} N_V^{\textup{tor}} \cdot \prod_{V\in V_{II}(\tilde{\Gamma})}\frac{(-1)^{w_{E_V}-1}}{w_{E_V}}\cdot \prod_{V\in V_{III}(\tilde{\Gamma})}\prod_{i=1}^{l_V}\frac{(-1)^{w_{V,i}-1}}{w_{V,i}}\right).
\end{eqnarray*}
Now by the construction of the map $\mathfrak{H}_d \rightarrow \tilde{\mathfrak{H}}_d$ in Construction \ref{con:tilde}, the two summations can be replaced by a summation over $\mathfrak{H}_d$. Note that for $\tilde{h}\in\tilde{\mathfrak{H}}_d$ we have
\[ \sum_{h\mapsto\tilde{h}}\frac{1}{|\textup{Aut}(h)|} = \frac{1}{|\textup{Aut}(\tilde{h})|}\sum_{(\textbf{w}_{V,+})_{V\in V_{III}(\tilde{\Gamma})}} \frac{1}{|\textup{Aut}(\textbf{w}_{V,+})|}, \]
where the sum is over all $h\in\mathfrak{H}_d$ giving $\tilde{h}$ via the map from Construction \ref{con:tilde}. This can be seen by multiplying both sides with $|\text{Aut}(\tilde{h})|$. 

Moreover, note that $V(\Gamma)=V_I(\tilde{\Gamma})\cup V_{III}(\tilde{\Gamma})$ and $E(\Gamma)=E(\tilde{\Gamma})\setminus\cup_{V\in V_{II}(\tilde{\Gamma})}\{E_V\}$, where,  for a vertex $V$ of type (II), $E_V$ is the unique edge containing the vertex $V$. Then
\begin{eqnarray*}
N_d &=& \sum_{h\in\mathfrak{H}_\beta}\left(\frac{1}{|\textup{Aut}(h)|} \cdot \prod_{E\in E(\Gamma)} w_E \cdot \prod_{V\in V(\Gamma)} N_V^{\textup{tor}} \right. \\
&& \textup{ } \hspace{2cm} \left. \cdot \prod_{V\in V_{II}(\tilde{\Gamma})}\frac{(-1)^{w_{E_V}-1}}{w_{E_V}}\cdot \prod_{V\in V_{III}(\tilde{\Gamma})}\prod_{i=1}^{l_V}\frac{(-1)^{w_{V,i}-1}}{w_{V,i}} \right).
\end{eqnarray*}
Using that
\[ \prod_{E\in L_\Delta(\Gamma)} \frac{(-1)^{w_E-1}}{w_E} = \prod_{V\in V_{II}(\tilde{\Gamma})}\frac{(-1)^{w_{E_V}-1}}{w_{E_V}}\cdot \prod_{V\in V_{III}(\tilde{\Gamma})}\prod_{i=1}^{l_V}\frac{(-1)^{w_{V,i}-1}}{w_{V,i}} \]
completes the proof.
\end{proof}

\begin{cor}
\[ N_d = \sum_{h\in\mathfrak{H}_d} N_h. \]
\end{cor}

\begin{proof}
This follows from Theorem \ref{thm:degmax} and $\mathfrak{H}_d = \coprod_{\substack{\underline{\beta}\in H_2^+(X,\mathbb{Z}) \\ D\cdot\underline{\beta}=dw_{\text{out}}}} \mathfrak{H}_\beta$.
\end{proof}

\begin{defi}
\label{defi:mult}
Let $h : \Gamma \rightarrow B$ be a tropical curve. For a trivalent vertex $V\in V(\Gamma)$ define $m_V=\lvert u_{(V,E_1)}\wedge u_{(V,E_2)}\rvert=\lvert\text{det}(u_{(V,E_1)}|u_{(V,E_2)})\rvert$, where $E_1,E_2$ are any two edges adjacent to $V$. For a vertex $V\in V(\Gamma)$ of valency $\nu_V>3$ let $h_V$ be the one-vertex tropical curve describing $h$ locally at $V$ and let $h'_V$ be a deformation of $h_V$ to a trivalent tropical curve. It has $\nu_V-2$ vertices. Define $m_V=\prod_{V'\in V(h'_V)}m_{V'}$. For a bounded leg $E\in L_\Delta(\Gamma)$ define $m_E=(-1)^{w_E+1}/w_E^2$. Then define the \textit{multiplicity} of $h$ to be
\[ \text{Mult}(h) = \frac{1}{|\text{Aut}(h)|} \cdot \prod_V m_V \cdot \prod_{E\in L_\Delta(\Gamma)} m_E. \]
\end{defi}

\begin{prop}
For a tropical curve $h : \Gamma \rightarrow B$ in $\mathfrak{H}_d$ we have
\[ N_h = \textup{Mult}(h) \]
\end{prop}

\begin{proof}
By the tropical correspondence theorem with point conditions on toric divisors (\cite{GPS}, Theorem 3.4) we have $m_V = \prod_{E\rightarrow V}w_E \cdot N_V^{\text{tor}}$. The product is over all eges of $\Gamma$ pointing towards $V$ with respect to the orientation of $\Gamma$ such that all edges point towards the root vertex $V_{\text{out}}$. Then $\prod_{V\in V(\Gamma)}m_V = \prod_{E\in E(\Gamma)} w_E \cdot \prod_{V\in V(\Gamma)} N_V^{\text{tor}}$ as each $E\in E(\Gamma)$ occurs exactly once. Plugging this and $m_E=(-1)^{w_E+1}/w_E^2$ into the definition of $\text{Mult}(h)$ we obtain $N_h$.
\end{proof}

Together with Theorem \ref{thm:degmax} this gives the tropical correspondence theorem (Theorem \ref{thm:trop}).

\subsection{Invariants with prescribed degree splitting}	%%

\begin{defi}
As in \S\ref{S:degdiv} fix a cyclic labelling of $D_t^0=D_1+\ldots+D_k$. For $[d_1,\ldots,d_k] \in \mathbb{N}^k$ define
\[ N_{[d_1,\ldots,d_k]} = \sum_{h\mapsto[d_1,\ldots,d_k]} N_h, \]
where the sum is over all $h \in \mathfrak{H}_d$ with degree splitting $[d_1,\ldots,d_k]$ (Definition \ref{defi:splitting}). 
\end{defi}

\begin{rem}
Barrott and Nabijou \cite{BN} define invariants with prescribed degree splittings by looking at the family $\mathfrak{X}_{t\neq 0} \rightarrow \mathbb{A}^1$ only degenerating the divisor and using torus localization. We conjecture that these coincide with the invariants defined above. This question will be investigated in future work.
\end{rem}

\section{Scattering calculations}								%%%
\label{S:scattering}

In this section we recall the notions of wall structures and scattering diagrams from \cite{GS11}, restricting to the $2$-dimensional case with trivial gluing data. We then explain how the dual intersection complex $(B,\mathscr{P},\varphi)$ of the toric degeneration $(\mathfrak{X},\mathfrak{D})\rightarrow\mathbb{A}^1$ of $(X,D)$ defines a consistent wall structure $\mathscr{S}_\infty$. For toric varieties the correspondence between scattering diagrams (wall structures with one vertex) and logarithmic Gromov-Witten invariants was shown in \cite{GPS}. We use the relation between scattering diagrams and tropical curves from \cite{GPS}, Theorem 2.8, (see Lemma \ref{lem:2.8}) to extend this correspondence to our non-toric case.

\subsection{Scattering diagrams}							%%
\label{S:scatteringdiag}

Let $M\simeq\mathbb{Z}^2$ be a lattice and write $M_{\mathbb{R}}=M\otimes_{\mathbb{Z}}\mathbb{R}$. Let $\Sigma$ be a fan in $M_{\mathbb{R}}$ and let $\varphi$ be an integral strictly convex piecewise affine function on $\Sigma$ with $\varphi(0)=0$. Note that the rays of $\Sigma$ form the corner locus of $\varphi$. Let $P_\varphi$ be the monoid of integral points in the upper convex hull of $\varphi$,
\[ P_\varphi = \{p=(\overline{p},h)\in M\oplus\mathbb{Z} \ | \ h \geq \varphi(\overline{p})\}. \]
Write $t:=z^{(0,1)}$ and let $R_\varphi$ be the $\mathbb{C}\llbracket t\rrbracket$-algebra obtained by completion of $\mathbb{C}[P_\varphi]$ with respect to $(t)$,
\[ R_\varphi = \varprojlim\mathbb{C}[P_\varphi]/(t^k). \]

\begin{defi}
A \textit{ray} for $\varphi$ is a half-line $\mathfrak{d} = \mathbb{R}_{\geq 0}\cdot m_{\mathfrak{d}} \subseteq M_{\mathbb{R}}$, with $m_{\mathfrak{d}}\in M$ primitive, together with an element $f_{\mathfrak{d}}\in R_\varphi$ such that
\begin{compactenum}[(1)]
\item each exponent $p=(\overline{p},h)$ in $f_{\mathfrak{d}}$ satisfies $\overline{p}\in\mathfrak{d}$ or $-\overline{p}\in\mathfrak{d}$. In the first case the ray is called \textit{incoming}, in the latter it is called \textit{outgoing};
\item if $m_{\mathfrak{d}}$, is a ray generator of $\Sigma$, then $f_{\mathfrak{d}} \equiv 1 \textup{ mod }(z^{m_{\mathfrak{d}}})$;
\item if $m_{\mathfrak{d}}$, is not a ray generator of $\Sigma$, then $f_{\mathfrak{d}} \equiv 1 \textup{ mod }(z^{m_{\mathfrak{d}}}t)$.
\end{compactenum}
A \textit{scattering diagram} for $\varphi$ is a set $\mathfrak{D}$ of rays for $\varphi$ such that for every power $k>0$ there are only a finite number of rays $(\mathfrak{d},f_{\mathfrak{d}})\in\mathfrak{D}$ with $f_{\mathfrak{d}} \not\equiv 1 \textup{ mod }(t^k)$.
\end{defi}

Let $\mathfrak{D}$ be a scattering diagram for $\varphi$. Let $\gamma : [0,1] \rightarrow M_{\mathbb{R}}$ be a closed immersion not meeting the origin and with endpoints not contained in any ray of $\mathfrak{D}$. Then, for each power $k>0$, we can find numbers $0<r_1\leq r_2\leq\ldots\leq r_s<1$ and rays $\mathfrak{d}_i=(\mathbb{R}_{\geq 0}m_i,f_i)\in\mathfrak{D}$ with $f_i \not\equiv 1 \textup{ mod }(t^k)$ such that (1) $\gamma(r_i)\in\mathfrak{d}_i$, (2) $\mathfrak{d}_i\neq\mathfrak{d}_j$ if $r_i=r_j$ and $i\neq j$, and (3) $s$ is taken as large as possible. 

For each ray $\rho\in\Sigma^{[1]}$ write $f_{\rho}=1+z^{(m_\rho,\varphi(m_\rho))}$ for $m_\rho\in M$ the primitive generator of $\rho$, and define
\[ \tilde{R}_\varphi^k = \left(R_\varphi/(t^{k+1})\right)_{\prod_{\rho}f_\rho}. \]
For each $i$, define a $\mathbb{C}\llbracket t\rrbracket$-algebra automorphism of $\tilde{R}_\varphi^k$ by $\theta_{\mathfrak{d}_i}=\text{exp}(\text{log}(f_i)\partial_{n_i})$ for $\partial_n(z^p):=\braket{n,\overline{p}}z^p$, i.e.,
\[ \theta_{\mathfrak{d}_i}^k(z^p) = f_i^{-\braket{n_i,\overline{p}}}z^p, \]
where $n_i\in N=\text{Hom}(M,\mathbb{Z})$ is the unique primitive vector satisfying $\braket{n_i,m_i} = 0$ and $\braket{n_i,\gamma'(r_i)} > 0$. Define
\[ \theta_{\gamma,\mathfrak{D}}^k = \theta^k_{\mathfrak{d}_s} \circ \ldots \circ \theta^k_{\mathfrak{d}_1}. \]
If $r_i=r_j$, then $\theta_{\mathfrak{d}_i}$ and $\theta_{\mathfrak{d}_j}$ commute. Hence, $\theta_{\gamma,\mathfrak{D}}^k$ is well-defined. Moreover, define
\[ \theta_{\gamma,\mathfrak{D}} = \lim_{k\rightarrow\infty}\theta_{\gamma,\mathfrak{D}}^k. \]

\begin{defi}
A scattering diagram $\mathfrak{D}$ is \textit{consistent to order $k$} if, for all $\gamma$ such that $\theta_{\gamma,\mathfrak{D}}^k$ is defined,
\[ \theta_{\gamma,\mathfrak{D}}^k \equiv 1 \textup{ mod }(t^{k+1}). \]
It is \textit{consistent to any order}, or simply \textit{consistent}, if $\theta_{\gamma,\mathfrak{D}} = 1$.
\end{defi}

\begin{prop}[\cite{GPS}, Theorem 1.4]
\label{prop:scatt}
Let $\mathfrak{D}$ be a scattering diagram such that $f_{\mathfrak{d}} \equiv 1 \text{ mod } (t)$ for all $\mathfrak{d}\in\mathfrak{D}$. Then there exists a consistent scattering diagram $\mathfrak{D}_\infty$ containing $\mathfrak{D}$ such that $\mathfrak{D}_\infty\setminus\mathfrak{D}$ consists only of outgoing rays.
\end{prop}

\begin{proof}
The proof is constructive, so we will give it here. 

Take $\mathfrak{D}_0=\mathfrak{D}$. We will show inductively that there exists a scattering diagram $\mathfrak{D}_k$ containing $\mathfrak{D}_{k-1}$ that is consistent to order $k$. Let $\mathfrak{D}'_{k-1}$ consist of those rays $\mathfrak{d}$ in $\mathfrak{D}_{k-1}$ with $f_{\mathfrak{d}} \not\equiv 1 \text{ mod } (t^{k+1})$. Let $\gamma$ be a closed simple loop around the origin. Then $\theta_{\gamma,\mathfrak{D}_{k-1}} \equiv \theta_{\gamma,\mathfrak{D}'_{k-1}} \text{ mod } (t^{k+1})$. By the induction hypothesis this can be uniquely written as
\[ \theta_{\gamma,\mathfrak{D}'_{k-1}} = \text{exp}\left(\sum_{i=1}^s c_iz^{m_i}\partial_{n_i}\right) \]
with $m_i\in M\setminus\{0\}$m $n_i\in m_i^\bot$ primitive and $c_i\in (t^k)$. Define
\[ \mathfrak{D}_k = \mathfrak{D}_{k-1} \cup \left\{(\mathbb{R}_{\geq 0}m_i,1\pm c_iz^{m_i}) \ | \ i=1,\ldots,s\right\}, \]
with sign chosen in each ray such that its contribution to $\theta_{\gamma,\mathfrak{D}_k}$ is $\text{exp}(c_iz^{m_i}\partial_{n_i})$ modulo $(t^{k+1})$. These contributions exactly cancel the contributions to $\theta_{\gamma,\mathfrak{D}_k}$ coming from $\mathfrak{D}_{k-1}$, so $\theta_{\gamma,\mathfrak{D}_k} \equiv 1 \text{ mod }(t^{k+1})$.

Then take $\mathfrak{D}_\infty$ to be the non-disjoint union of the $\mathfrak{D}_k$ for all $k\in\mathbb{N}$. The diagram $\mathfrak{D}_\infty$ will usually have infinitely many rays.
\end{proof}

\begin{defi}
Two scattering diagrams $\mathfrak{D}$, $\mathfrak{D}'$ are \textit{equivalent} if $\theta_{\gamma,\mathfrak{D}}=\theta_{\gamma,\mathfrak{D}'}$ for any closed immersion $\gamma$ for which both sides are defined.
\end{defi}

\begin{defi}
A scattering diagram $\mathfrak{D}$ is called \textit{minimal} if 
\begin{compactenum}[(1)]
\item any two rays $\mathfrak{d},\mathfrak{d}'$ in $\mathfrak{D}$ have distinct support , i.e., $m_{\mathfrak{d}}\neq m_{\mathfrak{d}'}$;
\item it contains no \textit{trivial} ray, i.e., with $f_{\mathfrak{d}}=1$.
\end{compactenum}
\end{defi}

\begin{rem}
Every scattering diagram $\mathfrak{D}$ is equivalent to a unique minimal scattering diagram. In fact, if $\mathfrak{d},\mathfrak{d}'$ have the same support, then we can replace these two rays with a single ray with the same support and attached function $f_{\mathfrak{d}}\cdot f_{\mathfrak{d}'}$. Moreover, we can remove any trivial ray without affecting $\theta_{\gamma,\mathfrak{D}}$.
\end{rem}

\subsection{Scattering diagrams and toric invariants}					%%

In \cite{GPS}, Theorem 2.4, a bijective correspondence between certain scattering diagrams and tropical curves is established, leading to an enumerative correspondence (\cite{GPS}, Theorem 2.8). Combining this result with the tropical correspondence theorem for torically transverse stable log maps with point conditions on the toric boundary (\cite{GPS}, Theorems 3.4, 4.4), we get the following.

\begin{lem}
\label{lem:2.8}
Let $\textbf{m}=(m_1,\ldots,m_n)$ be an $n$-tuple of (not necessarily distinct) primitive vectors of $M$. Let $\Sigma$ be a fan in $M_{\mathbb{R}}$ and let $\varphi$ be an integral strictly convex piecewise affine function on $\Sigma$ such that $\varphi(0)=0$. Let $\mathfrak{D}$ be a scattering diagram for $\varphi$ consisting of a number of lines\footnote{This means that the rays in $\mathfrak{D}$ come in pairs, one incoming and one outgoing, with the same attached function.}, one in each direction $m_i$ and with attached function $f_i\in\mathbb{C}[z^{(-m_i,0)}]\subseteq P_\varphi$. Write the logarithm of $f_i$ as
\[ \textup{log }f_i = \sum_{w=1}^\infty a_{iw}z^{(-wm_i,0)}, \quad a_{iw}\in\mathbb{C}. \]
Let $\mathfrak{D}_\infty$ be the associated minimal consistent scattering diagram and let $\mathfrak{d}\in \mathfrak{D}_\infty\setminus\mathfrak{D}$ be a ray in direction $m_{\mathfrak{d}}$ with attached function $f_{\mathfrak{d}}$. Then
\[ \textup{log }f_{\mathfrak{d}} = \sum_{w=1}^\infty\sum_{\textbf{w}}w\frac{N_{\textbf{m}}(\textup{\textbf{w}})}{|\textup{Aut}(\textup{\textbf{w}})|} \left(\prod_{\substack{1\leq i\leq n\\ 1\leq j\leq l_i}} a_{iw_{ij}}\right) z^{(-wm_{\mathfrak{d}},0)}, \]
where the sum is over all $n$-tuples of weight vectors $\textbf{w}=(\textbf{w}_1,\ldots,\textbf{w}_n)$ satisfying
\[ \sum_{i=1}^n |\textbf{w}_i|m_i = wm_{\mathfrak{d}}. \]
Here $N_{\textbf{m}}(\textbf{w})$ is the toric logarithmic Gromov-Witten invariant defined in \S\ref{S:toricinv} and $\text{Aut}(\textbf{w})$ is the subgroup of the permutation group $S_n$ stabilizing $(w_1,\ldots,w_n)$.

Moreover, let $m\in\mathbb{Q}_{>0}m_1+\ldots+\mathbb{Q}_{>0}m_n$ be a primitive vector. If there is no ray $\mathfrak{d}\in\mathfrak{D}_\infty$ in direction $m$, then $N_{\textbf{m}}(\textbf{w})=0$ for all $\textbf{w}$ satisfying $\sum_{i=1}^n |\textbf{w}_i|m_i = wm_{\mathfrak{d}}$.
\end{lem}

\begin{rem}
Note that \cite{GPS} deals with the case $\varphi\equiv 0$, where the $t$-order of an element $z^{(\overline{p},h)}$ is simply given by $h$. In our case, the $t$-order is $\varphi(-\overline{p})+h\geq 0$:
\[ z^{(\overline{p},h)} = \left(z^{(-\overline{p},\varphi(-\overline{p}))}\right)^{-1} t^{\varphi(-\overline{p})+h}. \]
The formula in \cite{GPS}, Theorem 2.8, contains some explicit $t$-factors. These are not visible in Lemma \ref{lem:2.8} due to this different notion of $t$-order.
\end{rem}

\subsection{Wall structures}								%%

Let $(B,\mathscr{P},\varphi)$ be a $2$-dimensional polarized polyhedral affine manifold. Note that for each $x\in B\setminus\Delta$ (where $\Delta$ is the singular locus of $B$), $\varphi$ defines an integral strictly convex piecewise affine function
\[ \varphi_x : \Lambda_{B,x} \simeq M \rightarrow \mathbb{R} \]
on $\Sigma_x$, the fan describing $(B,\mathscr{P})$ locally at $x$. If $\tau_x\in\mathscr{P}$ is the smallest cell containing $x$, then this is given by
\[ \Sigma_x = \{K_{\tau_x}\sigma \ | \ \tau_x\subseteq\sigma\in\mathscr{P}\}, \]
where $K_{\tau_x}\sigma$ is the cone generated by $\sigma$ relative to $\tau$, i.e.,
\[ K_{\tau_x}\sigma = \mathbb{R}_{\geq 0}(\sigma-\tau_x) = \{m \in M_{\mathbb{R}} \ | \ \exists m_0\in\tau_x, m_1\in\sigma,\lambda\in\mathbb{R}_{\geq 0} : m=\lambda(m_1-m_0)\}. \]
As in \S\ref{S:scatteringdiag} this defines a monoid by the integral points in the upper convex hull of $\varphi_x$,
\begin{equation} 
\label{eq:Px}
P_x := P_{\varphi_x} = \{p=(\overline{p},h)\in\Lambda_{B,x}\oplus\mathbb{Z} \ | \ h\geq\varphi_x(\overline{p})\}.
\end{equation}
Note that $\textup{Spec }\mathbb{C}[P_x]$ gives a local toric model for the toric degeneration defined by $(B,\mathscr{P},\varphi)$ at a point on the interior of the toric stratum corresponding to $\tau_x$.

\begin{defi}
\label{defi:structure}
Let $(B,\mathscr{P},\varphi)$ be a $2$-dimensional polarized polyhedral affine manifold such that $(B,\mathscr{P})$ is simple. For $x,x'\in B$ integral points, let $m_{xx'}\in\Lambda_B$ denote the primitive vector pointing from $x$ to $x'$. For a $1$-cell $\rho$ and $x\in\rho\setminus\Delta$ let $v[x]$ be the vertex of the irreducible component of $\rho\setminus\Delta$ containing $x$. This is unique by the construction of the discriminant locus $\Delta$ for polyhedral affine manifolds (\cite{GHS}, Construction 1.1).
\begin{compactenum}[(1)]
\item A \textit{slab} $\mathfrak{b}$ on $(B,\mathscr{P},\varphi)$ is a $1$-dimensional rational polyhedral subset of a $1$-cell $\rho_{\mathfrak{b}}\in\mathscr{P}^{[1]}$ together with elements $f_{\mathfrak{b},x}\in\mathbb{C}[P_x]$, one for each $x\in\mathfrak{b}\setminus\Delta$, satisfying the following conditions.
\begin{compactenum}[(1)]
\item $f_{\mathfrak{b},x} \equiv 1 \textup{ mod } (t)$ if $\rho_{\mathfrak{b}}$ does not contain an affine singularity;
\item $f_{\mathfrak{b},x} \equiv 1+z^{(m_{v[x]\delta},\varphi(m_{v[x]\delta}))} \textup{ mod } (t)$ if $\rho_{\mathfrak{b}}$ contains an affine singularity $\delta$;
\item $f_{\mathfrak{b},x} = z^{(m_{v[x]v[x']},\varphi(m_{v[x]v[x']}))}f_{\mathfrak{b},x'}$ for all $x,x'\in\mathfrak{b}\setminus\Delta$.
\end{compactenum}
Note that conditions (1) and (2) are compatible with (3).
\item A \textit{wall} $\mathfrak{p}$ on $(B,\mathscr{P},\varphi)$ is a $1$-dimensional rational polyhedral subset of a maximal cell $\sigma_{\mathfrak{p}}\in\mathscr{P}^{[2]}$ with $\mathfrak{p}\cap\textup{Int}(\sigma_{\mathfrak{p}})\neq\emptyset$ together with (i) a \textit{base point} \textup{Base}($\mathfrak{p})\in\mathfrak{p}\setminus\partial\mathfrak{p}$, (ii) an \textit{exponent} $p_{\mathfrak{p}} \in \Gamma(\sigma_{\mathfrak{p}},\Lambda \oplus \underline{\mathbb{Z}})$ such that $p_{\mathfrak{p},x}=(\overline{p}_{\mathfrak{p},x},h_{\mathfrak{p},x})\in P_x$ for all $\mathfrak{p}\setminus\Delta$ with $h_{\mathfrak{p},x}>\varphi(\overline{p}_{\mathfrak{p},x})$ for $x\neq\textup{Base}(\mathfrak{p})$, and (iii) $c_{\mathfrak{p}}\in\mathbb{C}$, such that
\[ \mathfrak{p} = (\textup{Base}(\mathfrak{p}) - \mathbb{R}_{\geq 0}\overline{p}_{\mathfrak{p}}) \cap \sigma_{\mathfrak{p}}. \]
For each $x\in\mathfrak{p}\setminus\Delta$ this defines a function
\[ f_{\mathfrak{p},x} = 1+c_{\mathfrak{p}}z^{p_{\mathfrak{p}}} \in \mathbb{C}[P_x]. \]
\item A \textit{wall structure} $\mathscr{S}$ on $(B,\mathscr{P},\varphi)$ is a locally finite set of slabs and walls with a polyhedral decomposition $\mathscr{P}_{\mathscr{S}}$ of its support $|\mathscr{S}|=\cup_{\mathfrak{b}\in\mathscr{S}}\mathfrak{b}$ such that
\begin{compactenum}[(1)]
\item The map sending a slab $\mathfrak{b}\in\mathscr{S}$ to its underlying $1$-cell of $\mathscr{P}$ is injective;
\item Each closure of a connected component of $B\setminus|\mathscr{S}|$ (\textit{chamber}) is convex and its interior is disjoint from any wall;
\item Any wall in $\mathscr{S}$ is a union of elements of $\mathscr{P}_{\mathscr{S}}$;
\item Any maximal cell of $\mathscr{P}$ contains only finitely many slabs or walls in $\mathscr{S}$.
\end{compactenum}
\item A \textit{joint} $\mathfrak{j}$ of a wall structure $\mathscr{S}$ on $(B,\mathscr{P},\varphi)$ is a vertex of $\mathscr{P}_{\mathscr{S}}$. At each joint $\mathfrak{j}$, the wall structure defines a scattering diagram $\mathfrak{D}_{\mathfrak{j}}$ for $\varphi_{\mathfrak{j}}$.
\end{compactenum}
\end{defi}

\begin{defi}
\label{defi:S0}
A polarized polyhedral affine manifold $(B,\mathscr{P},\varphi)$ induces an \textit{initial wall structure} $\mathscr{S}_0$ consisting only of slabs as follows. For each $1$-cell $\rho$ containing an affine singularity $\delta$ there is a slab $\mathfrak{b}$ with underlying polyhedral subset $\rho$ and with
\[ f_{\mathfrak{b},v} = 1+z^{(m_{v\delta},\varphi(m_{v\delta}))}, \]
where $m_{v\delta}\in \Lambda_{B,v}$ is the primitive vector pointing from $v$ to $\delta$.
\end{defi}

\begin{defi}[\cite{GS11}, Definition 2.28]
A wall structure $\mathscr{S}$ is \textit{consistent (to order $k$) at a joint $\mathfrak{j}$} if the associated scattering diagram $\mathfrak{D}_{\mathfrak{j}}$ is consistent (to order $k$). It is \textit{consistent (to order $k$)} if it is consistent (to order $k$) at any joint.
\end{defi}

\begin{defi}[\cite{GS11}, Definition 2.41]
Two wall structures $\mathscr{S}$, $\mathscr{S}'$ are \textit{compatible to order $k$} if the following conditions hold.
\begin{compactenum}[(1)]
\item If $\mathfrak{p}\in\mathscr{S}$ is a wall with $c_{\mathfrak{p}}\neq 0$ and $f_{\mathfrak{p},x} \not\equiv 1 \textup{ mod }(t^{k+1})$ for some $x\in\mathfrak{p}\setminus\Delta$, then $\mathfrak{p}\in\mathscr{S}'$, and vice versa.
\item If $x\in \textup{Int}(\mathfrak{b})\cap\textup{Int}(\mathfrak{b}')$ for slabs $\mathfrak{b}\in\mathscr{S},\mathfrak{b}'\in\mathscr{S}'$, then $f_{\mathfrak{b},x} \equiv f_{\mathfrak{b}',x} \textup{ mod }(t^{k+1})$.
\end{compactenum}
\end{defi}

\begin{prop}
\label{prop:Sk}
If $(B,\mathscr{P},\varphi)$ is a polarized polyhedral affine manifold with $(B,\mathscr{P})$ simple, then there exists as sequence of wall structures $(\mathscr{S}_k)_{k\in\mathbb{N}}$ such that
\begin{compactenum}[(1)]
\item $\mathscr{S}_0$ is the initial wall structure defined by $(B,\mathscr{P},\varphi)$;
\item $\mathscr{S}_k$ is consistent to order $k$;
\item $\mathscr{S}_k$ and $\mathscr{S}_{k+1}$ are compatible to order $k$.
\end{compactenum}
\end{prop}

\begin{proof}
By \cite{GS11}, Remark 1.29, if $(B,\mathscr{P})$ is simple, then the central fiber of the corresponding toric degeneration is locally rigid for any choice of open gluing data. In this case, the existence of a sequence of wall structures as claimed is the main part ({\S}3, {\S}4) of \cite{GS11}. Roughly speaking, the proof goes by induction as follows. To obtain $\mathscr{S}_k$ from $\mathscr{S}_{k-1}$, for each joint $\mathfrak{j}$ of $\mathscr{S}_{k-1}$ we calculate the scattering diagram $\mathfrak{D}_{\mathfrak{j},k}$ consistent to order $k$ from $\mathfrak{D}_{\mathfrak{j},k-1}$ as in Proposition \ref{prop:scatt}. Then we add walls corresponding to these rays to the scattering diagram $\mathscr{S}_{k-1}$. This will probably produce some new joint or complicate the scattering diagrams at other joints. However, it is shown in \cite{GS11} that this procedure after finitely many steps gibes a wall structure $\mathscr{S}_k$ consistent to order $k$.
\end{proof}

\subsection{Proof of the main theorem}							%%

Let $Q$ be a $2$-dimensional Fano polytope and let $(B,\mathscr{P},\varphi)$ be the dual intersection complex of the corresponding toric degeneration. Let $\mathscr{S}_\infty$ be the consistent wall structure defined by $(B,\mathscr{P},\varphi)$, i.e., the limit of the $\mathscr{S}_k$ in Proposition \ref{prop:Sk}. Let $\sigma_0$ be as in Definition \ref{defi:sigma0}.

\begin{lem}
\label{lem:global}
The support $|\mathscr{S}_\infty|$ is disjoint from the interior of $\sigma_0$.
\end{lem}

\begin{proof}
$(B,\mathscr{P},\varphi)$ defines an initial wall structure $\mathscr{S}_0$ as in Definition \ref{defi:S0}. The joints of $\mathscr{S}_0$ are the vertices of $\mathscr{P}$. Let $v$ be such a vertex. By \cite{GS11} Proposition 3.9, the walls in $\mathscr{S}_\infty$ with base point $v$ lie in the cone $v + \mathbb{R}_{\leq 0}m_{v\delta} + \mathbb{R}_{\leq 0}m_{v\delta'} \subseteq B$. Here $\delta,\delta'$ are the affine singularities on edges adjacent to $v$ and $m_{\delta v}$ is the primitive integral tangent vector on $B$ pointing from $v$ to $\delta$.

By inductively using \cite{GS11}, Proposition 3.9, all walls in $\mathscr{S}_\infty$ lie in the union of these cones, i.e.,
\[ |\mathscr{S}_\infty| \subseteq \bigcup_{v\in\mathscr{P}^{[0]}} v + \mathbb{R}_{\leq 0}m_{v\delta} + \mathbb{R}_{\leq 0}m_{v\delta'}. \]
In particular, there are no walls in $\mathscr{S}_\infty$ supported on the interior of $\sigma_0$.
\end{proof}

The unbounded walls in $\mathscr{S}_\infty$ are all parallel in the direction $m_{\textup{out}}\in \Lambda_B$. Let $f_{\textup{out}}$ be the product of all functions $f_{\mathfrak{p}}$ attached to unbounded walls $\mathfrak{p}$ in $\mathscr{S}_\infty$. Then $f_{\text{out}}$ can be regarded as an element of $\mathbb{C}\llbracket x\rrbracket$ for $x:=z^{(-m_{\textup{out}},0)}\in\mathbb{C}[\Lambda_B\oplus\mathbb{Z}]$.

\begin{thm}[Theorem \ref{thm:main}]
\label{thm:scattering}
\[ \textup{log }f_{\textup{out}} = \sum_{d=1}^\infty (D\cdot\underline{\beta}) \cdot N_d \cdot x^{D\cdot\underline{\beta}}. \]
\end{thm}

In fact, we will prove a more general statement, giving an enumerative meaning for the function attached to any wall in $\mathscr{S}_\infty$. For this we need the following.

\begin{defi}[\cite{GrP2}, Definition 2.2]
A \textit{tropical disk} $h : \Gamma \rightarrow B$ is a tropical curve with the choice of univalent vertex $V_\infty$, adjacent to a unique edge $E_\infty$, such that $h$ is balanced for all vertices $V \neq V_\infty$.
\end{defi}

\begin{defi}
Let $\mathfrak{p}\in\mathscr{S}_\infty$ be a wall and choose $x\in\textup{Int}(\mathfrak{p})$. Define $\mathfrak{H}_{\mathfrak{p},w}$ to be the set of all tropical disks $h : \Gamma \rightarrow B$ with $h(V_\infty)=x$ and $u_{(V_\infty,E_\infty)}=-w\cdot m_{\mathfrak{p}}$.
\end{defi}

\begin{defi}
For $h\in \mathfrak{H}_{\mathfrak{p},w}$ define, with $N_V^{\textup{tor}}$ as in Definition \ref{defi:Ntor},
\[ N_h = \frac{1}{|\textup{Aut}(h)|}\prod_{E\in E(\Gamma)}w_E\cdot\prod_{V\neq V_{\textup{out}}} N_V^{\textup{tor}} \cdot \prod_{E\in L_\Delta(\Gamma)}\frac{(-1)^{w_E-1}}{w_E}. \]
\end{defi}

\begin{rem}
Note that the sets $\mathfrak{H}_{\mathfrak{p},w}$ are in bijection for different choices of $x$. Moreover, for an unbounded wall $\mathfrak{p}$ the set $\mathfrak{H}_{\mathfrak{p},w}$ is empty for $w_{\text{out}}\nmid w$, and for each $d$ there is an injective map $\iota:\mathfrak{H}_{\mathfrak{p},dw_{\text{out}}} \hookrightarrow \mathfrak{H}_d$ by removing $V_\infty$ and extending $E_\infty$ to infinity, giving $E_{\text{out}}$. Hence, $N_{\iota(h)}=dw_{\text{out}}N_h$.
\end{rem}

\begin{prop}
\label{prop:scattering}
For a wall $\mathfrak{p}$ of $\mathscr{S}_\infty$ we have
\[ \textup{log }f_{\mathfrak{p}} = \sum_{w=1}^\infty\sum_{h\in\mathfrak{H}_{\mathfrak{p},w}} N_h z^{(wm_{\mathfrak{p}},0)}. \]
\end{prop}

\begin{proof}
We want to prove the claimed equality by induction, so we need a well-ordered set. For a wall $\mathfrak{p}$ of $\mathscr{S}_\infty$, define a set 
\[ \textup{Parents}(\mathfrak{p}) = \{\mathfrak{p}'\in\mathscr{S}_\infty \ | \ \mathfrak{p}\cap\mathfrak{p}'=\textup{Base}(\mathfrak{p}) \neq \textup{Base}(\mathfrak{p}')\}. \]
Here $\textup{Base}(\mathfrak{p})$ is the base point of $\mathfrak{p}$ (see Definition \ref{defi:structure}, (2)). Note that $\textup{Base}(\mathfrak{p}')$ is only defined for walls, so we define the condition $\textup{Base}(\mathfrak{p})\neq\textup{Base}(\mathfrak{p}')$ to be always satisfied when $\mathfrak{p}'$ is a slab. Then define inductively
\[ \textup{Ancestors}(\mathfrak{p}) = \{\mathfrak{p}\} \cup \bigcup_{\mathfrak{p}'\in\textup{Parents}(\mathfrak{p})} \textup{Ancestors}(\mathfrak{p}'). \]
For each $k$, the set of walls in $\mathscr{S}_k$ is finite and totally ordered by $\mathfrak{p}_1 \leq \mathfrak{p}_2$ if and only if $\mathfrak{p}_1 \in \textup{Ancestors}(\mathfrak{p}_2)$. Hence, it is well-ordered and we can use induction. The set of smallest elements with respect to this ordering is $\{\mathfrak{p}\in\mathscr{S}_\infty \textup{ wall}\ | \ \textup{Base}(\mathfrak{p})\in\mathscr{P}_{\mathscr{S}_0}^{[0]}=\mathscr{P}^{[0]}\}$. For such $\mathfrak{p}$, the set $\textup{Ancestors}(\mathfrak{p})$ consists of $\mathfrak{p}$ and two slabs $\mathfrak{b}_1,\mathfrak{b}_2$. This defines a scattering diagram at the joint $\textup{Base}(\mathfrak{p})$ that is equivalent to the scattering diagram obtained from two lines in the directions $m_1,m_2$ of the slabs $\mathfrak{b}_1,\mathfrak{b}_2$ with attached functions $f_1=1+z^{(-m_1,0)}$ and $f_2=1+z^{(-m_2,0)}$, respectively. Note that
\[ \textup{log }f_i = \sum_{w=1}^\infty \frac{(-1)^{w-1}}{w}z^{(wm_i,0)}. \]
Then Lemma \ref{lem:2.8} gives
\[ \textup{log }f_{\mathfrak{p}} = \sum_{w=1}^\infty\sum_{\textbf{w}}w\frac{N_{\textbf{m}}(\textup{\textbf{w}})}{|\textup{Aut}(\textup{\textbf{w}})|} \left(\prod_{\substack{1\leq i\leq n\\ 1\leq j\leq l_i}} \frac{(-1)^{w_{ij}-1}}{w_{ij}}\right) z^{(-wm_{\mathfrak{p}},0)}, \]
where the sum is over all $n$-tuples of weight vectors $\textbf{w}=(\textbf{w}_1,\ldots,\textbf{w}_n)$ satisfying
\[ \sum_{i=1}^n |\textbf{w}_i|m_i = wm_{\mathfrak{p}}. \]
Here $N_{\textbf{m}}(\textbf{w})$ is the toric logarithmic Gromov-Witten invariant defined in \S\ref{S:toricinv}. Tropical disks $h$ in $\mathfrak{H}_{\mathfrak{p},w}$ have a $1$-valent vertex $V_\infty$ mapping to the interior of $\mathfrak{p}$ and another vertex $V$ with one compact edge of weight $w$ and several bounded legs in directions $m_1$ or $m_2$ with weights $\textbf{w}_1=(w_{11},\ldots,w_{1l_1})$ and $\textbf{w}_2=(w_{21},\ldots,w_{2l_2})$, respectively, such that $\sum_{i=1}^n |\textbf{w}_i|m_i = wm_{\mathfrak{p}}$. Hence, the set $\mathfrak{H}_{\mathfrak{p},w}$ is in bijection with the set of pairs $\textbf{w}$ as above, with $w_{ij}$ being the weights of the bounded legs of the corresponding tropical disk $h$. Moreover, $|\textup{Aut}(h)|=|\textup{Aut}(\textbf{w})|$ and $N_V^{\textup{tor}}=N_{\textbf{m}}(\textbf{w})$ by definition. Hence, the summation over $\textbf{w}$ can be replaced by a summation over $\mathfrak{H}_{\mathfrak{p},w}$, and we obtain
\[ \textup{log }f_{\mathfrak{p}} = \sum_{w=1}^\infty\sum_{h\in\mathfrak{H}_{\mathfrak{p},w}} w\frac{N_V^{\textup{tor}}}{|\textup{Aut}(h)|}\left(\prod_{E\in L_\Delta(\Gamma)}\frac{(-1)^{w_E-1}}{w_E}\right)z^{(-wm_{\mathfrak{p}},0)}. \]
This is precisely the claimed formula for such $\mathfrak{p}$, completing the base case.

For the induction step, let $\mathfrak{p}$ be a wall of $\mathscr{S}_\infty$ and assume the claimed formula holds for all walls $\mathfrak{p}'\in\textup{Ancestors}(\mathfrak{p})\setminus\{\mathfrak{p}\}$. By Lemma \ref{lem:2.8} and using the induction hypothesis,
\[ \textup{log }f_{\mathfrak{p}} = \sum_{w=1}^n\sum_{\textbf{w}}\left(w\frac{N_{\textbf{m}}(\textup{\textbf{w}})}{|\textup{Aut}(\textup{\textbf{w}})|} \prod_{\substack{\mathfrak{p}'\in\textup{Parents}(\mathfrak{p})\\ 1\leq j\leq l_{\mathfrak{p}'}}} \sum_{h'\in\mathfrak{H}_{\mathfrak{p}',w_{\mathfrak{p}'j}}} N_{h'} \right) z^{(-wm_{\mathfrak{p}},0)}, \]
where the second sum is over all tuples $\textbf{w}=(\textbf{w}_{\mathfrak{p}'})_{\mathfrak{p}'\in\textup{Parents}(\mathfrak{p})}$ of weight vectors $\textbf{w}_{\mathfrak{p}'}=(w_{\mathfrak{p}'1},\ldots,w_{\mathfrak{p}'l_{\mathfrak{p}'}})$ with
\[ \sum_{\mathfrak{p}'\in\textup{Parents}(\mathfrak{p})}|\textbf{w}_{\mathfrak{p}'}|m_{\mathfrak{p}'} = wm_{\mathfrak{p}}. \]
Factoring out, we can replace the product over $\textup{Parents}(\mathfrak{p})$ and $1\leq j \leq l_{\mathfrak{p}'}$ and the sum over $h'\in\mathfrak{H}_{\mathfrak{p}',w_{\mathfrak{p}'j}}$ by a sum over tuples $(h')_{\mathfrak{p}',j}:=(h'\in\mathfrak{H}_{\mathfrak{p}',w_{\mathfrak{p}'j}})_{\substack{\mathfrak{p}'\in\textup{Parents}(\mathfrak{p}) \\ 1\leq j \leq l_{\mathfrak{p}'}}}$:
\[ \textup{log }f_{\mathfrak{p}} = \sum_{w=1}^n\sum_{\textbf{w}}\left(w\frac{N_{\textbf{m}}(\textup{\textbf{w}})}{|\textup{Aut}(\textup{\textbf{w}})|} \sum_{(h')_{\mathfrak{p}',j}} \prod_{h'\in (h')_{\mathfrak{p}',j}} N_{h'} \right) z^{(-wm_{\mathfrak{p}},0)}. \]
Further, we can replace the summations over $\textbf{w}$ and $(h')_{\mathfrak{p}',j}$ by a summation over $\mathfrak{H}_{\mathfrak{p},w}$, since this is precisely the data that determines a tropical disk in $\mathfrak{H}_{\mathfrak{p},w}$. We get the claimed formula by using that, for $h\in\mathfrak{H}_{\mathfrak{p},w}$ determined by $\textbf{w}$ and $(h')_{\mathfrak{p}',j}$, we have
\[ |\text{Aut}(h)| = |\text{Aut}(\textbf{w})| \cdot \prod_{(h')_{\mathfrak{p}',j}} |\text{Aut}(h')|, \]
so
\[ N_h = w\frac{N_{\textbf{m}}(\textbf{w})}{|\text{Aut}(\textbf{w})|} \cdot \prod_{h'\in (h')_{\mathfrak{p}',j}} N_{h'}. \]
This completes the proof.
\end{proof}
Now Theorem \ref{thm:scattering} follows by summation over all outgoing walls $\mathfrak{p}$ in $\mathscr{S}_\infty$.

\section{Torsion points}									%%%
\label{S:torsion}

In this section we consider $(\mathbb{P}^2,E)$. Choose a flex point $O$ on the elliptic curve $E$ and consider the group law on $E$ with $O$ the identity. An \textit{$m$-torsion point} on $E$ is a point $P$ such that $m \cdot P = O$. As a topological group, the elliptic curve is a torus $S^1\times S^1=\mathbb{R}^2/\mathbb{Z}^2$. The $m$-torsion points form a group $\mathbb{Z}_m \times \mathbb{Z}_m$.

\begin{lem}
\label{lem:torsion}
If $C$ is a rational degree $d$ curve intersecting $E$ in a single point $P$, then $P$ is a $3d$-torsion point.
\end{lem}

\begin{proof}
Let $C$ be a rational degree $d$ curve intersecting $E$ in a single point $P$ and let $L$ be the line tangent to $O$. Then the cycle $C-dL$ has degree $0$, so it is linearly equivalent to zero, since $\textup{Pic }\mathbb{P}^2 \cong \mathbb{Z}$ by the degree map. Moreover, it intersects $E$ in the cycle $3d(P-O)$ which in turn is linearly equivalent to zero.
\end{proof}

Let $T_d\simeq\mathbb{Z}_{3d}\times\mathbb{Z}_{3d}$ be the set of $3d$-torsion points on $E$ and let $\beta_d$ be the class of degree $d$ stable log maps (Definition \ref{defi:beta}). Then we have a decomposition
\[ \mathscr{M}(\mathbb{P}^2,\beta_d) = \coprod_{P\in T_d} \mathscr{M}(\mathbb{P}^2,\beta_d)_P, \]
where $\mathscr{M}(\mathbb{P}^2,\beta_d)_P$ is the subspace of $\mathscr{M}(\mathbb{P}^2,\beta_d)$ of maps intersecting $E$ in $P$. Let $\llbracket\mathscr{M}(\mathbb{P}^2,\beta_d)_P\rrbracket$ be the restriction of $\llbracket\mathscr{M}(\mathbb{P}^2,\beta_d)\rrbracket$ to $\mathscr{M}(\mathbb{P}^2,\beta_d)_P$ and define
\[ N_{d,P} := \int_{\llbracket\mathscr{M}(\mathbb{P}^2,\beta_d)_P\rrbracket}1. \]

\begin{defi}
For $P\in\cup_{d\geq 1}T_d$ denote by $k(P)$ the smallest integer $k\geq 1$ such that $P$ is $3k$-torsion.
\end{defi}

\begin{lem}
$N_{d,P}$ only depends on $P$ through $k(P)$.
\end{lem}

\begin{proof}
This was shown in \cite{Bou4}, Lemma 1.2, using ideas from \cite{CvGKT2}. The freedom of choice of $O$ and the fact that the monodromy of the family of smooth cubics in $\mathbb{P}^2$ maps surjectively to $SL(2,\mathbb{Z})$ acting on $T_k \simeq \mathbb{Z}_{3k} \times \mathbb{Z}_{3k}$ implies that two points $P,P'$ with $k(P)=k(P')$ are related to each other via a monodromy transformation. Then the deformation invariance of logarithmic Gromov-Witten invariants shows that $N_{d,P}=N_{d,P'}$ for all $d$.
\end{proof}

\begin{defi}
Write $N_{d,k}$ for $N_{d,P}$ with $P$ such that $k(P)=k$.
\end{defi}

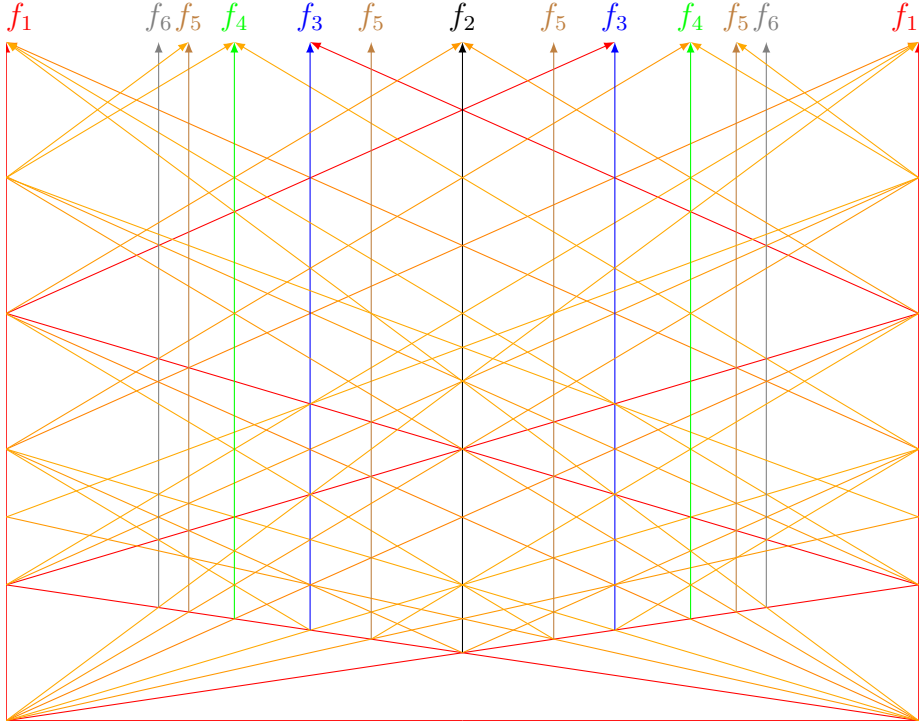
\begin{figure}[h!]
\centering
\begin{tikzpicture}[xscale=12,yscale=0.6,define rgb/.code={\definecolor{mycolor}{RGB}{#1}}, rgb color/.style={define rgb={#1},mycolor},rotate=90]
\clip (0,0) rectangle (16.2,1);
\draw[->,rgb color={255,0,0}] (0.000,0.500) -- (0.000,6.50);
\draw[->,rgb color={255,0,0}] (0.000,0.500) -- (0.000,-5.50);
\draw[->,rgb color={255,0,0}] (-1.50,-0.500) -- (15.0,5.00);
\draw[->,rgb color={255,0,0}] (-1.50,1.50) -- (15.0,-4.00);
\draw[->,rgb color={255,0,0}] (-1.50,1.50) -- (-16.5,6.50);
\draw[->,rgb color={255,0,0}] (-1.50,-0.500) -- (-16.5,-5.50);
\draw[->,rgb color={255,0,0}] (-6.00,-1.50) -- (15.0,2.00);
\draw[->,rgb color={255,0,0}] (-6.00,2.50) -- (15.0,-1.00);
\draw[->,rgb color={255,0,0}] (-6.00,2.50) -- (-30.0,6.50);
\draw[->,rgb color={255,0,0}] (-6.00,-1.50) -- (-30.0,-5.50);
\draw[->,rgb color={255,0,0}] (-13.5,-2.50) -- (15.0,0.667);
\draw[->,rgb color={255,0,0}] (-13.5,3.50) -- (15.0,0.333);
\draw[->,rgb color={255,0,0}] (-13.5,3.50) -- (-40.5,6.50);
\draw[->,rgb color={255,0,0}] (-13.5,-2.50) -- (-40.5,-5.50);
\draw[->,rgb color={255,0,0}] (-24.0,-3.50) -- (15.0,-0.250);
\draw[->,rgb color={255,0,0}] (-24.0,4.50) -- (15.0,1.25);
\draw[->,rgb color={255,0,0}] (-24.0,4.50) -- (-48.0,6.50);
\draw[->,rgb color={255,0,0}] (-24.0,-3.50) -- (-48.0,-5.50);
\draw[->,rgb color={255,0,0}] (-37.5,-4.50) -- (15.0,-1.00);
\draw[->,rgb color={255,0,0}] (-37.5,5.50) -- (15.0,2.00);
\draw[->,rgb color={255,0,0}] (-37.5,5.50) -- (-52.5,6.50);
\draw[->,rgb color={255,0,0}] (-37.5,-4.50) -- (-52.5,-5.50);
\draw[->,rgb color={255,0,0}] (-54.0,-5.50) -- (15.0,-1.67);
\draw[->,rgb color={255,0,0}] (-54.0,6.50) -- (15.0,2.67);
\draw[->,red] (0.000,0.000) -- (15.0,0.000) node[above left = 0 and -0.15]{$f_1$};
\draw[->,red] (-3.00,-1.00) -- (15.0,-1.00);
\draw[->,red] (0.000,1.00) -- (15.0,1.00) node[above right = 0 and -0.15]{$f_1$};
\draw[->,red] (-3.00,2.00) -- (15.0,2.00);
\draw[->,red] (-45.0,-5.00) -- (15.0,-5.00);
\draw[->,red] (-30.0,-4.00) -- (15.0,-4.00);
\draw[->,red] (-18.0,-3.00) -- (15.0,-3.00);
\draw[->,red] (-9.00,-2.00) -- (15.0,-2.00);
\draw[->,red] (-9.00,3.00) -- (15.0,3.00);
\draw[->,red] (-18.0,4.00) -- (15.0,4.00);
\draw[->,red] (-30.0,5.00) -- (15.0,5.00);
\draw[->,red] (-45.0,6.00) -- (15.0,6.00);
\draw[->,black] (1.50,0.500) -- (15.0,0.500) node[above]{$f_2$};
\draw[->,black] (-4.50,-1.50) -- (15.0,-1.50);
\draw[->,black] (0.000,-0.500) -- (15.0,-0.500);
\draw[->,black] (0.000,1.50) -- (15.0,1.50);
\draw[->,black] (-4.50,2.50) -- (15.0,2.50);
\draw[->,black] (-22.5,-3.50) -- (15.0,-3.50);
\draw[->,black] (-12.0,-2.50) -- (15.0,-2.50);
\draw[->,black] (-12.0,3.50) -- (15.0,3.50);
\draw[->,black] (-22.5,4.50) -- (15.0,4.50);
\draw[->,rgb color={255,169,0}] (0.000,0.000) -- (15.0,2.50);
\draw[->,rgb color={255,169,0}] (-18.0,-3.00) -- (15.0,-0.800);
\draw[->,rgb color={255,152,0}] (-9.00,-2.00) -- (15.0,0.000);
\draw[->,rgb color={255,132,0}] (-3.00,-1.00) -- (15.0,1.00);
\draw[->,rgb color={255,169,0}] (0.000,1.00) -- (15.0,6.00);
\draw[->,rgb color={255,132,0}] (-3.00,2.00) -- (-3.00,6.50);
\draw[->,rgb color={255,152,0}] (-9.00,3.00) -- (-19.5,6.50);
\draw[->,rgb color={255,169,0}] (-18.0,4.00) -- (-33.0,6.50);
\draw[->,rgb color={255,169,0}] (0.000,0.000) -- (15.0,-5.00);
\draw[->,rgb color={255,169,0}] (-18.0,-3.00) -- (-33.0,-5.50);
\draw[->,rgb color={255,152,0}] (-9.00,-2.00) -- (-19.5,-5.50);
\draw[->,rgb color={255,132,0}] (-3.00,-1.00) -- (-3.00,-5.50);
\draw[->,rgb color={255,169,0}] (0.000,1.00) -- (15.0,-1.50);
\draw[->,rgb color={255,132,0}] (-3.00,2.00) -- (15.0,0.000);
\draw[->,rgb color={255,152,0}] (-9.00,3.00) -- (15.0,1.00);
\draw[->,rgb color={255,169,0}] (-18.0,4.00) -- (15.0,1.80);
\draw[->,rgb color={255,152,0}] (3.00,0.000) -- (15.0,1.33);
\draw[->,rgb color={255,169,0}] (0.000,-1.00) -- (15.0,0.250);
\draw[->,rgb color={255,152,0}] (3.00,1.00) -- (15.0,3.00);
\draw[->,rgb color={255,169,0}] (0.000,2.00) -- (13.5,6.50);
\draw[->,rgb color={255,152,0}] (3.00,0.000) -- (15.0,-2.00);
\draw[->,rgb color={255,169,0}] (0.000,-1.00) -- (13.5,-5.50);
\draw[->,rgb color={255,152,0}] (3.00,1.00) -- (15.0,-0.333);
\draw[->,rgb color={255,169,0}] (0.000,2.00) -- (15.0,0.750);
\draw[->,rgb color={255,169,0}] (3.00,0.000) -- (15.0,-1.33);
\draw[->,rgb color={255,169,0}] (3.00,1.00) -- (15.0,0.000);
\draw[->,rgb color={255,152,0}] (1.50,0.500) -- (15.0,-1.00);
\draw[->,rgb color={255,169,0}] (0.000,-0.500) -- (15.0,-3.00);
\draw[->,rgb color={255,152,0}] (0.000,1.50) -- (15.0,0.250);
\draw[->,rgb color={255,152,0}] (1.50,0.500) -- (15.0,2.00);
\draw[->,rgb color={255,169,0}] (0.000,-0.500) -- (15.0,0.750);
\draw[->,rgb color={255,152,0}] (0.000,1.50) -- (15.0,4.00);
\draw[->,blue] (2.00,0.333) -- (15.0,0.333) node[above]{$f_3$};
\draw[->,blue] (2.00,0.667) -- (15.0,0.667) node[above]{$f_3$};
\draw[->,brown] (2.40,0.200) -- (15.0,0.200) node[above]{$f_5$};
\draw[->,brown] (0.000,-0.800) -- (15.0,-0.800);
\draw[->,brown] (1.80,1.20) -- (15.0,1.20);
\draw[->,brown] (2.40,0.800) -- (15.0,0.800) node[above]{$f_5$};
\draw[->,brown] (1.80,-0.200) -- (15.0,-0.200);
\draw[->,brown] (0.000,1.80) -- (15.0,1.80);
\draw[->,rgb color={255,169,0}] (3.00,0.000) -- (15.0,1.00);
\draw[->,rgb color={255,169,0}] (3.00,1.00) -- (15.0,2.33);
\draw[->,gray] (2.50,0.167) -- (15.0,0.167) node[above]{$f_6$};
\draw[->,gray] (2.50,0.833) -- (15.0,0.833) node[above]{$f_6$};
\draw[->,blue] (-13.0,-2.67) -- (15.0,-2.67);
\draw[->,blue] (-5.00,-1.67) -- (15.0,-1.67);
\draw[->,blue] (0.000,-0.667) -- (15.0,-0.667);
\draw[->,blue] (1.00,1.33) -- (15.0,1.33);
\draw[->,blue] (-3.00,2.33) -- (15.0,2.33);
\draw[->,blue] (-10.0,3.33) -- (15.0,3.33);
\draw[->,rgb color={255,169,0}] (2.00,0.333) -- (15.0,-0.750);
\draw[->,blue] (-10.0,-2.33) -- (15.0,-2.33);
\draw[->,blue] (-3.00,-1.33) -- (15.0,-1.33);
\draw[->,blue] (1.00,-0.333) -- (15.0,-0.333);
\draw[->,blue] (0.000,1.67) -- (15.0,1.67);
\draw[->,blue] (-5.00,2.67) -- (15.0,2.67);
\draw[->,blue] (-13.0,3.67) -- (15.0,3.67);
\draw[->,rgb color={255,169,0}] (2.00,0.667) -- (15.0,1.75);
\draw[->,rgb color={255,152,0}] (0.000,0.000) -- (15.0,3.33);
\draw[->,rgb color={255,169,0}] (-3.00,-1.00) -- (15.0,1.40);
\draw[->,rgb color={255,152,0}] (0.000,1.00) -- (8.25,6.50);
\draw[->,rgb color={255,169,0}] (-3.00,2.00) -- (-9.75,6.50);
\draw[->,rgb color={255,132,0}] (0.000,0.000) -- (15.0,1.67);
\draw[->,rgb color={255,169,0}] (-9.00,-2.00) -- (15.0,-0.400);
\draw[->,rgb color={255,152,0}] (-3.00,-1.00) -- (15.0,0.500);
\draw[->,rgb color={255,132,0}] (0.000,1.00) -- (15.0,3.50);
\draw[->,rgb color={255,152,0}] (-3.00,2.00) -- (10.5,6.50);
\draw[->,rgb color={255,169,0}] (-9.00,3.00) -- (-9.00,6.50);
\draw[->,rgb color={255,152,0}] (0.000,0.000) -- (8.25,-5.50);
\draw[->,rgb color={255,169,0}] (-3.00,-1.00) -- (-9.75,-5.50);
\draw[->,rgb color={255,152,0}] (0.000,1.00) -- (15.0,-2.33);
\draw[->,rgb color={255,169,0}] (-3.00,2.00) -- (15.0,-0.400);
\draw[->,rgb color={255,132,0}] (0.000,0.000) -- (15.0,-2.50);
\draw[->,rgb color={255,169,0}] (-9.00,-2.00) -- (-9.00,-5.50);
\draw[->,rgb color={255,152,0}] (-3.00,-1.00) -- (10.5,-5.50);
\draw[->,rgb color={255,132,0}] (0.000,1.00) -- (15.0,-0.667);
\draw[->,rgb color={255,152,0}] (-3.00,2.00) -- (15.0,0.500);
\draw[->,rgb color={255,169,0}] (-9.00,3.00) -- (15.0,1.40);
\draw[->,rgb color={255,169,0}] (1.80,0.400) -- (15.0,1.50);
\draw[->,rgb color={255,169,0}] (1.80,0.600) -- (15.0,-0.500);
\draw[->,green] (2.25,0.250) -- (15.0,0.250) node[above]{$f_4$};
\draw[->,green] (-5.25,-1.75) -- (15.0,-1.75);
\draw[->,green] (0.000,-0.750) -- (15.0,-0.750);
\draw[->,green] (1.50,1.25) -- (15.0,1.25);
\draw[->,green] (-2.25,2.25) -- (15.0,2.25);
\draw[->,green] (2.25,0.750) -- (15.0,0.750) node[above]{$f_4$};
\draw[->,green] (-2.25,-1.25) -- (15.0,-1.25);
\draw[->,green] (1.50,-0.250) -- (15.0,-0.250);
\draw[->,green] (0.000,1.75) -- (15.0,1.75);
\draw[->,green] (-5.25,2.75) -- (15.0,2.75);
\draw[->,brown] (1.80,0.400) -- (15.0,0.400) node[above]{$f_5$};
\draw[->,brown] (0.000,-0.600) -- (15.0,-0.600);
\draw[->,brown] (0.600,1.40) -- (15.0,1.40);
\draw[->,brown] (1.80,0.600) -- (15.0,0.600) node[above]{$f_5$};
\draw[->,brown] (0.600,-0.400) -- (15.0,-0.400);
\draw[->,brown] (0.000,1.60) -- (15.0,1.60);
\draw[->,rgb color={255,152,0}] (0.000,0.000) -- (15.0,1.25);
\draw[->,rgb color={255,169,0}] (-3.00,-1.00) -- (15.0,0.200);
\draw[->,rgb color={255,152,0}] (0.000,1.00) -- (15.0,2.67);
\draw[->,rgb color={255,169,0}] (-3.00,2.00) -- (15.0,5.00);
\draw[->,rgb color={255,152,0}] (0.000,0.000) -- (15.0,-1.67);
\draw[->,rgb color={255,169,0}] (-3.00,-1.00) -- (15.0,-4.00);
\draw[->,rgb color={255,152,0}] (0.000,1.00) -- (15.0,-0.250);
\draw[->,rgb color={255,169,0}] (-3.00,2.00) -- (15.0,0.800);
\draw[->,rgb color={255,169,0}] (0.000,0.000) -- (15.0,1.00);
\draw[->,rgb color={255,169,0}] (0.000,1.00) -- (15.0,2.25);
\draw[->,rgb color={255,169,0}] (0.000,0.000) -- (15.0,-1.25);
\draw[->,rgb color={255,169,0}] (0.000,1.00) -- (15.0,0.000);
\end{tikzpicture}
\caption{The wall structure $\mathscr{S}_6$ on one unbounded maximal cell of $\mathscr{P}$, showing the relevant attached functions $f_l$. }
\label{fig:sage2}
\end{figure}

Under the toric degeneration $\mathfrak{X} \rightarrow \mathbb{A}^1$ different $3d$-torsion points may map to the same point on the central fiber $X_0$, and even to a $0$-dimensional stratum. However, the limits of the $3$-torsion points all lie on the $1$-dimensional strata. The intersection points with $3d$-torsion correspond to the unbounded walls $\mathfrak{p}$ in $\mathscr{S}_\infty$ with non-zero $t^{3d}$-coefficient of $\textup{log }f_{\mathfrak{p}}$. Their number is exactly $3d$ and they are distributed as the $3l$-torsion points of a circle (see Figure \ref{fig:main} and Figure \ref{fig:sage2}).

This can be explained as follows. The $2$-dimensional torus is a $S^1$-fibration over $S^1$. In the SYZ limit, the $S^1$-fiber shrinks to a point and we only see the $S^1$-base, which is the tropicalization. In this limit, the $\mathbb{Z}_{3k}$-fibers of $T_k\simeq\mathbb{Z}_{3k}\times\mathbb{Z}_{3k}$ are identified and we only see the $\mathbb{Z}_{3k}$ in the base. The toric degeneration of divisors $\mathfrak{D}\rightarrow\mathbb{A}^1$ is an elliptic fibration. It contains a singular fiber that is a cycle of three $\mathbb{P}^1$, i.e., an $I_3$ fiber in Kodaira's classification of singular elliptic fibers. Then the monodromy acting on the first cohomology class of the general fiber is given by $M_1=\left(\begin{smallmatrix}1&3\\ 0&1\end{smallmatrix}\right)$ up to conjugation. Now the action of $M_1$ on $\mathbb{Z}_3\times\mathbb{Z}_3$ is trivial. This means that each $3$-torsion point really defines a section of the family, which will have a limit on the special fiber. Such limit is necessarily on the smooth part of the special fiber, i.e., on a $1$-dimensional toric stratum (see \cite{SS}, Theorem 6.3). The action of $M_1$ on $\mathbb{Z}_6\times\mathbb{Z}_6$ has some fixed points, corresponding to $6$-torsion points on $E$ that define sections with limit on the smooth part of the special fiber, i.e., on $1$-dimensional strata. The other $6$-torsion points are permuted by the action of $M_1$, so they only define multisections with limit on the singular part of the special fiber, i.e., on $0$-dimensional strata.

Now consider the refined degeneration $\tilde{\mathfrak{X}}_d\rightarrow\mathbb{A}^1$. By construction of the refinement, the limits of all $3d$-torsion points lie on $1$-dimensional strata of the central fiber $X:=Y$. Indeed, the central fiber of the elliptic fibration $\tilde{\mathfrak{D}}\rightarrow\mathbb{A}^1$ is a cycle of $3d$ lines, i.e., an $I_{3d}$-fiber. Then the monodromy acting on $\mathbb{Z}_{3k}\times\mathbb{Z}_{3k}$ is given by $\left(\begin{smallmatrix}1&3d\\ 0&1\end{smallmatrix}\right)$ which is the identity for all $k\mid d$.

\begin{defi}
Let $s_{k,l}$ be the number of points P on $E$ with $k(P)=k$, fixed by the action of $M_l:=\left(\begin{smallmatrix}1&3l\\ 0&1\end{smallmatrix}\right)$, but not fixed by the action of $M_{l'}$ for all $l'<l$. Note that $s_k:=\sum_{l\mid k}s_{k,l}$ is the number of points $P$ on $E$ with $k(P)=k$ and that $\sum_{k\mid d}s_k=(3d)^2$ is the number of $3d$-torsion points.
\end{defi}

\begin{table}[h!]
\begin{tabular}{|l|l|l|l|l|l|} 															\cline{1-1}
$s_{1,1}=9$ 	 															\\ \cline{1-2}
$s_{2,1}=9$		& $s_{2,2}=18$													\\ \cline{1-3}
$s_{3,1}=18$	& 			& $s_{3,3}=54$										\\ \cline{1-4}
$s_{4,1}=18$	& $s_{4,2}=18$	&			& $s_{4,4}=72$							\\ \cline{1-5}
$s_{5,1}=36$	&			&			&			& $s_{5,5}=180$				\\ \cline{1-6}
$s_{6,1}=18$	& $s_{6,2}=36$	& $s_{6,3}=54$	&			&			& $s_{6,6}=108$	\\ \hline
\end{tabular}
\vspace{7mm}
\label{tab:s}
\caption{The number $s_{k,l}$ of points $P$ on $E$ with $k(P)=k$, fixed by $M_l$, but not fixed by $M_{l'}$ for $l'<l$.}
\end{table}

\begin{defi}
For a wall $\mathfrak{p}\in\mathscr{S}_\infty$ let $l(\mathfrak{p})$ be the smallest number such that $\textup{log }f_{\mathfrak{p}}$ has non-trivial $t^{3l(\mathfrak{p})}$-coefficient. Let $r_l$ be the number of walls with $l(\mathfrak{p})=l$. 
\end{defi}

\begin{lem}
\label{lem:rl}
The number $r_l$ can be defined recursively by
\[ r_1 = 3, \quad r_l = 3l - \sum_{l'\mid l}r_{l'}. \]
\end{lem}

\begin{proof}
For a wall $\mathfrak{p}\in\mathscr{S}_\infty$, the condition $l(\mathfrak{p})=l$ means that the corresponding toric stratum of $X:=Y$ contains the limits of the points on $E$ that are fixed by the action of $M_l$ but not fixed by $M_{l'}$ for $l'<l$. Since in the SYZ limit we only see the base of the fibration, $r_l$ equals the number of points on a circle $S^1$ that are $3l$-torsion but not $3l'$-torsion. This number can be defined as above.
\end{proof}

\begin{table}[h!]
\begin{tabular}{|c|c|c|c|c|c|} 								\hline
$r_1=3$	& $r_2=3$	& $r_3=6$	& $r_4=6$	& $r_5=12$	& $r_6=6$	\\ \hline
\end{tabular}
\vspace{7mm}
\label{tab:r}
\caption{The number $r_l$ of walls $\mathfrak{p}$ with $l(\mathfrak{p})=l$.}
\end{table}

Note that $s_{k,l}/r_l$ is the number of points $P$ with $k(P)=k$ and with limit on the stratum corresponding to a particular wall $\mathfrak{p}$ with $l(\mathfrak{p})=l$. A direct consequence of Proposition \ref{prop:scattering} is the following.

\begin{cor}[Theorem \ref{thm:torsion}]
\label{cor:torsion}
Let $\mathfrak{p}$ be an unbounded wall of order $l$ in $\mathscr{S}_\infty$. Then
\[  \textup{log }f_{\mathfrak{p}} = \sum_{d=1}^\infty 3d \left(\sum_{k: l\mid k\mid d} \frac{s_{k,l}}{r_l} N_{d,k}\right) x^{3d}. \]
\end{cor}

\begin{defi}
Similarly to \S\ref{S:BPS} we can define integers $n_{d,k}$ recursively by
\[ n_{d,d} = N_{d,d}, \quad n_{d,k} = N_{d,k} - \sum_{d':k\mid d'\mid d} M_{d'}[d/d'] \cdot n_{d',k}. \]
\end{defi}

Some of the numbers $n_{d,k}$ have been calculated by Takahashi \cite{Ta1}. Their relation to local BPS numbers is studied in \cite{CvGKT}\cite{CvGKT2}\cite{CvGKT3}.

\begin{rem}
\label{rem:combine}
Unfortunately we are not able to apply the methods of this section to the refined situation of \S\ref{S:degdiv} in order to calculate the contributions to $N_{[d_1,\ldots,d_m]}$ with prescribed torsion: Let $N_{[d_1,\ldots,d_m],k}$ be the logarithmic Gromov-Witten invariant of stable log maps contributing to $N_{[d_1,\ldots,d_m]}$ and meeting $E$ in a fixed point of order $3k$, and let $n_{[d_1,\ldots,d_m],k}$ be the corresponding log BPS number. Define $l([d_1,\ldots,d_m]):=l(\mathfrak{p})$ for $\mathfrak{p}$ the unbounded wall for $[d_1,\ldots,d_m]$. Then
\begin{eqnarray*}
\sum_{k=1}^d \frac{s_{k,l}}{r_l} n_{[d_1,\ldots,d_m],k} &=& n_{[d_1,\ldots,d_m]} \quad \text{for } l=l([d_1,\ldots,d_m])\\
\sum_{l([d_1,\ldots,d_m])=l} n_{[d_1,\ldots,d_m],k} &=& n_{d,k} \qquad\quad\ \ \text{for all } l|k
\end{eqnarray*}
This gives a system of linear equations for the indeterminates $n_{[d_1,\ldots,d_m],k}$. In general the number of equations will be smaller than the number of indeterminates, so there will be no unique solution. However, for $d\leq 3$ we indeed have enough equations to determine the numbers $n_{[d_1,\ldots,d_m],k}$ as we will show in \S\ref{S:calcP2c}.
\end{rem}

\section{Higher genus and $q$-refined invariants}					%%%
\label{S:genus}

For an effective curve class $\underline{\beta}$ of $X$ let $\beta^g$ be the class of $1$-marked stable log maps to $X$ of genus $g$, class $\underline{\beta}$ and maximal tangency with $D$ at a single unspecified point. The moduli space $\mathscr{M}(X,\beta^g)$ of basic stable log maps of class $\beta^g$ has virtual dimension $g$. We can cut this dimension down to zero by inserting a \textit{lambda class}. Let $\pi : \mathcal{C} \rightarrow \mathscr{M}(X,\beta^g)$ be the universal curve, with relative dualizing sheaf $\omega_\pi$. Then $\mathbb{E}=\pi_\star\omega_\pi$ is a rank $g$ vector bundle over $\mathscr{M}(X,\beta^g)$, called the Hodge bundle. The lambda classes are the Chern classes of the Hodge bundle, $\lambda_j=c_j(\mathbb{E})$. We can define higher genus $1$-marked log Gromov-Witten invariants
\[ N_\beta^g = \int_{\llbracket\mathscr{M}(X,\beta^g)\rrbracket} (-1)^g\lambda_g \in \mathbb{Q}. \]

\begin{defi}
Let $h : \Gamma \rightarrow B$ be a tropical curve. For a trivalent vertex $V$ with multiplicity $m_V$ (Definition \ref{defi:mult}) define, with $q=e^{i\hbar}$,
\[ m_V(q) = \frac{1}{i\hbar}\left(q^{m_V/2}-q^{-m_V/2}\right). \]
For a vertex with higher valvency define $m_V(q) = \prod_{V'\in V(h'_V)} m_{V'}(q)$ with $h'_V$ as in Definition \ref{defi:mult}. For a bounded leg $E$ with weight $w_E$ define
\[ m_E(q) = \frac{(-1)^{w_E+1}}{w_E}\cdot \frac{i\hbar}{q^{w_E/2}-q^{-w_E/2}}. \]
Then define the \textit{$q$-refined multiplicity} of $h$ to be
\[ m_h(q) = \frac{1}{|\text{Aut}(h)|} \cdot \prod_{V\in V(\Gamma)} m_V(q) \cdot \prod_{E\in L_\Delta(\Gamma)} m_E(q). \]
\end{defi}

\begin{thm}
We have, with $q=e^{i\hbar}$,
\[ \sum_{g\geq 0} N_\beta^g\hbar^{2g} = \sum_{h\in\mathfrak{H}_\beta} m_h(q) \]
\end{thm}

\begin{proof}
Consider a stable log map in $\mathscr{M}(X,\beta^g)$ and let $h : \Gamma \rightarrow B$ be its tropicalization. In Definition \ref{defi:tropical} we defined the genus of $h$ to be $g_h = g_\Gamma + \sum_V g_V$. Using gluing and vanishing properties of lambda classes, Bousseau showed in \cite{Bou1} that $\Gamma$ is still a tree ($g_\Gamma=0$), i.e., all contributions to $g_h$ come from vertices. Hence, $h$ maps to an element of $\mathfrak{H}_\beta$ by forgetting genera at vertices $g_V$. So we can sum over $\mathfrak{H}_\beta$ but have to consider $q$-refined contributions of vertices. By \cite{Bou1}, Proposition 29, the contribution of a vertex $V$ with classical multiplicity $m_V$ is $m_V(q)$. By \cite{Bou2}, Lemma 5.9, the contribution of a bounded leg $L$ with weight $w_L$ is $m_L(q)$.
\end{proof}

\begin{expl}
For the tropical curve in Figure \ref{fig:balancing2} we have
\begin{eqnarray*}
m_h(q) &=& \frac{1}{i\hbar}\left(q^{9/2}-q^{-9/2}\right) \cdot \frac{1}{i\hbar}\left(q^3-q^{-3}\right) \cdot \frac{-1}{2} \cdot \frac{i\hbar}{q^1-q^{-1}} \cdot \left(\frac{i\hbar}{q^{1/2}-q^{-1/2}}\right)^2 \\
&=& -\frac{27}{2} + \frac{999}{16}\hbar^2 - \frac{137781}{1280}\hbar^4 + \ldots
\end{eqnarray*}
So this contributes $\frac{999}{16}$ to $N_3^1(\mathbb{P}^2)$.
\end{expl}

To obtain a higher genus version of Theorem \ref{thm:main}, we have to $q$-refine the slab functions in the initial wall structure $\mathscr{S}_0$. For $q$-refined wall structures it turns out to be more conventient to work with the logarithm of such functions. Define the $q$-refined initial wall structure $\mathscr{S}_0(q)$ to have the same slabs as $\mathscr{S}_0$ but with slab functions $f_{\mathfrak{b},v}=1+z^{(m_{v\delta},0)}$ replaced by $f_{\mathfrak{b},v}(q)$, where
\[ \text{log }f_{\mathfrak{b},v}(q) = \sum_{k\geq 1} \frac{(-1)^{k+1}i\hbar}{q^{k/2}-q^{-k/2}}z^{(km_{v\delta},0)}. \]
Note that the coefficient of $z^{(km_{v\delta},0)}$ is the $q$-multiplicity of a bounded leg of weight $k$. Let $\mathscr{S}_\infty(q)$ be the consistent $q$-refined wall structure obtained from $\mathscr{S}_0(q)$, let $f_{\text{out}}(q)$ be the product all of functions attached to unbounded walls, and write $x=z^{(-m_{\textup{out}},0)}$.

\begin{thm}
\[ \textup{log }f_{\textup{out}}(q) = \sum_{\underline{\beta}\in H_2^+(X,\mathbb{Z})} \sum_{g\geq 0} (D\cdot\underline{\beta}) N_\beta^g\hbar^{2g} x^{D\cdot\underline{\beta}}. \]
\end{thm}

\begin{proof}
The higher genus version of Lemma \ref{lem:2.8} is \cite{Bou2}, Proposition 6.2. Applying it inductively as in the proof of Proposition \ref{prop:scattering} we obtain the above formula.
\end{proof}

For $(\mathbb{P}^2,E)$, let $N_{d,k}^g$ be the log Gromov-Witten invariant of $1$-marked genus $g$, degree $d$ stable log maps meeting $E$ in a specified point of order $k$ with respect to the group law from \S\ref{S:torsion}. Everything we said respects the groups structure of $E$, so we get a $q$-refined version of Theorem \ref{thm:torsion}:

\begin{thm}
Let $\mathfrak{p}$ be an unbounded wall of order $l$ in $\mathscr{S}_\infty(q)$ for $(\mathbb{P}^2,E)$. Then
\[  \textup{log }f_{\mathfrak{p}}(q) = \sum_{d=1}^\infty \sum_{g\geq 0} 3d \left(\sum_{k: l\mid k\mid d} \frac{s_{k,l}}{r_l} N_{d,k}^g\right) \hbar^{2g}x^{3d}. \]
\end{thm}

\section{Explicit calculations}								%%%
\label{S:calc}

In this section we will calculate some logarithmic Gromov-Witten invariants and log BPS numbers explicitly. To this end, I wrote a sage code for calculating scattering diagrams and wall structures. It can be found on my webpage\footnote{\url{http://timgraefnitz.com/}}.

\subsection{$(\mathbb{P}^2,E)$}
\label{S:calcP2}

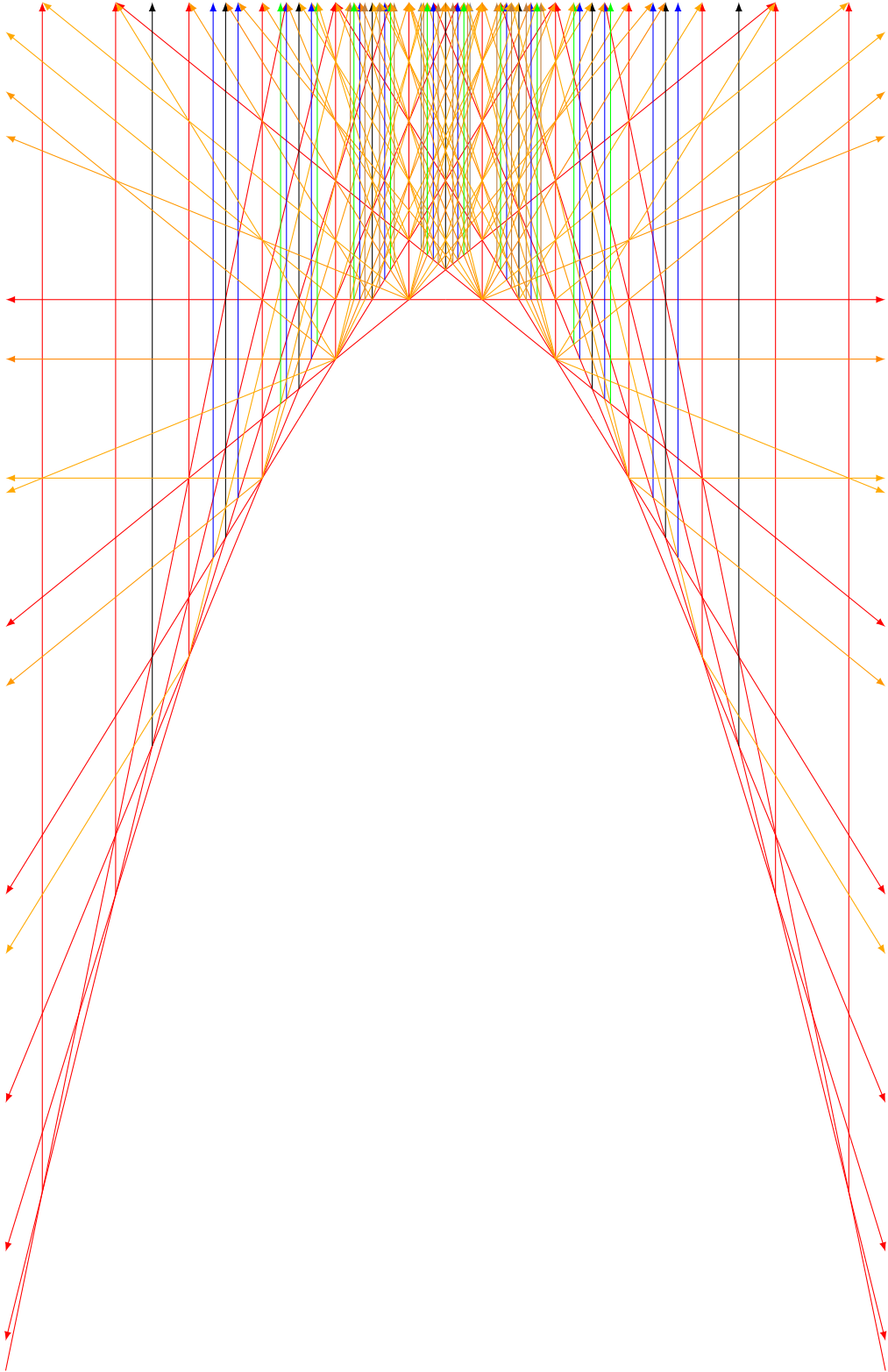
\begin{figure}[h!]
\centering
\begin{tikzpicture}[xscale=1.1,yscale=0.3,define rgb/.code={\definecolor{mycolor}{RGB}{#1}}, rgb color/.style={define rgb={#1},mycolor},rotate=90]
\draw[->,rgb color={255,0,0}] (0.000,0.500) -- (0.000,6.50);
\draw[->,rgb color={255,0,0}] (0.000,0.500) -- (0.000,-5.50);
\draw[->,rgb color={255,0,0}] (-1.50,-0.500) -- (15.0,5.00);
\draw[->,rgb color={255,0,0}] (-1.50,1.50) -- (15.0,-4.00);
\draw[->,rgb color={255,0,0}] (-1.50,1.50) -- (-16.5,6.50);
\draw[->,rgb color={255,0,0}] (-1.50,-0.500) -- (-16.5,-5.50);
\draw[->,rgb color={255,0,0}] (-6.00,-1.50) -- (15.0,2.00);
\draw[->,rgb color={255,0,0}] (-6.00,2.50) -- (15.0,-1.00);
\draw[->,rgb color={255,0,0}] (-6.00,2.50) -- (-30.0,6.50);
\draw[->,rgb color={255,0,0}] (-6.00,-1.50) -- (-30.0,-5.50);
\draw[->,rgb color={255,0,0}] (-13.5,-2.50) -- (15.0,0.667);
\draw[->,rgb color={255,0,0}] (-13.5,3.50) -- (15.0,0.333);
\draw[->,rgb color={255,0,0}] (-13.5,3.50) -- (-40.5,6.50);
\draw[->,rgb color={255,0,0}] (-13.5,-2.50) -- (-40.5,-5.50);
\draw[->,rgb color={255,0,0}] (-24.0,-3.50) -- (15.0,-0.250);
\draw[->,rgb color={255,0,0}] (-24.0,4.50) -- (15.0,1.25);
\draw[->,rgb color={255,0,0}] (-24.0,4.50) -- (-48.0,6.50);
\draw[->,rgb color={255,0,0}] (-24.0,-3.50) -- (-48.0,-5.50);
\draw[->,rgb color={255,0,0}] (-37.5,-4.50) -- (15.0,-1.00);
\draw[->,rgb color={255,0,0}] (-37.5,5.50) -- (15.0,2.00);
\draw[->,rgb color={255,0,0}] (-37.5,5.50) -- (-52.5,6.50);
\draw[->,rgb color={255,0,0}] (-37.5,-4.50) -- (-52.5,-5.50);
\draw[->,rgb color={255,0,0}] (-54.0,-5.50) -- (15.0,-1.67);
\draw[->,rgb color={255,0,0}] (-54.0,6.50) -- (15.0,2.67);
\draw[->,red] (0.000,0.000) -- (15.0,0.000);
\draw[->,red] (-3.00,-1.00) -- (15.0,-1.00);
\draw[->,red] (0.000,1.00) -- (15.0,1.00);
\draw[->,red] (-3.00,2.00) -- (15.0,2.00);
\draw[->,red] (-45.0,-5.00) -- (15.0,-5.00);
\draw[->,red] (-30.0,-4.00) -- (15.0,-4.00);
\draw[->,red] (-18.0,-3.00) -- (15.0,-3.00);
\draw[->,red] (-9.00,-2.00) -- (15.0,-2.00);
\draw[->,red] (-9.00,3.00) -- (15.0,3.00);
\draw[->,red] (-18.0,4.00) -- (15.0,4.00);
\draw[->,red] (-30.0,5.00) -- (15.0,5.00);
\draw[->,red] (-45.0,6.00) -- (15.0,6.00);
\draw[->,black] (1.50,0.500) -- (15.0,0.500);
\draw[->,black] (-4.50,-1.50) -- (15.0,-1.50);
\draw[->,black] (0.000,-0.500) -- (15.0,-0.500);
\draw[->,black] (0.000,1.50) -- (15.0,1.50);
\draw[->,black] (-4.50,2.50) -- (15.0,2.50);
\draw[->,black] (-22.5,-3.50) -- (15.0,-3.50);
\draw[->,black] (-12.0,-2.50) -- (15.0,-2.50);
\draw[->,black] (-12.0,3.50) -- (15.0,3.50);
\draw[->,black] (-22.5,4.50) -- (15.0,4.50);
\draw[->,rgb color={255,169,0}] (0.000,0.000) -- (15.0,2.50);
\draw[->,rgb color={255,169,0}] (-18.0,-3.00) -- (15.0,-0.800);
\draw[->,rgb color={255,152,0}] (-9.00,-2.00) -- (15.0,0.000);
\draw[->,rgb color={255,132,0}] (-3.00,-1.00) -- (15.0,1.00);
\draw[->,rgb color={255,169,0}] (0.000,1.00) -- (15.0,6.00);
\draw[->,rgb color={255,132,0}] (-3.00,2.00) -- (-3.00,6.50);
\draw[->,rgb color={255,152,0}] (-9.00,3.00) -- (-19.5,6.50);
\draw[->,rgb color={255,169,0}] (-18.0,4.00) -- (-33.0,6.50);
\draw[->,rgb color={255,169,0}] (0.000,0.000) -- (15.0,-5.00);
\draw[->,rgb color={255,169,0}] (-18.0,-3.00) -- (-33.0,-5.50);
\draw[->,rgb color={255,152,0}] (-9.00,-2.00) -- (-19.5,-5.50);
\draw[->,rgb color={255,132,0}] (-3.00,-1.00) -- (-3.00,-5.50);
\draw[->,rgb color={255,169,0}] (0.000,1.00) -- (15.0,-1.50);
\draw[->,rgb color={255,132,0}] (-3.00,2.00) -- (15.0,0.000);
\draw[->,rgb color={255,152,0}] (-9.00,3.00) -- (15.0,1.00);
\draw[->,rgb color={255,169,0}] (-18.0,4.00) -- (15.0,1.80);
\draw[->,rgb color={255,152,0}] (3.00,0.000) -- (15.0,1.33);
\draw[->,rgb color={255,169,0}] (0.000,-1.00) -- (15.0,0.250);
\draw[->,rgb color={255,152,0}] (3.00,1.00) -- (15.0,3.00);
\draw[->,rgb color={255,169,0}] (0.000,2.00) -- (13.5,6.50);
\draw[->,rgb color={255,152,0}] (3.00,0.000) -- (15.0,-2.00);
\draw[->,rgb color={255,169,0}] (0.000,-1.00) -- (13.5,-5.50);
\draw[->,rgb color={255,152,0}] (3.00,1.00) -- (15.0,-0.333);
\draw[->,rgb color={255,169,0}] (0.000,2.00) -- (15.0,0.750);
\draw[->,rgb color={255,169,0}] (3.00,0.000) -- (15.0,-1.33);
\draw[->,rgb color={255,169,0}] (3.00,1.00) -- (15.0,0.000);
\draw[->,rgb color={255,152,0}] (1.50,0.500) -- (15.0,-1.00);
\draw[->,rgb color={255,169,0}] (0.000,-0.500) -- (15.0,-3.00);
\draw[->,rgb color={255,152,0}] (0.000,1.50) -- (15.0,0.250);
\draw[->,rgb color={255,152,0}] (1.50,0.500) -- (15.0,2.00);
\draw[->,rgb color={255,169,0}] (0.000,-0.500) -- (15.0,0.750);
\draw[->,rgb color={255,152,0}] (0.000,1.50) -- (15.0,4.00);
\draw[->,blue] (2.00,0.333) -- (15.0,0.333);
\draw[->,blue] (2.00,0.667) -- (15.0,0.667);
\draw[->,brown] (2.40,0.200) -- (15.0,0.200);
\draw[->,brown] (0.000,-0.800) -- (15.0,-0.800);
\draw[->,brown] (1.80,1.20) -- (15.0,1.20);
\draw[->,brown] (2.40,0.800) -- (15.0,0.800);
\draw[->,brown] (1.80,-0.200) -- (15.0,-0.200);
\draw[->,brown] (0.000,1.80) -- (15.0,1.80);
\draw[->,rgb color={255,169,0}] (3.00,0.000) -- (15.0,1.00);
\draw[->,rgb color={255,169,0}] (3.00,1.00) -- (15.0,2.33);
\draw[->,gray] (2.50,0.167) -- (15.0,0.167);
\draw[->,gray] (2.50,0.833) -- (15.0,0.833);
\draw[->,blue] (-13.0,-2.67) -- (15.0,-2.67);
\draw[->,blue] (-5.00,-1.67) -- (15.0,-1.67);
\draw[->,blue] (0.000,-0.667) -- (15.0,-0.667);
\draw[->,blue] (1.00,1.33) -- (15.0,1.33);
\draw[->,blue] (-3.00,2.33) -- (15.0,2.33);
\draw[->,blue] (-10.0,3.33) -- (15.0,3.33);
\draw[->,rgb color={255,169,0}] (2.00,0.333) -- (15.0,-0.750);
\draw[->,blue] (-10.0,-2.33) -- (15.0,-2.33);
\draw[->,blue] (-3.00,-1.33) -- (15.0,-1.33);
\draw[->,blue] (1.00,-0.333) -- (15.0,-0.333);
\draw[->,blue] (0.000,1.67) -- (15.0,1.67);
\draw[->,blue] (-5.00,2.67) -- (15.0,2.67);
\draw[->,blue] (-13.0,3.67) -- (15.0,3.67);
\draw[->,rgb color={255,169,0}] (2.00,0.667) -- (15.0,1.75);
\draw[->,rgb color={255,152,0}] (0.000,0.000) -- (15.0,3.33);
\draw[->,rgb color={255,169,0}] (-3.00,-1.00) -- (15.0,1.40);
\draw[->,rgb color={255,152,0}] (0.000,1.00) -- (8.25,6.50);
\draw[->,rgb color={255,169,0}] (-3.00,2.00) -- (-9.75,6.50);
\draw[->,rgb color={255,132,0}] (0.000,0.000) -- (15.0,1.67);
\draw[->,rgb color={255,169,0}] (-9.00,-2.00) -- (15.0,-0.400);
\draw[->,rgb color={255,152,0}] (-3.00,-1.00) -- (15.0,0.500);
\draw[->,rgb color={255,132,0}] (0.000,1.00) -- (15.0,3.50);
\draw[->,rgb color={255,152,0}] (-3.00,2.00) -- (10.5,6.50);
\draw[->,rgb color={255,169,0}] (-9.00,3.00) -- (-9.00,6.50);
\draw[->,rgb color={255,152,0}] (0.000,0.000) -- (8.25,-5.50);
\draw[->,rgb color={255,169,0}] (-3.00,-1.00) -- (-9.75,-5.50);
\draw[->,rgb color={255,152,0}] (0.000,1.00) -- (15.0,-2.33);
\draw[->,rgb color={255,169,0}] (-3.00,2.00) -- (15.0,-0.400);
\draw[->,rgb color={255,132,0}] (0.000,0.000) -- (15.0,-2.50);
\draw[->,rgb color={255,169,0}] (-9.00,-2.00) -- (-9.00,-5.50);
\draw[->,rgb color={255,152,0}] (-3.00,-1.00) -- (10.5,-5.50);
\draw[->,rgb color={255,132,0}] (0.000,1.00) -- (15.0,-0.667);
\draw[->,rgb color={255,152,0}] (-3.00,2.00) -- (15.0,0.500);
\draw[->,rgb color={255,169,0}] (-9.00,3.00) -- (15.0,1.40);
\draw[->,rgb color={255,169,0}] (1.80,0.400) -- (15.0,1.50);
\draw[->,rgb color={255,169,0}] (1.80,0.600) -- (15.0,-0.500);
\draw[->,green] (2.25,0.250) -- (15.0,0.250);
\draw[->,green] (-5.25,-1.75) -- (15.0,-1.75);
\draw[->,green] (0.000,-0.750) -- (15.0,-0.750);
\draw[->,green] (1.50,1.25) -- (15.0,1.25);
\draw[->,green] (-2.25,2.25) -- (15.0,2.25);
\draw[->,green] (2.25,0.750) -- (15.0,0.750);
\draw[->,green] (-2.25,-1.25) -- (15.0,-1.25);
\draw[->,green] (1.50,-0.250) -- (15.0,-0.250);
\draw[->,green] (0.000,1.75) -- (15.0,1.75);
\draw[->,green] (-5.25,2.75) -- (15.0,2.75);
\draw[->,brown] (1.80,0.400) -- (15.0,0.400);
\draw[->,brown] (0.000,-0.600) -- (15.0,-0.600);
\draw[->,brown] (0.600,1.40) -- (15.0,1.40);
\draw[->,brown] (1.80,0.600) -- (15.0,0.600);
\draw[->,brown] (0.600,-0.400) -- (15.0,-0.400);
\draw[->,brown] (0.000,1.60) -- (15.0,1.60);
\draw[->,rgb color={255,152,0}] (0.000,0.000) -- (15.0,1.25);
\draw[->,rgb color={255,169,0}] (-3.00,-1.00) -- (15.0,0.200);
\draw[->,rgb color={255,152,0}] (0.000,1.00) -- (15.0,2.67);
\draw[->,rgb color={255,169,0}] (-3.00,2.00) -- (15.0,5.00);
\draw[->,rgb color={255,152,0}] (0.000,0.000) -- (15.0,-1.67);
\draw[->,rgb color={255,169,0}] (-3.00,-1.00) -- (15.0,-4.00);
\draw[->,rgb color={255,152,0}] (0.000,1.00) -- (15.0,-0.250);
\draw[->,rgb color={255,169,0}] (-3.00,2.00) -- (15.0,0.800);
\draw[->,rgb color={255,169,0}] (0.000,0.000) -- (15.0,1.00);
\draw[->,rgb color={255,169,0}] (0.000,1.00) -- (15.0,2.25);
\draw[->,rgb color={255,169,0}] (0.000,0.000) -- (15.0,-1.25);
\draw[->,rgb color={255,169,0}] (0.000,1.00) -- (15.0,0.000);
\end{tikzpicture}
\caption{The output of the sage code, giving the relevant part of $\bar{\mathscr{S}}_6$ to compute $N_{d,k}$ for $d\leq 6$. The colors correspond to different orders. The initial wall structure $\bar{\mathscr{S}}_0$ is red.}
\label{fig:sage}
\end{figure}

We want to calculate the numbers $N_{d,k}$ and $n_{d,k}$ for $d\leq 6$ as well as the numbers $n_{[d_1,\ldots,d_k]}$ for $d\leq 4$. Loading the code into a sage shell and typing
\begin{eqnarray*}
&&\verb+D = Diagram(case="P2",order=6)+ \\
&&\verb+D2 = D.scattering(order=6,case="P2")+\\
&&\verb+D2.tex(initial_diagram=D,print_directions=[(1,0)])+
\end{eqnarray*}
one can produce a TikZ code that under some small changes gives Figure \ref{fig:sage}. It shows the part of the wall structure $\bar{\mathscr{S}}_6$ on the discrete covering space $\bar{B}$ (see \S\ref{S:affinecharts}) that is relevant for computing the functions on the central maximal cell. The full $\bar{\mathscr{S}}_6$ would be symmetric, carrying much more walls on the outer area. We have
\begin{eqnarray*}
\text{log }f_{\text{out}} &=& 27x^3 + \frac{405}{2}x^6 + 2196x^9 + \frac{110997}{4}x^{12} \\
&& \text{ } \qquad\ \ + \frac{2051892}{5}x^{15} + 5527710x^{18} + \mathcal{O}(x^{20})
\end{eqnarray*}
This yields the following logarithmic Gromov-Witten invariants:

\renewcommand{\arraystretch}{1.2}
\begin{table}[h!]
\begin{tabular}{|c|c|c|c|c|c|} 								\hline
$N_1=9$	& $N_2=\frac{135}{4}$		& $N_3=244$	& $N_4=\frac{36999}{16}$	& $N_5=\frac{635634}{25}$	& $N_6=307095$	\\ \hline
\end{tabular}
\vspace{7mm}
\label{tab:N}
\caption{The invariants $N_d$ of $(\mathbb{P}^2,E)$ for $d \leq 6$.}
\end{table}

Subtracting multiple cover contributions we get the following log BPS numbers:

\begin{table}[h!]
\begin{tabular}{|c|c|c|c|c|c|} 								\hline
$n_1=9$	& $n_2=27$		& $n_3=234$	& $n_4=2232$	& $n_5=25380$	& $n_6=305829$	\\ \hline
\end{tabular}
\vspace{7mm}
\label{tab:n}
\caption{The log BPS numbers $n_d$ of $(\mathbb{P}^2,E)$ for $d \leq 6$.}
\end{table}

They are related to the local BPS numbers $n_d^{\text{loc}}$, shown in \cite{CKYZ}, Table 1, by Remark \ref{rem:local} and \cite{CKYZ}, (2.1). 

\subsubsection{Torsion points}										%
\label{S:calcP2a}

Write $f_l$ for the function $f_{\mathfrak{p}}$ attached to a wall $\mathfrak{p}$ with $l(\mathfrak{p})=l$.

The sage code gives the following:
\begin{align*}
\textup{log }f_1 &= 9x^{3}+\frac{63}{2}x^{6}+246x^{9}+\frac{9279}{4}x^{12}+\frac{175464}{5}x^{15}+307041x^{18} + \mathcal{O}(x^{21}) \\
\textup{log }f_2 &= 36x^{6} + 2322x^{12} + 307164x^{18} + \mathcal{O}(x^{21}) \\ 
\textup{log }f_3 &= 243x^9 + \frac{614061}{2}x^{18} + \mathcal{O}(x^{21}) \\
\textup{log }f_4 &= 2304x^{12} + \mathcal{O}(x^{21}) \\
\textup{log }f_5 &= 25425x^{15} + \mathcal{O}(x^{21}) \\
\textup{log }f_6 &= 307152x^{18} + \mathcal{O}(x^{21})
\end{align*}
From Corollary \ref{cor:torsion} we get:
\begin{compactenum}[(1)]
\item For $l=6$ we get
\[ 307041 = 18 \cdot \frac{108}{6} \cdot N_{6,6}, \]
hence $N_{6,6}=948$. There are no multiple cover contributions, as there are no curves of degree $<6$ meeting a point with $l=6$. This shows $n_{6,6}=948$.
\item For $l=5$ we get
\[ 25425 = 15 \cdot \frac{180}{12} \cdot N_{5,5}, \]
hence $N_{5,5}=113$. There are no multiple cover contributions, so $n_{5,5}=113$.
\item For $l=4$ we get
\[ 2304 = 12 \cdot \frac{72}{6} \cdot N_{4,4}, \]
hence $N_{4,4}=16$. There are no multiple cover contributions, so $n_{4,4}=16$.
\item For $l=3$ we get
\[ 243 = 9 \cdot \frac{54}{6} \cdot N_{3,3}, \]
hence $N_{3,3}=3$. There are no multiple cover contributions, so $n_{3,3}=3$. Moreover,
\[ \frac{614061}{2} = 18 \cdot \left(\frac{54}{6} \cdot N_{6,3} + \frac{54}{6} \cdot 948\right), \]
so $N_{6,3}=\frac{3789}{4}$. Subtracting $M_3[2]n_{3,3}=\frac{15}{4}\cdot 3$ we get $n_{6,3}=936$.
\item For $l=2$ we get
\[ 36 = 6 \cdot \frac{18}{3} \cdot N_{2,2}, \]
hence $N_{2,2}=1$. There are no multiple cover contributions, so $n_{2,2}=1$, and
\[ 2322 = 12\cdot\left(\frac{18}{3}\cdot N_{4,2}+\frac{18}{3}\cdot 16\right), \]
so $N_{4,2}=\frac{65}{4}$. Subtracting $M_1[4]n_{1,1}+M_2[2]n_{2,1}=\frac{35}{16}\cdot 1+\frac{9}{4}\cdot 0$ we get $n_{4,1}=14$. Moreover,
\[ 307164 = 18\cdot\left(\frac{18}{3}\cdot N_{6,2}+\frac{36}{3}\cdot 948\right), \]
hence $N_{6,2}=\frac{8533}{9}$. Subtracting $M_2[3]n_{2,2}=\frac{91}{9}\cdot 1$ we get $n_{6,2}=938$.
\item For $l=1$ we get
\[ 9 = 3 \cdot \frac{9}{3} \cdot N_{1,1}, \]
hence $N_{1,1}=1$. There are no multiple cover contributions, so $n_{1,1}=1$, and
\[ \frac{63}{2} = 6\cdot\left(\frac{9}{3}\cdot N_{2,1} + \frac{9}{3}\cdot 1\right), \]
so $N_{2,1}=\frac{3}{4}$. Subtracting $M_1[2]n_{1,1}=\frac{3}{4}\cdot 1$ we get $n_{2,1}=0$. Moreover,
\[ 246 = 9\cdot\left(\frac{9}{3}\cdot N_{3,1}+\frac{18}{3}\cdot 3\right), \]
so $N_{3,1}=\frac{28}{9}$. Subtracting $M_1[3]n_{1,1}=\frac{10}{9}\cdot 1$ we get $n_{3,1}=2$. Finally,
\[ 307041 = 18 \cdot \left(\frac{9}{3}\cdot N_{6,1}+\frac{9}{3}\cdot\frac{8533}{9}+\frac{18}{3}\cdot\frac{3789}{4}+\frac{18}{3}\cdot 948\right), \]
so $N_{6,1}=\frac{2842}{3}$. Subtracting $M_1[6]n_{1,1}+M_2[3]n_{2,1}+M_3[2]n_{3,1}=\frac{77}{6}\cdot 1+\frac{91}{9}\cdot 0+\frac{15}{4}\cdot 2$ we get $n_{6,1}=927$.
\end{compactenum}
In summary, the numbers $n_{d,k}$ for $d\leq 6$ are shown in Table \ref{tab:results}. The $n_{d,d}$ coincide with the $m_d$ in \cite{Ta2}, Theorem 1.4. The numbers $n_{d,k}$ for $d\leq 3$ are calculated in \cite{Ta1}. The sum $\sum_{k\mid d}s_kn_{d,k}$ is the log BPS number $n_d$ of $(\mathbb{P}^2,E)$. From $n_5$ and $n_{5,5}$ one also obtains $n_{5,1}$. To the best of my knowledge the numbers $n_{4,1}$, $n_{4,2}$, $n_{6,1}$, $n_{6,2}$ and $n_{6,3}$ are new.

\begin{table}[h!]
\begin{tabular}{|l|l|l|l|l|l|} 															\cline{1-1}
$n_{1,1}=1$ 	 															\\ \cline{1-2}
$n_{2,1}=0$	& $n_{2,2}=1$													\\ \cline{1-3}
$n_{3,1}=2$	& 			& $n_{3,3}=3$										\\ \cline{1-4}
$n_{4,1}=14$	& $n_{4,2}=14$	&			& $n_{4,4}=16$							\\ \cline{1-5}
$n_{5,1}=108$	&			&			&			& $n_{5,5}=113$				\\ \cline{1-6}
$n_{6,1}=927$	& $n_{6,2}=938$	& $n_{6,3}=936$	&			&			& $n_{6,6}=948$	\\ \hline
\end{tabular}
\vspace{7mm}
\label{tab:results}
\caption{The numbers $n_{d,k}$ calculated below.}
\end{table}

\subsubsection{Degenerating the divisor}						%
\label{S:calcP2b}

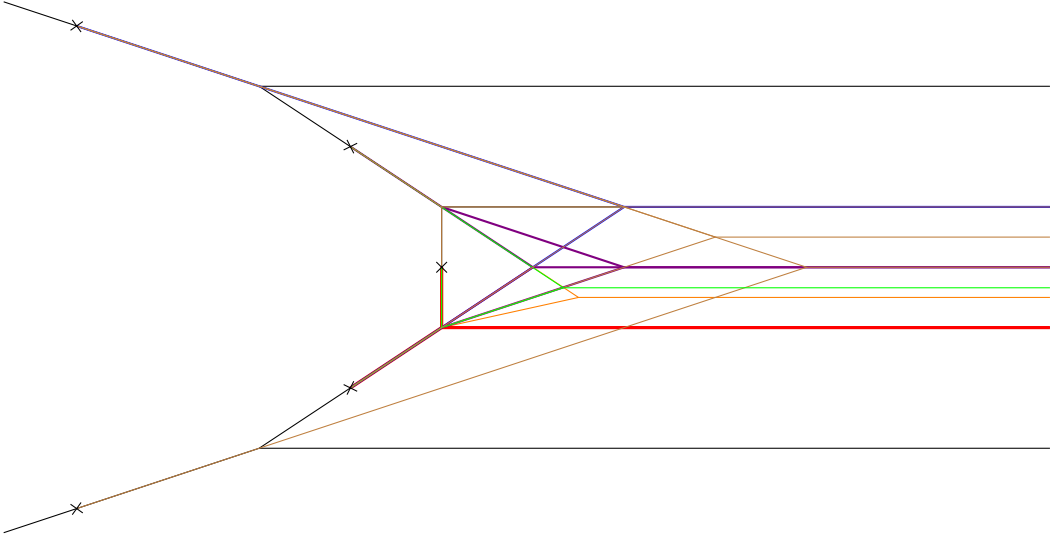
\begin{figure}[h!]
\centering
\begin{tikzpicture}[xscale=0.8,yscale=1.6]
\draw (-7.2,-1.7) -- (-3,-1) -- (0,0) -- (0,1) -- (-3,2) -- (-7.2,2.7);
\draw (-3,-1) -- (10,-1);
\draw (-3,2) -- (10,2);
\draw (0,0) -- (10,0);
\draw (0,1) -- (10,1);
\draw[orange] (-1.5,-0.5) -- (0,0) -- (0,0.5);
\draw[orange] (0,0) -- (2.25,0.25) -- (-1.5,1.5);
\draw[orange] (2.25,0.25) -- (10,0.25);
\draw[very thick,red] (-1.5,-0.5) -- (0,0) -- (0,0.5);
\draw[very thick,red] (0,0) -- (10,0);
\draw[thick, blue] (-6,2.5) -- (3,1) -- (-1.5,-0.5);
\draw[thick, blue] (3,1) -- (10,1);
\draw[thick,violet] (-1.5,-0.5) -- (1.5,0.5) -- (-1.5,1.5);
\draw[thick,violet] (1.5,0.5) -- (10,0.5);
\draw[thick,violet] (0,0) -- (3,0.5) -- (0,1);
\draw[thick,violet] (3,0.5) -- (10,0.5);
\draw[brown] (-1.5,-0.5) -- (0,0) -- (0,0.5);
\draw[brown] (0,0) -- (4.5,0.75) -- (-6,2.5);
\draw[brown] (4.5,0.75) -- (10,0.75);
\draw[brown] (-6,-1.5) -- (6,0.5) -- (-6,2.5);
\draw[brown] (6,0.5) -- (10,0.5);
\draw[green] (-1.5,-0.5) -- (0,0) -- (0,0.5);
\draw[green] (0,0) -- (2,0.33) -- (-1.5,1.5);
\draw[green] (2,0.33) -- (10,0.33);
\draw[brown] (0,0.5) -- (0,1) -- (-1.5,1.5);
\draw[brown] (-6,2.5) -- (3,1) -- (-1.5,-0.5);
\draw[brown] (0,1) -- (10,1);
\coordinate[fill,cross,inner sep=2pt] (0) at (0,0.5);
\coordinate[fill,cross,inner sep=2pt,rotate=18.4] (0) at (-1.5,-0.5);
\coordinate[fill,cross,inner sep=2pt,rotate=71.6] (0) at (-1.5,1.5);
\coordinate[fill,cross,inner sep=2pt,rotate=9.4] (0) at (-6,-1.5);
\coordinate[fill,cross,inner sep=2pt,rotate=80.6] (0) at (-6,2.5);
\end{tikzpicture}
\caption{Tropical curves contributing to $n_{[d_1,\ldots,d_k]}$ for $(\mathbb{P}^2,E)$.}
\label{fig:colors}
\end{figure}

The limit $s\rightarrow 0$ corresponds to a degeneration of the elliptic curve $E$ to a cycle of three lines $D_t^0=D_1+D_2+D_3$. In this section we consider invariants of $(\mathbb{P}^2,E)$ with given degrees over the components $D_i$, i.e., with given degree splitting $[d_1,\ldots,d_k]$ (Definition \ref{defi:splitting}). Since the components $D_i$ are isomorphic as divisors of $\mathbb{P}^2$ this does not depend on the labelling of the $D_i$ and we can omit zeros in $[d_1,\ldots,d_k]$. The tropical curves contributing to the invariants $N_{[d_1,\ldots,d_k]}$ for $d=\sum d_i \leq 4$ are shown in Figure \ref{fig:colors}. There is more than one tropical curve with degree splitting $[2,2]$ and $[2,1,1]$, respectively. We have:
\begin{eqnarray*}
f_{[1]} &=&1+9x^3 + \mathcal{O}(x^6) \\
f_{[2]} &=&(1+9x^3)(1+72x^6) + \mathcal{O}(t^9) \\
f_{[1,1]} &=& 1+36x^6 + \mathcal{O}(x^9) \\
f_{[3]} &=& (1+9x^3)(1+72x^6)(1-78x^9) + \mathcal{O}(x^{12}) \\
f_{[2,1]} &=& 1+243x^9 + \mathcal{O}(x^{12}) \\
f_{[1,1,1]} &=& 1+81x^9 + \mathcal{O}(x^{12}) \\
f_{[4]} &=& (1+9x^3)(1+72x^6)(1-78x^9)(1+5256x^{12}) + \mathcal{O}(x^{15}) \\
f_{[3,1]} &=& 1+1872x^{12} + \mathcal{O}(x^{15}) \\
f_{[2,2]} &=& (1+36x^6)(1+1296x^{12})(1+1530x^{12}) + \mathcal{O}(x^{15}) \\
f_{[2,1,1]} &=& (1+144x^{12})(1+432x^{12})(1+1296x^{12}) + \mathcal{O}(x^{15})
\end{eqnarray*}

The invariant $N_{[d]}$ has a $d$-fold cover contribution from $n_{[d]}$ and $N_{[2,2]}$ has a $2$-fold cover contribution from $n_{[1,1]}$. This gives the following log BPS numbers:

\begin{table}[h!]
\centering
\begin{tabular}{|l|l|l|l|} 																					\hhline{|-|~|~|~|}
\cellcolor{red!50}$n_{[1]}=3$ 	 																			\\ \hhline{|-|-|~|~|}
\cellcolor{red!50}$n_{[2]}=3$  	& \cellcolor{violet!50}$n_{[1,1]}=6$													\\ \hhline{|-|-|-|~|}
\cellcolor{red!50}$n_{[3]}=15$	& \cellcolor{green!50}$n_{[2,1]}=27$		& \cellcolor{blue!50}$n_{[1,1,1]}=9$						\\ \hhline{|-|-|-|-|}
\cellcolor{red!50}$n_{[4]}=72$	& \cellcolor{orange!50}$n_{[3,1]}=156$	& \cellcolor{violet!50}$n_{[2,2]}=168$	& \cellcolor{brown!70}$n_{[2,1,1]}=156$	\\ \hline
\end{tabular}
\vspace{7mm}
\label{tab:colors}
\caption{Log BPS numbers of $(\mathbb{P}^2,E)$ with given degrees over the $D_i$.}
\end{table}

\text{ } \\[-12mm] \text{ }

\subsubsection{Combining the methods}						%
\label{S:calcP2c}

As indicated in Remark \ref{rem:combine} in general it is not possible to combine the above methods and calculate invariants with given degree splitting and meeting a points on $E$ with prescribed torsion. However, for $d\leq 3$ this is indeed possible as we will show now.

For $d=1$ we have $\mathcal{G}_1=\{[1]\}$, so $n_{[1],1}=n_{1,1}=1$. For $d=2$ there are two types of tropical curves. One of them contributes to $N_{[1,1]}$ and corresponds to stable log maps meeting $E$ in a point of order $6$, so $n_{[1,1],1}=0$ and $n_{[1,1],2}=n_{2,2}=1$. The other one contributes to $N_{[2]}$, so $n_{[2],1}=n_{1,1}=0$.

For $d=3$ we have $l([3])=1$, $l([2,1])=3$ and $l([1,1,1])=1$. So $n_{[2,1],1}=0$ and the equations in Remark \ref{rem:combine} are the following:
\begin{align*}
3n_{[3],1} + 6n_{[3],3} &= 15  			& n_{[3],1} + n_{[1,1,1],1} &= 2 \\
9n_{[2,1],3} &= 27 				& n_{[3],3} + n_{[1,1,1],3} &= 3 \\
3n_{[1,1,1],1} + 6n_{[1,1,1],3} &= 9		& n_{[2,1],3} &= 3
\end{align*}
This system of linear equations has the unique solution:
\begin{table}[h!]
\centering
\begin{tabular}{|l|l|l|} 																					\hline
$n_{[3],1}=1$	& $n_{[2,1],1}=0$		& $n_{[1,1,1],1}=1$														\\ \hline
$n_{[3],3}=2$	& $n_{[2,1],3}=3$		& $n_{[1,1,1],3}=1$														\\ \hline
\end{tabular}
\vspace{7mm}
\label{tab:P1xP1}
\caption{Log BPS numbers of $(\mathbb{P}^2,E)$ for $d=3$ with prescribed degree splitting and torsion point.}
\end{table}
\text{ }\\
If we try the same method for $d>3$, we have more indeterminates than independent equations, so there is no unique solution.

\subsubsection{Higher genus}								%

Let us compute the higher genus numbers $N_d^g$ for $d\leq 2$ and $g\leq 4$. For $d=1$ we have three tropical curves, all isomorphic to each other. They have one vertex of classical multiplicity $m_V=3$ and two bounded legs of weight $w_E=1$ each. This gives, with $q=e^{i\hbar}$,
\begin{eqnarray*}
m_h(q) &=& \frac{1}{i\hbar}\left(q^{3/2}-q^{-3/2}\right) \cdot \left(\frac{i\hbar}{q^{1/2}-q^{-1/2}}\right)^2 \\
&=& 3 - \frac{7}{8}\hbar^2 + \frac{29}{640}\hbar^4 - \frac{137}{322560}\hbar^6 + \frac{41}{7372800}\hbar^8 + \mathcal{O}(\hbar^{10})
\end{eqnarray*}
Multiplying by $3$ gives
\begin{table}[h!]
\centering
\begin{tabular}{|c|c|c|c|c|} 								\hline
$N_1^0=9$		& $N_1^1=-\frac{21}{8}$	& $N_1^2=\frac{87}{640}$		& $N_1^3=-\frac{137}{107520}$	& $N_1^4=\frac{41}{2457600}$	\\ \hline
\end{tabular}
\vspace{7mm}
\label{tab:Ng}
\caption{Higher genus log Gromov-Witten invaraints $N_d^g$ for $(\mathbb{P}^2,E)$.}
\end{table}

\subsection{$\mathbb{P}^1\times\mathbb{P}^1$}					%%

Similarly, executing the code for $\mathbb{P}^1\times\mathbb{P}^1$ gives the following:
\begin{eqnarray*}
\text{log }f_{\text{out}} = 16x^2 + 72x^4 + 352x^6 + 3108x^8 + \frac{120016}{5}x^{10} + 198384x^{12} + \mathcal{O}(x^{14})
\end{eqnarray*}

\begin{table}[h!]
\centering
\begin{tabular}{|c|c|c|c|c|c|} 								\hline
$n_1=8$	& $n_2=16$		& $n_3=72$		& $n_4=368$	& $n_5=2400$	& $n_6=16320$	\\ \hline
\end{tabular}
\vspace{7mm}
\label{tab:n}
\caption{The log BPS numbers $n_d$ of $\mathbb{P}^1\times\mathbb{P}^1$ for $d \leq 6$.}
\end{table}

The functions and invariants with prescribed degree splitting are:
\begin{eqnarray*}
f_{[1]} &=& 1+4x^2 + \mathcal{O}(x^4) \\
f_{[2]} &=& (1+4x^2)(1+10x^4) + \mathcal{O}(x^6) \\
f_{[1,1]} &=& 1+16x^4 + \mathcal{O}(x^6) \\
f_{[3]} &=& (1+4x^2)(1+10x^4)(1-20x^6) + \mathcal{O}(x^8) \\
f_{[2,1]} &=&  1+36x^6 + \mathcal{O}(x^8) \\
f_{[1,1,1]} &=& 1+36x^6 + \mathcal{O}(x^8) \\
f_{[4]} &=& (1+4x^2)(1+10x^4)(1-20x^6)(1+115x^8) + \mathcal{O}(x^{10}) \\
f_{[3,1]} &=& 1+64x^8 + \mathcal{O}(x^{10}) \\
f_{[2,2]} &=& (1+16x^4)(1+64x^8)(1+264x^8) + \mathcal{O}(x^{10}) \\
f_{[2,1,1]} &=& 1+64x^8 + \mathcal{O}(x^{10}) \\
f_{[1,2,1]} &=& 1+256x^8 + \mathcal{O}(x^{10}) \\
f_{[1,1,1,1]} &=& 1+64x^8 + \mathcal{O}(x^{10})
\end{eqnarray*}
This gives the following log BPS numbers:

\begin{table}[h!]
\centering
\begin{tabular}{|l|l|l|l|l|l|} 																				\hhline{|-|~|~|~|~|~|}
$n_{[1]}=2$ 	 																				\\ \hhline{|-|-|~|~|~|~|}
$n_{[2]}=0$  	& $n_{[1,1]}=4$																		\\ \hhline{|-|-|-|~|~|~|}
$n_{[3]}=0$		& $n_{[2,1]}=6$		& $n_{[1,1,1]}=6$														\\ \hhline{|-|-|-|-|-|-|}
$n_{[4]}=0$		& $n_{[3,1]}=8$		& $n_{[2,2]}=24$		& $n_{[2,1,1]}=8$		& $n_{[1,2,1]}=32$	& $n_{[1,1,1,1]}=8$ 	\\ \hline
\end{tabular}
\vspace{7mm}
\label{tab:P1xP1}
\caption{The log BPS numbers of $\mathbb{P}^1\times\mathbb{P}^1$ with given degree splitting.}
\end{table}

From this we compute the log BPS numbers with given curve class (bidegree) as follows. The factors are the number of unbounded walls in the wall structure contributing to $n_{[d_1,\ldots,d_m]}$ for different labellings of $D=D_1+D_2+D_3+D_4$.
\begin{align*}
n_{(1,0)} &= 2 \cdot n_{[1]} &= 2 \\
n_{(2,0)} &= 2 \cdot n_{[2]} &= 0 \\
n_{(1,1)} &= 4 \cdot n_{[1,1]} &= 16 \\
n_{(3,0)} &= 2 \cdot n_{[3]} &= 0 \\
n_{(2,1)} &= 4 \cdot n_{[2,1]} + 2 \cdot n_{[1,1,1]} &= 36 \\
n_{(4,0)} &= 2 \cdot n_{[4]} &= 0 \\
n_{(3,1)} &= 4 \cdot n_{[3,1]} + 4 \cdot n_{[2,1,1]} &= 768 \\
n_{(2,2)} &= 4 \cdot n_{[2,2]} + 4 \cdot n_{[1,2,1]} + 4 \cdot n_{[1,1,1,1]} &= 1152
\end{align*}

\subsubsection{Deforming $\mathbb{F}_2$}						%
\label{S:calc8'a}

Consider the toric degeneration of $\mathbb{P}^1\times\mathbb{P}^1$ by deformation of the Hirzebruch surface $\mathbb{F}_2$ (case (8'a) in Figure \ref{fig:listb}) from Example \ref{expl:8'a2}.

\begin{figure}[h!]
\centering
\begin{tikzpicture}[scale=0.7]
% SOURCE
\draw[green,thick] (-8,3.5) -- (-9,4) node[fill,circle,inner sep=1pt]{} -- (-10,4.25);
\draw[green,thick] (-9,4) -- (-8.25,4) node[above]{$4$} -- (-7.5,4);
\draw[red,thick] (-10,2.5) -- (-9.5,2.25) node[above]{$2$} -- (-9,2) node[fill,circle,inner sep=1pt]{} -- (-8.9,1.5);
\draw[red,thick] (-9,2) -- (-9.1,1.5) node[fill,circle,inner sep=1pt]{} -- (-9.2,1);
\draw[red,thick] (-9,2) -- (-8.25,2) node[above]{$4$} -- (-7.5,2);
\draw[red,thick] (-10,0.4) -- (-9,0) node[fill,circle,inner sep=1pt]{} -- (-8.9,-0.5);
\draw[red,thick] (-9,0) -- (-9.1,-0.5) node[fill,circle,inner sep=1pt]{} -- (-9.2,-1);
\draw[red,thick] (-9,0) -- (-8.25,0) node[above]{$4$} -- (-7.5,0);
\draw[red,thick] (-10,0.6) -- (-9,0);
\draw[blue,thick] (-10,-2.5) -- (-9,-2) node[fill,circle,inner sep=1pt]{} -- (-10,-1.5);
\draw[blue,thick] (-9,-2) -- (-8.25,-2) node[above]{$4$} -- (-7.5,-2);
% ARROW
\draw[->,thick] (-7,1) -- (-6,1);
% TARGET
\draw[dashed] (-5.6,3.6) -- (-2,3);
\draw (-2,3) -- (0,2) -- (0,0) -- (-2,-1);
\draw[dashed] (-2,-1) -- (-5.6,-1.6);
\draw (-2,3) -- (6,3) node[right]{$L_1+L_2$};
\draw (0,2) -- (6,2) node[right]{$L_2$};
\draw (0,1) -- (6,1) node[right]{$L_1-L_2$};
\draw (0,0) -- (6,0) node[right]{$L_2$};
\draw[dashed] (-2,-1) -- (6,-1);
\coordinate[fill,cross,inner sep=2pt,rotate=80.54] (0) at (-5,3.5);
\coordinate[fill,cross,inner sep=2pt,rotate=63.43] (0) at (-1,2.5);
\coordinate[fill,cross,inner sep=2pt] (0) at (0,1.5);
\coordinate[fill,cross,inner sep=2pt] (0) at (0,0.5);
\coordinate[fill,cross,inner sep=2pt,rotate=26.57] (0) at (-1,-0.5);
\coordinate[fill,cross,inner sep=2pt,rotate=9.46] (0) at (-5,-1.5);
% CURVES
\draw[green,thick] (-1,2.5) -- (-2,3) -- (2,3) node[above]{$4$} -- (6,3);
\draw[green,thick,dashed] (-2,3) -- (-5,3.5);
\draw[red,very thick] (-1,2.5) -- (0,2) node[left]{$2$} -- (0,1) node[left]{$1$} -- (0,0.5);
\draw[red,very thick] (0,2) -- (3,2) node[above]{$4$} -- (6,2);
\draw[blue,thick] (-1,2.5) -- (2,1) -- (-1,-0.5);
\draw[blue,thick] (2,1) -- (4,1) node[above]{$4$} -- (6,1);
\end{tikzpicture}
\caption{Tropical curves corresponding to $\underline{\beta}=L_1+L_2$.}
\label{fig:calc8'a}
\end{figure}
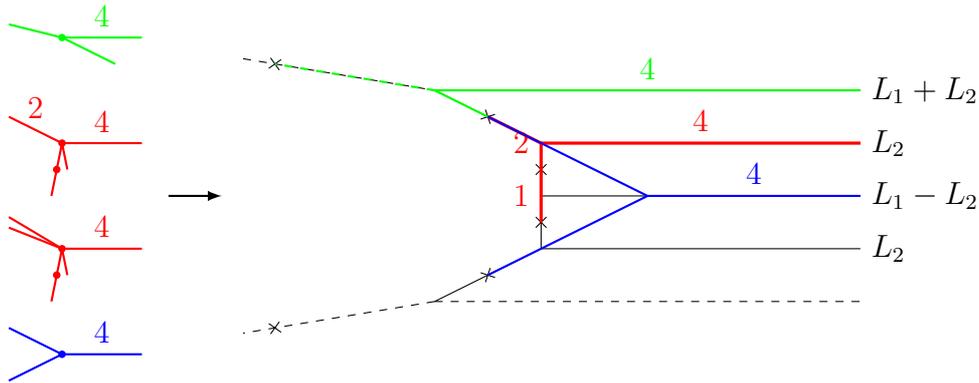

For $d=1$ it is clear by the symmetry $L_1\leftrightarrow L_2$ that $n_{(1,0)}=n_{L_1}=n_{L_2}=4$. For $d=2$ we have $n_{(2,0)}=n_{2L_1}=n_{2L_2}$ and $n_{(1,1)}=n_{L_1+L_2}$. Figure \ref{fig:calc8'a} shows the tropical curves corresponding to stable log maps of class $\underline{\beta}=L_1+L_2$. The first one has multiplicity $4$, the second one $-4$, the third one $8$ and the last one $4$. By symmetry there are two tropical curves similar to the red ones at the lower vertex, again with multiplicities $-4$ and $8$. This gives $n_{(1,1)}=n_{L_1+L_2}=16$ and in turn $n_{(2,0)}=0$. One can proceed similarly for higher degrees.

\subsection{Cubic surface}								%%
\label{S:calccubic}

The dual intersection complex of the cubic surface (case (3a) in Figure \ref{fig:listb}) is quite similar to the one of $(\mathbb{P}^2,E)$. The only differences are that for each vertex the determinant of primitive generators of adjacent bounded edges is $1$ instead of $3$ and that there are three affine singularities on each bounded edge. As a consequence, by the change of lattice trick (\cite{GHK2}, Proposition C.13), the wall structure of $(X,D)$ is in bijection with the wall structure of $(\mathbb{P}^2,E)$, and the functions attached to walls in direction $m_{\text{out}}$, in particular the unbounded walls, coincide. This immediately implies $N_d(X,D) = 3 \cdot N_d(\mathbb{P}^2,E)$. Subtracting multiple covers we get (note that $w_{\text{out}}=1$ and the multiple cover contributions of of degree $1$ curves are $M_1[k]=\frac{1}{k^2}\binom{-1}{k-1}=\frac{(-1)^{k-1}}{k^2}$):

\begin{table}[h!]
\begin{tabular}{|c|c|c|c|c|c|} 								\hline
$n_1=27$	& $n_2=108$	& $n_3=729$	& $n_4=6912$	& $n_5=76275$	& $n_6=920727$	\\ \hline
\end{tabular}
\vspace{7mm}
\label{tab:n}
\caption{The log BPS numbers $n_d$ of the cubic surface for $d \leq 6$.}
\end{table}

\subsubsection{Curve classes}								%
A smooth cubic surface $X$ can be given by blowing up six general points on $\mathbb{P}^2$. Let $e_1,\ldots,e_1$ be the classes of the exceptional divisors and let $\ell$ be the pullback of the class of a line in $\mathbb{P}^2$. Then $\ell,e_1,\ldots,e_6$ generate $\text{Pic}(X)\simeq H_2^+(X,\mathbb{Z})\simeq\mathbb{Z}^7$.

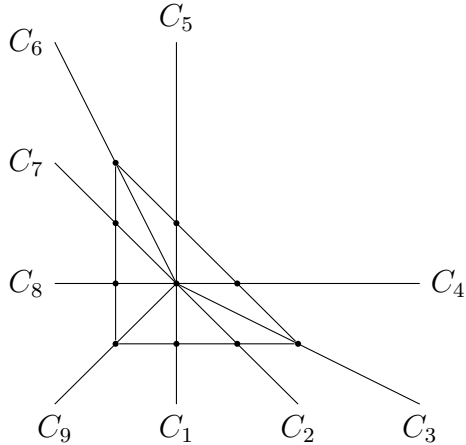
\begin{figure}[h!]
\centering
\begin{tikzpicture}[scale=0.8]
\coordinate[fill,circle,inner sep=0.8pt] (0) at (0,0);
\coordinate[fill,circle,inner sep=0.8pt] (1) at (-1,-1);
\coordinate[fill,circle,inner sep=0.8pt] (2) at (2,-1);
\coordinate[fill,circle,inner sep=0.8pt] (3) at (-1,2);
\coordinate (1a) at (-2,-2);
\coordinate (2a) at (4,-2);
\coordinate (3a) at (-2,4);
\coordinate (12) at (-0.5,-1);
\coordinate (12a) at (-1.5,-2);
\coordinate (12b) at (-0.5,-2);
\coordinate (12') at (0.5,-1);
\coordinate (12'a) at (0.5,-2);
\coordinate (12'b) at (1.5,-2);
\coordinate (12'') at (1.5,-1);
\coordinate (12''a) at (2.5,-2);
\coordinate (12''b) at (3.5,-2);
\coordinate (23) at (-0.5,1.5);
\coordinate (23a) at (-0.5,4);
\coordinate (23b) at (-1.75,4);
\coordinate (23') at (0.5,0.5);
\coordinate (23'a) at (4,0.5);
\coordinate (23'b) at (0.5,4);
\coordinate (23'') at (1.5,-0.5);
\coordinate (23''a) at (4,-0.5);
\coordinate (23''b) at (4,-1.75);
\coordinate (31) at (-1,-0.5);
\coordinate (31a) at (-2,-1.5);
\coordinate (31b) at (-2,-0.5);
\coordinate (31') at (-1,0.5);
\coordinate (31'a) at (-2,0.5);
\coordinate (31'b) at (-2,1.5);
\coordinate (31'') at (-1,1.5);
\coordinate (31''a) at (-2,2.5);
\coordinate (31''b) at (-2,3.5);
\draw (1) -- (2) -- (3) -- (1);
\draw (0,0) -- (1a) node[below]{$C_9$};
\draw (0,0) -- (0,-2) node[below]{$C_1$};
\draw (0,0) -- (2,-2) node[below]{$C_2$};
\draw (0,0) -- (2a) node[below]{$C_3$};
\draw (0,0) -- (4,0) node[right]{$C_4$};
\draw (0,0) -- (0,4) node[above]{$C_5$};
\draw (0,0) -- (3a) node[left]{$C_6$};
\draw (0,0) -- (-2,2) node[left]{$C_7$};
\draw (0,0) -- (-2,0) node[left]{$C_8$};
\coordinate[fill,circle,inner sep=0.8pt] (b) at (-1,0);
\coordinate[fill,circle,inner sep=0.8pt] (c) at (-1,1);
\coordinate[fill,circle,inner sep=0.8pt] (d) at (0,-1);
\coordinate[fill,circle,inner sep=0.8pt] (e) at (1,-1);
\coordinate[fill,circle,inner sep=0.8pt] (f) at (1,0);
\coordinate[fill,circle,inner sep=0.8pt] (g) at (0,1);
\end{tikzpicture}
\caption{The dual intersection complex of a smooth nef toric surface $X^0$ deforming to a smooth cubic surface $X$.}
\label{fig:fan}
\end{figure}

The dual intersection complex of a smooth nef toric surface $X^0$ deforming to $X$ is shown in Figure \ref{fig:fan}. Its asymptotic fan is the fan of $X^0$. Denote the curves corresponding to the rays of this fan by $C_1,\ldots,C_9$, labelled as in Figure \ref{fig:fan}. Then (see \cite{KM}, \S4) an isomorphism $\text{Pic}(X^0) \simeq \text{Pic}(X)$ is given as follows:
\begin{align*}
[C_1] &\mapsto e_2-e_5 & [C_2] &\mapsto \ell-e_2-e_3-e_6 & [C_3] &\mapsto e_6 \\
[C_4] &\mapsto e_3-e_6 & [C_5] &\mapsto \ell-e_1-e_3-e_4 & [C_6] &\mapsto e_4 \\
[C_7] &\mapsto e_1-e_4 & [C_8] &\mapsto \ell-e_1-e_2-e_5 & [C_9] &\mapsto e_5 \\
\end{align*}
\text{ }\\[-14mm]
Now we know the curve classes of the cubic surface $X$ corresponding to the unbounded edges in the dual intersection complex $(B,\mathscr{P},\varphi)$. In turn, we are able to associate to each tropical curve in $(B,\mathscr{P},\varphi)$ the curve class of the corresponding stable log maps, by composition of the maps from Constructions \ref{con:HG} and \ref{con:GH}.

As shown in \cite{Hos}, the Weyl group $W_{E_6}$ of type $E_6$ acts on $\text{Pic}(X)$ as symmetries of configurations of the $27$ lines and this action preserves the local Gromov-Witten invariants $N_\beta^{\text{loc}}$ of $X$. Hence, by the log-local correspondence \cite{vGGR}, it preserves the logarithmic Gromov-Witten invariants $N_\beta$ of $X$ considered here. The curve classes $\underline{\beta}$ of the cubic $X$ giving a nonzero contribution $N_\beta^{\text{loc}}$, up to action of $W_{E_6}$, are given in \cite{KM}, Table 1, along with the corresponding local BPS number $n_\beta^{\text{loc}}$. 

For $d=1$ and $d=2$ there is, up to the action of $W_{E_6}$, only one curve class giving a nonzero contribution, so this is trivial. For $d=1$ this is $\underline{\beta}=e_6$ and the length of its orbit is $27$, so $n_\beta=1$. For $d=2$ it is $\underline{\beta}=\ell-e_1$, with orbit length $27$, so $n_\beta=4$. For $d=3$ there are two equivalence classes giving a nonzero contribution, with representatives $\ell$ and $3\ell-\sum_{i=1}^6 e_i$, respectively.

\begin{figure}[h!]
\centering
\begin{tikzpicture}[scale=0.65]
\draw[dashed] (-5,-4) -- (-3,-3);
\draw (-3,-3) -- (0,0) -- (0,3) -- (-3,6);
\draw[dashed] (-3,6) -- (-5,7);
\draw (-3,6) -- (10,6) node[right]{$e_5$};
\draw (-2,5) -- (10,5) node[right]{$e_2-e_5$};
\draw (-1,4) -- (10,4) node[right]{$\ell-e_2-e_3-e_6$};
\draw (0,3) -- (10,3) node[right]{$e_6$};
\draw (0,2) -- (10,2) node[right]{$e_3-e_6$};
\draw (0,1) -- (10,1) node[right]{$\ell-e_1-e_3-e_4$};
\draw (0,0) -- (10,0) node[right]{$e_4$};
\draw (-1,-1) -- (10,-1) node[right]{$e_1-e_4$};
\draw (-2,-2) -- (10,-2) node[right]{$\ell-e_1-e_2-e_5$};
\draw[dashed] (-3,-3) -- (10,-3);
\coordinate[fill,cross,inner sep=2pt,rotate=63.43] (0) at (-4,6.5);
\coordinate[fill,cross,inner sep=2pt,rotate=45] (0) at (-0.5,3.5);
\coordinate[fill,cross,inner sep=2pt,rotate=45] (0) at (-1.5,4.5);
\coordinate[fill,cross,inner sep=2pt,rotate=45] (0) at (-2.5,5.5);
\coordinate[fill,cross,inner sep=2pt] (0) at (0,0.5);
\coordinate[fill,cross,inner sep=2pt] (0) at (0,1.5);
\coordinate[fill,cross,inner sep=2pt] (0) at (0,2.5);
\coordinate[fill,cross,inner sep=2pt,rotate=45] (0) at (-0.5,-0.5);
\coordinate[fill,cross,inner sep=2pt,rotate=45] (0) at (-1.5,-1.5);
\coordinate[fill,cross,inner sep=2pt,rotate=45] (0) at (-2.5,-2.5);
\coordinate[fill,cross,inner sep=2pt,rotate=26.57] (0) at (-4,-3.5);
% CURVES
\draw[red,very thick] (-2.5,5.5) -- (-2,5) node[left]{$1$} -- (-1,4) node[left]{$2$} -- (0,3) node[left]{$3$} -- (0,2) node[left]{$2$} -- (0,1) node[left]{$1$} -- (0,0.5);
\draw[red,very thick] (0,3) -- (5,3) node[above]{$3$} -- (10,3);
\draw[blue,thick] (-4,6.5) -- (3,3) -- (-2.5,-2.5);
\draw[blue,thick] (3,3) -- (7,3) node[above]{$3$} -- (10,3);
\end{tikzpicture}
\caption{Tropical curves corresponding to $\underline{\beta}=3\ell-\sum_{i=1}^6e_i$.}
\label{fig:calccubic}
\end{figure}
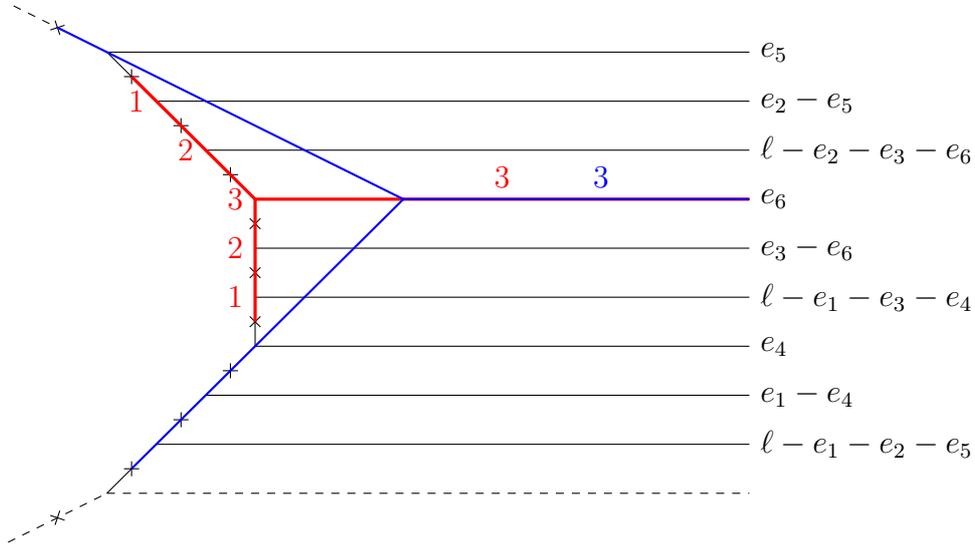

The red tropical curve in Figure \ref{fig:calccubic} corresponds to the class
\[ 1 \cdot (e_2-e_5) + 2 \cdot (\ell-e_2-e_3-e_6) + 3 \cdot e_6 + 2 \cdot (e_3-e_6) + 1 \cdot (\ell-e_1-e_3-e_4) = 3\ell-\sum_{i=1}^6e_i \]
and similarly for the blue tropical curve. Changing the affine singularities in which the bounded legs end may change the curve class. It turns out that for the red tropical curve any change leads to the class $\underline{\beta}=\ell$ or to a class giving a nonzero contribution. Its multiplicity is $18$, so together with the choice of outgoing edge this gives a contribution of $54$ to $n_{3\ell-\sum_{i=1}^6e_i}$. For the blue tropical curve there are two changes leading again to $\underline{\beta}=3\ell-\sum_{i=1}^6e_i$ and six changes leading to $\underline{\beta}=\ell$. The multiplicity of any of these tropical curves is $3$. Together with the choice of outgoing edge this gives a contribution of $3 \cdot 3 \cdot 3 = 27$ to $n_{3\ell-\sum_{i=1}^6e_i}$. The orbit length is $1$, so $n_{3\ell-\sum_{i=1}^6e_i}=81$. The orbit length of $\ell$ is $72$, so $n_\ell = (729-81)/72 = 9$. So, in agreement with \cite{KM}, Table 1, we have:
\begin{table}[h!]
\begin{tabular}{|c|c|c|c|} 											\hline
$n_{e_i}=1$	& $n_{\ell-e_i}=4$	& $n_\ell=9$	& $n_{3\ell-\sum_{i=1}^6e_i}=81$ 	\\ \hline
\end{tabular}
\vspace{7mm}
\label{tab:n}
\caption{The log BPS numbers $n_\beta$ of the cubic surface for $d \leq 3$.}
\end{table}

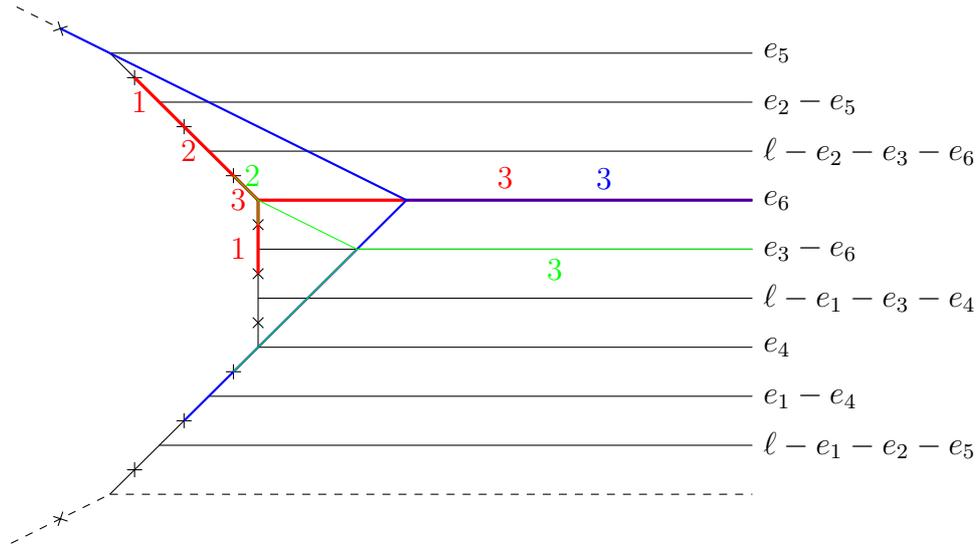
\begin{figure}[h!]
\centering
\begin{tikzpicture}[scale=0.65]
\draw[dashed] (-5,-4) -- (-3,-3);
\draw (-3,-3) -- (0,0) -- (0,3) -- (-3,6);
\draw[dashed] (-3,6) -- (-5,7);
\draw (-3,6) -- (10,6) node[right]{$e_5$};
\draw (-2,5) -- (10,5) node[right]{$e_2-e_5$};
\draw (-1,4) -- (10,4) node[right]{$\ell-e_2-e_3-e_6$};
\draw (0,3) -- (10,3) node[right]{$e_6$};
\draw (0,2) -- (10,2) node[right]{$e_3-e_6$};
\draw (0,1) -- (10,1) node[right]{$\ell-e_1-e_3-e_4$};
\draw (0,0) -- (10,0) node[right]{$e_4$};
\draw (-1,-1) -- (10,-1) node[right]{$e_1-e_4$};
\draw (-2,-2) -- (10,-2) node[right]{$\ell-e_1-e_2-e_5$};
\draw[dashed] (-3,-3) -- (10,-3);
\coordinate[fill,cross,inner sep=2pt,rotate=63.43] (0) at (-4,6.5);
\coordinate[fill,cross,inner sep=2pt,rotate=45] (0) at (-0.5,3.5);
\coordinate[fill,cross,inner sep=2pt,rotate=45] (0) at (-1.5,4.5);
\coordinate[fill,cross,inner sep=2pt,rotate=45] (0) at (-2.5,5.5);
\coordinate[fill,cross,inner sep=2pt] (0) at (0,0.5);
\coordinate[fill,cross,inner sep=2pt] (0) at (0,1.5);
\coordinate[fill,cross,inner sep=2pt] (0) at (0,2.5);
\coordinate[fill,cross,inner sep=2pt,rotate=45] (0) at (-0.5,-0.5);
\coordinate[fill,cross,inner sep=2pt,rotate=45] (0) at (-1.5,-1.5);
\coordinate[fill,cross,inner sep=2pt,rotate=45] (0) at (-2.5,-2.5);
\coordinate[fill,cross,inner sep=2pt,rotate=26.57] (0) at (-4,-3.5);
% CURVES
\draw[red,very thick] (-2.5,5.5) -- (-2,5) node[left]{$1$} -- (-1,4) node[left]{$2$} -- (0,3) node[left]{$3$} -- (0,2) node[left]{$1$} -- (0,1.5);
\draw[red,very thick] (0,3) -- (5,3) node[above]{$3$} -- (10,3);
\draw[blue,thick] (-4,6.5) -- (3,3) -- (-1.5,-1.5);
\draw[blue,thick] (3,3) -- (7,3) node[above]{$3$} -- (10,3);
\draw[green] (-0.5,3.5) node[right]{$2$} -- (0,3) -- (0,2.5);
\draw[green] (0,3) -- (2,2) -- (-0.5,-0.5);
\draw[green] (2,2) -- (6,2) node[below]{$3$} -- (10,2);
\end{tikzpicture}
\caption{Tropical curves corresponding to $\underline{\beta}=\ell$.}
\label{fig:calccubic2}
\end{figure}

Figure \ref{fig:calccubic2} shows tropical curves corresponding to stable log maps of class $\underline{\beta}=\ell$. In particular, any change of bounded legs of the green tropical curves still leads to class $\underline{\beta}=\ell$. For instance, the red tropical curve has class
\[ 1 \cdot (e_2-e_5) + 2 \cdot (\ell-e_2-e_3-e_6) + 3 \cdot e_6 + 1 \cdot (e_3-e_6) = 2\ell-e_2-e_3-e_5 \]
Under the action of $W_{E_6}$ this is equivalent to $\underline{\beta}=2\ell-e_1-e_2-e_3$ and in turn (via the map $s_6$ from \cite{KM}, \S3) to
\[ 2 \cdot (2\ell-e_1-e_2-e_3)-(\ell-e_2-e_3)-(\ell-e_1-e_3)-(\ell-e_1-e_2) = \ell. \]
Similarly, one computes the classes of the other tropical curves.

\appendix

\section{Artin fans and logarithmic modifications}					%%%
\label{A:artin}

An \textit{Artin fan} is a logarithmic algebraic stack that is logarithmically \'etale over a point. Artin fans were introduced in \cite{AW} to prove the invariance of logarithmic Gromov-Witten invariants under \textit{logarithmic modifications}, that is, proper birational logarithmically \'etale morphisms. We will briefly summarize this subject.

To any fine saturated log smooth scheme $X$ one can define an associated Artin fan $\mathcal{A}_X$. It has an \'etale cover by finitely many \textit{Artin cones} -- stacks of the form $[V/T]$, where $V$ is a toric variety and $T$ its dense torus. In this way, $\mathcal{A}_X$ encodes the combinatorial structure of $X$. A \textit{subdivision} of the Artin fan $\mathcal{A}_X$ induces a logarithmic modification of $X$ via pull-back. Moreover, all logarithmic modifications of $X$ arise this way. This ultimately leads to a proof of the birational invariance in logarithmic Gromov-Witten theory \cite{AW}.

Olsson \cite{Ol} showed that a logarithmic structure on a given underlying scheme $\underline{X}$ is equivalent to a morphism $\underline{X}\rightarrow\underline{\textbf{Log}}$, where $\underline{\textbf{Log}}$ is a zero-dimensional algebraic stack -- the moduli stack of logarithmic structures. It carries a universal logarithmic structure whose associated logarithmic algebraic stack we denote by $\textbf{Log}$ -- providing a universal family of logarithmic structures $\textbf{Log}\rightarrow\underline{\textbf{Log}}$. As shown in \cite{AW}, if $X$ is a fine saturated log smooth scheme, then the morphism $X\rightarrow\textbf{Log}$ factors through an initial morphism $X\rightarrow\mathcal{A}_X$, where $\mathcal{A}_X$ is an Artin fan and $\mathcal{A}_X\rightarrow\textbf{Log}$ is \'etale and representable. While this serves as a definition of the associated Artin fan $\mathcal{A}_X$, there is a more explicit description of $\mathcal{A}_X$ in terms of the \textit{tropicalization} of $X$, given below.

Let $S$ be a fine saturated log scheme.

\begin{defi}
\label{defi:Log}
The \textit{moduli stack of log structures over $S$} is the category $\underline{\textbf{Log}_S}$ fibered over the category of $\underline{S}$-schemes defined as follows. The objects over a scheme morphism $\underline{X}\rightarrow\underline{S}$ are the log morphisms $X\rightarrow S$ over $\underline{X}\rightarrow\underline{S}$. The morphisms from $X\rightarrow S$ to $X'\rightarrow S$ are the log morphisms $h : X \rightarrow X'$ over $S$ for which $h^\star\mathcal{M}_{X'}\rightarrow\mathcal{M}_X$ is an isomorphism.
\end{defi}

\begin{prop}[\cite{Ol} Theorem 1.1]
$\underline{\textup{\textbf{Log}}_S}$ is an algebraic stack locally of finite presentation over $\underline{S}$.
\end{prop}

\begin{defi}
An \textit{Artin fan} is a logarithmic algebraic stack that is logarithmically \'etale over a point. An \textit{Artin cone} is a logarithmic algebraic stack isomorphic to $[V/T]$, where $V$ is a toric variety and $T$ its dense torus.
\end{defi}

\begin{rem}
\label{rem:artin}
If a logarithmic algebraic stack has a strict representable \'etale cover by Artin cones, then it is an Artin fan. In fact, in \cite{AW} Artin fans were defined this way. Later the definition was generalized to the one above.
\end{rem}

\begin{lem}[\cite{AW}, Lemma 2.3.1]
\label{lem:artin}
An algebraic stack that is representable and \'etale over $\textup{\textbf{Log}}$ has a strict \'etale cover by Artin cones.
\end{lem}

\begin{prop}[\cite{ACMW}, Proposition 3.1.1]
\label{prop:artin}
Let $X$ be a logarithmic algebraic stack that is locally connected in the smooth topology. Then there is an initial factorization of $X\rightarrow \textbf{Log}$ through a strict \'etale  morphism $\mathcal{A}_X \rightarrow \textbf{Log}$ which is representable by algebraic spaces.
\end{prop}

\begin{defi}
Let $X$ be a fine saturated log smooth scheme. The \textit{Artin fan of $X$} is the stack $\mathcal{A}_X$ from Proposition \ref{prop:artin}. Indeed, this is an Artin fan by Lemma \ref{lem:artin} and Remark \ref{rem:artin}.
\end{defi}

We now give a more explicit description of the Artin fan $\mathcal{A}_X$ of a fine saturated log smooth scheme $X$. By Lemma \ref{lem:artin} and Proposition \ref{prop:artin}, $\mathcal{A}_X$ has a strict \'etale cover by Artin cones. In fact, $\mathcal{A}_X$ is a colimit of Artin cones $\mathcal{A}_\sigma$ corresponding to the cones $\sigma$ in the tropicalization $\Sigma(X)$ of $X$.

\begin{defi}
For a cone $\sigma\subseteq N_{\mathbb{R}}$, let $P=\sigma^\vee\cap M$ be the corresponding monoid. The \textit{Artin cone defined by $\sigma$} is the logarithmic algebraic stack
\[ \mathcal{A}_\sigma = \left[\faktor{\textup{Spec }\mathbb{C}[P]}{\textup{Spec }\mathbb{C}[P^{\textup{gp}}]}\right] \]
with the toric log structure coming from the global chart $P \rightarrow \mathbb{C}[P]$.
\end{defi}

\begin{defi}
Let $\Sigma$ be a generalized cone complex (Definition \ref{defi:Cones}) that is a colimit of a diagram of cones $s : I \rightarrow \textup{\textbf{Cones}}$. Then define $\mathcal{A}_\Sigma$ to be the colimit as sheaves over $\textup{Log}$ of the corresponding diagram of sheaves given by $I \ni i \mapsto \mathcal{A}_{s(i)}$.
\end{defi}

\begin{prop}[\cite{ACGS1}, Proposition 2.2.2]
Let $X$ be a fine saturated log smooth scheme with tropicalization (Definition \ref{defi:trop}) a generalized cone complex $\Sigma(X)$. Then 
\[ \mathcal{A}_X \cong \mathcal{A}_{\Sigma(X)}. \]
\end{prop}

\begin{defi}
A \textit{subdivision} of an Artin fan $\mathcal{X}$ is a morphism of Artin fans $\mathcal{Y}\rightarrow\mathcal{X}$ whose base change via any map $\mathcal{A}_\sigma\rightarrow\mathcal{X}$ is isomorphic to $\mathcal{A}_\Sigma$ for some subdivision $\Sigma$ of $\sigma$.
\end{defi}

\begin{defi}
A \textit{logarithmic modification} of fine saturated log smooth schemes is a proper surjective logarithmically \'etale morphism.
\end{defi}

Let $X$ be a fine saturated log smooth scheme with tropicalization $\Sigma(X)$. Then a subdivision $\tilde{\Sigma}(X)$ of $\Sigma(X)$ gives a subdivision $\mathcal{A}_{\tilde{\Sigma}(X)} \rightarrow \mathcal{A}_X$ of the Artin fan of $X$. The pull back $\tilde{X} := \mathcal{A}_{\tilde{\Sigma}(X)} \times_{\mathcal{A}_X} X \rightarrow X$ is a logarithmic modification. Moreover, all logarithmic modifications of $X$ arise this way:

\begin{prop}[\cite{AW}, Corollary 2.6.7]
If $Y\rightarrow X$ is a logarithmic modification of fine saturated log smooth schemes, then $Y\rightarrow\mathcal{A}_Y \times_{\mathcal{A}_X} X$ is an isomorphism.
\end{prop}

\begin{thm}[\cite{AW}, Theorem 1.1]
\label{thm:AW}
Let $h : Y \rightarrow X$ be a logarithmic modification of log smooth schemes. This induces a projection $\pi : \bar{\mathscr{M}}(Y) \rightarrow \bar{\mathscr{M}}(X)$ with
\[ \pi_\star\llbracket\bar{\mathscr{M}}(Y)\rrbracket = \llbracket\bar{\mathscr{M}}(X)\rrbracket, \]
where $\bar{\mathscr{M}}(X)$ is the stack of stable log maps to $X$.
\end{thm}

\begin{cor}
Logarithmic Gromov-Witten invariants are invariant under logarithmic modifications.
\end{cor}

%%%% Literature


\begin{thebibliography}{9}
\bibitem[AC]{AC} D. Abramovich, Q. Chen, \textit{Stable logarithmic maps to Deligne-Faltings pairs II}, Asian J. Math. 18 (3), 465--488, 2014.
\bibitem[ACC+]{ACC+} M. Akhtar, T. Coates, A. Corti, L. Heuberger, A. Kasprzyk, A. Oneto, A. Petracci, T.  Prince,  K.  Tveiten, \textit{Mirror symmetry and the classification of orbifold del Pezzo surfaces}, Proc. Amer. Math. Soc. 144 (2), 513--527, 2016.
\bibitem[ACGS1]{ACGS1} D. Abramovich, Q. Chen, M. Gross, B. Siebert, \textit{Decomposition of degenerate Gromov-Witten invariants}, Compos. Math. 156 (10), 2020--2075, 2020.
\bibitem[ACGS2]{ACGS2} D. Abramovich, Q. Chen, M. Gross, B. Siebert, \textit{Punctured logarithmic maps}, \href{https://arxiv.org/abs/2009.07720}{arXiv:2009.07720}, 2020.
\bibitem[ACMW]{ACMW} D. Abramovich, Q. Chen, S. Marcus, and J. Wise, \textit{Boundedness of the space of stable logarithmic maps}, J. Eur. Math. Soc. 19 (9), 2783--2809, 2017.
\bibitem[AF]{AF} D. Abramovich and B. Fantechi, \textit{Orbifold techniques in degeneration formulas}, Ann. Sc. Norm. Super. Pisa Cl. Sci. 16 (5), 519--579, 2016.
\bibitem[AMW]{AMW} D. Abramovich, S. Marcus, J, Wise, \textit{Comparison theorems for Gromov-Witten invariants of smooth pairs and of degenerations}, Ann. Inst. Fourier 64 (4), 1611--1667, 2014.
\bibitem[AW]{AW} D. Abramovich, J. Wise, \textit{Birational invariance in logarithmic Gromov-Witten theory}, Compos. Math. 154 (3), 595--620, 2018.
\bibitem[Arg]{Arg} H. Arg\"uz, \textit{Mirror symmetry for the Tate curve via tropical and log corals}, \href{https://arxiv.org/abs/1712.10260}{arXiv:1712.10260}, 2017.
\bibitem[BN]{BN} L. J. Barrott, N. Nabijou, \textit{Tangent curves to degenerating hypersurfaces}, \href{https://arxiv.org/abs/2007.05016}{arXiv:2007.05016}, 2020.
\bibitem[Bou1]{Bou1} P. Bousseau, \textit{Tropical refined curve counting from higher genera and lambda classes}, Invent. Math. 215 (1), 1--79, 2019.
\bibitem[Bou2]{Bou2} P. Bousseau, \textit{The quantum tropical vertex}, Geom. Topol. 24, 1297--1379, 2020.
\bibitem[Bou3]{Bou3} P. Bousseau, \textit{Scattering diagrams, stability conditions, and coherent sheaves on $\mathbb{P}^2$}, \href{https://arxiv.org/abs/1909.02985}{arXiv:1909.02985}, 2019.
\bibitem[Bou4]{Bou4} P. Bousseau, \textit{A proof of N. Takahashi's conjecture on genus zero Gromov-Witten theory of $(\mathbb{P}^2,E)$},\href{https://arxiv.org/abs/1909.02992}{arXiv:1909.02992}, 2019.
\bibitem[BP]{BP} J. Bryan, R. Pandharipande, \textit{Curves in Calabi-Yau threefolds and Topological Quantum Field Theory}, Duke Math. J. 126 (2), 369--396, 2005.
\bibitem[CGG+]{CGG+} T. Coates, A. Corti, S. Galkin, V. Golyshev, A. Kasprzyk, \textit{Mirror symmetry and Fano manifolds}, in European Congress of Mathematics Krak\'ow, 2--7 July, 2012 (Eur. Math. Soc., Z\"urich, 2014), 285--300.
\bibitem[CPS]{CPS} M. Carl, M. Pumperla, B. Siebert, \textit{A tropical view on Landau-Ginzburg models}, \href{https://www.math.uni-hamburg.de/home/siebert/preprints/LGtrop.pdf}{preprint}, 2010.
\bibitem[Che1]{Che1} Q. Chen, \textit{Stable logarithmic maps to Deligne-Faltings pairs I}, Ann. of Math. 180 (2), 455--521, 2014.
\bibitem[Che2]{Che2} Q. Chen, \textit{The degeneration formula for logarithmic expanded degenerations}, J. Alg. Geom. 23 (2), 341-392, 2014.
\bibitem[CGKT1]{CvGKT} J. Choi, M. van Garrel, S. Katz, N. Takahashi, \textit{Local BPS invariants: enumerative aspects and wall-crossing}, Int. Math. Res. Notices, rny 171, 2018
\bibitem[CGKT2]{CvGKT2} J. Choi, M. van Garrel, S. Katz, N. Takahashi, \textit{Log BPS numbers of log Calabi-Yau surfaces}, \href{https://arxiv.org/abs/1810.02377}{arXiv:1810.02377}, 2018.
\bibitem[CGKT3]{CvGKT3} J. Choi, M. van Garrel, S. Katz, N. Takahashi, \textit{Sheaves of maximal intersection and multiplicities of stable log maps}, \href{https://arxiv.org/abs/1908.10906}{arXiv:1908.10906}, 2019.
\bibitem[CKYZ]{CKYZ} T.-M. Chiang, A. Klemm, S.-T. Yau, E. Zaslow, \textit{Local mirror symmetry: calculations and interpretations}, Adv. Theor. Math. Phys. 3 (3), 495--565, 1999.
\bibitem[Fuk]{Fuk} K. Fukaya, \textit{Multivalued Morse theory, asymptotic analysis and mirror symmetry}, Proc. Sympos. Pure Math. 73, 205--278, 2005.
\bibitem[Ful]{Ful} W. Fulton, \textit{Intersection theory}, Springer, New York, ISBN: 978-0-387-98549-7, 1998.
\bibitem[GGR]{vGGR} M. van Garrel, T. Graber, H. Ruddat, \textit{Local Gromov-Witten invariants are log invariants}, Adv. Math. 350 (9), 860--876, 2019.
\bibitem[GWZ]{vGWZ} M. van Garrel, T. Wong, G. Zaimi, \textit{Integrality of relative BPS state counts of toric del Pezzo surfaces}, Comm. in Number Theory and Physics 7 (4), 671--687, 2013.
\bibitem[Gat]{Ga} A. Gathmann, \textit{Relative Gromov-Witten invariants and the mirror formula}, Math. Ann. 325 (2), 393--412, 2003.
\bibitem[GV]{GV} R. Gopakumar, C. Vafa, \textit{M-theory and topological strings II}, \href{https://arxiv.org/abs/hep-th/9812127}{arXiv:9812127}, 1998.
\bibitem[Gra]{Gra} T. Gr\"afnitz, \textit{Theta functions, broken lines and $2$-marked log Gromov-Witten invariants}, \href{https://arxiv.org/abs/2204.12257}{arXiv:2204.12257}, 2022.
\bibitem[GRZ]{GRZ} T. Gr\"afnitz, H. Ruddat, E. Zaslow, \textit{The proper Landau-Ginzburg potential is the open mirror map}, \href{https://arxiv.org/abs/2204.12249}{arXiv:2204.12249}, 2022.
\bibitem[Gro1]{GrP2} M. Gross, \textit{Mirror symmetry for $\mathbb{P}^2$ and tropical geometry}, Adv. in Math. 224 (1), 169--245, 2010.
\bibitem[Gro2]{Gr10} M. Gross, \textit{Mirror Symmetry and Tropical Geometry}, Springer, Regional Conference Series in Mathematics 144, ISBN: 978-0-8218-5232-3, 2011.
\bibitem[GHK]{GHK1} M. Gross, P. Hacking, S. Keel, \textit{Mirror symmetry for log Calabi-Yau surfaces I}, Publ. Math. Inst. Hautes \'Etudes Sci. 122, 65--168, 2015.
\bibitem[GHKK]{GHK2} M. Gross, P. Hacking, S. Keel, M. Kontsevich, \textit{Canonical bases for cluster algebras}, J. Amer. Math. Soc. 31, 497--608, 2018.
\bibitem[GHS]{GHS} M. Gross, P. Hacking, B. Siebert, \textit{Theta functions on varieties with effective anti-canonical class}, \href{https://arxiv.org/abs/1601.07081}{arXiv:1601.07081}, 2016.
\bibitem[GPS]{GPS} M. Gross, R. Pandharipande, B. Siebert, \textit{The tropical vertex}, Duke Math. J. 153 (2), 297--362, 2010.
\bibitem[GS1]{DataI} M. Gross, B. Siebert, \textit{Mirror symmetry via logarithmic degeneration data I}, J. Diff. Geom. 72 (2), 169--338, 2006.
\bibitem[GS2]{DataII} M. Gross, B. Siebert, \textit{Mirror symmetry via logarithmic degeneration data II}, J. Alg. Geom. 19 (4), 679--780, 2010.
\bibitem[GS3]{GS11} M. Gross, B. Siebert, \textit{From real affine geometry to complex geometry}, Ann. of Math. 174 (3), 1301--1428, 2011.
\bibitem[GS4]{Inv} M. Gross, B. Siebert, \textit{An Invitation to toric degenerations},  Surv. Differ. Geom. 16, 43--78, 2011. 
\bibitem[GS5]{LogGW} M. Gross, B. Siebert, \textit{Logarithmic Gromov-Witten invariants}, J. Amer. Math. Soc. 26, 451--510, 2013.
\bibitem[GS6]{Intr1} M. Gross, B. Siebert, \textit{Intrinsic mirror symmetry and punctured Gromov-Witten invariants}, \href{https://arxiv.org/abs/1609.00624}{arXiv:1609.00624}, 2016.
\bibitem[GS7]{Intr2} M. Gross, B. Siebert, \textit{Intrinsic mirror symmetry}, \href{https://arxiv.org/abs/1909.07649}{arXiv:1909.07649}, 2019.
\bibitem[Hos]{Hos} S. Hosono, \textit{Counting BPS states via Holomorphic Anomaly Equations}, Fields Inst. Commun. 38, 57--86, 2003.
\bibitem[IP]{IoPa} E.-N. Ionel, T. H. Parker, \textit{The symplectic sum formula for Gromov-Witten invariants}, Ann. of Math. 159 (3), 935--1025, 2004.
\bibitem[Kat1]{Kat1} K. Kato, \textit{Logarithmic structures of Fontaine-Illusie}, Algebraic analysis, geometry, and number theory (J.-I. Igusa, ed.), Johns Hopkins University Press, 191--224, 1989.
\bibitem[Kat2]{Kat2} F. Kato, \textit{Log smooth deformation and moduli of log smooth curves}, Int. J. Math. 11 (2), 215--232, 2000.
\bibitem[KLR]{KLR} B. Kim, H. Lho, H. Ruddat, \textit{The degeneration formula for stable log maps}, \href{https://arxiv.org/abs/1803.04210}{arXiv:1803.04210}, 2018.
\bibitem[KM]{KM} Y. Konishi, S. Minabe, \textit{Local Gromov-Witten invariants of cubic surfaces via nef toric degeneration}, Arkiv f\"or Matematik 47, 345--360, 2009.
\bibitem[KS]{KS} M. Kontsevich, Y. Soibelman, \textit{Affine structures and non-Archimedean analytic spaces}. Progr. Math. 244, 321--385, 2006.
\bibitem[LR]{LiRu} A.-M. Li and Y. Ruan, \textit{Symplectic surgery and Gromov-Witten invariants of Calabi-Yau 3-folds}, Invent. Math. 145 (1), 151--218, 2001.
\bibitem[Li1]{Li} J. Li, \textit{Stable morphisms to singular schemes and relative stable morphisms}, J. Diff. Geom. 57 (3), 509--578, 2001.
\bibitem[Li2]{Li2} J. Li, \textit{A degeneration formula of GW-invariants}, J. Diff. Geom., 60 (2), 199--293, 2002.
\bibitem[Lin]{Lin} Y.-S. Lin, \textit{Enumerative geometry of del Pezzo surfaces}, \href{https://arxiv.org/abs/2005.08681}{arXiv:2005.08681}, 2020.
\bibitem[MR]{MR} T. Mandel, H. Ruddat, \textit{Descendant log Gromov-Witten invariants for toric varieties and tropical curves}, Trans. Amer. Math. Soc. 373, 1109--1152, 2020.
\bibitem[Mum]{Mum} D. Mumford, \textit{An analytic construction of degenerating abelian varieties over complete rings}, Compos. Math. 24 (3), 1-51, 1972.
\bibitem[NS]{NS} T. Nishinou, B. Siebert, \textit{Toric degenerations of toric varieties and tropical curves}, Duke Math. J. 135 (1), 1-51, 2006.
\bibitem[MFO]{Oberwolfach} D. Abramovich, M. van Garrel, H. Ruddat (organizers), \textit{Logarithmic enumerative geometry and mirror symmetry}, Oberwolfach Report No. 27/2019, \href{https://www.mfo.de/occasion/1925a/www_view}{workshop webpage}, 2019.
\bibitem[Ols]{Ol} M. Olsson, \textit{Logarithmic geometry and algebraic stacks}, Ann. Sci. \'Ecole Norm. Sup. 36 (5), 747--791, 2003.
\bibitem[Pri1]{Pri1} T. Prince, \textit{Smoothing toric Fano surfaces using the Gross-Siebert algorithm}, Proc. London Math. Soc 117 (3), 617--660, 2018
\bibitem[Pri2]{Pri2} T. Prince, \textit{The tropical superpotential for $\mathbb{P}^2$}, Algebr. Geom. 7 (1), 30--58, 2020.
\bibitem[Pum]{Pum} M. Pumperla, \textit{Unifying constructions in toric mirror symmetry}, PhD thesis, Hamburg University, 2011.
\bibitem[SS]{SS} M. Sch\"utt, T. Shioda, \textit{Elliptic surfaces}, \href{https://arxiv.org/abs/0907.0298}{arXiv:0907.0298}, 2010.
\bibitem[SYZ]{SYZ} A. Strominger, S.-T. Yau, E. Zaslow, \textit{Mirror symmetry is T-duality}, Nucl. Phys. B 479 (1-2): 243--259, 1996.
\bibitem[Tak1]{Ta1} N. Takahashi, \textit{Curves in the complement of a smooth plane cubic whose normalizations are $\mathbb{A}^1$}, \href{https://arxiv.org/abs/alg-geom/9605007}{arXiv:9605007}, 1996.
\bibitem[Tak2]{Ta2} N. Takahashi, \textit{Log mirror symmetry and local mirror symmetry}, Comm. Math. Phys. 220 (2), 293--299, 2001.
\end{thebibliography}
\end{document}